\magnification=1000
\hsize=11.7cm
\vsize=18.9cm
\lineskip2pt \lineskiplimit2pt
\nopagenumbers

\hoffset=-1truein
\voffset=-1truein

\advance\voffset by 4truecm
\advance\hoffset by 4.5truecm

\newif\ifentete

\headline{\ifentete\ifodd	\count0 
      \rlap{\head}\hfill\tenrm\llap{\the\count0}\relax
    \else
        \tenrm\rlap{\the\count0}\hfill\llap{\head} \relax
    \fi\else
\global\entetetrue\fi}

\def\entete#1{\entetefalse\gdef\head{#1}} 
\entete{}

\input amssym.def
\input amssym.tex

\def\-{\hbox{-}}
\def\.{{\cdot}}

\def\K{{\cal K}}
\def\F{{\cal F}}

\def\L{{\cal L}}
\def\M{{\cal M}}
\def\N{{\cal N}}
\def\P{{\cal P}}

\def\T{{\cal T}}

\def\I{{\cal I}}

\def\H{{\cal H}}
\def\s{{\cal S}}

\def\X{{\cal X}}

\def\ch{\frak c\frak h}
\def\ad{\frak a\frak c}
\def\ab{\frak a\frak b}
\def\Ab{\frak A\frak b}

\def\Gr{\frak G\frak r}

\def\Fct{\frak F\frak c\frak t}
\def\Nat{\frak N\frak a\frak t}
\def\Ker{\frak K\frak e\frak r}

\def\int{\frak i\frak n\frak t}

\def\qq{\quad{\rm and}\quad}

\def\lv{\frak l\frak v}

\def\res{\frak r\frak e\frak s}

\def\too{\longrightarrow}
\def\aut{\frak a\frak u\frak t}

\def\Set{\frak S\frak e\frak t}
\def\Loc{\frak L\frak o\frak c}
\def\loc{\frak l\frak o\frak c}

 3
 2
\font\large=cmr10  scaled \magstep 2
 2
\font\larti=cmti10  scaled \magstep 2
 1
 2

\font\cds=cmr7
\font\cdt=cmti7
\font\cdy=cmsy7
\font\cdi=cmmi7

\count0=1

\centerline{\large Existence, uniqueness and functoriality}
\medskip
\centerline{\large of the perfect locality over a Frobenius {\larti P}-category}
\medskip
\centerline{\bf Lluis Puig }
\medskip
\noindent 
\centerline{\cds CNRS, Institut de Math\'ematiques de Jussieu, lluis.puig@imj-prg.fr}
\par
\noindent
\centerline{\cds 6 Av Bizet, 94340 Joinville-le-Pont, France}

\medskip
\noindent
{\bf Abstract:} {\cds  Let {\cdt p} be a prime, {\cdt P} a finite {\cdt p}-group and {\cdy F} a {\cdt Frobenius P-category\/}. The question on the existence of  a suitable category {\cdy L}{${}^{_{\rm sc}}$} extending the 
{\cdt full\/} subcategory of  {\cdy F} over the set of   {\cdt {\cdy F}-selfcentralizing\/} subgroups of {\cdt P} goes back to Dave Benson in 1994. In 2002  Carles Broto, Ran Levi and Bob Oliver formulate the existence and the uniqueness of the category {\cdy L}{${}^{_{\rm sc}}$} in terms of the annulation of an {\cdt obstruction {\cds 3}-cohomology element\/} 
 and of the vanishing of a  {\cdt {\cds 2}-cohomology group\/}, and they state a sufficient condition for the vanishing of these   {\cdt n-cohomology groups\/}. Recently, Andrew Chermak has proved the existence and the uniqueness of 
 {\cdy L}{${}^{_{\rm sc}}$} {\cdt via\/} his {\cdt objective partial groups\/}, and Bob Oliver, following some of Chermak's methods, has also proved the vanishing of those {\cds n-cohomology groups\/} for n {\cdi  >} 1, both applying the 
 {\cdt Classification of the finite   simple groups\/}. Here we  give   {\cdt direct} proofs of the existence and the uniqueness of {\cdy L}{${}^{_{\rm sc}}$}; moreover, in [11] we already show that 
 {\cdy L}{${}^{_{\rm sc}}$} can be completed  in a suitable category {\cdy L} extending  {\cdy F} and we also prove some {\cdt functoriality\/} of this correspondence.  }

\bigskip
\noindent
{\bf £1. Introduction  }
\medskip
£1.1. Let $p$ be a prime, $P$ a finite $p\-$group and $\F$ a 
{\it Frobenius  $P\-$cate-gory\/}~[10]. The question on the existence of  a suitable category  $\P^{^{\rm sc}}$ extending the {\it full\/} subcategory 
$\F^{^{\rm sc}}$ of  $\F$ over the set of   {\it $\F\-$selfcentralizing\/} subgroups of~$P$~[10,~\S3] goes back to Dave Benson in 1994 [1]. Indeed, considering our  suggestion of constructing a {\it topological space\/} from the family of {\it classi-fying spaces\/} of the  {\it $\F\-$localizers\/} --- a family of finite groups indexed by the {\it $\F\-$selfcentralizing\/} subgroups of $P$ we had just introduced at that time [8] --- Benson, in his tentative construction,  had foreseen the interest of this extension, actually as a generalization for Frobenius $P\-$categories of our old {\it $O\-$locality\/} for finite groups in [7].

\medskip
£1.2. In [2]  Carles Broto, Ran Levi and Bob Oliver formulate the existence and the uniqueness of the category $\P^{^{\rm sc}}$ in terms of the annulation of an {\it obstruction $3\-$cohomology element\/} and of
the vanishing of a  {\it $2\-$cohomology group\/}, respectively. They actually state a sufficient condition for the va-nishing of the corresponding {\it $n\-$cohomology groups\/} and moreover, assuming the existence of $\P^{^{\rm sc}}\,,$ they succeed on the construction of a {\it classifying space\/}.

\medskip
£1.3. As a matter of fact, if $G$ is a finite group and $P$  a Sylow $p\-$subgroup of $G\,,$ the corresponding {\it Frobenius $P\-$category\/} $\F_{\!G}$ introduced in [7] admits an extension $\P_G$ defined over {\it all\/} the subgroups of~$P$ where,   for any pair of subgroups $Q$ and $R$ of~$P\,,$ the set of {\it morphisms\/} from $R$ to $Q$  is the following {\it quotient set\/} of the {\it $G\-$transporter\/}
$$\P_G (Q,R) = T_G (R,Q)\big/{\Bbb O}^p\big(C_G(R)\big)
\eqno £1.3.1.$$ 
\eject
\noindent
Analogously, in the general setting, if we are interested in some {\it functoriality\/}
 for our constructions, we need not only the existence of $\P^{^{\rm sc}}$ but the existence of a suitable category $\P$ 
 extending $\F$ and containing  $\P^{^{\rm sc}}$ as a {\it full\/} subcategory. Soon after~[2], on the one hand we showed that the {\it contravariant\/} functor from $\F_G$ mapping $Q$ on $C_G (Q)\big/{\Bbb O}^p\big(C_G(Q)\big)$ can be indeed generalized to a {\it contravariant\/} functor $\frak c^\frak h_\F$ from any {\it Frobenius $P\-$category\/} $\F$ (see £2.4 below); on the other hand, we already prove in [9] that the existence of the so-called {\it perfect $\F^{^{\rm sc}}\-$locality\/} 
 $\P^{^{\rm sc}}$ forces the existence  of a {\it unique extension\/} $\P$ 
 of~$\F$  by~$\frak c^\frak h_\F\,,$  called the {\it perfect $\F\-$locality\/}.

\medskip 
£1.4. Recently, Andrew Chermak [3] has proved the existence and the uniqueness of $\P^{^{\rm sc}}\!$ {\it via\/} his {\it objective partial groups\/}, and Bob Oliver [6], following some of Chermak's methods, has also proved for $n\ge 2$ the vanishing of the {\it $n\-$cohomology groups\/} mentioned above. 
In reading their preprints, we were disappointed not only because their proofs depend on the so-called  {\it Classification of the finite simple groups\/} (CFSG), but because in their arguments they need strong properties of  the finite groups.  Indeed, since [7] we are convinced that a previous classification of the so-called ``local structures'' will be the way to clarify CFSG in future versions; thus, our effort in creating the {\it Frobenius $P\-$categories\/} was  directed to provide a precise formal support to the vague notion of ``local structures'', independently of
``environmental'' finite groups and of most of finite group properties.

\medskip
£1.5. Here we will show that till now our intuition was correct, namely that there is a {\it direct\/} proof of the existence and the uniqueness of $\P^{^{\rm sc}}\,;$ that is to say, a proof that can be qualified of {\it inner\/} or 
{\it tautological\/} in the sense that only pushes far enough the initial axioms of  {\it Frobenius $P\-$categories\/}. Moreover, as~we mention above,
 the existence and the uniqueness of the {\it perfect $\F^{^{\rm sc}}\-$locality\/} $\P^{^{\rm sc}}$ will guarantee  the existence and the uniqueness of the 
 {\it perfect $\F\-$locality\/} $\P$ defined over all the subgroups of~$P$ --- which no longer can be described in terms of Chermak's {\it objective partial groups\/} --- and then it makes sense to discuss the {\it functorial\/}
 nature of the correspondence mapping $\F$ on~$\P\,.$

\medskip
£1.6. Let us explain how our method works. In [11, Chap.~18] we  introduce the 
{\it $\F\-$localizers\/} mentioned above and, as a matter of fact, we already introduce the $\F\-$localizer $L_\F (Q)$ for {\it any\/} subgroup $Q$ of $P$ 
(see Theorem 2.10 below), which is indeed an extension of the group $\F (Q)$ of 
$\F\-$automorphisms of $Q$  by the $p\-$group $\frak c^\frak h_\F (Q)$ 
(cf.~£1.3). More generally, in [11, Chap.~17] we introduce the {\it $\F\-$localities\/} as a wider framework where to look for the {\it perfect $\F\-$locality\/}. Namely, considering the category $\T_P$ --- where the objects are all the subgroups of $P\,,$  the set of morphisms from $R$ to $Q$ is the {\it $P\-$transporter\/} 
$T_P (R,Q)\,,$ and  the composition is induced by the product in $P$ --- 
we call {\it $\F\-$locality\/} any extension $\pi\,\colon \L\to \F$ of the category 
$\F\,,$ endowed with a functor $\tau\,\colon \T_P\to \L$ 
such that the composition $\pi\circ \tau\,\colon \T_P\to \F$ is the canonical functor defined by the conjugation in~$P\,;$ of course, we add some suitable conditions as {\it divisibility\/}  and  {\it $p\-$coherence\/} (see £2.8 below). As a matter of fact, a {\it perfect $\F\-$locality\/} is just a {\it divisible $\F\-$locality\/} $\L$ where the group $\L (Q)$ of $\L\-$automorphisms of any subgroup $Q$ of $P$ coincides with the {\it $\F\-$localizer\/} of $Q$ (see £2.13 below).

\medskip
£1.7. It turns out that there are indeed other $\F\-$localities --- easier to construct --- which deserve consideration; their construction depends on the existence of the 
{\it $\F\-$basic $P\times P\-$sets\/} $\Omega$ (see section 3 below) introduced in [11,~Chap.~21], which allows the realization of $\F$ inside  the 
{\it symmetric group\/} of~$\Omega$ and then it allows the possibility of 
considering {\it localities\/}  as defined in~[7].  In [11,~Chap.~22] we introduce the so-called {\it basic $\F\-$locality\/} $\L^{^{\rm b}}$ which is canonically associated with $\F$ (see section 4 below). More precisely,  in [11, Chap.~24] we show that  the very structure of a {\it perfect 
$\F^{^{\rm sc}}\-$locality\/}~$\P^{^{\rm sc}}$  supplies particular 
{\it  $\F\-$basic  $P\times P\-$sets\/}.

\medskip
£1.8.  From anyone of these {\it  $\F\-$basic  $P\times P\-$sets\/} $\Omega$ we can  construct a  par-ticular {\it  $\F\-$locality\/} $\L^{^{\Omega}}$ (see £4.7 below) in such a way that  the {\it full\/} subcategory $\L^{^{\Omega,\rm sc}}$ of $\L^{^{\Omega}}$ over the set of {\it $\F\-$selfcentralizing\/} subgroups 
of~$P$   admits a  quotient $\bar\L^{^{\rm n,sc}}\!$ --- independent of our choice --- containing $\P^{^{\rm sc}}\!$ (see~Corollary~£5.20 below). The point is that these particular {\it $\F\-$basic $P\times P\-$sets\/} can be described {\it directly\/}, without assuming the existence of $\P^{^{\rm sc}}$ (see Proposition~3.4 below); hence, we can introduce  the so-called {\it natural $\F^{^{\rm sc}}\-$locality\/}~$\bar\L^{^{\rm n,sc}}\,;$ then, the
$\F^{^{\rm sc}}\!\-$locality $\bar\L^{^{\rm n,sc}}$ supplies a {\it support\/} for the proof of the existence and  the uniqueness of the {\it perfect 
$\F^{^{\rm sc}}\-$locality\/}. Moreover, the  {\it basic $\F\-$locality\/} 
$\L^{^{\rm b}}$ also admits a quotient $\bar\L^{^{\rm b}}$ in such a way
 that   the {\it full\/} subcategory  $\bar\L^{^{\rm b,sc}}$ of $\bar\L^{^{\rm b}}$ over the set of {\it $\F\-$selfcentralizing\/} subgroups of~$P$  contains $\bar\L^{^{\rm n,sc}}$ (see~£4.13.3 below); this fact
 is the key point in order to discuss {\it functoriality\/}.

\medskip
£1.9. More explicitly, we replace the whole set of $\F\-$selfcentralizing subgroups
of $P$ by a nonempty set $\frak X$ of them, {\it containing any subgroup of $P$ 
which admits an $\F\-$morphism from some subgroup in $\frak X\,,$\/} and replace the categories $\F^{^{\rm sc}}\!$ and 
$\bar\L^{^{\rm n,sc}}$ by their respective {\it full\/} subcategories $\F^{^{\frak X}}\!$ and $\bar\L^{^{\rm n,\frak X}}\!$ 
over $\frak X\,;$ then, we prove the existence and the uniqueness  of the  {\it perfect $\F^{^{\frak X}}\!\-$sublocality\/} 
$\P^{^{\frak X}}$ in  $\bar\L^{^{\rm n,\frak X}}\!$ arguing by induction on $\vert\frak X\vert\,.$ Our  proof depends on the annulation of the  cohomology classes of  suitable $2\-$ and $1\-${\it cocycles\/}, respectively; although the corresponding 
cohomology groups need not vanish,{\footnote{\dag}{\cds We thank Bob Oliver who  showed us explicit examples.}} some vanishing result in cohomology allows us to push our {\it cocycles\/} towards an easier situation where the existence and the uniqueness of a suitable {\it localizer\/} completes the proof (see Theorem~6.22 below). The arguments in order to prove this vanishing  result  together with the lifting result mentioned below admit a wider framework that we develop  in [12].

\medskip
£1.10. At this point, we have the  {\it perfect $\F^{^{\rm sc}}\!\-$locality\/} 
$\P^{^{\rm sc}}$ as  an  $\F^{^{\rm sc}}\!\-$sublo-cality of $\bar\L^{^{\rm n,sc}}$ and therefore,  as  an  $\F^{^{\rm sc}}\!\-$sublocality of $\bar\L^{^{\rm b,sc}}\!$ (cf.~£1.8). But, as mentioned in~£1.3 above,  the existence of the 
 {\it perfect $\F^{^{\rm sc}}\-$locality\/} $\P^{^{\rm sc}}$ forces the existence  of the {\it perfect $\F\-$locality\/} $\P$  [11, Theorem~20.24]; more generally, for any {\it $p\-$coherent $\F\-$locality\/} $\hat\L\,,$ denoting by 
 $\hat\L^{^{\rm sc}}$  the {\it full\/} subcategory of $\hat\L$ over the set of 
 $\F\-$selfcentralizing subgroups of $P\,,$ any {\it $\F\-$locality\/} functor from $\P^{^{\rm sc}}$ to $\hat\L^{^{\rm sc}}$ can be extended to a unique  
 {\it $\F\-$locality\/} functor from $\P$ to~$\hat \L$ (see section~£7 below, where we  give a slightly different and more correct proof{\footnote{\dag}{\cds Our argument in [11, 20.16] has been scratched.  }}  of the existence of $\P$). In particular, we get an  {\it $\F\-$locality\/} functor from the
{\it perfect $\F\-$locality\/} $\P$ to $\bar\L^{^{\rm b}}\,;$ once again, a {\it homotopically trivial\/} situation exhibited in~[12] allows us to show that this functor can be lifted to an essentially unique  functor from $\P$ to the {\it basic $\F\-$locality\/} $\L^{^{\rm b}}\!$
(see Theorem~£8.10 below).

\medskip
£1.11. Let $P'$ be another finite $p\-$group, $\F'$ a Frobenius $P'\-$category 
and $\alpha\,\colon P\to P'$ an {\it $(\F,\F')\-$functorial\/} group homomorphism [11, 12.1], so that it determines a so-called 
{\it Frobenius functor\/}
$$\frak f_\alpha : \F\too \F'
\eqno £1.11.1;$$
once we know the existence and the uniqueness of the respective 
{\it perfect $\F\-$\/} and {\it $\F'\-$localities\/} $\P$ and $\P'\,,$ it is reasonable to ask for the existence and the uniqueness of a suitable isomorphism class of functors
$$\frak l_\alpha : \P\too \P'
\eqno £1.11.2\phantom{.}$$
lifting $\frak f_\alpha\,;$ but here we only get a positive answer for the quotients
$$\bar\frak l_\alpha : \bar\P = \P/[\frak c_\F^{\frak h},\frak c_\F^{\frak h}]\too 
\bar\P' = \P'/[\frak c_{\F'}^{\frak h},\frak c_{\F'}^{\frak h}]
\eqno £1.11.3;$$
then,  for a third finite $p\-$group $P''$ together with a Frobenius $P''\-$category 
$\F''$ and an 
{\it $(\F',\F'')\-$functorial\/} group homomorphism $\alpha'\,\colon P'\to P''\,,$ 
{\it the functors $\bar\frak l_{\alpha'}\circ\bar\frak l_{\alpha}$ and 
$\bar\frak l_{\alpha'\circ\alpha}$
are naturally isomorphic\/} (see section 9 below).

\medskip
£1.12. Actually, if $\alpha$ is surjective and $\F' = \F/{\rm Ker}(\alpha)$ then the existence of $\frak l_\alpha$ follows from [11, Theorem~17.18] where, assuming the existence of the~{\it perfect $\F\-$locality\/}~$ \P\,,$ we exhibit a 
{\it perfect $\F'\-$locality\/} $\P'$ and a functor  
$\frak l_\alpha\,\colon  \P\to \P'$ lifting $\frak f_\alpha\,.$ Thus, here we actually may assume that  $\alpha$ is injective; in this case, we start by getting a relationship between  the 
{\it natural $\F^{^{\rm sc}}\-$locality\/} $\bar\L^{^{\rm n,sc}}$ and 
{\it the basic $\F'\-$locality\/}~$\L'^{^{\rm b}}\,;$  more explicitly, the converse image  ${\rm Res}_{\F^{^{\rm sc}}}(\L'^{^{\rm b}})$ of~$\F^{^{\rm sc}}$
 in~$\L'^{^{\rm b}}$ is clearly a {\it $p\-$coherent $\F^{^{\rm sc}}\-$locality\/} and we will exhibit  a {\it $\F^{^{\rm sc}}\-$locality\/} functor  from 
 $\bar\L^{^{\rm n,sc}}$ to a suitable {\it quotient $\F^{^{\rm sc}}\-$locality\/}  
 $\overline{{\rm Res}_{\F^{^{\rm sc}}}(\L'^{^{\rm b}})}$ (see Theorem~£9.10 below). Finally,
 from the  {\it $\F\-$locality\/} functors 
 $$\bar\P^{^{\rm sc}}\too \L^{^{\rm b,sc}}\too 
 \overline{{\rm Res}_{\F^{^{\rm sc}}}(\L'^{^{\rm b}})}
 \longleftarrow {\rm Res}_{\F^{^{\rm sc}}}(\L'^{^{\rm b}})  
 \longleftarrow {\rm Res}_{\F^{^{\rm sc}}}(\bar\P')
 \eqno £1.12.1$$
 we will obtain an  {\it $\F\-$locality\/} functor $\P^{^{\rm sc}}\too {\rm Res}_{\F^{^{\rm sc}}}(\bar\P')$ which can be extended to  an  {\it $\F\-$locality\/} functor $\bar\P\too {\rm Res}_{\F}(\bar\P')$ (see Theorem~£9.15 below).
 \eject

\bigskip
\noindent
{\bf £2. Frobenius $P\-$categories and coherent $\F\-$localities}

\medskip
£2.1. Denote  by  $\frak i\Gr$ the category formed by the finite groups and by the injective group  homomorphisms. Recall that, for any category $\frak C\,,$ we denote by $\frak C^\circ$  the {\it opposite\/} category and, for any $\frak C\-$object 
$C\,,$ by $\frak C_C$  (or $(\frak C)_C$ to avoid confusion)  the category of ``$\frak C\-$morphisms to $C$'' 
[11, 1.7]. If any $\frak C\-$object  admits {\it inner\/} automorphisms then we denote by $\tilde\frak C$
the corresponding quotient and call it the {\it exterior\/} quotient of $\frak C$ [11,~1.3]. Let $p$ be a prime; 
for any finite $p\-$group $P$ we denote   by  $\F_{\!P}$ the subcategory of  $\frak i\Gr$ where the objects are all the  subgroups of $P$ and where the morphisms are the group homomorphisms induced by conjugation by elements of $P\,.$

 \medskip
£2.2. A {\it Frobenius  $P\-$category\/} $\F$ is a subcategory 
of $\frak i\Gr$ containing $\F_{\!P}$ where the objects are all the  subgroups of $P$
and the morphisms fulfill the following three conditions [11, 2.8 and Proposition~2.11]
\medskip
\noindent
£2.2.1\quad {\it For any subgroup $Q$ of $P\,,$ the inclusion functor $(\F)_Q\to 
(\frak i\Gr)_Q$~is~full.\/}
\smallskip
\noindent
£2.2.2\quad {\it $\F_P (P)$ is a Sylow $p\-$subgroup of $\F (P)\,.$\/}
\smallskip
\noindent
£2.2.3\quad {\it Let $Q$ be a subgroup of $P$ fulfilling $\xi \big(C_P (Q)\big)
= C_P\big(\xi (Q)\big)$ for any $\F\-$morphism $\xi\,\colon Q\.C_P(Q)\to P\,,$
let $\varphi\,\colon Q\to P$ be an $\F\-$morphism  and let $R$ be a subgroup of $N_P\big(\varphi(Q)\big)$ containing $\varphi (Q)$ such that $\F_P(Q)$ contains the action of $\F_R \big(\varphi(Q)\big)$ over $Q$ via $\varphi\,.$ Then there is an $\F\-$morphism 
$\zeta\,\colon R\to P$ fulfilling $\zeta\big(\varphi (u)\big) = u$ for any $u\in Q\,.$\/}
\medskip
\noindent
As in [11,~1.2], for any pair of subgroups $Q$ and $R$ of $P\,,$ we denote by $\F (Q,R)$ the set of $\F\-$morphisms from $R$ to $Q$ and set $\F (Q) = \F (Q,Q)\,.$ If $G$ is a finite subgroup admitting $P$ as a Sylow $p\-$subgroup, we denote by 
$\F_G$ the {\it Frobenius $P\-$category\/} where  the morphisms
 are the group homomorphisms induced by the conjugation by elements of $G\,.$

\medskip
£2.3. Fix a Frobenius $P\-$category $\F\,;$ for any subgroup $Q$ of $P$ and
any subgroup $K$ of the group ${\rm Aut}(Q)$ of automorphisms of $Q\,,$
we say that $Q$ is {\it fully $K\-$normalized\/} in $\F$ if, for any $\F\-$morphism 
$\xi\,\colon Q\.N_P^K (Q)\to P\,,$ we have [11, 2.6]
$$\xi \big(N_P^K (Q)\big) = N_P^{\,{}^\xi \!K}\big(\xi (Q)\big)
\eqno £2.3.1,$$
 where $N_P^K (Q)$ is the converse image of  $K$ in $N_P (Q)$ {\it via\/} the canonical group homomorphism 
 $N_P (Q)\to {\rm Aut}(Q)\,,$ and ${}^\xi \! K$ denotes the image in  ${\rm Aut}\big(\xi (Q)\big)$ of $K$ {\it via\/} 
$\xi\,.$  Recall that if $Q$ is fully $K\-$normalized in $\F$ then we have a new Frobenius $N_P^K(Q)\-$category $N_\F^K(Q)$  [11, 2.14 and Proposition~2.16] where, for any pair of subgroups $R$ and $T$ of $N_P^K (Q)\,,$ the set of morphisms $\big(N_\F^K(Q)\big)(R,T)$ is the set of group homomorphisms from 
$T$ to $R$ induced by the $\F\-$morphisms  $\psi\,\colon Q\.T\to Q\.R$ which stabilize $Q$ and induce on it an\break
\eject
\noindent element of $K\,.$ Note that from [11,~statement~2.13.2 and Corollary~5.14] it is not difficult to prove that if $T$ contains $Q\.C_P (Q)$ then the canonical map
$$\big(N_\F^K(Q)\big)(R,T)\too T_{K\cap \F(Q)}\big(\F_T(Q),\F_R(Q)\big)
\eqno £2.3.2\phantom{.}$$
is surjective.

\medskip
£2.4. We denote by $H_\F$ the {\it $\F\-$hyperfocal\/} subgroup of $P\,,$ which is the subgroup generated by the sets $\{u^{-1}\sigma (u)\}_{u\in Q}$ where $Q$ runs over the set of subgroups of $P$ and $\sigma$ over the set of $p'\-$elements of 
$\F (Q)$ [11,~13.2]. As above, for any subgroup $Q$ of $P$ {\it fully centralized\/} in $\F$ --- namely, {\it fully $\{1\}\-$normalized\/} in~$\F$ --- we have the Frobenius 
$C_P (Q)\-$category $C_\F (Q) = N_\F^{\{1\}}(Q)$ and therefore we can consider the $C_\F (Q)\-$hyperfocal subgroup $H_{C_\F (Q)}$ of $C_P (Q)\,;$
 then,~in [11,~Proposition~13.14] we exhibit a unique {\it contravariant\/} functor
$$\frak c^\frak h_\F : \F\too \widetilde\Gr
\eqno£2.4.1,$$
where $\widetilde\Gr$ denotes the {\it exterior\/} quotient of the category $\Gr$ of finite groups (cf.~£2.1), mapping any  subgroup $Q$ of $P$ fully centralized in $\F$ on the quotient $C_P (Q)/H_{C_\F (Q)}$ 
and any $\F\-$morphism $\varphi\,\colon R\to Q$ from a subgroup $R$ of $P$ {\it fully centralized\/} in $\F$ 
on a $\widetilde\Gr\-$morphism induced by an $\F\-$morphism
$$\varphi (R)\.C_P (Q)\too R\.C_P (R)
\eqno £2.4.2\phantom{.}$$
sending $\varphi (v)$ to $v$ for any $v\in R$ (cf.~condition~£2.2.3).

\medskip 
£2.5. We say that a subgroup $U$ of $P$ is {\it $\F\-$stable\/} if we have $\varphi (Q\cap U)\i U$
for any subgroup $Q$ of $P$ and any $\F\-$morphism $\varphi\,\colon Q\to P\,;$ then, setting $\bar P = P/U\,,$
 there is  a Frobenius $\bar P\-$category $\bar \F = \F/U$ such that the canonical homomorphism 
 $\varpi\,\colon P\to \bar P$ is {\it $(\F,\bar\F)\-$functorial\/} and that the corresponding {\it Frobenius functor\/}
 $\frak f_\varpi\,\colon \F\to \bar \F$ is {\it full\/} over the subgroups of $P$ containing~$U$ [11,~Proposition~12.3]. In particular, if $Q$ is a subgroup of~$P$ {\it fully normalized\/} in $\F\,,$ it follows from 
 [11,~Proposition~13.9] that $H_{C_\F (Q)}$ is an {\it $N_\F (Q)\-$stable\/} subgroup of $N_P (Q)$
 and therefore we can consider the {\it quotients\/}
 $$\overline{N_P (Q)}^\frak h = N_P (Q)/H_{C_\F (Q)}\qq 
 \overline{N_\F (Q)}^\frak h =N_\F (Q)/H_{C_\F (Q)}
 \eqno £2.5.1.$$

\medskip
£2.6. We say that a subgroup $Q$ of $P$  is  {\it $\F\-$selfcentralizing\/} if we have
$$C_P\big(\varphi (Q))\i \varphi (Q)
\eqno £2.6.1\phantom{.}$$
 for any $\varphi \in \F (P,Q)\,;$  we denote by $\F^{^{\rm sc}}$ the full subcategory of $\F$ over the set of $\F\-$selfcentralizing subgroups of $P\,.$ More generally, as mentioned above we consider  a nonempty set $\frak X$ of  subgroups of $P$ containing any subgroup of~$P$ admitting an $\F\-$morphism from some subgroup in $\frak X\,,$ and then we denote by~$\F^{^{\frak X}}$ the {\it full\/} subcategory of $\F$ over the set $\frak X$ of objects; in most situations, the subgroups in $\frak X$ will be
$\F\-$selfcentralizing and if $\frak X$ is the set of all the $\F\-$selfcentralizing subgroups of $P\,,$ we write $\rm sc$ instead of 
$\frak X\,.$
\eject

\medskip
£2.7.  Denote by $\T_P^{^{\frak X}}$ the {\it full\/} subcategory of $\T_P$ (cf.~£1.6) over the set 
$\frak X$ and  by $\kappa^{_{\frak X}}\,\colon\T_{P}^{^{\frak X}}\to \F^{^{\frak X}}$ the canonical functor 
determined by the conjugation. An  $\F^{^{\frak X}}\-${\it locality\/} $\L^{^{\frak X}}$ is a category, where
$\frak X$ is the set of objects, endowed with two functors
$$\tau^{_{\frak X}} : \T_{ P}^{^{\frak X}}\too \L^{^{\frak X}}\quad{\rm and}\quad 
\pi^{_{\frak X}} : \L^{^{\frak X}}\too \F^{^{\frak X}}
\eqno £2.7.1\phantom{.}$$
which are the identity on the set of objects and fulfill $\pi^{_{\frak X}}\circ\tau^{_{\frak X}} =\kappa^{_{\frak X}}\,,$ 
$\pi^{_{\frak X}}$~being {\it full\/}; as above, for any pair of subgroups $Q$ and~$R$ in~$\frak X\,,$ we denote 
by~$\L^{^{\frak X}}\! (Q,R)$ the set of $\L^{^{\frak X}}\-$morphisms from $R$ to~$Q$ and by 
$$\tau_{_{Q,R}}^{_{\frak X}}\,\colon\T_{P}^{^{\frak X}} (Q,R)\too \L^{^{\frak X}} \!(Q,R)\qq
\pi_{_{Q,R}}^{_{\frak X}} \,\colon  \L^{^{\frak X}}\! (Q,R)\too \F^{^{\frak X}} \!(Q,R)
\eqno £2.7.2\phantom{.}$$
the corresponding maps; we write $Q$ only once if $Q= R\,.$

\medskip
£2.8. We say that $\L^{^{\frak X}}$ is {\it divisible\/} if, for any pair of subgroups 
$Q$ and $R$ in~$\frak X\,,$ ${\rm Ker}(\pi_{_R}^{_{\frak X}})$~acts  regularly on the ``fibers'' of 
$\pi_{_{Q,R}}^{_{\frak X}}\,,$ and that $\L^{^{\frak X}}$ is {\it coherent\/} if moreover for any 
$x\in\L^{^{\frak X}} (Q,R)$ and any  $v\in R$ we have [11,~17.8 and~17.9]
$$ x\.\tau_{_R}^{_{\frak X}} (v) = \tau_{_Q}^{_{\frak X}}\Big(\big(\pi_{_{Q,R}}^{_{\frak X}}(x)\big) (v)\Big)\.x 
\eqno £2.8.1;$$ 
in this case, it actually follows from [11,~Remark~17.11] that, if $\tau^{_{\frak X}}$ is {\it faithful\/} and $N_P(T)$
contains $Q$ and $R$ for some $T\in \frak X$ {\it fully normalized\/} in $\F\,,$ $\tau^{_{\frak X}}$ induces a bijection
$$\big(N_{\F^{^{\frak X}}} (T)\big) (Q,R) \cong 
\F_{\L^{^{\frak X}}\! (T)} \big(\tau^{_{\frak X}} (Q),\tau^{_{\frak X}}(R)\big)
\eqno £2.8.2;$$
moreover, the {\it divisibility\/} determines a functor
$$\Ker (\pi^{_{\frak X}}) : \tilde\F^{^{\frak X}}\too \Ab
\eqno £2.8.3\phantom{.}$$
sending any $Q\in \frak X$ to ${\rm Ker}(\pi_{_Q}^{_{\frak X}})\,;$  indeed, any  
$z\in {\rm Ker}(\pi_{_Q}^{_{\frak X}})$ induces the identity on $Q$ and therefore it stabilizes the image of $R$ for any 
$x\in\L^{^{\frak X}} (Q,R)\,,$  so that there is a unique $w\in {\rm Ker}(\pi_{_R}^{_{\frak X}})$ fulfilling $x\.w = z\.x\,;$
then,  the correspondence sending $z$ to $w$ is a group homomorphism which only depends on the class
of $x$ in $\tilde\F^{^{\frak X}}\!(Q,R)$ [11,~Proposition~17.10].  On the other hand, we say that $\L^{^{\frak X}}$ is {\it $p\-$coherent\/}  if it is {\it coherent\/} and, for any subgroup $Q$ in  $\frak X\,,$ the kernel ${\rm Ker} (\pi_{_Q}^{_{\frak X}})$ is a $p\-$group; in this case, it follows from  [11,~17.13] that if $Q$ is fully centralized in $\F$ then we have
$$H_{C_\F (Q)} \i {\rm Ker}(\tau_{_Q}^{_{\frak X}})
\eqno £2.8.4;$$
we say that $\L^{^{\frak X}}$ is {\it $P\-$bounded\/} if it is {\it coherent\/} and, for any subgroup $Q$
 in~$\frak X$ fully normalized in $\F\,,$ we have ${\rm Ker} (\pi_{_Q}^{_{\frak X}})\i 
 \tau_{_Q}^{_{\frak X}}\big(N_P(Q)\big)\,.$ 
 Finally, we say that $\L^{^{\frak X}}$ is {\it perfect\/} if it is {\it $P\-$bounded\/} and
  for any subgroup $Q$ in $\frak X$   fully centralized in $\F$ we have [11,~17.13]
$$H_{C_\F (Q)} = {\rm Ker}(\tau_{_Q}^{_{\frak X}})
\eqno £2.8.5.$$
\eject

\medskip
£2.9. If $\L'^{^{\frak X}}$ is a second $\F^{^{\frak X}}\!\-${\it locality\/} with structural functors $\tau'^{_{\frak X}}$ 
and $\pi'^{_{\frak X}}\,,$  we call {\it $\F^{^{\frak X}}\!\-$locality functor\/} from $\L^{^{\frak X}}$ to $\L'^{^{\frak X}}$ 
any functor $\frak l^{^{\frak X}}\colon \L^{^{\frak X}}\!\to \L'^{^{\frak X}}$ fulfilling 
$$\tau'^{_{\frak X}} = \frak l^{^{\frak X}}\!\circ\tau^{_{\frak X}}\qq 
\pi'^{_{\frak X}} \!\circ\frak l^{^{\frak X}} = \pi^{_{\frak X}}
\eqno £2.9.1;$$
the composition of two {\it $\F^{^{\frak X}}\!\-$locality functors\/} is obviously an {\it $\F^{^{\frak X}}\!\-$locality functor\/}; we say that two  {\it $\F^{^{\frak X}}\!\-$locality functors\/} $\frak l^{^{\frak X}}$ and 
$\bar\frak l^{^{\frak X}}$ from $\L^{^{\frak X}}$ to $\L'^{^{\frak X}}$ are {\it natu-rally $\F^{^{\frak X}}\!\-$isomorphic\/} whenever we have a {\it natural isomorphism\/} $\lambda^{^{\!\frak X}}\colon \frak l^{^{\frak X}}
\cong \bar\frak l^{^{\frak X}}$ such that 
$\pi'^{_{\frak X}} \! * \lambda^{^{\!\frak X}} = {\rm id}_{\pi^{_{\!\frak X}} }\,;$ in this case,  $(\lambda^{^{\!\frak X}})_Q$ belongs to $\Ker (\pi'^{_{\frak X}}_{_Q})$ for any $Q\in \frak X$ 
and, since 
$$\frak l^{^{\frak X}}\big(\tau^{_{\!\frak X}}_{_{P,Q}} (1)\big) = 
\tau'^{_{\!\frak X}} _{_{P,Q}} (1) = \bar\frak l^{^{\frak X}}\big(\tau^{_{\!\frak X}}_{_{P,Q}} (1)\big) 
\eqno £2.9.2,$$
if $\L'^{^{\frak X}}$ is {\it divisible\/} then $\lambda^{^{\!\frak X}}$ is uniquely determined by $(\lambda^{^{\!\frak X}})_P\,;$ indeed, we have
$$(\lambda^{^{\!\frak X}})_P\.\tau'^{_{\!\frak X}} _{_{P,Q}} (1) = \tau'^{_{\!\frak X}} _{_{P,Q}} (1)\.(\lambda^{^{\!\frak X}})_Q
\eqno £2.9.3.$$

\medskip
£2.10. Once again, the composition of a {\it natural $\F^{^{\frak X}}\!\-$isomorphism\/} between 
{\it $\F^{^{\frak X}}\!\-$locality functors\/} with an {\it $\F^{^{\frak X}}\!\-$locality functor\/} or with another such a {\it natural $\F^{^{\frak X}}\!\-$isomorphism\/}
is  a {\it natural $\F^{^{\frak X}}\!\-$isomorphism\/} between {\it $\F^{^{\frak X}}\!\-$locality 
functors\/}. Note that  if $\L^{^{\frak X}}$ and $\L'^{^{\frak X}}$ are {\it coherent\/} then 
any {\it $\F^{^{\frak X}}\!\-$locality functor\/} 
$\frak l^{^{\frak X}}\colon \L^{^{\frak X}}\!\to \L'^{^{\frak X}}$ determines a {\it natural map\/}
$$\nu_{\frak l^{^{\frak X}}} : \Ker (\pi^{_{\frak X}})\too \Ker (\pi'^{_{\frak X}})
\eqno £2.10.1\phantom{.}$$
which is clearly {\it compatible\/} with the restrictions to 
$\Ker (\kappa^{_{\frak X}})$ of $\tau^{_{\frak X}}$ and $\tau'^{_{\frak X}}\,;$ in this case, it is quite clear that any subfunctor $\frak k^{^{\frak X}}$ of $\Ker (\pi^{_{\frak X}})$ determines a {\it quotient $\F\-$locality\/} $\L^{^\frak X}/\frak k^{^{\frak X}}$ defined, for any pair of subgroups $Q$ and $R$ in $\frak X\,,$ by the quotient set
$$(\L^{^\frak X}/\frak k^{^{\frak X}})(Q,R) = \L^{^\frak X} (Q,R)/\frak k^{^{\frak X}}(R)
\eqno £2.10.2,$$
 and by the corresponding induced composition.

 \medskip
£2.11. With the notation in~£2.5.1, we are interested in the  
$\overline{N_\F (Q)}^\frak h\!\-$locality $\overline{N_{\F,Q} (Q)}^\frak h$ where the morphisms are the pairs formed by an $\overline{N_\F (Q)}^\frak h\!\-$mor-phism and by an automorphism of $Q\,,$ both  determined by the {\it same\/} 
$\F\-$mor-phism [11,~18.3], and where the composition and the structural
 functors are the obvious ones. Similarly, if $L$ is a finite group acting on $Q\,,$ we are interested in the 
 $\F_L\-$locality $\F_{L,Q}$ where the morphisms  are the pairs formed by an $\F_L\-$morphism and by an automorphism of $Q\,,$ both  determined by the {\it same\/} element of $L\,.$ We are ready to describe
the {\it $\F\-$localizer\/} of $Q$ [11,~Theorem~18.6].
\eject

\bigskip
\noindent
{\bf Theorem~£2.12.}{\footnote{\dag}{\cds We thank John Rognes who pointed out that the proof of [11,~Lemma~18.8] only works whenever the normal subgroup {\cdt Q} of {\cdt M} is Abelian. To complete the proof of [11,~Theorem~18.6], we replace the application of this lemma in page~342,
by quoting  [5, Proposition~4.9]; indeed, [5,~condition~4.9.1]  follows from [11, conditions 18.6.2 and 18.6.3],
and [5, condition~4.9.2]  follows from [11, condition 18.8.1].}}
 {\it For any subgroup $Q$ of $P$ fully normalized in $\F$ there is a triple
 formed by a finite  group $L_\F (Q)$ and by two group homomorphisms
$$\tau_{_Q} : N_P (Q)\too L_\F (Q)\quad and\quad  \pi_{_Q} : L_\F (Q)\too \F (Q)
\eqno £2.12.1\phantom{.}$$
such that $\pi_{_Q}\circ \tau_{_Q}$ is induced by the $N_P (Q)\-$conjugation, that we have the exact sequence
$$1\too H_{C_\F (Q)}\too C_P (Q)\buildrel \tau_{_Q}\over\too  L_\F (Q) 
\buildrel \pi_{_Q}\over\too \F (Q)\too 1
\eqno £2.12.2\phantom{.}$$
and that $\pi_{_Q}$ and $\tau_{_Q}$ induce an equivalence of categories
$$\overline{N_{\F,Q} (Q)}^\frak h\cong \F_{L_\F (Q),Q}
\eqno £2.12.3.$$
\smallskip
\noindent
Moreover, for another such a triple $L'\,,$ $\tau'_{_Q}$ and $\pi'_{_Q}\,,$ there is
 a group isomorphism $\lambda\,\colon L_\F (Q)\cong L'\,,$ unique up to $\frak c^\frak
h_\F (Q)\-$conjugation, fulfilling $\lambda\circ\tau_{_Q} = \tau'_{_Q}$ and
$\pi'_{_Q}\circ\lambda = \pi_{_Q}\,.$\/}

\medskip
£2.13. For any subgroup $Q$ of $P$ fully normalized in $\F\,,$ we call 
{\it $\F\-$localizer\/} of $Q$  any finite group $L$ endowed with two group homomorphisms as in~£2.12.1 fulfilling the conditions~£2.12.2 and~£2.12.3. Note that, if $\L^{^{\frak X}}$ is an $\F^{^{\frak X}}\-$locality  then, for any 
$Q\in \frak X\,,$ the structural functors $\tau^{_{\frak X}}$ and $\pi^{_{\frak X}}$ determine two group homomorphisms (cf.~£2.7.2)
$$\tau_{_Q}^{_{\frak X}} : N_P (Q)\too \L^{^{\frak X}}\! (Q)\qq  
\pi_{_Q}^{_{\frak X}} : \L^{^{\frak X}}\! (Q)\too \F (Q)
\eqno £2.13.1\phantom{.}$$
and $\pi_{_Q}^{_{\frak X}}$ is surjective; in particular, if $Q$ is fully normalized in 
$\F$ then, since $\F_P (Q)$ is a Sylow $p\-$subgroup of $\F(Q)$ [11,~Proposition~2.11],  $\tau_{_Q}^{_{\frak X}}\big(N_P (Q)\big)$ is a Sylow 
$p\-$subgroup of $\L^{^{\frak X}}\!(Q)$ if and only if it contains a Sylow 
$p\-$subgroup of~${\rm Ker}(\pi_{_Q}^{_{\frak X}})\,.$ 
Hence, if $\L^{^{\frak X}}$ is {\it divisible\/} and for any $Q\in \frak X$ fully normalized in~$\F$ the group $\L^{^{\frak X}} (Q)$ endowed with 
$\tau_{_Q}^{_{\frak X}}$ and $\pi_{_Q}^{_{\frak X}}$ is an $\F\-$localizer 
of~$Q\,,$
 it is easily checked from [11,~Proposition~17.10] that $\L^{^{\frak X}}$
is {\it $p\-$coherent\/} and therefore that it is a 
{\it perfect $\F^{^{\frak X}}\!\-$locality\/} (cf.~£2.8). Actually, the
converse statement is true and it is easily checked from [11,~Proposition~18.4].

\medskip
£2.14.  We also need the {\it $\F\-$localizing functor\/} $\loc_\F$ --- defined in [11,~18.12.1 and~Proposition~18.19]
--- from the {\it proper category of chains\/} 
$\ch^*(\F)$ of $\F\,.$ Explicitly, for any {\it $\F\-$chain\/}
$\frak q\,\colon \Delta_n\to \F$  [11,~A2.8], we denote by $\F (\frak q)$ the subgroup~of elements in $\F\big(\frak q(n)\big)$ which can be lifted to an automorphism of~$\frak q\,;$ then, if $\frak q$ is {\it fully normalized\/} in $\F$ [11,~2.18], we consider the {\it $\F\-$localizer\/} $L_\F \big(\frak q (n)\big)$ 
of $\frak q (n)$ and denote by $L_\F (\frak q)$ and by $N_P (\frak q)$ the respective converse images of $\F (\frak q)$ in $L_\F \big(\frak q(n)\big)$ and in 
$N_P\big(\frak q (n)\big)\,,$ endowed with the suitable group homomorphisms
$$\tau_{\frak q} : N_P (\frak q)\too L_\F (\frak q)\qq  
\pi_{\frak q} : L_\F (\frak q)\too \F (\frak q)
\eqno £2.14.1;$$
actually, in the construction of $\loc_\F$ we only can consider the Abelian part of the extension $\pi_{\frak q}\,,$ namely the quotient
$$\bar\pi_{\frak q} : \bar L_\F (\frak q) = L_\F (\frak q)/[{\rm Ker}(\pi_{\frak q}),{\rm Ker}(\pi_{\frak q})]\too \F (\frak q)
\eqno £2.14.2.$$

\medskip
£2.15. Moreover, we denote by $\frak L\frak o\frak c$ the category
where the objects are the pairs $(L,Q)$ formed by a finite group $L$ and a normal
$p\-$subgroup $Q$ of~$L\,,$ and where the morphisms from  $(L,Q)$ to $(L',Q')$ are the group homomorphisms $f\,\colon L\to L'$ fulfilling $f(Q)\i Q'$ [11,~18.12]; this category has an obvious {\it inner structure\/} [11,~1.3] mapping any object 
$(L,Q)$ on the subgroup of the group of automorphisms of $L$  determined by the 
$Q\-$conjugation, and we~denote by $\widetilde{\Loc}$ the corresponding 
{\it exterior quotient\/} and by $\lv\,\colon \widetilde{\Loc}\to \Gr$ 
the functor sending $(L,Q)$ to $L/Q\,.$ Then, it follows 
from [11,~Proposition~18.19] that there is a suitable functor 
$$\loc_\F : \ch^*(\F)\too \widetilde\Loc
\eqno £2.15.1$$
mapping any {\it $\F\-$chain\/} $\frak q\,\colon \Delta_n\to \F$  {\it fully normalized\/} in $\F$ [11,~2.18]
on the pair~$\big(\bar L_\F (\frak q),{\rm Ker}(\bar\pi_\frak q)\big)\,.$

\medskip
£2.16. But, for a  {\it $p\-$coherent $\F^{^\frak X}\!\-$locality\/} $\L^{^\frak X}$ it follows from  [11,~Proposition~A2.10] that we have a functor
$$\aut_{\L^{^\frak X}} : \ch^*(\L^{^\frak X})\too \Gr
\eqno £2.16.1$$
mapping  any {\it $\L^{^\frak X}\-$chain\/} $\hat\frak q\,\colon 
\Delta_n\to \L^{^\frak X}$ on the group $\L^{^\frak X}\!(\hat\frak q)$ of all the automorphisms of $\hat\frak q\,;$ moreover,
$\aut_{\L^{^\frak X}}$ induces an obvious functor 
$$\loc_{\L^{^\frak X}}: \ch^*(\F^{^\frak X})\too \widetilde\Loc
\eqno £2.16.2$$
mapping  any {\it $\F^{^\frak X}\-$chain\/} 
$\frak q\,\colon \Delta_n\to \F^{^\frak X}$ on the pair 
$\big(\L^{^\frak X}\!(\hat\frak q), {\rm Ker}(\pi^{_\frak X}_{\hat\frak q})\big)$ for a  
{\it $\L^{^\frak X}\-$chain\/} $\hat\frak q\,\colon \Delta_n\to \L^{^\frak X}$ lifting $\frak q\,;$ here, we are interested in the following $\frak X\-$relative version of [11,~Proposition~18.21].

\bigskip
\noindent
{\bf Proposition~£2.17.} {\it If $\L^{^\frak X}\!$ is a  $p\-$coherent 
$\F^{^\frak X}\!\-$locality such  that ${\rm Ker}(\pi^{_\frak X}_{_Q})$ is Abelian for any $Q\in \frak X\,,$ then there is a unique natural map
$$\lambda_{\L^{^\frak X}} : \loc_{\F^{^\frak X}}\too \loc_{\L^{^\frak X}}$$
such that $\lv * \lambda_{\L^{^\frak X}} = {\rm id}_{\aut_{\F^{^\frak X}}}$ and that, for any  $\F^{^\frak X}\-$chain $\frak q\,\colon \Delta_n\to \F^{^\frak X}$  fully normalized in $\F^{^\frak X}\,,$ we have $(\lambda_{\L^{^\frak X}})_\frak q \circ \bar\tau^{_\frak X}_\frak q = \tau^{_{\L^{^\frak X}}}_\frak q\,.$\/}
\eject

\bigskip
 \bigskip
\noindent
{\bf £3. The natural $\F\-$basic $P\times P\-$sets}
\medskip
£3.1. Recall that a {\it basic $P\times P\-$set\/} [11,~21,4] is a finite nonempty $P\times P\-$set $\Omega$ such that $\{1\}\times P$ acts {\it freely\/} on $\Omega\,,$ that we have
$$\Omega^\circ \cong \Omega\qq \vert\Omega\vert/\vert P\vert \not\equiv 0 \bmod{p}
\eqno £3.1.1\phantom{.}$$
where we  denote by $\Omega^\circ$ the $P\times P\-$set obtained by exchanging both factors,
and that, for any subgroup $Q$ of $P$ and any injective  group homomorphism $\varphi\,\colon Q\to P$ such that $\Omega$ contains a $P\times P\-$subset isomorphic to $(P\times P)/\Delta_\varphi (Q)\,,$ we have a 
$Q\times P\-$set isomorphism
$${\rm Res}_{\varphi\times {\rm id}_P} (\Omega)\cong {\rm Res}_{\iota_Q^P\times {\rm id}_P} (\Omega)
\eqno £3.1.2\phantom{.}$$
where,  for any pair of group homomorphisms $\varphi$ and $\varphi'$ from $Q$ to $P\,,$ we set 
$$\Delta_{\varphi,\varphi'} (Q)  =\{(\varphi (u),\varphi' (u))\}_{u\in Q}\qq 
\Delta_\varphi (Q) 
= \Delta_{\varphi,\iota_Q^P} (Q)
\eqno £3.1.3,$$ 
and denote by $\iota_Q^P$ the corresponding inclusion map.

\medskip
£3.2. Then, for any pair of subgroups $Q$ and $R$ of $P\,,$ denoting by $\F^\Omega (Q,R)$  the set of 
injective group homomorphisms $\varphi\,\colon R\to P$ such that 
$$\varphi (R)\i Q\qq {\rm Res}_{\,\varphi\times {\rm id}_P}(\Omega)\cong  
{\rm Res}_{\,\iota_R^P\times {\rm id}_P}(\Omega)
\eqno £3.2.1,$$
it follows from [11,~Proposition~21.9] that $\F^\Omega$ is a {\it Frobenius $P\-$category\/}.
Moreover, if $\F$ is a  {\it Frobenius $P\-$category\/}, let us say that
$\Omega$ is a {\it $\F\-$basic $P\times P\-$set\/} whenever $\F^\Omega = \F\,;$ then,
it follows from [11,~Proposition~21.12] that any {\it Frobenius $P\-$category\/} $\F$ admits 
an {\it $\F\-$basic $P\times P\-$set.\/}

\medskip
£3.3. From now on, we fix  a {\it Frobenius $P\-$category\/} $\F$ and a nonempty set~$\frak X$ of 
subgroups of $P$ as in~£2.6 above;  more generally, we say that a $P\times P\-$set $\Omega^{^\frak X}$ is
{\it $\F^{^\frak X}\!\-$basic\/} if  it fulfills condition~£3.1.1 and the statement [11,~21.7]
\smallskip
\noindent
£3.3.1\quad {\it The stabilizer of any element of $\,\Omega^{^\frak X}\!$ coincides with
 $\Delta_{\psi,\psi'} (R)$ for some $R\in \frak X$ and suitable
$\psi,\psi'\in \F (P,R)\,,$ and we have 
$$\big\vert(\Omega^{^\frak X})^{\Delta_{\varphi,\varphi'} (Q)}\big\vert =
 \big\vert (\Omega^{^\frak X})^{\Delta (Q)}\big\vert$$
 for any~$Q\in \frak X$ and any $\varphi,\varphi'\in \F (P,Q)\,.$\/}
\smallskip
\noindent
 Recall that, according to [11,~Proposition~21.12], for any  $\F^{^\frak X}\!\-$basic 
 $P\times P\-$set $\Omega^{^\frak X}$ there is an $\F\-$basic $P\times P\-$set 
$\Omega$ containing $\Omega^{^\frak X}$  and fulfilling
$$\Omega^{\Delta_\varphi (Q)} = (\Omega^{^\frak X})^{\Delta_\varphi (Q)}
\eqno £3.3.2\phantom{.}$$
for any $Q\in \frak X$ and any $\varphi\in \F (P,Q)\,.$ In order to describe the $\F\-$basic 
 $P\times P\-$set~$\Omega$ announced in~£1.8 above, we need the notation of [11,~Chap.~6],
 which we actually recall in section~£5 below (cf.~£5.3.1).
 \eject

\bigskip
\noindent
{\bf Proposition~£3.4.} {\it Assume that any element of $\,\frak X$ is $\F\-$selfcentralizing. Then, the $P\times P\-$set
$$\Omega^{^\frak X} = \bigsqcup_{Q}\, \bigsqcup_{\tilde\varphi} 
(P\times P)/\Delta_{\varphi} (Q) 
\eqno £3.4.1,$$
where $Q$ runs over a set of representatives for the set of  $P\-$conjugacy classes in~$\frak X$ and
 $\tilde\varphi$ runs over  a set of representatives for the set of $\tilde\F_P(Q)\-$orbits 
in~$\tilde\F (P,Q)_{\tilde\iota_{Q}^P}\,,$  is an $\F^{^\frak X}\-$basic $P\times P\-$set which fulfills 
$\vert(\Omega^{^\frak X})^{\Delta (Q)}\vert = \vert Z(Q)\vert$ for any $Q\in \frak X\,.$\/}

\medskip
\noindent
{\bf Proof:} Since we clearly have 
$$(\Omega^{^\frak X})^\circ\cong \Omega^{^\frak X}\qq
\vert \Omega^{^\frak X}\!/P\vert \equiv \vert\tilde\F(P)\vert \bmod{p}
\eqno £3.4.2,$$
 it suffices to check that, for any $R\in \frak X$ and any $\psi\in \F (P,R)\,,$ we have 
$$\vert(\Omega^{^\frak X})^{\Delta_\psi (R)}\vert = \vert Z(R)\vert
\eqno £3.4.3;$$ 
but, for any subgroup $Q$ of $P$ and any $\varphi\in \F (P,R)\,,$ $\Delta_\psi (R)$  fixes  the class of $(u,v)\in P\times P$ in $(P\times P)/\Delta_\varphi (Q) $ if and only if it is contained in $\Delta_\varphi (Q)^{(u,v)}$ or, equivalently, we have 
 $$vRv^{-1}\i Q\qq \varphi (vwv^{-1}) = u\psi (w)u^{-1} \hbox{ for any $w\in R$}
 \eqno £3.4.4,$$
 which amounts to saying that the following $\tilde\F\-$diagram is commutative
$$\matrix{\hskip-25ptP&&&&P\hskip-15pt\cr
\hskip-25pt\Vert&\hskip-10pt\nwarrow\hskip-5pt{\tilde\varphi\atop}&\phantom{\Big\uparrow}&\hskip-5pt{\tilde\iota_Q^P\atop}\hskip-5pt\nearrow&\Vert\hskip-15pt\cr 
&&Q&&\cr
\hskip-25pt P&&&&P\hskip-15pt\cr
&&\hskip-25pt{\scriptstyle \tilde\kappa_{_{Q,R}}(v)}\hskip-1pt\big\uparrow\cr
\phantom{\Big\uparrow}&\hskip-25pt{\atop \tilde \psi}\hskip-5pt\nwarrow
&\hskip-27pt&\hskip-15pt\nearrow\hskip-5pt{\atop \tilde \iota_R^P}\hskip-30pt\cr
&&R&&\cr}
\eqno £3.4.5\phantom{.}$$
where $\kappa_{_{Q,R}}(v)\,\colon R\to Q$ is the group homomorphism determined by  the conjugation by~$v\,.$

\smallskip
Since $\tilde\varphi$ belongs to $\tilde\F (P,Q)_{\tilde\iota_Q^P}$ (see [11,~6.4.1] or £5.1.1 below), it follows from [11, Proposition~6.7]   that
the pair $(\tilde\psi,\tilde\iota_R^P)$ determines  the isomorphism class of the  $(\tilde\F^\circ )_R\-$object
(cf.~£2.1)
$$\tilde\kappa_{_{Q,R}}(v) :  R\too Q
\eqno £3.4.6;$$
 that is to say, if $(u',v')\in P\times P$ is another element  such that 
$\Delta_\varphi (Q)^{(u',v')}$ contains $\Delta_\psi (R)\,,$ we have $v' = sv$ for some $s\in Q$ and therefore we get
$$\psi (w) = u'^{-1}\varphi (svwv^{-1}s^{-1})u' = 
\varphi (vwv^{-1})^{\varphi(s)^{-1}u'}
\eqno £3.4.7\phantom{.}$$
 for any $w\in R\,;$  at this point, it follows from [11,~Proposition~4.6] that,  for a suitable  $z\in Z(R)\,,$ we have $\varphi(s)^{-1}u' = uz\,,$ which proves our claim.
 \eject

\medskip
£3.5. If any element of $\,\frak X$ is $\F\-$selfcentralizing, we call 
$\Omega^{^\frak X}$ the {\it natural $\F^{^\frak X}\-$basic $P\times P\-$set\/}.
Recall that we say that an $\F\-$basic $P\times P\-$set $\Omega$ is {\it thick\/}
 if the multiplicity of the indecomposable $P\times P\-$set 
 $(P\times P)/\Delta_\varphi (Q)$  is at least two for any subgroup $Q$ of $P$ and any $\varphi\in \F (P,Q)$
 [11,~21.4]. Let us call {\it natural\/} any $\F\-$basic 
 $P\times P\-$set $\Omega$ which  fulfills
$$\vert \Omega^{\Delta_\varphi (Q)}\vert = \vert Z (Q)\vert
\eqno £3.5.1\phantom{.}$$
 for any $\F\-$selfcentralizing subgroup $Q$ of $P$ and any $\varphi\in \F (P,Q)\,,$ and is thick outside of the set  of $\F\-$selfcentralizing subgroups of $P\,,$  namely  the multiplicity of  $(P\times P)/\Delta_\psi (R)$  is at least two if $R$ is not $\F\-$selfcentralzing; the existence of {\it natural $\F\-$basic $P\times P\-$sets\/} follows from Proposition~£3.4 together with [11,~Proposition~21.12].

\medskip
£3.6. Let $\Omega$ be an $\F\-$basic $P\times P\-$set and  $Q$ a subgroup of 
$P\,;$ it follows from our definition in~£3.2 that any $Q\times P\-$orbit in 
${\rm Res}_{Q\times P}(\Omega)$ is isomorphic to the quotient set 
$(Q\times P)/\Delta_\eta (T)$ (cf.~£3.1.3) for some subgroup $T$ of
$P$ and~some $\eta\in \F (Q, T)\,;$ note that  the isomorphism class of this 
$Q\times P\-$set $(Q\times P)/\Delta_\eta (T)$ only depends on the conjugacy class of~$T$ in $P$ and on the
class $\tilde\eta$ of $\eta$ in~$\tilde\F (Q,T)\,;$ moreover, it is quite clear that 
$\bar N_{Q\times P}\big(\Delta_\eta (T)\big)$  acts regularly on 
$\big((Q\times P)/\Delta_\eta (T)\big)^{\Delta_\eta (T)}$ and that we have a group isomorphism
$${\rm Aut}_{Q\times P}\big((Q\times P)/\Delta_\eta (T)\big)\cong \bar N_{Q\times P}\big(\Delta_\eta (T)\big)
\eqno £3.6.1\phantom{.}$$

\bigskip
\noindent
{\bf Proposition~£3.7.} {\it Let $\Omega$ be a natural $\F\-$basic 
$P\times P\-$set, $Q$ and $T$ a pair of   $\F\-$selfcentralizing subgroups of $P$  and $\eta$ an element of~$\F (Q,T).$ Then, the multiplicity of 
$(Q\times P)/\Delta_\eta (T)$ in ${\rm Res}_{Q\times P}(\Omega)$ is at most one, and it is one if and only if $\tilde\eta$ belongs to $\tilde\F(Q,T)_{\tilde\iota_T^P}\,.$ Moreover, in this case we have 
$${\rm Aut}_{Q\times P}\big((Q\times P)/\Delta_\eta (T)\big)\cong Z(T)
\eqno £3.7.1\phantom{.}$$
and the multiplicity of $(Q\times P)/\Delta_\eta (T)$ in any $\F\-$basic $P\times P\-$set $\,\Omega'$ is at least one. \/}

\medskip
\noindent
{\bf Proof:} According to our definition (cf.~£3.5.1), we have
$$\vert \Omega^{\Delta_\eta (T)}\vert = \vert Z (T)\vert
\eqno £3.7.2;$$
hence, if  the
multiplicity of $(Q\times P)/\Delta_\eta (T)$ in ${\rm Res}_{Q\times P}(\Omega)$ is not zero,
then it is one and we have  (cf.~£3.6)
$$\vert \bar N_{Q\times P}\big(\Delta_\eta (T)\big)\vert\le \vert Z (T)\vert
\eqno £3.7.3\phantom{.}$$
which forces isomorphism~£3.7.1. In this case, since 
$N_{Q\times P}\big(\Delta_\eta (T)\big)$ {\it covers\/} the intersection 
$ \F_Q\big(\eta (T)\big)\cap {}^{\eta}\F_{\! P}  (T)$ where 
${}^{\eta}\F_{\! P}  (T)$ is the  image of $\F_{\! P}  (T)$ in  ${\rm Aut}\big(\eta (T)\big)$ {\it via\/} $\eta$ (cf.~£2.3), it follows 
from~[11,~6.5] that $\tilde\eta$ belongs to~$\tilde\F(P,T)_{\tilde\iota_T^Q}\,.$
\eject
\smallskip
 Moreover, for  any  $\F\-$basic $P\times P\-$set $\,\Omega'\,,$ denoting by  
 $\frak T$ the set of subgroups $R$ of $Q$ such that $\F(R,T)\not= \emptyset$ and $\vert R\vert \not= \vert T\vert\,,$ and by $\Omega'^{\frak T}$ 
the subset of $\omega'\in \Omega'$ such that the stabilizer 
$(Q\times P)_{\omega'}$ coincides with $\Delta_\psi (R)$ for some $R\in \frak T\,,$
it is quite clear that $\tilde\eta\in \tilde\F(Q,T)_{\tilde\iota_T^P}$ forces 
$(\Omega'^{\frak T})^{\Delta_\eta (T)} = \emptyset\,;$ since 
$\Omega'^{\Delta (T)}$ is clearly not empty, this proves the last statement.

\bigskip 
 \bigskip
\noindent
{\bf £4. Construction of $\F\-$localities from $\F\-$basic $P\times P\-$sets}

\medskip
£4.1. Let $\Omega$ be an $\F\-$basic $P\times P\-$set  and denote by $G$ the group of $\{1\}\times P\-$set automorphisms  
of~${\rm Res}_{\{1\}\times P}(\Omega)\,;$ it is clear that we have an injective~map from $P\times \{1\}$ into $G\,;$  we identify this image with the $p\-$group $P$ so that, from now on, $P$ is contained in $G$ and acts freely on $\Omega\,.$ Recall that, for any pair of subgroups $Q$ and $R$ of $P\,,$ we have (cf.~£3.2)
$$T_G (R,Q)/C_G(R)\cong \F (Q,R)
\eqno £4.1.1\phantom{.}$$
where $T_G (R,Q)$ is the {\it $G\-$transporter\/} from $R$ to $Q\,.$

\medskip
£4.2.   Let $Q$ be a subgroup of $P\,;$ clearly, the centralizer $C_G (Q)$ coincides with the group of 
$Q\times P\-$set automorphisms of ${\rm Res}_{Q\times P}(\Omega)$ and therefore, denoting by~$\frak O_{\Omega,Q}$  
the set of isomorphism classes  of~$Q\times P\-$orbits of $\Omega\,,$ by~$k_{\tilde O}$ the number of 
$Q\times P\-$orbits  of isomorphism class $\tilde O\in \frak O_{\Omega,Q}$ in 
$\Omega\,,$ by~$\frak S_{k_{\tilde O}}$ the corresponding 
{\it $k_{\tilde O}\-$symmetric\/} group and by ${\rm Aut}(O)$ the group of 
$Q\times P\-$set automorphisms of $O\in \tilde O\,,$  
 it is easily checked that we have a {\it canonical\/} $\widetilde\Gr\-$isomorphism [11,~22.5.1]
$$\tilde\omega_{_Q} : C_G (Q)\cong \prod_{\tilde O\in \frak O_{\Omega,Q}} {\rm Aut}(O)\wr  \frak S_{k_{\tilde O}}
\eqno £4.2.1.$$
More precisely, as in  [11,~Proposition~22.11], for any subgroup $R$ of $Q$ we have  
a commutative $\widetilde\Gr\-$diagram
 $$\matrix{C_G (Q)&\too & C_G (R)\cr
\big\uparrow&\phantom{\Big\uparrow}&\big\uparrow\cr
{\displaystyle\prod_{\tilde O\in \frak O_{\Omega,Q}}}\!\frak S_{k_{\tilde O}}&\too
&{\displaystyle\prod_{\tilde M\in  \frak O_{\Omega,R}}}\! \frak S_{k_{\tilde M}}\cr}
\eqno £4.2.2.$$

\medskip
£4.3. As in [11,~Proposition~22.7], let us denote by $\frak S_\Omega (Q)$ the minimal normal subgroup of~$C_G (Q)$ containing  $(\omega_{_Q})^{-1}
\big(\prod_{\tilde O\in \frak O_{\Omega,Q}}\frak S_{k_{\tilde O}}\big)$ for a 
representative~$\omega_{_Q}$ of $\tilde\omega_{_Q}\,;$ then, denoting by 
$\frak O^1_{\Omega,Q}$ the subset of isomorphism classes 
$\tilde O\in \frak O_{\Omega,Q}$ with multiplicity one in $\Omega$ and 
by~$\ab \big({\rm Aut}(O)\big)$ the maximal Abelian quotient
of  ${\rm Aut}(O)\,,$ it follows from [11,~Lemma~22.8] that
$$C_G(Q)/\frak S_\Omega (Q)\cong 
\prod_{\tilde O\in \frak O^1_{\Omega,Q}} {\rm Aut}(O)\times  
\prod_{\tilde O\in \frak O_{\Omega,Q} - \frak O^1_{\Omega,Q}} 
\ab \big({\rm Aut}(O)\big)
\eqno £4.3.1;$$
\eject
\noindent
let us denote by $\frak S_\Omega^{1} (Q)$ the converse image in $C_G(Q)$ of  the {\it commutator\/} subgroup of this quotient, so that we have
$$C_G(Q)/\frak S_\Omega^1 (Q)\cong \prod_{\tilde O\in 
\frak O_{\Omega,Q}} \ab \big({\rm Aut}(O)\big)
\eqno £4.3.2.$$

\medskip
£4.4. Although in [11,~Chap.~22] we assume that $\Omega$ is {\it thick\/} (cf.~£3.5), it is easily checked that the elementary arguments in [11,~Proposition~22.11] still prove that, for any subgroup $R$ of $Q\,,$ we have
$$\frak S_\Omega (Q)\i \frak S_\Omega (R)
\eqno £4.4.1\phantom{.}$$
and therefore we still get
$$\frak S_\Omega^1 (Q)\i \frak S_\Omega^1 (R)
\eqno £4.4.2;$$
as a matter of fact, these arguments do not depend on the conditions~£3.1.1.
Thus, as in [11,~22.13], we get a {\it contravariant\/} functor 
$$\tilde\frak c^{^\Omega} : \tilde\F \too \Ab
\eqno £4.4.3\phantom{.}$$
mapping any subgroup $Q$ of $P$ on the Abelian group
$$\tilde\frak c^{^\Omega}(Q) = C_G(Q)/\frak S_\Omega^1 (Q)\cong 
\prod_{\tilde O\in \frak O_{\Omega,Q}} \ab \big({\rm Aut}(O)\big)
\eqno £4.4.4\phantom{.}$$
and any $\tilde\F\-$morphism $\tilde\varphi\,\colon R\to Q$ on the group homomorphism
$$\tilde\frak c^{^\Omega}(\tilde\varphi) :  C_G(Q)/\frak S_\Omega^1 (Q)\too  C_G(R)/\frak S_\Omega^1 (R)
\eqno £4.4.5\phantom{.}$$
induced by  conjugation in $G$ by any element $x\in T_G (R,Q)$ lifting $\tilde\varphi$ (cf.~£4.1.1).

\medskip
£4.5. More precisely, for any $\tilde\F\-$morphism $\tilde\varphi\,\colon R\to Q\,,$ any  $\tilde O\in \frak O_{\Omega,Q}$ and any $\tilde M\in \frak O_{\Omega,R}\,,$ we consider the (possibly empty) set of the injective $R\times P\-$set homomorphisms
$$f : M\too {\rm Res}_{\varphi\times {\rm id}_P}(O)
\eqno £4.5.1,$$
for $M\in \tilde M\,,$ $O\in \tilde O$ and $\varphi\in \tilde\varphi\,;$ it is clear that ${\rm Aut}(M)\times {\rm Aut}(O)$ acts on this set by left- and right-hand composition and let us denote by $\I_{\tilde M}^{\tilde O}(\tilde\varphi)$ a set of representatives for the set of  ${\rm Aut}(M)\times {\rm Aut}(O)\-$orbits. Then, 
if $f$ is such an  injective $R\times P\-$set homomorphism, denoting by
 ${\rm Aut}(O)_f$ the stabilizer of $f(M)$ in ${\rm Aut}(O)\,,$ we get an obvious group homomorphism
 $$\delta_f : {\rm Aut}(O)_f\too {\rm Aut}(M)
 \eqno £4.5.2\phantom{.}$$
 and we denote by $\varepsilon_f \,\colon {\rm Aut}(O)_f\to {\rm Aut}(O)$ the corresponding inclusion group homomorphism; we are interested in the maximal Abelian quotients of these groups;  explicitly, we denote by
 $$\matrix{\ab\big({\rm Aut}(O)\big)&&\ab\big({\rm Aut}(M)\big)\cr
{\atop \ab^\circ (\varepsilon_f)}\hskip-5pt\searrow\hskip-30pt
&&\hskip-30pt\nearrow\hskip-5pt{\atop\ab(\delta_f) }\cr
&\ab\big({\rm Aut}(O)_f\big) \cr}
\eqno £4.5.3\phantom{.}$$
\eject
\noindent
the group homomorphisms respectively determined by the {\it transfert\/} induced by $\varepsilon_f\,,$ and by
$\delta_f\,.$ With all this notation, from [11,~Proposition~22.17] we get the following description of
 $\tilde\frak c^{^{\Omega}} (\tilde\varphi)$

 \bigskip
\noindent
{\bf Proposition~£4.6.} {\it For any $\tilde\F\-$morphism 
$\tilde\varphi\,\colon R\to Q\,,$ we have
$$\tilde\frak c^{^{\Omega}} (\tilde\varphi) = 
\sum_{\tilde O\in \frak O_{\Omega,Q} }\,\sum_{\tilde M\in \frak O_{\Omega,R} }
\,\sum_{f\in \I_{\tilde M}^{\tilde O}(\tilde\varphi)}\ab(\delta_f)\circ \ab^\circ (\varepsilon_f)
\eqno £4.6.1.$$\/}

\medskip
£4.7. Now, {\it the correspondence sending any pair of subgroups $Q$ and $R$ of~$P$ to the quotient set
$$\L^{^{\Omega}}(Q,R) = T_G (R,Q)/\frak S_\Omega^1 (R)
\eqno £4.7.1,$$
endowed with the canonical maps
$$\tau^{_{\Omega}}_{_{Q,R}} : \T_{\!P} (Q,R)\to \L^{^{\Omega}}(Q,R) \quad
and\quad\pi^{_{\Omega}}_{_{Q,R}} : \L^{^{\Omega}}(Q,R)\to \F (Q,R)
\eqno £4.7.2,$$
defines a $p\-$coherent $\F\-$locality $(\tau^{_{\Omega}},\L^{^{\Omega}},\pi^{_{\Omega}})\,.$\/}
Indeed,  from inclusion~£4.4.2 it is not difficult to check that, for
any triple of subgroups $Q\,,$ $R$ and $T$ of $P\,,$ the product in $G$
induces a map
$$\L^{^{\Omega}}(Q,R)\times \L^{^{\Omega}}(R,T)\too \L^{^{\Omega}}(Q,T)
\eqno £4.7.3\,;$$
then, it is quite clear that these maps determine a {\it composition\/} in the correspondence
$\L^{^{\!\Omega}}$ defined above and that the canonical maps in~£4.7.2 define structural functors
$$\tau^{_{\Omega}} : \T_{\!P}\too \L^{^{\Omega}} \qq \pi^{_{\Omega }} : \L^{^{\Omega}}\too \F 
\eqno £4.7.4;$$
moreover, the {\it divisibility\/} and the {\it coherence\/}  of $\L^{^{\Omega}}$ 
(cf.~£2.8) are easy consequences of the fact that $G$ is a group, whereas the 
{\it $p\-$coherence\/} follows from isomorphisms~£3.6.1 and~£4.3.2.

\medskip
£4.8.  Let us denote by $P\-\Set$ the category of finite $P\-$sets endowed with the {\it disjoint union\/} and with the {\it inner\/} direct product mapping any pair of finite $P\-$sets $X$ and $Y$ on the {\it $P\-$set\/} --- still noted $X \times Y$ --- obtained from the restriction of the $P\times P\-$set $X\times Y$ through the 
{\it diagonal\/} map $\Delta \,\colon P\to P\times P\,;$ as in [13,~2.5], let us consider the functor 
$$\frak f_\Omega :  P\-\Set\too P\-\Set
\eqno £4.8.1$$
 mapping any $P\-$set $X$ on the $P\-$set noted $\Omega \times_P X $
 and, as in [13,~3.5],  we consider the positive integer 
$$\ell (\frak f_\Omega) = {\vert\frak f_\Omega (P)\vert\over\vert P\vert} 
= \vert\Omega/P\vert
\eqno £4.8.2;$$ 
recall that  if $\alpha\,\colon P'\to P$ is a group homomorphism then the {\it restriction\/} defines a functor [13,~2.3]
$$\res_\alpha : P\-\Set\too  P'\-\Set
\eqno £4.8.3.$$
\eject
\noindent
As a matter of fact, for any pair of subgroups $Q$ and $R$ of~$P\,,$ we have a canonical bijection [13,~2.5.2]
$$T_G^\varphi (R,Q)\cong  \Nat^* (\res_{\iota_R^P}\circ\frak f_\Omega,\res_{\iota_Q^P\circ\varphi}\circ\frak f_\Omega)
\eqno £4.8.4\phantom{.}$$
where $T_G^\varphi (R,Q)$ denotes the converse image of $\varphi\in \F (Q,R)$  in $T_G (R,Q)$ {\it via\/} the bijection~£4.1.1 and,
for any pair of functors $\frak f$ and $\frak g$ from  $R\-\Set$ to $R\-\Set\,,$
$\Nat^*(\frak f,\frak g)$ denotes the set  of {\it natural isomorphisms\/} from $\frak f$ to $\frak g\,;$ for short, we set $\Nat^*(\frak f) = \Nat^*(\frak f,\frak f)\,.$
 In particular, we have 
$$C_G (R) \cong \Nat^*(\res_{\iota_R^P}\circ\frak f_\Omega)
\eqno £4.8.5\phantom{.}$$
and the image of $\frak S_\Omega (R)$ is easy to describe.

\smallskip
£4.9. Recall that,  for any $P\-$set $X\,,$ the correspondence sending any $P\-$set~$Y$ to the $P\-$set $X\times Y$ and any $P\-$set map
$g\,\colon Y\to Y'$ to the $P\-$set map
$${\rm id}_X\times g : X\times Y\too X\times Y'
\eqno £4.9.1\phantom{.}$$
defines a functor  preserving disjoint unions [13,~2.9]
$$\frak m_X : P\-\Set \too P\-\Set
\eqno £4.9.2;$$
recall that $(\frak m_X)^\circ = \frak m_X$ [13,~5.3 and~6.1] and that 
 $\ell (\frak m_X) = \vert X\vert$ [13,~5.3]; we say that $X$ is {\it $\F\-$stable\/}
 if we have
$${\res}_\varphi(X) ={\res}_{\iota_Q^P}(X)
\eqno £4.9.3\phantom{.}$$
 for any subgroup $Q$ of $P$ and any $\varphi\in \F (P,Q)\,.$

\bigskip
\noindent
{\bf Proposition~£4.10.} {\it For any $\F\-$stable $P\-$set $X$ such that $p$
does not divide~$\vert X\vert\,,$ there are an $\F\-$basic $P\times P\-$set 
$\Omega'$ containing $\Omega$ and fulfilling $\frak f_{\Omega'} = 
\frak f_\Omega\circ \frak m_X\,,$ an $\F\-$locality functor
$\L^{^\Omega}\to \L^{^{\Omega'}}$ and an injective natural map
$\tilde\frak c^{^\Omega}\to \tilde\frak c^{^{\Omega'}}\,.$\/}

\medskip
\noindent
{\bf Proof:} From the very definition of a {\it $\F\-$basic $P\times P\-$set\/},
it is easily checked that $\frak f_\Omega$ is {\it $\F\-$stable\/} [13,~6.3]
and that $p$ does not divide $\ell (\frak f_\Omega)\,;$ hence, it follows from
[13,~Theorem~6.6] that $\frak f_\Omega$ and $\frak m_X$ centralizes each
other and, with the notation there, that we have $\frak m_X = \frak k + \frak h$ for some $\frak k\in \K_\H$ and $\frak h\in \H^\F\,;$ then, it follows from [13,~2.8 and~Proposition~6.5]
that we get
$$\frak f_\Omega\circ \frak m_X 
= \frak m_X\circ \frak f_\Omega  = \frak h\circ \frak f_\Omega\qq 
(\frak f_\Omega\circ \frak m_X)^\circ = \frak f_\Omega\circ \frak m_X 
\eqno £4.10.1,$$
so that $\frak f_\Omega\circ \frak m_X$ is also {\it $\F\-$stable\/}; moreover,
since [13,~5.3]
$$\ell (\frak f_\Omega\circ \frak m_X) = \ell (\frak f_\Omega)\ell (\frak m_X)
\eqno £4.10.2,$$
$p$ does not divide $\ell (\frak f_\Omega\circ \frak m_X)\,.$
Consequently, if $\Omega'$ is the unique  $P\times P\-$set $\Omega'$ 
 fulfilling $\frak f_{\Omega'} = \frak f_\Omega\circ \frak m_X$ [13,~2.8.1 and~3.1.1] then it is easily checked that $\Omega'$ fulfills conditions~£3.1.1
 and~£3.1.2 above, so that it is a {\it basic $P\times P\-$set\/}; moreover, according
 to [13,~Proposition~7.6], it is a {\it $\F\-$basic $P\times P\-$set\/}.
 \eject

 \smallskip
Denote by $G'$ the group of $\{1\}\times P\-$set automorphisms  of~${\rm Res}_{\{1\}\times P}(\Omega')\,;$ now,
for any pair of subgroups $Q$ and $R$ of $P\,,$ and any $\varphi\in \F (Q,R)\,,$  the composition with~$\frak m_X$
determines a map
$$\Nat^* (\res_{\iota_R^P}\circ\frak f_\Omega,\res_{\iota_Q^P\circ\varphi}\circ\frak f_\Omega)\too  \Nat^* (\res_{\iota_R^P}\circ\frak f_{\Omega'},\res_{\iota_Q^P\circ\varphi}\circ\frak f_{\Omega'})
\eqno £4.10.3\phantom{.}$$
sending any {\it natural isomorphism\/}
$$\nu : \res_{\iota_R^P}\circ\frak f_\Omega\cong \res_{\iota_Q^P\circ\varphi}\circ\frak f_\Omega
\eqno £4.10.4\phantom{.}$$
to the {\it natural isomorphism\/}
$$\nu * \frak m_X : \res_{\iota_R^P}\circ\frak f_\Omega\circ \frak m_X\cong \res_{\iota_Q^P\circ\varphi}\circ\frak f_\Omega\circ \frak m_X
\eqno £4.10.5;$$
consequently, from the canonical bijections~£4.8.4 we get a canonical map
$$\frak t_{_{R,Q}}^\varphi : T_G^\varphi (R,Q)\too T_{G'}^\varphi (R,Q)
\eqno £4.10.6.$$

\smallskip
Moreover, for any third subgroup $T$ of $P\,,$  any  $\F\-$morphism
$\psi\,\colon T\to R$ and any {\it natural isomorphism\/}
$$\eta : \res_{\iota_T^P}\circ\frak f_\Omega\cong \res_{\iota_R^P\circ\psi}\circ\frak f_\Omega
\eqno £4.10.7\phantom{.}$$
we have the {\it natural isomorphism\/}
$$(\res_\psi * \nu)\circ\eta : \res_{\iota_T^P}\circ\frak f_\Omega\cong \res_{\iota_R^P\circ\psi}\circ\frak f_\Omega\circ \res_{\iota_Q^P\circ\varphi\circ\psi}\circ\frak f_\Omega
\eqno £4.10.8\phantom{.}$$
and thus we get
$$\big((\res_\psi * \nu)\circ\eta\big) * \frak m_X = 
\big(\res_\psi * (\nu * \frak m_X)\big)\circ (\eta * \frak m_X)
\eqno £4.10.9;$$
that is to say, we obtain $\,\frak t_{_{T,Q}}^{\varphi\circ\psi} = 
\frak t_{_{R,Q}}^\varphi \circ \,\frak t_{_{T,R}}^\psi \,.$ In particular, we get 
a group homomorphism $\frak t_{_R}^{{\rm id}_R} \,\colon C_G (R)\to C_{G'} (R)$
and it is easily checked that $\frak S_{\Omega'} (R)$ contains  $\frak t_{_R}^{{\rm id}_R}
\big(\frak S_\Omega (R)\big)\,,$  which forces 
$$\frak t_{_R}^{{\rm id}_R}\big(\frak S_\Omega^1 (R)\big)\i \frak S_{\Omega'}^1 (R)
\eqno £4.10.10.$$

\smallskip
 Consequently, since we have
$$T_G (R,Q) = \bigsqcup_{\varphi\in \F (Q,R)} T_G^\varphi (R,Q)
\eqno £4.10.11,$$
the family of maps $\frak t_{_{R,Q}}^\varphi$ induces a canonical functor and 
a {\it natural map\/}
$$\frak l : \L^{^\Omega}\too \L^{^{\Omega'}}\qq 
\lambda :\tilde\frak c^{^\Omega}\too \tilde\frak c^{^{\Omega'}}
\eqno £4.10.12;$$
it is easily checked that $\frak l$ is compatible with the structural functors
$\tau^{_\Omega}$ and~$\tau^{_{\Omega'}}\,,$ and with the structural functors 
$\pi^{_\Omega}$ and $\pi^{_{\Omega'}}\,,$ so that it is an {\it $\F\-$locality functor\/} (cf.~2.9). Moreover, since $p$ does not divide $\vert X\vert\,,$ $X$ contains the trivial $P\-$set with {\it multiplicity\/} $k$ prime to $p$ and therefore $\Omega'$ contains $\Omega$ [13,~3.1.1]; more precisely, we claim that the group homomorphism
$$\lambda_R :\tilde\frak c^{^\Omega}(R)\too \tilde\frak c^{^{\Omega'}}(R)
\eqno £4.10.13\phantom{.}$$
\eject
\noindent
is injective. Indeed, otherwise choose a nontrivial element 
$a = (a_{\tilde O})_{\tilde O\in \frak O_{\Omega,R}}$~in
$${\rm Ker}(\lambda_R)\i \prod_{\tilde O\in \frak O_{\Omega,R}} 
\ab \big({\rm Aut}(O)\big)
\eqno £4.10.14\phantom{.}$$
and an element $\tilde O_\circ\in \frak O_{\Omega,Q}$ with $\vert O_\circ\vert$ minimal in such a way that 
$a_{\tilde O_\circ}\not= 0\,;$ then, it is easily checked that the component of $\lambda_R (a)$ in the factor 
$\tilde O_\circ$ of
$$\tilde\frak c^{^{\Omega'}}(R)\cong \prod_{\tilde O\in \frak O_{\Omega',R}} 
\ab \big({\rm Aut}(O)\big)
\eqno £4.10.15$$
coincides with $k\.a_{\tilde O_\circ}\not= 0\,,$ a contradiction. We are done.

\bigskip
\noindent
{\bf Corollary~£4.11.\/}{\footnote{\dag}{\cds The proof of the uniqueness of the
basic {\cdy F}-locality in [12,~Proposition~22.12] is not correct.}} {\it If $\Omega$ is thick and $\Omega'$ is an $\F\-$basic
$P\times P\-$set then we have  an $\F\-$locality functor
$\L^{^{\Omega'}}\!\to \L^{^{\Omega}}$ and an injective natural map
$\tilde\frak c^{^{\Omega'}}\!\to \tilde\frak c^{^{\Omega}}\,.$
In particular, if $\Omega'$ is thick then we have an $\F\-$locality isomorphism
$\L^{^{\Omega'}}\cong \L^{^{\Omega}}\,.$ \/}

\medskip
\noindent
{\bf Proof:} Denote by $X$ and by $X'$ the respective images by $\frak f_\Omega$
and by $\frak f_{\Omega'}$ of the trivial $P\-$set; since $X$ and  $X'$ are
$\F\-$stable [13,~Proposition~6.5] and $p$ does not divide
$$\vert\Omega/P\vert = \vert X\vert\qq \vert\Omega'/P\vert = \vert X'\vert
\eqno £4.11.1,$$ 
it follows from Proposition~£4.10 above that  there are two
 $\F\-$basic $P\times P\-$sets $\Omega''\j \Omega$ and $\Omega'''\j \Omega'$   fulfilling 
 $$\frak f_{\Omega''} = \frak f_\Omega\circ \frak m_{X'}\qq 
 \frak f_{\Omega'''} = \frak f_{\Omega'}\circ \frak m_{X}
 \eqno £4.11.2,$$
 two {\it $\F\-$locality\/} functors
$$\L^{^\Omega}\too \L^{^{\Omega''}}\qq 
\L^{^{\Omega'}}\too \L^{^{\Omega'''}}
\eqno £4.11.3\phantom{.}$$
 and two injective natural maps
$$\tilde\frak c^{^\Omega}\too \tilde\frak c^{^{\Omega''}}\qq
\tilde\frak c^{^{\Omega'}}\too \tilde\frak c^{^{\Omega'''}}
\eqno £4.11.4.$$

\smallskip
Now, we claim that $\frak f_{\Omega''} = \frak f_{\Omega'''}\,;$ indeed, 
it follows from [13,~Theorem~6.6] that, since $\frak f_{\Omega''}$ and 
$\frak f_{\Omega'''}$ are {\it $\F\-$stable\/}, it suffices to prove that the
images by $\frak f_{\Omega''}$ and by $\frak f_{\Omega'''}$ of the trivial $P\-$set
coincide with each other; but, we clearly have
$$\eqalign{(\frak f_\Omega\circ \frak m_{X'})(1) &= (\frak m_{X'}\circ\frak f_\Omega)(1) = X'\times X\cr
 &= X\times  X' = (\frak m_X\circ\frak f_{\Omega'})(1) 
 = (\frak f_{\Omega'}\circ \frak m_X)(1)\cr}
\eqno £4.11.5.$$
Moreover, if $\Omega$ is {\it thick\/} then $\Omega''$ is also {\it thick\/}
and isomorphism~£4.4.4 implies that the left-hand {\it natural map\/} in~£4.11.4
is a {\it natural isomorphism\/} and therefore that the left-hand {\it $\F\-$locality functor\/} 
is an isomorphism, so that we have the {\it $\F\-$locality functors\/}  and the injective {\it natural maps\/}
$$\L^{^{\Omega'}}\!\too \L^{^{\Omega'''}}\! = \L^{^{\Omega''}}\!
\cong \L^{^\Omega}\qq 
\tilde\frak c^{^{\Omega'}}\!\too \tilde\frak c^{^{\Omega'''}} \!
= \tilde\frak c^{^{\Omega''}}\!\cong \tilde\frak c^{^\Omega}
\eqno £4.11.6.$$
We are done.
\eject

\medskip
£4.12. When $\Omega$ is {\it thick\/} we call $\L^{^\Omega}$ the 
{\it basic $\F\-$locality\/} and set $\L^{^\Omega} = \L^{^{\rm b}}\,,$ 
$\tau^{_\Omega}  = \tau^{_{\rm b}}\,,$ $\pi^{_\Omega}  = \pi^{_{\rm b}}$  
and $\tilde\frak c^{^\Omega} = \tilde\frak c^{^{\rm b}}\,;$ according to Corollary~£4.11, $\L^{^{\rm b}}$ and $\tilde\frak c^{^{\rm b}}$ do not depend 
on the choice of the thick $\F\-$basic $P\times P\-$set.  Let us denote by 
$\L^{^{\rm b,sc}}$  the {\it full\/} subcategory of $\L^{^{\rm b}}$ over the set 
of $\F\-$selfcentralizing subgroups of~$P$ and by  $\tilde\frak c^{^{\rm b,sc}}$ the restriction to $\tilde\F^{^{\rm sc}}$ of the natural map $\tilde\frak c^{^{\rm b}}\,;$ 
in general, denote by  $\L^{^{\Omega,\rm sc}}$  the {\it full\/} subcategory of 
$\L^{^{\Omega}}$ over the set of $\F\-$selfcentralizing subgroups of~$P$
 and by  $\tilde\frak c^{^{\Omega,\rm sc}}$ the restriction to 
 $\tilde\F^{^{\rm sc}}$ of the natural map $\tilde\frak c^{^{\Omega}}\,,$
so that we have an {\it $\F\-$locality  functor\/} and an injective {\it natural map\/}
$$\L^{^{\Omega,\rm sc}}\too \L^{^{\rm b,sc}}\qq
\tilde\frak c^{^{\Omega,\rm sc}}\too \tilde\frak c^{^{\rm b,sc}}
\eqno £4.12.1.$$

\medskip
£4.13. Moreover, for any subgroup $Q$ of $P\,,$ consider the set of isomorphism classes $\frak O^{^{\rm nsc}}_Q$ of indecomposable $Q\times P\-$sets  
$(Q\times P)/\Delta_\theta (U)$  where $U$ is a {\it $\F\-$nonselfcentralizing\/} subgroup of~$P$ and $\theta$ belongs to $\F (Q,U)\,;$  according to our arguments in [11,~23.2], it is  easily checked that the correspondence mapping any  
$\F\-$selfcentralizing subgroup $Q$ of $P$ on 
$$\tilde\frak c^{^{\rm nsc}}\! (Q) = \prod_{\tilde O\in \frak O^{^{\rm nsc}}_Q} \ab \big({\rm Aut}(O)\big)
\eqno £4.13.1\phantom{.}$$
defines a  {\it contravariant\/} subfunctor 
$\tilde\frak c^{^{\rm nsc}} \colon\tilde\F^{^{\rm sc}}\to  \Ab$ of
$\tilde\frak c^{^{\rm b,sc}}.$
Then, let us consider the {\it quotient $\F^{^{\rm sc}}\-$localities\/} (cf.~£2.10.2)
$$\bar\L^{^{\rm b,sc}}\! = \L^{^{\rm b,sc}}\!/\tilde\frak c^{^{\rm nsc}}\qq
\bar\L^{^{\Omega,\rm sc}}\! = \L^{^{\Omega,\rm sc}}\!/
(\tilde\frak c^{^{\Omega,\rm sc}}\cap \tilde\frak c^{^{\rm nsc}})
\eqno £4.13.2;$$
we still have  a {\it faithful $\F^{^{\rm sc}}\-$locality functor\/} and an injective 
{\it natural map\/}
 $$ \bar\L^{^{\Omega,\rm sc}}\too \bar\L^{^{\rm b,sc}}\qq 
\tilde\frak c^{^{\Omega,\rm sc}}/ (\tilde\frak c^{^{\Omega,\rm sc}}\cap 
\tilde\frak c^{^{\rm nsc}}) \too \tilde\frak c^{^{\rm b,sc}}/
\tilde\frak c^{^{\rm nsc}}
 \eqno £4.13.3.$$

\medskip
£4.14. If~$\Omega$  is {\it natural\/} (cf.~£3.5), we claim that 
$\bar\L^{^{\Omega,\rm sc}}$ does not depend on the choice of $\Omega$ and  set $\bar\L^{^{\Omega,\rm sc}} = \bar\L^{^{\rm n,sc}},$ 
$\bar\tau^{_{\Omega,\rm sc}}  = \bar\tau^{_{\rm n,sc}}\,,$  
$\bar\pi^{_{\Omega,\rm sc}}  = \bar\pi^{_{\rm n,sc}}$  and  
$\tilde\frak c^{^{\Omega,\rm sc}} = \tilde\frak c^{^{\rm n,sc}}\,;$  explicitly,
  it follows from Proposition~£3.7 and from [11,~Lemma~22.8] that for any 
 $\F\-$selfcentralizing subgroup $Q$ of $P$ we have
 $$\tilde\frak c^{^{\Omega,\rm sc}}\!(Q) \cong 
 Z(Q\cap^{\tilde\F^{^{\rm sc}}}\!\! \! P)\times \tilde\frak c^{^{\rm nsc}}\! (Q)
 \eqno £4.14.1\phantom{.}$$
where,  with the notation in~£5.7.2 below, we set
$$Z(Q\cap^{\tilde\F^{^{\rm sc}}}\!\! \! P) = 
\big(\prod_T\,\prod_{\tilde\gamma\in \tilde\F(P,T)_{\tilde\iota_T^Q}}\, 
Z(T)\big)^{Q\times \F_Q(T)}
\eqno £4.14.2,$$
$T$ running over the set of $\F\-$selfcentralizing subgroups of $Q\,.$
\eject

\medskip
£4.15.  In order to prove our claim, for any pair of $\F\-$selfcentralizing subgroups 
$Q$ and $T$ of $P\,,$  and any $\F\-$morphism $\eta\,\colon T\to Q\,,$ consider the $Q\times P\-$set $O = (Q\times P)/\Delta_\eta (T)\,,$  set 
$$\eqalign{Z_O = Z(T)&\qq A_O = {}^{\eta^*}\!\tilde\F_Q\big(\eta (T)\big)\cap \tilde\F_P(Q)\cr
\bar Z_O = Z_O/[A_O,Z_O]&\qq\bar A_O = A_O/[A_O,A_O]\cr}
\eqno £4.15.1,$$
where ${}^{\eta^*}\!\tilde\F_Q\big(\eta (T)\big)$ is the corresponding image of 
$\tilde\F_Q\big(\eta (T)\big)$ in $\widetilde{\rm Aut}(T)\,,$
and recall that we have [11,~23.8.5]
$$\ab\big({\rm Aut}(O)\big) \cong \bar Z_O\times \bar A_O
\eqno £4.15.2.$$
Then, denoting by $\frak O^{^{\rm sc}}_Q$ the set of isomorphism classes  of such indecomposable $Q\times P\-$sets, in [11,~Proposition~23.10] we exhibit two {\it contravariant\/} functors 
$$\bar\frak z^{^{\rm sc}} : \tilde\F^{^{\rm sc}}\too\Ab\qq 
\bar\frak a^{^{\rm sc}} : \tilde\F^{^{\rm sc}}\too\Ab
\eqno £4.15.3\phantom{.}$$
respectively mapping any $\F\-$selfcentralizing subgroup $Q$ of $P$ on
$$\bar\frak z^{^{\rm sc}} (Q) = \prod_{\tilde O\in \frak O^{^{\rm sc}}_Q}\bar Z_O\qq \bar\frak a^{^{\rm sc}} (Q) = \prod_{\tilde O\in \frak O^{^{\rm sc}}_Q}
\bar A_O
\eqno £4.15.4\phantom{.}$$
 and fulfilling $\tilde\frak c^{^{\rm b,sc}}/\tilde\frak c^{^{\rm nsc}}\cong \bar\frak z^{^{\rm sc}}\times \bar\frak a^{^{\rm sc}}\,.$

 \medskip
 £4.16. In particular, with an obvious notation, it is quite clear that we get (cf.~£2.9)
  $$\bar\L^{^{\rm b,sc}}\!\cong (\bar\L^{^{\rm b,sc}}\!/
 \,\bar\frak z^{^{\rm sc}}\, ) 
 \times^{\F^{^{\rm sc}}}  (\bar\L^{^{\rm b,sc}}\!/\,\bar\frak a^{^{\rm sc}}\, ) 
 \eqno £4.16.1\phantom{.}$$
 and then we set  $\L^{^{\rm c,sc}}\! = \bar\L^{^{\rm b,sc}}\!/\,\bar\frak a^{^{\rm sc}}$ [11,~23.10.2];
  if~$\Omega$  is {\it natural\/}, it is easily checked that the {\it  $\F^{^{\rm sc}}\-$locality functor\/} 
 $\bar\L^{^{\Omega,\rm sc}}\to \bar\L^{^{\rm b,sc}}$ still determines a 
 {\it faithful $\F^{^{\rm sc}}\-$locality functor\/} 
$\bar\L^{^{\Omega,\rm sc}}\!\to \L^{^{\rm c,sc}}\,.$ Moreover, in [11,~Corollaries~23.24 and~23.28] we successively exhibit a sequence of  
$\F^{^{\rm sc}}\-$sublocalities
$$\L^{^{\rm r,sc}}\i \L^{^{\rm d,sc}}\i \L^{^{\rm c, sc}}
\eqno £4.16.2\phantom{.}$$ 
and it is not difficult to check that we may assume that $\L^{^{\rm r,sc}}$ 
contains the image of $\bar\L^{^{\Omega,\rm sc}};$ finally, 
 it follows from equality~[11,~£23.28.1] and 
 from isomorphism~£4.14.1 above that the functors $\Ker (\pi^{_{\rm r,sc}})$ and 
 $\tilde\frak c^{^{\Omega,\rm sc}}/\tilde\frak c^{^{\rm nsc}}$ coincide with each other.  Consequently, we get an {\it $\F^{^{\rm sc}}\-$locality isomorphism\/}
$$\bar\L^{^{\Omega,\rm sc}}\!\cong \L^{^{\rm r,sc}}
\eqno £4.16.3,$$
which proves our claim.

\bigskip
 \bigskip
\noindent
{\bf £5. Construction of an $\F^{^{\frak X}}\-$basic $P\times P\-$set from a perfect 
$\F^{^{\frak X}}\-$locality}

\medskip
£5.1.  Let $\F$ be a Frobenius $P\-$category and $\frak X$  a nonempty set  of 
$\F\-$self-centralizing subgroups of $P$ which contains any subgroup of $P$ admitting an $\F\-$morphism from some subgroup in~$\frak X\,;$ let us denote 
by $\F^{^\frak X}$ and by $(\bar\L^{^{\rm n,sc}})^{^\frak X}\!$ the respective
{\it full\/} subcategories of $\F$ and of $\bar\L^{^{\rm n,sc}}\!$ (cf.~£4.14)  over $\frak X\,;$   
it is quite clear that the correspondence mapping any $Q\in \frak X$ on (cf.~£4.14)
$$\tilde\frak k^{^{\frak X}}\! (Q) = 
\big(\prod_T\,\prod_{\tilde\gamma\in \tilde\F(P,T)_{\tilde\iota_T^Q}}\, Z(T)\big)^{Q\times \F_Q(T)}
\eqno £5.1.1,$$
where $T$ runs over the set of $\F\-$selfcentralizing subgroups of $Q$ which does not belong to $\frak X\,,$ determines a subfunctor  $\tilde\frak k^{^{\frak X}} \colon \tilde\F^{^\frak X}\to \Ab$ of 
$\big(\Ker (\bar\pi^{^{\rm n,sc}})\big)^{\!^\frak X}\,;$
let us consider the quotient $\F^{^{\frak X}} \-$locality  (cf.~£2.9)
$$\bar\L^{^{\rm n,\frak X}}\! = (\bar\L^{^{\rm n,sc}})^{^\frak X}\!/\tilde\frak k^{^{\frak X}}
\eqno £5.1.2,$$
 where we denote the structural functors by
$$\bar\tau^{^{\rm n,\frak X}} : \T^{^\frak X}_{\!P} \too\bar\L^{^{\rm n,\frak X}}\qq
\bar\pi^{^{\rm n,\frak X}} : \bar\L^{^{\rm n,\frak X}}\too \F^{^\frak X}
\eqno £5.1.3.$$
In this section we show  that  any possible {\it perfect $\F^{^\frak X}\-$locality\/} 
$(\tau^{_{\frak X}},\P^{^{\frak X}},\pi^{_{\frak X}})$  is contained in $(\bar\tau^{^{\rm n,\frak X}}\!,\bar\L^{^{\rm n,\frak X}}\!,\bar\pi^{^{\rm n,\frak X}})$ --- called  the {\it natural $\F^{^\frak X}\-$locality\/} (cf.~£1.8). Actually, when $\frak X$ is the set of all the $\F\-$selfcentralizing  subgroups of $P\,,$ this is already proved
in [11,~Corollary~24.18]; but, although the arguments there still hold for $\frak X\,,$ we state below the main
steps of the proof for $\frak X$ since we explicitly need the present context for our inductive argument.

\medskip
£5.2. First of all, let us recall the {\it distributive direct product\/} in $\ad (\tilde\F^{^\frak X})$ 
[11,~Proposition~6.14]. The {\it additive cover\/} 
$\ad(\tilde\F^{^{\frak X}})$ of $\tilde\F^{^{\frak X}}$ is the category
where the objects   are the finite sequences  $\bigoplus_{i\in I} Q_i$  of subgroups $Q_i$ in~$\frak X\,,$ and where a morphism from  another object  $R = \bigoplus_{j\in J} R_j$ to 
$Q = \bigoplus_{i\in I} Q_i$ is a pair  $(\tilde\alpha,f)$ formed by a map $f\,\colon J\to I$ and by a family $\tilde\alpha
=\{\tilde\alpha_j\}_{j\in J}$ of\/ $\tilde\F^{^{\frak X}}\-$morphisms
$\tilde\alpha_j\,\colon R_j\to Q_{f (j)}\,.$   The composition of $(\tilde\alpha,f)$ with
another $\ad (\tilde\F^{^{\frak X}})\-$morphism
$$(\tilde\beta,g) : T = \bigoplus_{\ell\in L} T_\ell \too R= \bigoplus_{j\in J} R_j
\eqno £5.2.1,$$
 formed by a map $g\,\colon L\to J$ and by a family  $\tilde\beta = \{\tilde\beta_\ell\}_{\ell\in L}\,,$ is the pair
 formed by $f\circ g$ and by the family $\{\tilde\alpha_{g  (\ell)}\circ\tilde\beta_\ell\}_{\ell\in L}$ of composed morphisms 
$$\tilde\alpha_{g  (\ell)}\circ \tilde\beta_\ell : T_\ell\too R_{g(\ell)} \too
 Q_{(f\circ g)(\ell)}
\eqno £5.2.2.$$
\eject

\medskip
£5.3. It follows from [11,~Corollary~4.9] that, for any triple of
subgroups $Q\,,$ $R$ and~$T$ in $\frak X\,,$ any  $\tilde\F\-$morphism\/ $\tilde\alpha\,\colon Q\to R$
 induces an injective map from $\tilde\F (T,R)$ to~$\tilde\F (T,Q)$
 and then  we set
$$\tilde\F (T,Q)_{\tilde\alpha} = \tilde\F (T,Q) - \bigcup_{\tilde\theta'}
\tilde\F (T,Q')\circ\tilde\theta'
\eqno £5.3.1,$$
where $\tilde\theta'$ runs over the set of {\it $\tilde\F\-$nonisomorphisms\/} $\tilde\theta'\,\colon Q\to  Q'$ 
from $Q$ such that
$\tilde\alpha'\circ\tilde\theta' = \tilde\alpha$ for~some $\tilde\alpha'\in \tilde\F (R,Q')\,;$ 
in this case, according to  [11,~Corollary~4.9], $\tilde\alpha'$  is uniquely determined, and we simply say that $\tilde\theta'$ {\it divides\/}
$\tilde\alpha$ setting $\tilde\alpha' = \tilde\alpha/\tilde\theta'\,.$ Note that we have 
$\tilde\F (T,Q)_{\tilde\alpha} = \tilde\F (T,Q)$ if and only if $\tilde\alpha$ is an isomorphism.

\medskip
£5.4. Actually, an element $\tilde\beta\in \tilde\F (T,Q)$ which can be
extended to $Q'$ {\it via\/}~$\tilde\theta'\,,$ {\it a fortiori\/}  can be extended to
$N_{Q'}\big(\theta' (Q)\big)$ for a representative $\theta'\in \tilde\theta'\,;$ 
hence, it follows from condition~£2.2.3 above that $\tilde\beta$ belongs to 
$\tilde\F (T,Q)_{\tilde\alpha}$  if and only if for some representatives 
$\alpha\in \tilde\alpha$ and  $\beta\in \tilde\beta$ we have
$${}^{\alpha^*}\! \F_{\! R} \big(\alpha (Q)\big)\cap {}^{\beta^*}\!\F_T\big(\beta (Q)\big) = \F_{Q}(Q)
\eqno £5.4.1\phantom{.}$$
where $\alpha^*\,\colon \alpha (Q)\cong Q$ and $\beta^*\,\colon \beta (Q)\cong Q$ denote
the inverse of the isomorphisms respectively induced by $\alpha$ and $\beta\,;$  in particular, we get [11,~6.5.2
and~6.6.4]
\smallskip
\noindent
£5.4.2\quad $\tilde\beta\in \tilde\F(T,Q)_{\tilde\alpha}$ 
{\it is equivalent to $\tilde\alpha\in \tilde\F(R,Q)_{\tilde\beta}\,,$ and moreover $\bar N_R \big(\alpha (Q)\big)$ acts freely on\/ $\tilde\F
(T,Q)_{\tilde\alpha}\,.$\/}
\smallskip
\noindent
The next result follows from [11,~Proposition~6.7].

\bigskip
\noindent
{\bf Proposition~£5.5.}\phantom{.}  {\it For any triple of subgroups $Q\,,$
$R$ and~$T$ in $\frak X$ and any $\tilde\alpha\in \tilde\F (R,Q)\,,$ we have
$$\tilde\F (T,Q) = \bigsqcup_{ \tilde\theta'} \tilde\F (T,Q')_{\tilde\alpha/
\tilde\theta'} \circ \tilde\theta'
\eqno £5.5.1\phantom{.}$$
where $\tilde\theta'\,\colon Q\to Q'$ runs over a set of representatives 
for the isomorphism classes of $\tilde\F^{^{\frak X}}\-$morphism from $Q$ dividing 
$\tilde\alpha\,.$ In particular, $p$ does not divide~$\vert \tilde\F (P,Q)\vert\,.$\/}

\medskip
£5.6.  At this point, Proposition~£5.5 allows us to define a {\it distributive direct product\/} in 
$\ad(\tilde\F^{^{\frak X}})$ (see [11,~Chap.~6] and also [12,~Proposition~4.5]). If~$R$ and~$T$ 
are two subgroups in~$\frak X\,,$ let us consider the set~$\tilde\frak T_{R,T}^{^{\frak X}}$ of {\it strict triples\/}  
$(\tilde\alpha, Q,\tilde\beta)$ where $Q$ is a subgroup in~$\frak X\,,$
$\tilde\alpha$ and $\tilde\beta$ respectively belong to $\tilde\F (R,Q)$ and 
to $\tilde\F (T,Q)\,,$ and we have $\tilde\alpha\in \tilde\F (R,Q)_{\tilde\beta}$
or, equivalently, $\tilde\beta\in \tilde\F (T,Q)_{\tilde\alpha}\,.$
We say that two {\it strict triples\/}  $(\tilde\alpha, Q,\tilde\beta)$ and $(\tilde\alpha', Q',\tilde\beta')$ are {\it equivalent\/} 
if there is an $\tilde\F\-$isomorphism $\tilde\theta\,\colon Q\cong Q'$ fulfilling
$$\tilde\alpha'\circ\tilde\theta = \tilde\alpha\qq \tilde\beta'\circ\tilde\theta =
\tilde\beta
\eqno £5.6.1;$$
\eject
\noindent
then, $\tilde\theta$ is unique since, assuming that the triples coincide each other and choosing $\alpha\in
\tilde\alpha\,,$ $\beta\in \tilde\beta$ and~$\theta\in \tilde\theta\,,$ it is easily
checked that $\theta$ belongs to (cf.~£5.4.1)
$${}^{\alpha^*}\! \F_{\! R} \big(\alpha (Q)\big)\cap {}^{\beta^*}\!\F_T\big(\beta (Q)\big) = \F_Q(Q)
\eqno £5.6.2.$$

\medskip
£5.7. Denoting by $\check{\tilde\frak T}_{R,T}^{_{\frak X}}$ a set of representatives for the set of equivalence classes 
in  $\tilde\frak T_{R,T}^{^{\frak X}}\,,$ we call {\it $\tilde\F^{^{\frak X}}\!\-$intersection\/}  of $R$ and $T$ the 
$\ad(\tilde\F^{^{\frak X}})\-$object
$$R\cap^{\tilde\F^{^{\frak X}}}\!\! T = \bigoplus_{(\tilde\alpha, Q,\tilde\beta)\in 
\check{\tilde\frak T}_{R,T}^{_{\frak X}}} Q
\eqno £5.7.1;$$
note that, if we choose another set of representatives, then the uniqueness of the isomorphism in~£5.6.1 above yields a unique 
$\ad(\tilde\F^{^{\frak X}})\-$isomorphism between both $\ad(\tilde\F^{^{\frak X}})\-$objects; in particular, 
 we have
$$R\cap^{\tilde\F^{^{\frak X}}}\!\! T \cong 
\bigoplus_Q\,\bigoplus_{\tilde\gamma}\, Q
\eqno £5.7.2.$$
where $Q$ runs over a set of representatives for the set of $R\-$conjugacy classes of elements in $\frak X$ contained in $R$ and, for such a $Q\,,$ $\tilde\gamma$ runs over a set of representatives for the $\F_R(Q)\-$classes in 
$\tilde\F(T,Q)_{\tilde\iota_Q^R}\,.$ Finally, for any pair of $\ad(\tilde\F^{^{\frak X}})\-$objects 
$R = \bigoplus_{j\in J} R_j$ and $T = \bigoplus_{\ell\in L} T_\ell\,,$  we define
$$R\cap^{\tilde\F^{^{\frak X}}}\!\! T = \bigoplus_{(j,\ell)\in J\times L} R_j\cap^{\tilde\F^{^{\frak X}}}\!\!T_\ell
\eqno £5.7.3.$$
The argument in [11,~Proposition~6.14] still shows that the  {\it $\tilde\F^{^{\frak X}}\!\-$intersection\/} defines  a 
{\it distributive direct product\/} in $\ad(\tilde\F^{^{\frak X}})$ (see also [12,~Proposition~4.5]).

\medskip
£5.8. Analogously, the existence of a {\it perfect $\F^{^{\frak X}}\!\-$locality\/} 
$(\tau^{^{\frak X}}\!,\P^{^{\frak X}}\!,\pi^{^{\frak X}})$ actually determines a   {\it distributive direct product\/} in the 
{\it additive cover\/} $\ad(\P^{^{\frak X}})$ of $\,\P^{^{\frak X}}\,;$ as we show in [12,~Proposition~4.5],
this fact depends on Lemma~£5.9 and on Proposition~£5.11 below which admit the same proofs as the proofs of [11,~Proposition~24.2] and~[11,~Proposition~£24.4].

\bigskip
\noindent
{\bf Lemma~£5.9.} {\it Any $\P^{^{\frak X}}\!\-$morphism $x\,\colon R\to Q$
 is a monomorphism and an epimorphism.\/}

\medskip 
£5.10.  Thus, for any triple of subgroups $Q\,,$ $R$ and~$T$ in $\frak X\,,$ as in £5.3 above any $\P^{^{\frak X}}\!\-$morphism  $\,x\in\P^{^{\frak X}}\! (R,Q)$  induces an injective map from 
$\P^{^{\frak X}}\! (T,R)$ to~$\P^{^{\frak X}}\! (T,Q)$ and then, as in £5.3.1, we set
$$\P^{^{\frak X}} (T,Q)_x = \P^{^{\frak X}} (T,Q) - \bigcup_{z'} \P^{^{\frak X}} (T,Q')\. z'
\eqno £5.10.1\phantom{.}$$
where $z'$ runs over the set of {\it $\P^{^{\frak X}}\!\-$nonisomorphisms\/} 
$z'\,\colon Q\to Q'$ from $Q$ such\break
\eject
\noindent 
that $x'. z' = x$ for some $x'\in \P^{^{\frak X}}\! (R,Q')\,;$ then, $x'$ is uniquely determined by  this equality and we simply say that $z'$ {\it divides\/}~$x$ setting 
$x' = x/z'\,.$  Note that the existence of $x'$ for some $z'\in\P^{^{\frak X}}\! (Q',Q)$ is equivalent to the existence of a subgroup of~$R$ which is 
$\F\-$isomorphic to~$Q'$ and contains $\big(\pi_{_{R,Q}}(x)\big) (Q)\,;$
thus, it is quite clear that
\smallskip
\noindent
£5.10.2\quad {\it $\P^{^{\frak X}} \!(T,Q)_x$ is the converse image of $\tilde\F^{^{\frak X}}
(T,Q)_{\widetilde{\pi_{_{R,Q}}(x)}}$ in $\P^{^{\frak X}} (T,Q)\,.$\/}
Then, Proposition~£5.5 implies the following result.

\bigskip
\noindent
{\bf Proposition~£5.11.}  {\it For any triple of elements $Q\,,$ $R$ and~$T$ in $\frak X\,,$ and any 
$x\in \P^{^{\frak X}}\! (R,Q)\,,$ we have
$$\P^{^{\frak X}} \!(T,Q) = \bigsqcup_{z'} \P^{^{\frak X}} \!(T,Q')_{x/z'}\. z'
\eqno £5.11.1\phantom{.}$$
where $z'\,\colon Q\to Q'$ runs over a set of representatives 
for the isomorphism classes of $\P^{^{\frak X}}\!\-$morphism from $Q$ dividing 
$x\,.$ \/}

\medskip
£5.12. As above, if~$R$ and~$T$  are two  subgroups in  $\frak X\,,$ we consider the 
set~$\frak T^{^\frak X}_{R,T}$ of {\it strict $\P^{^{\frak X}}\!\-$triples\/}  $(x,Q,y)$ where 
$Q$ belongs to $\frak X\,,$  $x$ and $y$  respectively belong  to 
$\P^{^{\frak X}}\! (R,Q)$ and  to $\P^{^{\frak X}} \!(T,Q)\,,$ and we have
$x\in \P^{^{\frak X}} (R,Q)_y$ or, equivalently, $y\in \P^{^{\frak X}} (T,Q)_x\,.$
Note that, for any $v\in R$ and any $w\in T\,,$ the $\P^{^{\frak X}}\!\-$triple
$$v\.(x,Q,y)\.w^{-1} = \big(\tau^{_\frak X}_{_R}(v)\.x,Q,\tau^{_\frak X}_{_T}(w)\.y\big)
\eqno £5.12.1\phantom{.}$$
 still belongs to $\frak T^{^\frak X}_{R,T}$ and therefore the quotient set 
 $(R\times T)\backslash \frak T^{^\frak X}_{R,T}$
 clearly coincides with~$\tilde\frak T^{^\frak X}_{R,T}\,.$
 Similarly, we say that two {\it strict $\P^{^{\frak X}}\-$triples\/}  $(x,Q,y)$ and $(x',Q',y')$ are {\it equivalent\/} 
if there exists a $\P^{^{\frak X}}\-$isomorphism $z\, \,\colon Q\cong Q'$ fulfilling 
$$x'\. z = x\qq y'\. z = y
\eqno £5.12.2;$$
since $\P^{^{\frak X}}$ is divisible, such a $\P^{^{\frak X}}\-$isomorphism $z$ is {\it unique\/};
 in particular, in any equivalent class we may find a unique element fulfilling 
 $$Q\i R\qq x= \tau^{_\frak X}_{_{R,Q}}(1)
 \eqno £5.12.3.$$

\medskip
£5.13. Coherently, for any $Q\in \frak X$ denoting by $\s_Q^{^{\frak X}}$ the set of subgroups~of~$Q$ belonging 
 to $\frak X\,,$  we call 
{\it $\P^{^{\frak X}}\!\-$intersection\/}  of $R$ and $T$ the 
$\ad(\P^{^{\frak X}})\-$ object
$$R\cap^{\P^{^{\frak X}}}\! T = \bigoplus_{Q\in \s_R^{^{\frak X}}}\, \bigoplus_{y\in \P^{^{\frak X}}(T,Q)_{\tau^{_\frak X}_{_{R,Q}}(1)}} Q
\eqno £5.13.1\phantom{.}$$
and we clearly have canonical  $\ad (\P^{^{\frak X}})\-$morphisms
$$R\longleftarrow R \cap^{\P^{^{\frak X}}}\! T\too T
\eqno £5.13.2\phantom{.}$$
\eject
\noindent
respectively determined by $\tau^{_\frak X}_{_{R,Q}}(1)$ and $y\,.$ Note that, for any other choice of a set of  
representatives for the set of equivalence classes in $\frak T^{^\frak X}_{R,T}\,,$ we get an isomorphic object 
and a {\it unique\/} $\ad (\P^{^{\frak X}})\-$isomorphism which is compatible with the canonical
morphisms. Then, either the arguments in [11,~Proposition~£24.8] or  [12,~Proposition~4.5] prove the following.

\bigskip
\noindent
{\bf Proposition~£5.14.} {\it The category  $\ad (\P^{^{\frak X}})$ admits a distributive direct product 
mapping any pair of elements $R$ and~$T$ of $\frak X$ on their $\,\P^{^{\frak X}}\!\!\-$intersection
$R\,\cap^{\P^{^{\frak X}}}\! T $.\/}

\medskip
£5.15. Here, we are particularly interested in the $\P^{^\frak X}\!\-$intersection of $P$ with itself; 
more explicitly, denoting by $\Omega^{^\frak X}$ the set of pairs $(Q,y)$ formed by  $Q\in \frak X$ and by $y\in \P^{^\frak X}\! (P,Q)_{\tau^{_\frak X}_{_{P,Q}}(1)}\,,$  we have 
$$P \cap^{\P^{^\frak X}}\!\! P = \bigoplus_{(Q,y)\in \Omega^{^\frak X}}  Q
\eqno £5.15.1;$$
moreover, since $P\times P$ acts on the set  $\frak T^{^\frak X}_{P,P}$ (cf.~£5.12.1) preserving the equi-valence 
classes, this group acts on~$\Omega^{^\frak X}$ and it is easily checked that [11,~24.9]
\smallskip
\noindent
£5.15.2\quad {\it $(u,v)\in P\times P$ maps $(Q,y)\in \Omega^{^\frak X}$ on
$\big(Q^{u^{-1}},\tau^{_\frak X}_{_P}(v)\.y\.\tau^{_\frak X}_{_{Q,Q^{u^{-1}}}} (u^{-1}) \big)\,.$\/}
\smallskip
\noindent
In particular, $\{1\}\times P$ acts {\it freely\/} on $\Omega^{^\frak X}\,.$ On the other hand, it is clear that the map sending a {\it strict $\P^{^{\frak X}}\!\-$triple\/} $(x,Q,y)\in \frak T^{^\frak X}_{P,P}$ to $(y,Q,x)$ induces a 
$P\times P\-$set isomorphism $\Omega^{^\frak X}\cong (\Omega^{^\frak X})^\circ\,.$ The point is that from
 [11,~Proposition~24.10 and Corollary~24.11] and from Proposition~£3.4 above we get (cf.~£3.5).
 \smallskip
\noindent
£5.15.3\quad {\it   $\Omega^{^\frak X}$ is the natural $\F^{^\frak X}\-$basic $P\times P\-$set.\/}

\medskip
£5.16. Consequently, we may assume that $\Omega^{^\frak X}$ is contained in a {\it natural $\F\-$basic $P\times P\-$set\/} $\Omega$ (cf.~£3.5) and our purpose is to show that the perfect $\F^{^\frak X}\!\-$locality  $\P^{^\frak X}\!$ is contained in the {\it natural $\F^{^\frak X}\!\-$locality\/}~$\bar\L^{^{\rm n,\frak X}}\!$ (cf.~£5.1.3). 
First of all, it follows from Proposition~£5.14 that for any $Q\in \frak X$ the inclusion $Q\i P$ determines an  $\ad(\P^{^\frak X})\-$morphism 
$$\tau^{_\frak X}_{_{P,Q}}(1)\cap^{\P^{^\frak X}}\!\!\! \tau^{_\frak X}_{_{P}}(1): 
Q \cap^{\P^{^\frak X}}\!\!  P \too P \cap^{\P^{^\frak X}}\! \! P
\eqno £5.16.1;$$
actually, according to~£5.13.1 and denoting by $\Omega^{^\frak X}_Q$ the set of pairs $(T,z)$ formed by a subgroup  $T$ in $\frak X$ contained in~$Q$ and by an element $z$ of 
$\P^{^\frak X}\!(P,T)_{\tau^{_\frak X}_{_{Q,T}}(1)}\,,$ we have
$$Q \cap^{\P^{^\frak X}}\!\!  P = \bigoplus_{(T,z)\in \Omega^{^\frak X}_Q} T
\eqno £5.16.2,$$
the group $Q\times P$ acts on $\Omega^{^\frak X}_{Q}\,,$ and the $\ad(\P^{^\frak X})\-$morphism~£5.16.1  determines a $Q\times P\-$set homomorphism
$$f^{^\frak X}_Q :\Omega^{^\frak X}_{Q}\too {\rm Res}_{Q\times P}(\Omega^{^\frak X})\i 
{\rm Res}_{Q\times P}(\Omega)
\eqno £5.16.3.$$ 
From the arguments in  [11,~Proposition~24.15] we get the following result.

\bigskip
\noindent
{\bf Proposition~£5.17.} {\it For any $Q\in \frak X\,,$ the map $f^{^\frak X}_Q\,\colon \Omega^{^\frak X}_{Q} \to \Omega^{^\frak X}$  sends an element $(T,z)\in \Omega^{^\frak X}_Q$ to~$(R,y)\in \Omega^{^\frak X}$ 
if and only if we have  $T = Q\cap R$ and~$z = y\.\tau^{_\frak X}_{_{R,T}}(1)\,.$ In particular, this map is injective.\/}

 \medskip
£5.18. Thus, according to this proposition,  the image of $\,\Omega^{^\frak X}_{Q}$ in the  {\it na-tural  $\F\-$basic 
$P\times P\-$set} $\Omega$ coincides  with the union of all the  $Q\times P\-$orbits isomorphic to 
$(Q\times P)/\Delta_\eta (T)$ for some $T\in \frak X$ and some $\tilde\eta\in \tilde\F (Q,T)_{\tilde\iota_T^P}\,.$ 
On the other hand, for any $\P^{^\frak X}\!\-$isomorphism  $x\,\colon Q\cong Q'\,,$ it follows again from 
Proposition~£5.14 that we have an $\ad (\P^{^\frak X})\-$isomorphism 
$$x\cap^{\P^{^\frak X}}\! \!\tau^{_\frak X}_{_{P}}(1) : 
Q\cap^{\P^{^\frak X}}\!\! P \cong Q'\cap^{\P^{^\frak X}}\!\! P
\eqno £5.18.1\phantom{.}$$
and therefore we get a bijection between the sets of indices $\Omega^{^\frak X}_Q$ and $\Omega^{^\frak X}_{Q'}\,,$ which is compatible {\it via\/} $\pi^{_\frak X}_{_{Q',Q}}(x)$ with the respective actions of  $Q\times P$ and~$Q'\times P\,;$
 that is to say, we get a $Q\times P\-$set isomorphism
$$f^{^\frak X}_x : \Omega^{^\frak X}_{Q}\cong {\rm Res}_{\pi_{_{Q',Q}}(x)\times {\rm id}_P}(\Omega^{^\frak X}_{Q'})
\eqno £5.18.2.$$
As above, we set $G = {\rm Aut}_{\{1\}\times P}(\Omega)\,.$

\bigskip
\noindent
{\bf Proposition~£5.19.} {\it For any $\P^{^\frak X}\!\-$isomorphism $x\,\colon Q\cong Q'\,,$ the $Q\times P\-$set isomorphism 
$$f^{^\frak X}_x : \Omega^{^\frak X}_{Q} \cong {\rm Res}_{\pi^{_\frak X}_{_{Q',Q}}(x)\times {\rm id}_P}
(\Omega^{^\frak X}_{Q'})
\eqno £5.19.1\phantom{.}$$
 can be extended to an element $f_x$ of~$T_{G} (Q,Q')$ and the image of $f_x$ in $\bar\L^{^{\rm n,\frak X}}\!(Q',Q)$ is uniquely determined by $x\,.$\/}

\medskip
\noindent
{\bf Proof:} Since the $Q\times P\-$sets ${\rm Res}_{\,Q\times P} (\Omega)$
and ${\rm Res\,}_{\pi^{_\frak X}_{_{Q',Q}}(x)\times {\rm id}_P}\big({\rm Res}_{\,Q'\times P}
(\Omega)\big)$ are isomorphic (cf.~£3.1.2), and the $Q\times P\-$
and $Q'\times P\-$set homomorphisms
$$f^{^\frak X}_Q : \Omega^{^\frak X}_{Q}\too {\rm Res}_{\,Q\times P}(\Omega)\qq  f^{^\frak X}_{Q'} : \Omega^{^\frak X}_{Q'}\too {\rm Res}_{\,Q'\times P}(\Omega)
\eqno £5.19.2\phantom{.}$$
are injective (cf.~Proposition~£5.17), identifying $\Omega^{^\frak X}_{Q}$ and $\Omega^{^\frak X}_{Q'}$
with their images in $\Omega\,,$ $f^{^\frak X}_x$ can be extended to a $Q\times P\-$set isomorphism
$$f_x : {\rm Res}_{\,Q\times P}(\Omega)\cong {\rm Res\,}_{\pi^{_\frak X}_{_{Q',Q}}(x)\times 
{\rm id}_P}\big({\rm Res}_{\,Q'\times P}(\Omega)\big)
\eqno £5.19.3;$$
that is to say, we get an element $f_x$ of $T_{G}(Q,Q')$ (cf.~£3.1).
\eject

\smallskip
 Then, we claim that the image of $f_x$ in $\bar\L^{^{\rm n,\frak X}} (Q',Q)$ is independent of our choices;  indeed, for another choice $g_x\in \T_{G} (Q',Q)$ fulfilling the above conditions, the composed map $(f_x )^{-1}\circ g_x$ belongs to $C_G (Q)$ and induces the identity on~$\Omega^{^\frak X}_Q\,;$ since we know that (cf.~£4.4.4 and~£4.14.1)
$$C_G (Q)/\frak S_\Omega^1 (Q) \cong Z(Q\cap^{\tilde\F^{^{\rm sc}}}\!\!  P)\times \tilde\frak c^{^{\rm nsc}}\! (Q)
\eqno £5.19.4\phantom{.}$$
and $(f_x )^{-1}\circ g_x$ induces the identity on~$\Omega^{^\frak X}_Q\,,$ it is clear from~£5.18 that
the image in this quotient  of $(f_x )^{-1}\circ g_x$ belongs~to $\tilde\frak k^{^{\frak X}}\! (Q)$
and therefore it has a trivial image in $\bar\L^{^{\rm n,\frak X}} (Q)\,,$ so that $f_x$ and $g_x$
have the same image in $\bar\L^{^{\rm n,\frak X}} (Q',Q)\,.$
We are done.

\bigskip
\noindent
{\bf Corollary~£5.20.} {\it There is a faithful $\F^{^\frak X}\-$locality functor 
$\lambda\!^{^\frak X}\,\colon \P^{^\frak X}\to \bar\L^{^{\rm n,\frak X}}$  sending any $\P^{^\frak X}\!\-$isomorphism $x\,\colon Q\cong Q'$ to the image of $f_x$ in~$\bar\L^{^{\rm n, \frak X}}\! (Q',Q)\,.$ Moreover, any 
$\F^{^{\frak X}}\!\-$locality functor  $\mu^{_{\frak X}} \,\colon \P^{^{\frak X}}\to \bar\L^{^{\rm n,\frak X}}$ is naturally $\F^{^{\frak X}}\!\-$isomorphic 
to~$\lambda^{^{\!\frak X}}\,.$\/}

\medskip
\noindent
{\bf Proof:} Let us denote by $\lambda\!^{^\frak X} (x)$  the image of $f_x$ in 
$\bar\L^{^{\rm n, \frak X}} (Q',Q)\,;$  first of all, let  $x'\,\colon Q'\cong Q''$ be a second  $\P^{^\frak X}\!\-$isomorphism; it is clear that the automorphism 
${\rm Res}_{\pi^{_\frak X}_{_{Q',Q}}(x)\times {\rm id}_P} (f_{x'}) \circ f_{x}$ 
of ${\rm Res}_{\,Q\times P}(\Omega)$ extends ${\rm Res}_{\pi^{_\frak X}_{_{Q',Q}}(x)\times {\rm id}_P}(f^{^\frak X}_{x'})\circ f^{^\frak X}_{x}\,;$ consequently, by the proposition
above, we get
$$\lambda\!^{^\frak X}(x'\.x) = \lambda\!^{^\frak X} (x')\.  \lambda\!^{^\frak X} (x)
\eqno £5.20.1.$$

\smallskip
On the other hand, by the {\it divisibility\/} of $\P^{^\frak X}\!$, any 
$\P^{^\frak X}\!\-$morphism $z\,\colon T\to Q$ is the composition of 
$\tau^{_\frak X}_{_{Q,T'}}(1)$ with a  $\P^{^\frak X}\!\-$isomorphism  
$z_*\,\colon T\cong T_* $ where we set  $T_* = \big(\pi^{_\frak X}_{_{Q,T}}(z)\big)(T)\,;$ then, we simply define
$$\lambda\!^{^\frak X} (z) = \bar\tau_{_{Q,T_*}}^{_{\rm n, \frak X}}(1)\.\lambda\!^{^\frak X} (z_*)
\eqno £5.20.2.$$
Now, in order to prove that this correspondence defines a functor, it suffices to show that, for any $\P^{^\frak X}\!\-$isomorphism $x\,\colon Q\cong Q'$ and any subgroup $R$ of $Q\,,$ setting $R' = \big(\pi^{_\frak X}_{_{Q',Q}}(x )\big)(R)$ and denoting by $y\,\colon R\cong R'$ the $\P^{^\frak X}\!\-$isomorphism induced by $x$ (cf.~£2.8), we still have
$$\lambda\!^{^\frak X} (x)\.\bar\tau_{_{Q,R}}^{_{\rm n, \frak X}}(1) = 
\bar\tau_{_{Q',R'}}^{_{\rm n, \frak X}}(1) \.\lambda\!^{^\frak X} (y)
\eqno £5.20.3.$$

\smallskip
But, it is quite clear that the commutative $\ad (\P^{^\frak X})\-$diagram
(cf.~Proposition~£5.14)
$$\hskip-20pt\matrix{ R \cap^{\P^{^\frak X}}\!\!  P
&\buildrel \tau^{_\frak X}_{_{Q,R}}(1)\,\cap^{\P^{^\frak X}}\!\! \tau^{_\frak X}_{_{P}}(1)\over{\hbox to
75pt{\rightarrowfill}}& Q \cap^{\P^{^\frak X}}\!\! P\cr
\hskip-50pt{\scriptstyle y\,\cap^{\P^{^\frak X}}\!\!\tau^{_\frak X}_{_{P}}(1)}
\big\downarrow&\phantom{\Big\uparrow}
&\big\downarrow {\scriptstyle x\,\cap^{\P^{^\frak X}}\!\! \tau^{_\frak X}_{_{P}}(1)}\hskip-40pt\cr 
 R' \cap^{\P^{^\frak X}}\!\!  P
&\buildrel \tau^{_\frak X}_{_{Q',R'}}(1)\,\cap^{\P^{^\frak X}}\!\! 
\tau^{_\frak X}_{_{P}}(1)\over{\hbox to 80pt{\rightarrowfill}} 
& Q' \cap^{\P^{^\frak X}}\! \! P\cr}
\eqno £5.20.4\phantom{.}$$
determines a commutative diagram of $R\times P\-$sets (cf.~£5.16)
$$\matrix{\Omega^{^\frak X}_R\hskip-15pt&\too &\hskip-15pt {\rm Res}_{R\times P}^{Q\times P} 
(\Omega^{^\frak X}_Q)\cr 
{\scriptstyle  f^{^\frak X}_{y}}\big\downarrow &\phantom{\Big\uparrow} &\big\downarrow
{\scriptstyle {\rm Res}_{R\times P}^{Q\times P}( f^{^\frak X}_{x})}\hskip-30pt\cr
{\rm Res}_{\pi_y\times {\rm id}_P}(\Omega^{^\frak X}_{R'})&\too & {\rm Res}_{R\times P}^{Q\times
P}\big({\rm Res}_{\pi_x\times {\rm id}_P}(\Omega^{^\frak X}_{Q'})\big)\cr}
\eqno £5.20.5.$$
Consequently, the element $f_x$ of $T_G (Q,Q')$ extending $ f^{^\frak X}_{x}$ also extends 
$ f^{^\frak X}_{y}$ and we can choose $f_y = f_x\,.$
On the other hand, since  $\bar\tau^{_{\rm n,\frak X}}$ is  {\it faithful\/}, it is easily checked that $\lambda^{^{\!\frak X}}$ induces an injective group homomorphism
$\P^{^\frak X}\! (Q)\to \bar\L^{^{\rm n,\frak X}}\! (Q)$ for any $Q\in \frak X$ and therefore  this functor is faithful too.

\smallskip
Moreover, if $\mu^{_{\frak X}} \,\colon \P^{^{\frak X}}\!\to 
\bar\L^{^{\rm n,\frak X}}$ is another {\it $\F^{^{\frak X}}\!\-$locality functor\/}
then, for any $\P^{^{\frak X}}\!\-$morphism $x\,\colon R\to Q\,,$ we have
$\mu^{_{\frak X}}\!(x) = \lambda^{^{\!\frak X}}\!(x)\.c_x$ for some 
$c_x\in {\rm Ker}(\bar\pi^{^{\rm n,\frak X}}_{_R})\,;$ since 
$\mu^{_{\frak X}}\!\circ\tau^{_{\frak X}}\! = \bar\tau^{^{\rm n,\frak X}} = 
\lambda^{^{\!\frak X}}\!\circ\tau^{_{\frak X}}\! $ (cf.~£2.9), $c_x$ only
depends on the class $\tilde x$ of $x$ in~$\tilde\F^{^{\frak X}}\! (Q,R)\,;$
then, for another $\P^{^{\frak X}}\!\-$morphism $y\,\colon T\to R\,,$ we get
 $$\eqalign{\lambda^{^{\!\frak X}}\!(x\.y)\.c_{x\.y} &= \mu^{_{\frak X}}\!(x\.y) = 
 \mu^{_{\frak X}}\!(x)\.\mu^{_{\frak X}}\!(y) = \big(\lambda^{^{\!\frak X}}\!(x)\.c_x\big)\.\big(\lambda^{^{\!\frak X}}\!(y)\.c_y\big)\cr
&= \lambda^{^{\!\frak X}}\!(x\.y)\.\big(\Ker(\bar\pi^{^{\rm n,\frak X}})(\tilde y)\big)(c_x)\.c_y\cr}
 \eqno £5.20.6;$$
 thus, employing additive notation in the Abelian group $\big(\Ker(\bar\pi^{^{\rm n,\frak X}})\big)(T)\,,$ we still get
 $$0 = \big(\Ker(\bar\pi^{^{\rm n,\frak X}}) (\tilde y)\big)(c_x) 
 - c_{x\.y} + c_y
 \eqno £5.20.7.$$

 \smallskip
 That is to say,  setting
$$\Bbb C^n \big(\tilde\F^{^{\frak X}}\!,\Ker (\bar\pi^{^{\rm n,\frak X}})\big) = 
\prod_{\tilde\frak q\in \Fct(\Delta_n,\tilde\F^{^{\frak X}})}{\rm Ker}(\bar\pi^{^{\rm n,\frak X}}_{\tilde\frak q (0)})
\eqno £5.20.8\phantom{.}$$
 for any $n\in \Bbb N\,,$ the family $c = (c_x)_{\tilde x}$ where $\tilde x$ runs over 
 the set of $\tilde\F^{^{\frak X}}\!\-$mor-phisms is an element of 
 $\Bbb C^1 \big(\tilde\F^{^{\frak X}}\!,\Ker (\bar\pi^{^{\rm n,\frak X}})\big)$
 and equality~£5.20.7 shows that this element belongs to the kernel of the usual differential map
 $$d^{^{\frak X}}_1 : \Bbb C^1 \big(\tilde\F^{^{\frak X}}\!,\Ker (\bar\pi^{^{\rm n,\frak X}})\big)
\too \Bbb C^2 \big(\tilde\F^{^{\frak X}}\!,\Ker (\bar\pi^{^{\rm n,\frak X}})\big)
\eqno £5.20.9.$$
But, according to [12,~4.2] and~£5.7 above, the main result of [12,~\S4]
can be applied to the category $\tilde\F^{^{\frak X}}$ and  to the functor $\Ker (\bar\pi^{^{\rm n,\frak X}})$ since for any $Q\in \frak X$ we have
(cf.~£5.7.1)
$$\Ker (\bar\pi^{^{\rm n,\frak X}})(Q) =  Z(Q\cap^{\tilde\F^{^{\frak X}}}\!\!  P)
\eqno £5.20.10.$$
 Consequently, for any $n\ge 1$ we have
 $$\Bbb H^n \big(\tilde\F^{^{\frak X}}\!,\Ker (\bar\pi^{^{\rm n,\frak X}})\big) = \{0\}
 \eqno £5.20.11\phantom{.}$$
and, in particular, we have $c = d^{^{\frak X}}(z)$ for a suitable element
 $z = (z_Q)_{Q\in \frak X}$ in~$\Bbb C^0 \big(\tilde\F^{^{\frak X}},\Ker (\bar\pi^{^{\rm n,\frak X}})\big)\,;$ 
 hence,  for any $\P^{^{\frak X}}\!\-$morphism $x\,\colon R\to Q\,,$ in the additive notation we get
 $$c_x = \Big(\big(\Ker(\bar\pi^{^{\rm n,\frak X}})\big) (\tilde x)\Big)(z_Q) - z_R
 \eqno £5.20.12\phantom{.}$$
  \eject
\noindent
 and therefore we still get $\mu^{_{\frak X}}\!(x)\.z_R = 
 z_Q\.\lambda^{^{\!\frak X}}\!(x)\,,$ so that the family of $\bar\L^{^{\rm n, \frak X}}\!\-$iso-morphisms 
 $z_Q\,\colon Q\cong Q$ where $Q$ runs over $\frak X$ defines a {\it natural  $\tilde\F^{^{\frak X}}\-$isomor-phism\/}
 between $\lambda^{^{\!\frak X}}\!$ and $\mu^{^{\!\frak X}}$ (cf.~£2.9).

\bigskip
\noindent
{\bf Corollary~£5.21.} {\it Let $\P^{^{\frak X}}\!$ and $\P'^{^{\frak X}}\!$ be 
perfect  $\F^{^{\frak X}}\!\-$localities and assume that they are 
$\F^{^{\frak X}}\!\-$locality isomorphic. Then, there is an $\F^{^{\frak X}}\!\-$locality isomorphim
$\rho^{_{\frak X}} : \P^{^{\frak X}}\cong \P'^{^{\frak X}}$ such that we have
the commutative diagram
$$\matrix{\P^{^{\frak X}}&\buildrel \rho^{_{\frak X}}\over
\cong &\P'^{^{\frak X}}\cr
{\atop \lambda^{^{\!\frak X}}}\searrow&
&\hskip-10pt\swarrow{\atop \lambda'^{^\frak X}}\cr
&\bar\L^{^{\rm n,\frak X}}\cr}
\eqno £5.21.1.$$\/}
 
 \par
 \noindent
 {\bf Proof:} Considering  the set $\Omega'^{^\frak X}$ of pairs $(Q,y)$ formed by  $Q\in \frak X$ and by $y\in \P'^{^\frak X}\! 
 (P,Q)_{\tau'^{_\frak X}_{_{P,Q}}(1)}\,,$  it follows from £5.15.3 that the 
 $P\times P\-$sets $\Omega'^{^\frak X}\!$ and~$\Omega^{^\frak X}$ are mutually isomorphic; hence, up to suitable identifications, it follows from Corollary~£5.20 that  there is also  a faithful $\F^{^{\frak X}}\!\-$locality  functor 
 $\lambda'^{^\frak X}\,\colon\P'^{^\frak X}\to \bar\L^{^{\rm n,\frak X}}$ 
 as above;  thus, if $\rho^{_{\frak X}} : \P^{^{\frak X}}\cong \P'^{^{\frak X}}$ is an {\it $\F^{^{\frak X}}\!\-$locality 
 isomorphism\/},  we have a new $\F^{^{\frak X}}\!\-$locality  functor 
 $$\lambda'^{^\frak X}\!\circ \rho^{_{\frak X}}\,\colon \P^{^\frak X}\too 
 \bar\L^{^{\rm n,\frak X}}
 \eqno £5.21.2\phantom{.}$$
 and therefore, according to Corollary~£5.20 and~£2.9, it suffices~to modify
 the identification between $\Omega'^{^\frak X}\!$ and~$\Omega^{^\frak X}$
 with a suitable element of~$C_G (P)$ to get 
 $$\lambda'^{^\frak X}\!\circ \rho^{_{\frak X}} = \lambda^{^\frak X}
 \eqno £5.21.3.$$

 \bigskip
\noindent
{\bf £6. The perfect $\F^{^{\frak X}}\-$locality contained in the natural $\F^{^{\frak X}}\-$locality }

\medskip
£6.1. Let $\F$ be a Frobenius $P\-$category and $\frak X$  a nonempty set  of 
$\F\-$self-centralizing  subgroups of $P$ which contains any subgroup of $P$ admitting an $\F\-$morphism from some subgroup in~$\frak X\,;$ we keep our notation in~£5.1 above. In this section we prove the existence and the uniqueness of a 
{\it perfect $\F^{^\frak X}\-$locality\/} $(\tau^{_{\frak X}},\P^{^{\frak X}},\pi^{_{\frak X}})\,.$ 
The existence and the uniqueness of the {\it localizer\/} $(\tau_{_P},L_\F(P),\pi_{_P})$ (cf.~Theorem~£2.12) 
proves the existence and the uniqueness of the  {\it perfect $\F^{^\frak X}\-$locality\/} whenever $\frak X = \{P\}\,;$  thus,  assume that 
$\frak X\not= \{P\}\,,$  choose a minimal element $U$ in $\frak X$ {\it fully normalized\/} in $\F$ and set 
$$\frak Y = \frak X - \{\theta(U)\mid \theta\in \F(P,U)\}
\eqno £6.1.1;$$
then, arguing by induction on $\vert\frak X\vert$ we may assume that  there is a
{\it  perfect $\F^{^{\frak Y}}\-$locality\/}~$(\tau^{_{\frak Y}},\P^{^{\frak Y}},
\pi^{_{\frak Y}})$ which is unique up to $\F^{^\frak Y}\-$locality isomorphisms.
 At this point, according to Corollary~£5.20, we may assume that 
 $(\tau^{_{\frak Y}},\P^{^{\frak Y}},\pi^{_{\frak Y}})$ 
is an {\it $\F^{^{\frak Y}}\-$sublocality\/} of the 
{\it natural $\F^{^{\frak Y}}\-$locality\/} $(\bar\tau^{^{\rm n,\frak Y}}\!,
\bar\L^{^{\rm n,\frak Y}}\!,\bar\pi^{^{\rm n,\frak Y}})$ (cf.~£5.1);\break
\eject
\noindent
then, denoting by $(\bar\L^{^{\rm n,\frak X}}\!)^{^\frak Y}$ the {\it full\/}
subcategory of $\bar\L^{^{\rm n,\frak X}}\!$ over $\frak Y\,,$ we clearly
have an obvious functor $(\bar\L^{^{\rm n,\frak X}}\!)^{^\frak Y}\too 
\bar\L^{^{\rm n,\frak Y}}\!$ and we look to the {\it pull-back\/}
$$\matrix{\P^{^{\frak Y}}&\i &\bar\L^{^{\rm n,\frak Y}}\!\cr
\uparrow&\phantom{\big\uparrow}&\uparrow\cr
\M^{^{\frak Y}}\!&\i &(\bar\L^{^{\rm n,\frak X}}\!)^{^\frak Y}\cr}
\eqno £6.1.2,$$
so that we get a  {\it $p\-$coherent $\F^{^{\frak Y}}\-$locality\/} 
$(\upsilon^{_{\frak Y}}\!,\M^{^{\frak Y}}\!,\rho^{_{\frak Y}})\,;$  more explicitly, it is easily checked from~£5.1 that, for any $Q\in \frak Y\,,$ we have the exact sequence
$$1\too \prod_V\,\prod_{\tilde\theta }
\, Z(V)\too \M^{^{\frak Y}}(Q)\too \P^{^{\frak Y}}(Q)\too 1
\eqno £6.1.3\phantom{.}$$
where $V$ runs over a set of representatives for the set of $Q\-$conjugacy classes   of elements of $\frak X -\frak Y$ contained in $Q\,,$ and $\tilde\theta$ over a set of representatives for the set of $\F_Q (V)\-$classes in   $\tilde\F(P,V)_{\tilde\iota_V^Q}\,.$

\medskip
£6.2. Then, let us consider the quotient $\F^{^{\frak Y}} \-$locality   
$(\bar\upsilon^{_{\frak Y}}\!,\bar\M^{^{\frak Y}}\!,\bar\rho^{_{\frak Y}})$ of 
$\M^{^{\frak Y}}$ defined by
$$\bar\M^{^{\frak Y}}\!(Q,R) = \M^{^{\frak Y}}(Q,R)\big/\upsilon^{_{\frak Y}}_{_R} \big(Z (R)\big)
\eqno £6.2.1\phantom{.}$$
together  with the induced natural maps
$$\bar\upsilon^{_{\frak Y}}_{_{Q,R}} : \T_{\!P} (Q,R)\to \bar\M^{^{\frak Y}}(Q,R) \qq \bar\rho^{_{\frak Y}}_{_{Q,R}} : \bar\M^{^{\frak Y}}(Q,R)\to \F (Q,R)
\eqno £6.2.2\phantom{.}$$
for any $Q,R\in \frak Y\,;$ in order to show the existence of a {\it perfect $\F^{^\frak X}\-$locality\/} 
$(\tau^{_{\frak X}},\P^{^{\frak X}},\pi^{_{\frak X}})$  it suffices to prove that $\bar\rho^{_{\frak Y}}$ admits a {\it functorial\/}
 section 
 $$\sigma^{_{\frak Y}} : \F^{^{\frak Y}}\too  \bar\M^{^{\frak Y}}
 \eqno £6.2.3\phantom{.}$$
 such that the image of $\sigma^{_{\frak Y}}$ contains the image of 
 $\upsilon^{_{\frak Y}}\,.$ Note that the {\it exterior\/} quotients of 
 $\M^{^{\frak Y}}$ and $\bar\M^{^{\frak Y}}$ coincide with each other.

 \medskip
£6.3.   Indeed, in this case it is clear that the converse image~$\widehat\P^{^{\frak Y}}$ in  
$\M^{^{\frak Y}}$ of the ``image'' of $\sigma^{_{\frak Y}}$ in  
$\bar\M^{^{\frak Y}}$ is an 
$\F^{^{\frak Y}}\-$sublocality  isomorphic to $\P^{^{\frak Y}}\,,$
 so that it is a {\it perfect $\F^{^{\frak Y}}\-$locality\/}; at this point,  we consider
 the $\F^{^{\frak X}}\-$sublocality  $\widehat\P^{^{\frak X}}\!$ of 
$\bar\L^{^{\rm n,\frak X}}\!$ containing $\widehat\P^{^{\frak Y}}\!$
as a {\it full\/} subcategory over $\frak Y$ and fulfilling 
$$\widehat\P^{^{\frak X}}\!(Q,V) = \bar\L^{^{\rm n,\frak X}}\!(Q,V)
\eqno £6.3.1\phantom{.}$$
for any $Q\in \frak X$ and any $V\in \frak X -\frak Y\,.$ On the other hand, by the very definition 
of~$\bar\L^{^{\rm n,\frak X}}\!$
(cf.~£5.1), we have
$${\rm Ker}(\bar\pi^{^{\rm n,\frak X}}_{_V})  = \tilde\frak c^{^{\rm n,\frak X}}(V)/
\tilde\frak k^{^{\frak X}}\! (V)   =  \prod_{\tilde\theta\in \tilde\F(P,V)}Z(V)
\eqno £6.3.2\phantom{.}$$
\eject
\noindent
and therefore, since $p$ does not divide $\vert \tilde\F(P,V)\vert$ 
(cf.~Proposition~£5.5), we have  a surjective group homomorphism
 $$\nabla^{^{\frak X}}_V :{\rm Ker}(\bar\pi^{^{\rm n,\frak X}}_{_V}) \too Z(V)
 \eqno £6.3.3\phantom{.}$$
mapping $z = (z_{\tilde\theta})_{\tilde\theta\in \tilde\F(P,V)}$ on 
$$\nabla^{^{\frak X}}_V (z) = {1\over\vert \tilde\F(P,V)\vert}\.
\sum_{\tilde\theta\in \tilde\F(P,V)} z_{\tilde\theta}
\eqno £6.3.4,$$
 which is a section of the restriction to $Z(V)$ and  ${\rm Ker}(\bar\pi^{^{\rm n,\frak X}}_{_V}) $ of (cf.~£5.1.3)
$$\bar\tau^{^{\rm n,\frak X}}_{_V}\! : N_P(V)\too \bar\L^{^{\rm n,\frak X}}\!(V)
\eqno £6.3.5.$$

\medskip
£6.4. Finally, considering  the {\it contravariant Dirac functor\/}  
$\frak d^{^{\!\frak X}}\,\colon \widehat\P^{^{\frak X}}\to \Ab$
mapping any $Q\in \frak Y$ on $\{0\}$ and any $V\in \frak X -\frak Y$ on 
${\rm Ker} (\nabla^{^{\frak X}}_V)\,,$  the quotient $\F^{^{\frak X}}\! \-$locality  
$\P^{^{\frak X}}\! = \widehat\P^{^{\frak X}}\!/\frak d^{^{\!\frak X}}$ (cf.~£2.10), 
where we denote the structural functors by
$$\tau^{_{\frak X}} : \T_{\!P}\too \P^{^{\frak X}} \qq
\pi^{_{\frak X}} : \P^{^{\frak X}}\too \F 
\eqno £6.4.1,$$
is actually a {\it perfect  $\F^{^{\frak X}}\!\-$locality\/}; indeed, it is quite clear that the $\F^{^{\frak X}}\!\-$locality 
$\widehat\P^{^{\frak X}}$ and therefore the 
$\F^{^{\frak X}}\!\-$locality $\P^{^{\frak X}}$ are both {\it $p\-$coherent\/} 
(cf.~£2.8); moreover, for any $Q\in \frak Y$ fully normalized in $\F\,,$ we already know that 
$\P^{^{\frak X}}(Q) = \widehat\P^{^{\frak X}}(Q)$ is an {\it $\F\-$localizer\/} of $Q$ and, for  any $V\in \frak X -\frak Y\,,$ we have the exact sequence
$$1\too Z(V)\too \P^{^{\frak X}}\!(V)\too \F(V)\too 1
\eqno £6.4.2$$
which, together with the group homomorphisms
$$\tau^{_{\frak X}}_{_{V}} : N_P (V)\too \P^{^{\frak X}}\!(V) \qq
\pi^{_{\frak X}}_{_{V}} : \P^{^{\frak X}}\!(V)\too \F (V)
\eqno £6.4.3,$$
is an {\it $\F^{^{\frak X}}\-$localizer\/} of $V$ whenever $V$ is fully normalized in $\F\,;$ now, our claim follows from
£2.13 above.

\medskip
£6.5. Similarly,  in order to show the uniqueness of $(\tau^{_{\frak X}},\P^{^{\frak X}},\pi^{_{\frak X}})\,,$  in the general setting  we have to consider the $\F^{^{\frak X}}\-$sublocality  $(\upsilon^{_{\frak X}}\!,\M^{^{\frak X}}\!,\rho^{_{\frak X}})$ of 
$\bar\L^{^{\rm n,\frak X}}\!$ containing $\M^{^{\frak Y}}\!$ as a {\it full\/} subcategory over $\frak Y$ and fulfilling 
$$\M^{^{\frak X}}\!(Q,V) = \bar\L^{^{\rm n,\frak X}}\!(Q,V)
\eqno £6.5.1\phantom{.}$$
for any $Q\in \frak X$ and any $V\in \frak X -\frak Y\,,$ and, as above, the quotient 
$\F^{^{\frak X}} \-$locality   $(\bar\upsilon^{_{\frak X}}\!,\bar\M^{^{\frak X}}\!,\bar\rho^{_{\frak X}})$ of~$\M^{^{\frak X}}$ defined by
$$\bar\M^{^{\frak X}}\!(Q,R) = \M^{^{\frak X}}(Q,R)\big/\upsilon^{_{\rm n,\frak X}}_{_R} \big(Z (R)\big)
\eqno £6.5.2\phantom{.}$$
together  with the induced natural maps
$$\bar\upsilon^{_{\frak X}}_{_{Q,R}} : \T_{\!P} (Q,R)\to \bar\M^{^{\frak X}}(Q,R) \qq \bar\rho^{_{\frak X}}_{_{Q,R}} : \bar\M^{^{\frak X}}(Q,R)\to \F (Q,R)
\eqno £6.5.3\phantom{.}$$
\eject
\noindent
for any $Q,R\in \frak X\,;$ thus,  for any $Q\in \frak Y\,,$ from the exact 
sequence~£6.1.3 we obtain the exact sequence
$$1\too \prod_W\,\prod_{\tilde\theta} Z(W)\too \bar\M^{^{\frak X}}(Q)
\too \F^{^{\frak X}}(Q)\too 1
\eqno £6.5.4,$$
where $W$ runs over a set of representatives for the set of $Q\-$conjugacy classes   of elements of $\frak X -\frak Y$ contained in $Q$  and $\tilde\theta$ over a set of representatives for the $\F_Q (W)\-$classes in $\tilde\F(P,W)_{\tilde\iota_W^Q}\,,$ whereas for any 
$V\in \frak X-\frak Y$ it follows again from~£5.1 that we have the exact sequence
$$1\too {\rm Ker} (\nabla^{^{\frak X}}_V)\too \bar\M^{^{\frak X}}(V)\too \F^{^{\frak X}}(V)\too 1
\eqno £6.5.5.$$

\medskip
£6.6. At this point, if  $(\tau'^{^{\frak X}}\!,\P'^{^{\frak X}}\!,\pi'^{^{\frak X}})$ is another  
{\it perfect  $\F^{^{\frak X}}\!\-$locality\/}, once again  it follows from Corollary~£5.20 that we may assume that   
$(\tau'^{_{\frak X}},\P'^{^{\frak X}},\pi'^{_{\frak X}})$ is a {\it $\F^{^{\frak X}}\!\-$sublocality\/} of the   
{\it natural $\F^{^{\frak X}}\-$locality\/} $(\bar\tau^{^{\rm n,\frak X}}\!,\bar\L^{^{\rm n,\frak X}}\!,\bar\pi^{^{\rm n,\frak X}})$ (cf.~£5.1); then, for the corresponding {\it full\/} subcategories over $\frak Y\,,$ we have  
$(\P'^{^{\frak X}})^{^\frak Y} \i (\bar\L^{^{\rm n,\frak X}}\!)^{^\frak Y}\,;$ moreover, from the induction hypothesis and from  Corollary~£5.21, we may assume that the image $\P'^{^{\frak Y}}$ of $(\P'^{^{\frak X}})^{^\frak Y}$ in 
$\bar\L^{^{\rm n,\frak Y}}\!$ coincides with $\P^{^\frak Y}\,.$ Consequently, $(\P'^{^{\frak X}})^{^\frak Y} \!$ is contained 
in $\M^{^\frak Y}$ and  therefore  $\P'^{^{\frak X}}\!$ is contained in $\M^{^{\frak X}}\,;$ thus, the image of 
$\P'^{^{\frak X}}\!$ in $\bar\M^{^{\frak X}}\!,$ determines a  {\it functorial\/} section
$$\sigma'^{_{\frak X}} : \F^{^{\frak X}}\too  \bar\M^{^{\frak X}}
 \eqno £6.6.1\phantom{.}$$
 which induces a {\it functorial\/} section $\sigma'^{_{\frak Y}}$ as in~£6.2.3 above and, according to Theorem £6.22 below, this {\it functorial\/} section is 
 {\it naturally $\F^{^\frak Y}\-$isomorphic\/} to~$\sigma ^{_{\frak Y}}\,,$ 
 so that  we may assume that both coincide with each other (cf.~£2.9); that is to say, we get an $\F^{^{\frak X}}\-$locality isomorphism between 
 $(\tau'^{^{\frak X}}\!,\P'^{^{\frak X}}\!,\pi'^{^{\frak X}})$ and    the quotient 
 $\F^{^{\frak X}} \-$locality   $(\tau^{^{\frak X}}\!,\P^{^{\frak X}}\!,
 \pi^{^{\frak X}})$ in~£6.4 above. In conclusion, the existence and the uniqueness of the  {\it perfect  
 $\F^{^{\frak X}}\!\-$locality\/} depends on the existence and the uniqueness of a {\it functorial section\/}
 of $\bar\rho^{_{\frak X}}\,,$ proved in Theorem~£6.22  below.

\medskip
£6.7. In order to prove this theorem, we apply the methods developed in~[12,~\S4];  as  in~[12,~4.2],   $\tilde\F^{^{\frak X}}\!$ can be considered as an {\it $\tilde\F^{^{\frak X}}_{\! P}\-$category\/} and clearly it fulfills condition~[12,~4.2.1]; 
then, through the  {\it $\tilde\F^{^{\frak X}}\!\-$intersection\/} in~£5.7 above, $\tilde\F^{^{\frak X}}\!$ becomes a
{\it multiplicative $\tilde\F^{^{\frak X}}_{\! P}\-$category\/} and, as in [12,~4.6], we consider the functor and the 
{\it natural map\/} defined there
$$\frak m^{_\frak X}_P : \tilde\F^{^{\frak X}}\! \too \ad (\tilde\F^{^{\frak X}})\qq 
\omega^{_\frak X} : \frak m^{_\frak X}_P\too \frak  j^{_\frak X}
\eqno £6.7.1\phantom{.}$$
where $\frak j^{_\frak X}\,\colon \tilde\F^{^{\frak X}}\! \to \ad (\tilde\F^{^{\frak X}} )$ denotes the canonical functor.
 In our situation, we have to consider  the   {\it contravariant\/} functor
$$\frak z^{_{\frak X}}_U  : \tilde\F^{^{\frak X}}\too \Ab
\eqno £6.7.2\phantom{.}$$
\eject
\noindent
mapping any $Q\in \frak Y$ on $\{0\}$ and any $V\in \frak X -\frak Y$ on $Z(V)\,;$ as in~[12,~4.6], it is clear that $\frak z^{_{\frak X}}_U$ determines an {\it additive contravariant\/} functor
$$\frak z^{_{\frak X,\ad}}_U  : \ad (\tilde\F^{^{\frak X}})\too \Ab
\eqno £6.7.3.$$
Now, considering the {\it contravariant\/} functor  (cf.~£2.8)
$$\Ker (\bar\rho^{_{\frak X}}) : \tilde\F^{^{\frak X}}\too \Ab
\eqno £6.7.4\phantom{.}$$
which maps any $Q\in \frak X$ on ${\rm Ker} (\bar\rho^{_{\frak X}}_{_Q})$ (cf.~£6.2),  from  the functor 
$\frak m^{_\frak X}_P$ and the {\it natural map\/} $\omega^{_\frak X}$ we easily get  the following exact sequence of {\it contravariant\/} functors
$$\{0\}\too \frak z^{_{\frak X}}_U\buildrel 
\frak z^{_{\frak X,\ad}}_U \!\!* \omega^{_\frak X} \over{\hbox to 35pt{\rightarrowfill}} \frak z^{_{\frak X,\ad}}_U \!\circ \frak  m^{_\frak X}_P\;
{\hbox to 25pt{\rightarrowfill}}\; \Ker (\bar\rho^{_{\frak X}})\too \{0\}
\eqno £6.7.5.$$

\medskip
£6.8. Moreover,  since $U$ is  {\it fully normalized\/} in $\F\,,$ we have the Frobenius  $N_P(U)\-$category 
$N_\F (U)$ and we can consider the analogous {\it contravariant\/} functors; more explicitly,  set $N = N_P (U)$ 
and $\N = N_\F (U)\,,$ denote by $\frak N$ the set of $Q\in \frak X$ contained in $N\,,$ and consider the corresponding  
{\it contravariant\/} functors
$$\frak z^{_{\frak N}}_U  : \tilde\N^{^{\frak N}}\too \Ab\qq
\frak z^{_{\frak N,\ad}}_U  : \ad (\tilde\N^{^{\frak N}})\too \Ab
\eqno £6.8.1;$$
we clearly have an {\it inclusion\/} functor $\frak i_U\,\colon 
\N^{^{\frak N}}\!\!\to \F^{^{\frak X}}\,,$ which induces a functor 
$\tilde\frak i_U\,\colon \tilde\N^{^{\frak N}}\!\!\to \tilde\F^{^{\frak X}}$ 
fulfilling $\frak z^{_{\frak N}}_U = 
\frak z^{_{\frak X}}_U\circ \tilde\frak i_U\,.$ Our proof of Theorem~£6.22 below starts by showing that
the question can be reduced to the analogous question over the Frobenius $N\-$category $\N$ and the set $\frak N$
of subgroups of $N\,;$ the next proposition is a key point in this reduction.

 \bigskip
 \noindent
 {\bf Proposition~£6.9.} {\it For any $n\in \Bbb N$ the inclusion functor 
  $\,\tilde\frak i_U\,\colon \tilde\N^{^{\frak N}}\!\!\to \tilde\F^{^{\frak X}}$ induces 
 a  group isomorphism
 $$\Bbb H^n (\tilde\F^{^{\frak X}},\frak z^{_ {\frak X}}_U)\cong 
 \Bbb H^n \big(\tilde\N^{^{\frak N}}\!,\frak z^{_ {\frak N}}_U\big)
 \eqno £6.9.1.$$\/}

\par
\noindent
{\bf Proof:}  It follows from [2,~Proposition~3.2] that  both members of this isomorphism are  canonically isomorphic to $\Bbb H^n (\tilde\T_{\tilde\F (U)},\frak d_{\{1\}})$ where $\tilde\T_{\tilde\F (U)}$ denotes the {\it exterior quotient\/} (cf.~£2.1)
of the category over the set of $p\-$subgroups of~$\tilde\F (U)$ determined by the {\it transporter\/} in 
this group --- that is to say, for any pair of subgroups $Q$ and $R$ of $\tilde\F (U)\,,$ we set 
$$\tilde\T_{\tilde\F (U)} (Q,R) = Q\backslash T_{\tilde\F (U)}(R,Q)
\eqno £6.9.2,$$
the composition in $\tilde\T_{\tilde\F (U)}$ being defined by the product in $\tilde\F (U)$ --- and 
the {\it contravariant\/} functor $\frak d_{\{1\}}\,\colon \tilde\T_{\tilde\F (U)}\to \Ab$ maps any nontrivial $p\-$subgroup of $\tilde\F (U)$ on $\{0\}$ and $\{1\}$ on $Z(U)\,.$ Then, it is easily checked that these canonical isomorphisms are
compatible with the group homomorphism induced by~$\tilde\frak i_U\,.$
\eject

\medskip
£6.10. {\it Mutatis mutandis\/}, as in~£6.7 above,   $\tilde\N^{^{\frak N}}\!$ is actually an {\it $\tilde\F^{^{\frak N}}_{\! N}\-$category\/} 
and it fulfills the conditions in~[12,~2.1]; thus,  we can consider the corresponding functor and the  corresponding
{\it natural map\/} defined  in [12,~4.6]
$$\frak m^{_\frak N}_N : \tilde\N^{^{\frak N}}\! \too \ad (\tilde\N^{^{\frak N}})
\qq \omega^{_\frak N} : \frak m^{_\frak N}_N\too \frak  j^{_\frak N}
\eqno £6.10.1\phantom{.}$$
where $\frak j^{_\frak N}\,\colon \tilde\N^{^{\frak N}}\! \to 
\ad (\tilde\N^{^{\frak N}} )$ denotes the canonical functor. Then, it is quite clear that
we have the commutative diagram
$$\matrix{\tilde\F^{^{\frak X}}&\buildrel \frak j^{_{\frak X}}\over\too 
& \ad (\tilde\F^{^{\frak X}} )\cr
\hskip-20pt{\scriptstyle \tilde\frak i_U}\uparrow&\phantom{\big\uparrow}
&\uparrow{\scriptstyle \ad (\tilde\frak i_U)}\hskip-30pt\cr
\tilde\N^{^{\frak N}}&\buildrel \frak j^{_{\frak N}}\over\too 
&\ad (\tilde\N^{^{\frak N}} )\cr}
\eqno £6.10.2;$$
explicitly, for any $Q\in \frak N$ we may assume that
$$\frak m^{_\frak N}_N (Q) = \bigoplus_T\,\bigoplus_{\tilde\delta}\, T
\eqno £6.10.3,$$
where $T$ runs over a set of representatives for the set of $Q\-$conjugacy classes   of elements of $\frak N$ contained in $Q$
and then $\tilde\delta$ runs over a set of representatives in $\tilde\N (N,T)_{\tilde\iota_T^Q}$ for the set of 
$\tilde\F_Q (T)\-$orbits; note that
$$\tilde\iota_N^P\circ\tilde\N (N,T)_{\tilde\iota_T^Q} = 
\tilde\F (P,T)_{\tilde\iota_T^Q}\cap \big(\tilde\iota_N^P\circ\tilde\N (N,T)\big)
\eqno £6.10.4;$$
indeed, for any $\tilde\delta\in \tilde\N (N,T)$ and any subgroup $R$ of $Q$ strictly containing~$T\,,$ 
the existence of an $\tilde\F\-$morphism $\tilde\psi\,\colon R\to P$ fulfilling $\tilde\psi\circ\tilde\iota_T^R = 
\tilde\iota_N^P\circ\tilde\delta$ forces $\tilde\psi = \tilde\iota_N^P\circ\tilde\eta$
for some $\tilde\eta\in \tilde\N (N,R)\,.$

\medskip
£6.11. In particular, we clearly have the {\it natural map\/}
$$\theta_U : \ad (\tilde\frak i_U)\circ \frak m^{_\frak N}_N \too 
\frak m^{_\frak X}_P\circ \tilde\frak i_U 
\eqno £6.11.1,$$
 between the functors $\ad (\tilde\frak i_U)\circ \frak m^{_\frak N}_N$
and $\frak m^{_\frak X}_P\circ \tilde\frak i_U$  from 
$\tilde\N^{^{\frak N}}$ to~$\ad (\tilde\F^{^{\frak X}} )\,,$
sending $Q\in \frak N$ to the $\ad (\tilde\F^{^{\frak X}} )\-$morphism
$\frak m^{_\frak N}_N (Q)\to \frak m^{_\frak X}_P (Q)$ determined by
the identity on the subgroups $T$ of $Q$ and by the maps  
$$\tilde\N (N,T)_{\tilde\iota_T^Q}\big/\tilde\F_Q (T)\too \tilde\F (P,T)_{\tilde\iota_T^Q}\big/\tilde\F_Q (T)
\eqno £6.11.2\phantom{.}$$
sending the class of $\tilde\delta\in \tilde\N (N,T)_{\tilde\iota_T^Q}$ to  the class of $\tilde\iota_N^P\circ\tilde\delta\,;$ consequently, we still have  the 
{\it natural map\/}
$$\frak z^{_{\frak X,\ad}}_U * \theta_U :  
(\frak z^{_{\frak X,\ad}}_U\circ  \frak m^{_\frak X}_P)\circ \tilde\frak i_U
\too \frak z^{_{\frak N,\ad}}_U \circ \frak m^{_\frak N}_N 
\eqno £6.11.3\phantom{.}$$
 between the {\it $\Ab\-$valued contravariant\/} functors 
$(\frak z^{_{\frak X,\ad}}_U\circ  \frak m^{_\frak X}_P)\circ \tilde\frak i_U$ and $\frak z^{_{\frak N,\ad}}_U \circ \frak m^{_\frak N}_N $ from $\tilde\N^{^{\frak N}}\,;$ then, from the exact sequence~£6.7.5 we easily  get the commutative diagram of  {\it $\Ab\-$valued contravariant\/} functors  from 
$\tilde\N^{^{\frak N}}$
$$\matrix{0\too& \frak z^{_{\frak N}}_U&\too 
&(\frak z^{_{\frak X,\ad}}_U \!\circ \frak  m^{_\frak X}_P)\circ \tilde\frak i_U
&\too &\Ker (\bar\rho^{_{\frak X}})\circ \tilde\frak i_U &\too 0\cr
&\Vert&&\hskip-30pt{\scriptstyle \frak z^{_{\frak X,\ad}}_U * \theta_U}\big\downarrow&\phantom{\Big\downarrow}&
\hskip-30pt{\scriptstyle \overline{\frak z^{_{\frak X,\ad}}_U * \theta_U\!}\,}\big\downarrow\cr
0\too& \frak z^{_{\frak N}}_U&\too &\frak z^{_{\frak N,\ad}}_U \circ \frak m^{_\frak N}_N &\too 
&\overline{\frak z^{_{\frak N,\ad}}_U \circ \frak m^{_\frak N}_N \!}\, &\too 0\cr}
\eqno £6.11.4.$$
\eject

\medskip
£6.12. In conclusion, from the exact sequence~£6.7.5 and the diagram above
 we obtain the commutative diagram of {\it long exact sequences\/} of  cohomology groups
$$\hskip-2pt\matrix{\to\hskip-7pt&\Bbb H^n (\tilde\F^{^{\frak X}}\!\!,
\frak z^{_{\frak X}}_U)\!&\hskip-10pt\to\! 
&\hskip-5pt\Bbb H^n (\tilde\F^{^{\frak X}}\!\!,
\frak z^{_{\frak X,\ad}}_U \!\circ \frak  m^{_\frak X}_P)\!&\hskip-9pt\to 
&\hskip-10pt\Bbb H^n (\tilde\F^{^{\frak X}}\!\!,\Ker (\bar\rho^{_{\frak X}}))
&\hskip-10pt\to\cr
&\downarrow&\phantom{\big\downarrow}&\downarrow&&\downarrow\cr
\to\hskip-7pt&\Bbb H^n (\tilde\N^{^{\frak N}}\!\!,\frak z^{_{\frak N}}_U)\!
&\hskip-10pt\to\! &\hskip-10pt\Bbb H^n (\tilde\N^{^{\frak N}}\!\!,
\frak z^{_{\frak X,\ad}}_U \!\!\circ \frak  m^{_\frak X}_P\circ \tilde\frak i_U)
&\hskip-10pt\to \!&\hskip-7pt\Bbb H^n (\tilde\N^{^{\frak N}}\!\!,
\Ker (\bar\rho^{_{\frak X}})\circ \tilde\frak i_U)& \hskip-8pt\to\cr
&\Vert&\phantom{\Big\downarrow}&\downarrow&&\downarrow\cr
\to\hskip-7pt&\Bbb H^n (\tilde\N^{^{\frak N}}\!\!,\frak z^{_{\frak N}}_U)\!
&\hskip-10pt\to\! &\hskip-5pt\Bbb H^n (\tilde\N^{^{\frak N}}\!\!,
\frak z^{_{\frak N,\ad}}_U \!\circ \frak  m^{_\frak N}_N)
& \hskip-9pt\to& \hskip-5pt\Bbb H^n (\tilde\N^{^{\frak N}}\!\!,
\overline{\frak z^{_{\frak N,\ad}}_U \circ \frak m^{_\frak N}_N \!})
&\hskip-10pt\to\cr}
\eqno £6.12.1;$$
but, it follows from [12,~Proposition~4.16] that for any $n\ge 1$ we have
$$\Bbb H^n (\tilde\F^{^{\frak X}}\!\!, \frak z^{_{\frak X,\ad}}_U \!\circ \frak  m^{_\frak X}_P) = \{0\} = 
\Bbb H^n (\tilde\N^{^{\frak N}}\!\!, \frak z^{_{\frak N,\ad}}_U \!\circ \frak  m^{_\frak N}_N)
\eqno £6.12.2;$$
hence,  for any $n\ge 1\,,$ from the commutative diagram above of {\it long exact sequences\/}, we get
the commutative diagram
$$\matrix{\Bbb H^n (\tilde\F^{^{\frak X}}\!\!,\Ker (\bar\rho^{_{\frak X}}))
&\cong &\Bbb H^{n+1} (\tilde\F^{^{\frak X}}\!\!,\frak z^{_{\frak X}}_U)\cr
\hskip-60pt{\scriptstyle \Bbb H^{n} (\tilde\frak i_U,\Ker (\bar\rho^{_{\frak X}}))}\hskip4pt\downarrow&\phantom{\Big\downarrow}\cr
\Bbb H^n (\tilde\N^{^{\frak N}}\!\!,\Ker (\bar\rho^{_{\frak X}})\circ 
\tilde\frak i_U) &&\Big\downarrow\hskip4pt
{\scriptstyle \Bbb H^{n+1} (\tilde\frak i_U,\frak z^{_{\frak X}}_U)}\hskip-50pt\cr
\hskip-70pt{\scriptstyle \Bbb H^{n}(\tilde\N^{^{\frak N}}\!,\,\overline{\frak z^{_{\frak X,\ad}}_U * \theta_U\!}\,)}\hskip4pt\downarrow&\phantom{\Big\downarrow}\cr
\Bbb H^n (\tilde\N^{^{\frak N}}\!\!,
\overline{\frak z^{_{\frak N,\ad}}_U \circ \frak m^{_\frak N}_N \!})&\cong &
\Bbb H^{n+1} (\tilde\N^{^{\frak N}}\!\!,\frak z^{_{\frak N}}_U)\cr}
\eqno £6.12.3.$$

\medskip
£6.13. At this point,  for any $n\ge 1\,,$ it follows from Proposition~£6.9 that in order to show that
an {\it $n\-$cocycle\/} in $\Bbb C^n \big(\tilde\F^{^{\frak X}}\!\!,\Ker (\bar\rho^{_{\frak X}})\big)$ is 
an {\it $n\-$coboundary\/}, it suffices to prove  that the image in 
$\Bbb C^n (\tilde\N^{^{\frak N}}\!\!,\overline{\frak z^{_{\frak N,\ad}}_U \circ \frak m^{_\frak N}_N \!}\,)$ of
 its restriction to~$\Bbb C^n (\tilde\N^{^{\frak N}}\!\!,\Ker (\bar\rho^{_{\frak X}})\circ \tilde\frak i_U)$ is 
 an {\it $n\-$coboundary\/}; but, the announced existence of a {\it functorial section\/} of~$\bar\rho^{_{\frak X}}$ 
 depends on the proof that a suitable {\it $2\-$cocycle\/} 
 in~$\Bbb C^2 \big(\tilde\F^{^{\frak X}}\!\!,\Ker (\bar\rho^{_{\frak X}})\big)$ is a {\it $2\-$coboundary\/}; hence,
it suffices to prove that the corresponding image of this {\it $2\-$cocycle\/} in 
$\Bbb C^2 (\tilde\N^{^{\frak N}}\!\!,\overline{\frak z^{_{\frak N,\ad}}_U \circ \frak m^{_\frak N}_N \!}\,)$
is a {\it $2\-$coboundary\/}. This goes towards our announced reduction (cf.~£6.8), provided we show that 
this image coincides with the corresponding  {\it $2\-$cocycle\/} in the reduced situation, which demands an explicit connexion between 
both situations.

\medskip
£6.14. Recall that the {\it $\F^{^\frak X}\!\-$locality\/} $\bar\L^{^{\rm n,\frak X}}\!$ comes from a {\it natural $\F\-$basic 
$P\times P\-$set\/} $\Omega$ (cf.~£3.5); more explicitly,  denoting by $G$ the group of $\{1\}\times P\-$set automorphisms~of 
${\rm Res}_{\{1\}\times P}(\Omega)$ and identifying $P$ with the image of~$P\times \{1\}$ into~$G$ (cf.~£4.1), 
$\bar\L^{^{\rm n,\frak X}}\!$ is a suitable quotient (cf.~£4.13 and~£4.14) of the~{\it $\F^{^\frak X}\!\-$locality\/} defined by the {\it transporter\/} over $\frak X$ in $G$ (cf.~£4.1). Then, it follows from [11,~Proposition~21.11] that the subset
$$\Gamma = \bigcup_{\chi\in \F(U)} \Omega^{\Delta_\chi (U)}
\eqno £6.14.1\phantom{.}$$
\eject
\noindent
is actually an {\it $\N\-$basic $N\times N\-$set\/}; {\it mutatis mutandis\/},  denote by $H$ the group of $\{1\}\times N\-$set automorphisms of ${\rm Res}_{\{1\}\times N}(\Gamma)$ and, to avoid confusion, denote by $\bar N$  
the image of $N\times \{1\}$ into~$H\,;$ since $N_G (U)/C_G (U) \cong \F (U)$ (cf.~£4.1.1), it is clear that $N_G (U)$ stabilizes $\Gamma$ and therefore we have a canonical group homomorphism from $N_G (U)$ to $H\,;$ 
moreover, denote by 
$$\iota_{_U} : N\too N_H (\bar U)\qq \kappa_{_U} : N_H (\bar U)\too \N (U)  =\F (U)
\eqno £6.14.2\phantom{.}$$ 
the obvious maps.

\bigskip
\noindent
{\bf Proposition~£6.15.} {\it With the notation above, for any pair of  subgroups $Q$ and~$R$ of $N$  containing~$U$ 
and any element $\psi$ in $\N (Q,R)\,,$ there exists at most one $Q\times N\-$orbit in $\Gamma$ isomorphic to 
$(Q\times N)/\Delta_\psi (R)$ and then $\psi$ belongs 
to~$\tilde\N (N,R)_{\tilde\iota_R^Q}\,.$ In particular, $\Gamma$ is a natural 
$\N\-$basic $N\times N\-$set and $C_H (\bar Q)$ is an Abelian $p\-$group.\/}

\medskip
\noindent
{\bf Proof:} Since $\Gamma$ is  an {\it $\N\-$basic $N\times N\-$set\/}, for any $\omega\in\Gamma$ we already  know that 
$(Q\times N)_\omega = \Delta_\psi (R)$ for some  subgroup $R$ of~$P$
and some $\psi\in\F (P,R)$ (cf.~£3.1), and, by the very definition of $\Gamma\,,$ $\Delta_\psi (R)$
contains $\Delta_\chi (U)$ for some $\chi\in \F (U)\,;$ that is to say, $R$ contains $U$ and $\psi$ extends
$\chi\,,$ so that $\omega$ determines $\chi\,.$

\smallskip
We claim that we still have $(Q\times P)_\omega = \Delta_\psi (R)\,;$ indeed, if $(v,w)\in Q\times P$ fixes $\omega$
then $\Delta_{\chi^v} (U)$ fixes $\omega\.w^{-1}\,;$ that is to say, for any $u\in U$ we have
$$\omega\.w^{-1} = \chi^v (u)\.\omega\. w^{-1}u^{-1} = \omega\.(\chi\circ \chi^v) (u)w^{-1}u^{-1} 
\eqno £6.15.1$$
which forces $w^{-1} = (\chi\circ \chi^v) (u)w^{-1}u^{-1} \,,$ so that $u^w =  (\chi\circ \chi^v) (u)\,;$
hence, $w$ belongs to $N$ which proves our claim. Consequently, any $Q\times N\-$orbit in $\Gamma$ isomorphic 
to~$(Q\times N)/\Delta_\psi (R)$ is the intersection of $\Gamma$ with a $Q\times P\-$orbit in~$\Omega$ isomorphic 
to~$(Q\times P)/\Delta_\psi (R)$ and it follows from Proposition~£3.7 and from equality~£6.10.4 that such a $Q\times P\-$orbit is unique and that  $\psi$ belongs  to~$\tilde\N (N,R)_{\tilde\iota_R^Q}\,.$

\smallskip
In particular, the last statement follows from isomorphism~£4.2.1 applied to the {\it $\N\-$basic $N\times N\-$set\/} $\Gamma$ together with 
isomorphism~£3.7.1. We are done.

\medskip
£6.16. With the notation above, denote by $\T_H^{^\frak N}$ the category over the set~$\frak N$ determined by the 
{\it transporter\/} in $N_H (\bar U)\,,$ namely for any pair of subgroups $Q$ and $R$ in $\frak N$ we set 
$$\T_H^{^\frak N} (Q,R) = T_{N_H (\bar U)}(\bar R,\bar Q)
\eqno £6.16.1,$$
the composition in $\T_H^{^\frak N}$ being defined by the product in $H\,;$ recall that in~£4.7 above we have defined an 
{\it $\N\-$locality\/} $\L^{^\Gamma}$ and, since $\Gamma$   is a {\it natural\/} $\N\-$basic $N\times N\-$set and 
$C_H (\bar Q)$ is an Abelian $p\-$group,
it follows easily from the very definitions that the {\it full\/} subcategory $(\L^{^\Gamma})^{^\frak N}$ of~$\L^{^\Gamma}$ 
over $\frak N$ coincides with $\T_H^{^\frak N}\,;$ in particular, we denote by $\bar\T_H^{^\frak N}$
the corresponding quotient defined in~£4.13.2 above. 
\eject

\medskip
£6.17. Actually,  it is easily checked from this proposition and from equalities~£6.10.4   that the canonical group homomorphism
from $N_G (U)$ to $H$ induces a {\it full\/} functor
$$\frak g_U : N_{\bar\L^{^{\rm n,\frak X}}\!} (U) \too \bar\T_H^{^\frak N}  
\eqno £6.17.1;$$
but, recall that $\bar\L^{^{\rm n,\frak Y}}\!$ contains  the {\it perfect $\F^{^\frak Y}\-$locality\/} $\P^{^\frak Y}$
and note that,  setting $\frak M = \frak Y\cap \frak N\,,$  the normalizer $\P_U^{^\frak M} =N_{\P^{^\frak Y}}(U)$  
is nothing but a {\it perfect $\N^{^\frak M}\-$locality\/}; consequently, denoting by $\bar \T_H^{^\frak M}$
the corresponding quotient  defined in~£4.13.2 above of the  category $\T_H^{^\frak M}$ over the set~$\frak M$ determined by the 
{\it transporter\/} in~$N_H (\bar U)\,,$ the functor $\frak g_U$ induces a functor 
$$N_{\bar\L^{^{\rm n,\frak Y}}\!} (U) \too \bar\T_H^{^\frak M}
\eqno £6.17.2$$
which maps $\P_U^{^\frak M}$ onto an {\it equivalent $\N^{^\frak M}\-$locality\/}.

\medskip
£6.18. On the other hand, let $(\tau_{_U},L,\pi_{_U})$  be~an {\it $\F\-$localizer\/} of $U$ 
 (cf.~£2.13); thus, $L$ is a finite group, we have the injective and the surjective group homomorphisms
$$\tau_{_U} : N\too L \qq \pi_{_U} : L\too \F (U)
\eqno £6.18.1,$$
$\tau_{_U} (N)$ is a Sylow $p\-$subgroup of $L\,,$ and we also have the exact sequence
$$1\too Z(U)\buildrel \tau_{_U}\over\too L\buildrel \pi_{_U}\over\too \F(U)\too 1
\eqno £6.18.2;$$
note that, denoting by~$\T_L^{^\frak M}$ the {\it $\N^{^\frak M}\!\!\-$locality\/} determined by $\tau_{_U}$ and by the 
{\it transporter\/} over the set~$\frak M$ in the group $L\,,$ it is easily checked that $\T_L^{^\frak M}$ is also a 
{\it perfect $\N^{^\frak M}\-$locality\/}; consequently, it follows from the induction hypothesis that the 
{\it $\N^{^\frak M}\!\!\-$localities\/} $\P_U^{^\frak M}$ and  $\T_L^{^\frak M}$ are isomorphic.  Going towards our announced reduction (cf.~£6.8), we need the following result,  essentially proved in [11,~Proposition~18.20]. 

\bigskip
\noindent
{\bf Proposition~£6.19.} {\it With the notation above, there is a unique $C_H (\bar U)\-$con-jugacy class of group homomorphisms
$$\lambda_{_U} : L\too N_H (\bar U)
\eqno £6.19.1$$
fulfilling $\lambda_{_U}\circ\tau_{_U} = \iota_{_U} $ and 
$\kappa_{_U}\circ \lambda_{_U} = \pi_{_U}\,.$\/}

\medskip
\noindent
{\bf Proof:} We will apply [11,~Lemma~18.8] to the groups $L$ and 
$N_H(\bar U)\,,$ to the quotient $\F(U)$ of both, and to the group homomorphism 
$$\eta : \tau_{_U} (N)\too N_H(\bar U)
\eqno £6.19.2\phantom{.}$$
mapping $\tau_{_U}(u)$ on $\bar u$ for any $u\in N\,,$ which makes sense since $\tau_{_U}$ is injective; we already know that 
$\tau_{_U}(N)$ is a Sylow $p\-$subgroup of $L$  and, according to Proposition~£6.15, that $C_H (\bar U)$ is an 
Abelian $p\-$group; moreover, it is clear that the restriction of $\pi_{_U}$ to
$\tau_{_U}(N)$ coincides with  $\kappa_{_U}\circ\eta$ (cf.~£6.14.2).
\eject

\smallskip
Let~$Q$~be~a subgroup of $N$ and $x\in L$ such that $\tau_{_U}(Q)\i \tau_{_U}(N)^x\,;$ since we still have 
$$\tau_{_U}(U\.Q)\i \tau_{_U}(N)^x
\eqno £6.19.3,$$ 
it follows from the very definition of the $\F\-$localizer [11, Theorem~18.6] that there is an $\F\-$morphism 
$\varphi\,\colon U\.Q\to N$ such that we have
$$\varphi (u) = \big(\pi_{_U} (x)\big) (u)\qq \tau_{_U}\big(\varphi (v)\big) 
= x\tau_{_U}(v) x^{-1}
\eqno £6.19.4\phantom{.}$$
for any $u\in U$ and any $v\in U\.Q\,;$ thus, $\varphi$ is also an $\N\-$morphism from $U\.Q$ to $N$
and therefore there is $\bar x\in T_H (\bar U\.\bar Q,\bar N)$ fulfilling [11,~Proposition~21.9]
$$\varphi (u) = \big(\kappa_{_U}({\bar x})\big) (u) \qq \overline{\varphi (v)} = \bar x \bar v\bar x^{-1}
\eqno £6.19.5\phantom{.}$$
for any $u\in U$ and any $v\in U\.Q\,.$ Consequently,  we get
$$\kappa_{_U}(\bar x) = \pi_{_U}( x)\qq
\eta\big(x\tau_{_U}(v) x^{-1}\big) =\bar x\bar v\bar x^{-1} = \bar x\eta\big(\tau_{_U} (v)\big)\bar x^{-1}
\eqno £6.19.6\phantom{.}$$
for any $v\in U\.Q\,;$ that is to say, the condition in~[11,~Lemma~18.8] is fulfilled
and therefore, since $C_H (\bar U)$ is an Abelian $p\-$group,  it follows from [11,~Lemma~18.8] that there is a unique
$C_H (\bar U)\-$conjugacy class of group homomorphisms as
announced.

\medskip
£6.20.  As above,  denote by~$\T_L^{^\frak N}$ the {\it $\N^{^\frak N}\!\!\-$locality\/} determined by $\tau_{_U}$ and by the 
{\it transporter\/} over the set~$\frak N$ in the group $L\,;$ since $\bar\T_H^{^\frak N}$ is a quotient of the category
$\T_H^{^\frak N}$ determined by the  {\it transporter\/}  over the set~$\frak N$ in $N_H (\bar U)\,,$  any group homomorphism
$\lambda_{_U}\,\colon  L\to N_H (\bar U)$ in~£6.19.1 above determines an {\it $\N^{^\frak N}\!\-$locality functor\/}
$$\frak l_{_U} : \T_L^{^\frak N}\too \bar\T_H^{^\frak N}
\eqno £6.20.1\phantom{.}$$
and two of them are {\it naturally $\N^{^\frak N}\!\-$isomorphic\/} (cf.~£2.9). Such a functor induces an obvious 
{\it $\N^{^\frak N}\!\-$locality functor\/} from $\T_L^{^\frak M}$ to $\bar\T_H^{^\frak M}$
and then, since  the {\it $\N^{^\frak M}\!\!\-$localities\/} $\P_U^{^\frak M}$ and  $\T_L^{^\frak M}$ are isomorphic,
 it follows from Corollary~£5.21 that we may assume that they have the same image in $\bar\T_H^{^\frak M}\,.$
 At this point, denoting by $(\bar\T_H^{^\frak N})^{^\frak M}$ the {\it full\/} subcategory of  $\bar\T_H^{^\frak N}$
 over $\frak M\,,$ as in 6.2.1 above  we can look to the {\it pull-back\/}
$$\matrix{\P^{^{\frak M}}_U&\i &\bar\T_H^{^\frak M}\cr
\uparrow&\phantom{\big\uparrow}&\uparrow\cr
\M^{^{\frak M}}\!&\i &(\bar\T_H^{^\frak N})^{^\frak M}\cr}
\eqno £6.20.2\phantom{.}$$
and it is quite clear that $\M^{^{\frak M}}\!$ contains the image of $\T_L^{^\frak M}\,.$

\medskip
£6.21. Moreover, as in £6.5, we can extend $\M^{^{\frak M}}\!$ to an {\it $\N^{^\frak N}\!\!\-$sublocality\/} 
$(\upsilon^{_\frak N}\!,\M^{^\frak N}\!\!,\rho^{_\frak N})$ of $\bar\T_H^{^\frak N}$ and  $\frak l_{_U}$ actually
becomes a {\it $\N^{^\frak N}\!\-$locality functor\/} from $\T_L^{^\frak N}$ to $\M^{^\frak N}\,;$ moreover,
 it is easily checked that $\frak g_U$ induces a {\it full\/} functor
$$\frak h_U : N_{\M^{^{\frak X}}} (U) \too \M^{^{\frak N}} 
\eqno £6.21.1.$$
As in £6.5, we also consider  the quotient $\N^{^{\frak N}} \-$locality   
$(\bar\upsilon^{_{\frak N}}\!,\bar\M^{^{\frak N}}\!,\bar\rho^{_{\frak N}})$ of~$\M^{^{\frak N}}$ defined by
$$\bar\M^{^{\frak N}}\!(Q,R) = \M^{^{\frak N}}(Q,R)\big/\upsilon^{_{\frak N}}_{_R} \big(Z (R)\big)
\eqno £6.21.2\phantom{.}$$
together  with the induced natural maps
$$\bar\upsilon^{_{\frak N}}_{_{Q,R}} : \T_{\!P} (Q,R)\to \bar\M^{^{\frak N}}(Q,R) \qq 
\bar\rho^{_{\frak N}}_{_{Q,R}} : \bar\M^{^{\frak N}}(Q,R)\to \F (Q,R)
\eqno £6.21.3\phantom{.}$$
for any $Q,R\in \frak N\,,$ and we denote by $\bar\frak h_U\,\colon N_{\bar\M^{^{\frak X}}} (U) \too \bar\M^{^{\frak N}}$ the functor induced by $\frak h_U\,;$ note that the analogous quotient $\N^{^{\frak N}} \-$locality $\bar \T_L^{^\frak N}$ of
$\T_L^{^\frak N}$ coincides with $\N^{^{\frak N}}$ and similarly we denote by 
$$\bar\frak l_{_U} : \bar\T_L^{^\frak N} = \N^{^{\frak N}}\too \bar\M^{^{\frak N}}
\eqno £6.21.4\phantom{.}$$
the functor induced by $\frak l_{_U}\,.$ Finally, as in~£6.7 above,  considering the {\it contravariant\/} functor  (cf.~£2.8)
$$\Ker (\bar\rho^{_{\frak N}}) : \tilde\N^{^{\frak N}}\too \Ab
\eqno £6.21.5,$$
 from  the functor $\frak m^{_\frak N}_N$ and from the {\it natural map\/} $\omega^{_\frak N}$ (cf.~£6.10) we easily get  the following exact sequence of {\it contravariant\/} functors
$$\{0\}\too \frak z^{_{\frak N}}_U\buildrel 
\frak z^{_{\frak N,\ad}}_U \!\!* \omega^{_\frak N} \over{\hbox to 35pt{\rightarrowfill}} \frak z^{_{\frak N,\ad}}_U \!\circ \frak  m^{_\frak N}_N\;
{\hbox to 25pt{\rightarrowfill}}\; \Ker (\bar\rho^{_{\frak N}})\too \{0\}
\eqno £6.21.6\phantom{.}$$
and therefore we  have a {\it natural isomorphism\/} (cf. £6.11.4)
$$\overline{\frak z^{_{\frak N,\ad}}_U \!\circ \frak m^{_\frak N}_N \!}\, 
\cong \Ker (\bar\rho^{_\frak N})
\eqno £6.21.7.$$

\bigskip
\noindent
{\bf Theorem~£6.22.} {\it With the notation above, the structural functor $\bar\rho^{^{\frak X}}$ admits a unique natural 
$\F^{^{\frak X}}\!\-$isomorphism class of $\F^{^{\frak X}}\!\-$locality functorial sections
 $$\sigma^{_{\frak X}} : \F^{^{\frak X}}\too  \bar\M^{^{\frak X}}
 \eqno £6.22.1.$$\/} 

\par
\noindent
{\bf Proof:} For any $\F^{^{\frak X}}\!\-$morphism $\varphi\,\colon R\to Q\,,$ choose a lifting $x_\varphi$ 
in~$\M^{^{\frak X}}\! (Q,R)$ (cf.~£6.1 and~£6.5) and denote by $\bar x_\varphi$ the image of $x_\varphi$
in $\bar\M^{^{\frak X}}\! (Q,R)\,;$ actually, for any element $w$ in the {\it $P\-$transporter\/} 
$\T^{^{\frak X}}_P\!(Q,R)$ we can assume that
$$ x_{\kappa^{_{\frak X}}_{_{Q,R}}(w)}= \upsilon^{_{\frak X}}_{_{Q,R}}(w)
\eqno £6.22.2;$$
moreover, we can do our choice in such a way that we have
$$\bar x_{\kappa^{_{\frak X}}_{_Q}(u)\circ\varphi} = \bar\upsilon^{_{\frak X}}_{_Q} (u)\.\bar x_\varphi
\eqno £6.22.3\phantom{.}$$
for any $u\in Q\,;$ indeed, if we have $\kappa^{_{\frak X}}_{_Q}(u)\circ\varphi =Ê\varphi$ then we get $u = \varphi (z)$
for a suitable $z\in Z(R)\,;$ but,  since $\M^{^{\frak X}}$ is a $\F^{^{\frak X}}\!\-$sublocality of $\bar\L^{^{\rm n,\frak X}}\!$,
$\bar\M^{^{\frak X}}\!$~is~{\it coherent\/}; hence, in this case we obtain
$$\bar\upsilon^{_{\frak X}}_{_Q} (u)\.\bar x_\varphi = \bar\upsilon^{_{\frak X}}_{_Q} \big(\varphi (z)\big)\.\bar x_\varphi
= \bar  x_\varphi\.\bar\upsilon^{_{\frak X}}_{_R} (z) = \bar x_\varphi
\eqno £6.22.4.$$
\eject
\noindent
More precisely, if $\varphi\,\colon R\to Q$ is an $\N^{^{\frak N}}\!\-$morphism, we may assume that $x_\varphi$ belongs to
$\big(N_{\M^{^{\frak X}}} (U)\big)(Q,R)$ and then that $\frak h_U (x_\varphi)$  belongs to the image 
of~$\T_L^{^\frak N}(Q,R)\,;$ in particular, it follows from~£6.21 above that acually we have
$$\bar\frak h_U (\bar x_\varphi) = \bar \frak l_U (\varphi)
\eqno £6.22.5.$$

\smallskip
Then, for any triple of subgroups $Q\,,$ $R$ and $T$ in~$\frak X\,,$ 
and any pair of $\F\-$morphisms $\psi\,\colon T\to R$ and $\varphi\,\colon R\to Q\,,$ since 
$ x_\varphi\. x_\psi$ and $ x_{\varphi\circ\psi}$ have the same image $\varphi\circ\psi$ in $\F(Q,T)\,,$ 
the {\it divisibility\/} of $\M^{^{\frak X}}$ guarantees the existence and the uniqueness of
$k_{\varphi,\psi}\in {\rm Ker}(\rho^{_{\frak X}}_T)$ fulfilling
$$ x_\varphi\. x_\psi =  x_{\varphi\circ\psi}\. k_{\varphi,\psi}
\eqno £6.22.6.$$
Denote by $\bar k_{\varphi,\psi}$ the image of $k_{\varphi,\psi}$ in ${\rm Ker}(\bar\rho^{_{\frak X}}_T)\,;$
 since $\bar\M^{^{\frak X}}\!$ is  {\it coherent\/}, on the one hand for any $u\in Q$ and any $v\in R$ we get (cf.~£6.22.2)
$$\eqalign{\bar x_{\kappa^{_{\frak X}}_{_Q}(u)\circ\varphi}\. \bar x_{\kappa^{_{\frak X}}_{_R}(v)\circ\psi}
&= \big(\bar\upsilon^{_{\frak X}}_{_Q} (u)\. \bar x_\varphi\big)\. \big(\bar\upsilon^{_{\frak X}}_{_R} (v)\. \bar x_\psi\big)
= \bar\upsilon^{_{\frak X}}_{_Q} \big(u\varphi (v)\big)\.\bar x_\varphi\.\bar x_\psi\cr
\bar  x_{(\kappa^{_{\frak X}}_{_Q}(u)\circ\varphi)\circ(\kappa^{_{\frak X}}_{_R}(v)\psi)} &=
\bar x_{\kappa^{_{\frak X}}_{_Q} (u\varphi (v))\circ\varphi\circ\psi} = 
\bar\upsilon^{_{\frak X}}_{_Q} \big(u\varphi (v)\big)\. \bar x_{\varphi\circ\psi}\cr}
\eqno £6.22.7;$$
hence, from the {\it divisibility\/} of $\bar\M^{^{\frak X}}$ we obtain 
$$\bar k_{\kappa^{_{\frak X}}_{_Q}(u)\circ\varphi,
\kappa^{_{\frak X}}_{_R}(v)\circ\psi} = \bar k_{\varphi,\psi}
\eqno £6.22.8;$$
on the other hand,  we have a {\it contravariant\/} functor (cf.~£2.8)
$$\Ker (\bar\rho^{_{\frak X}}) : \tilde\F^{^{\frak X}}\too \Ab
\eqno £6.22.9\phantom{.}$$
mapping $T$ on ${\rm Ker}(\bar\rho^{_{\frak X}}_T)\,.$ That is to say,  for any $n\in \Bbb N$ setting
$$\Bbb C^n \big(\tilde\F^{^{\frak X}}\!,\Ker (\bar\rho^{_{\frak X}})\big) = 
\prod_{\tilde\frak q\in \Fct(\Delta_n,\tilde\F^{^{\frak X}})}
{\rm Ker}(\bar\rho^{^{\frak X}}_{\tilde\frak q (0)})
\eqno £6.22.10,$$
 we have obtained an element 
 $\bar k = \{\bar k_{\tilde\frak q}\}_{\tilde\frak q\in \Fct(\Delta_2,\tilde\F^{^{\frak X}})}$ in 
 $\Bbb C^2 \big(\tilde\F^{^{\frak X}}\!,\Ker (\bar\rho^{_{\frak X}})\big)$ where we set 
$\bar k_{\tilde\frak q} = \bar  k_{\frak q(1\bullet 2),\frak q (0\bullet 1)}$ for some representative 
$\frak q\,\colon \Delta_2\to \F^{^{\frak X}}$ of $\tilde\frak q\,.$

\smallskip
We claim that $\bar k$ is actually a {\it $2\-$cocycle.\/}; explicitly, considering the usual differential map [11,~A13.11]
$$\bar d^{^{\frak X}}_2 : \Bbb C^2 \big(\tilde\F^{^{\frak X}}\!\!,
\Ker (\bar\rho^{^{\frak X}})\big)\too 
\Bbb C^3 \big(\tilde\F^{^{\frak X}}\!\!,\Ker (\bar\rho^{^{\frak X}})\big)
\eqno £6.22.11,$$
we claim that $\bar d^{^{\frak X}}_2 (\bar k) = 0\,;$ indeed, with the notation above and denoting $\bar k_{\tilde\frak q}$
 by $\bar  k_{\tilde\frak q(1\bullet 2),\tilde\frak q (0\bullet 1)}$ , for a third 
$\F^{^{\frak X}}\-$morphism $\eta\,\colon W\to T$ we get
$$\eqalign{(\bar x_\varphi\.\bar x_\psi)\.\bar x_\eta &= (\bar x_{\varphi\circ\psi}
\.\bar k_{\tilde\varphi,\tilde\psi})\.\bar x_\eta
= (\bar x_{\varphi\circ\psi}\.\bar x_\eta)\.\big(\Ker (\bar\rho^{^{\frak X}})
(\tilde\eta)\big) (\bar k_{\tilde\varphi,\tilde\psi})\cr
&= \bar x_{\varphi\circ\psi\circ\eta}\.\bar k_{\tilde\varphi\circ\tilde\psi,\tilde\eta}
\.\big(\Ker (\bar\rho^{^{\frak X}})(\tilde\eta)\big)
(\bar k_{\tilde\varphi,\tilde\psi})\cr
\bar x_\varphi\.(\bar x_\psi\.\bar x_\eta) &= \bar x_\varphi\.(\bar x_{\psi\circ\eta}\.\bar k_{\tilde\psi,\tilde\eta}) = 
\bar x_{\varphi\circ\psi\circ\eta}\. \bar k_{\tilde\varphi,\tilde\psi\circ\tilde\eta}\.
\bar k_{\tilde\psi,\tilde\eta}\cr}
\eqno £6.22.12\phantom{.}$$
\eject
\noindent
and   the {\it divisibility\/} of $\bar\M^{^{\frak X}}$ forces
$$\bar k_{\tilde\varphi\circ\tilde\psi,\tilde\eta}
\.\big(\Ker (\bar\rho^{^{\frak X}})(\tilde\eta)\big)
(\bar k_{\tilde\varphi,\tilde\psi}) = \bar  k_{\tilde\varphi,\tilde\psi\circ\tilde\eta}
\. \bar k_{\tilde\psi,\tilde\eta}
\eqno £6.22.13;$$
since ${\rm Ker}(\bar\rho^{^{\frak X}}_W)$ is Abelian, in the additive notation we obtain
$$0 = \big(\Ker (\bar\rho^{^{\frak X}})(\tilde\eta)\big)(\bar k_{\tilde\varphi,\tilde\psi}) -  \bar k_{\tilde\varphi,\tilde\psi\circ\tilde\eta} 
+ \bar k_{\tilde\varphi\circ\tilde\psi,\tilde\eta} - \bar  k_{\tilde\psi,\tilde\eta}
\eqno £6.22.14,$$
proving our claim.

\smallskip
Then, it is well-known that we can modify our choice of the family of liftings 
$\{x_\varphi\}_\varphi$ getting a {\it functorial\/} section of 
$\bar\rho^{^{\frak X}}\,,$ provided we prove that $\bar k$ is a 
{\it $2\-$coboundary\/}; indeed, if $\bar k = \bar d^{^{\frak X}}_1 (\bar\ell)$ for some element 
$\bar\ell = (\bar\ell_{\tilde\frak r})_{\tilde\frak r\in \Fct(\Delta_1,\tilde\F^{^{\frak X}})}$ in 
$\Bbb C^1 \big(\tilde\F^{^{\frak X}}\!,\Ker (\bar\rho^{^{\frak X}})\big)\,,$ with the notation above we get
$$\bar k_{\tilde\varphi,\tilde\psi} = \big(\Ker (\bar\rho^{^{\frak X}})
(\tilde\psi)\big) (\bar\ell_{\tilde\varphi})\.
(\bar\ell_{\tilde\varphi\circ\tilde\psi})^{-1}\.\bar\ell_{\tilde\psi}
\eqno £6.22.15\phantom{.}$$
where we  identify any $\tilde\F^{^{\frak X}}\!\-$morphism with the obvious 
{\it $\tilde\F^{^{\frak X}}\!\-$chain\/} $\Delta_1\to \tilde\F^{^{\frak X}}\,;$ hence, from equality~£6.22.6 we obtain
$$\eqalign{\big(\bar x_\varphi\.(\bar\ell_{\tilde\varphi})^{-1}\big)\.
\big(\bar x_\psi \.(\bar\ell_{\tilde\psi})^{-1}\big)
&=    (\bar x_\varphi\.\bar x_\psi)\. \Big(\big(\Ker (\bar\rho^{^{\frak X}})
(\tilde\psi)\big)(\bar\ell_{\tilde\varphi})\.\bar\ell_{\tilde\psi}\Big)^{-1}\cr
&= \bar x_{\varphi\circ\psi}\.(\bar\ell_{\tilde\varphi\circ\tilde\psi})^{-1}\cr}
\eqno £6.22.16\phantom{.}$$
and it suffices to replace $x_\varphi$ by $x_\varphi\.(\ell_{\tilde\varphi})^{-1}$ for any lifting $\ell_{\tilde\varphi}$
of $\bar\ell_{\tilde\varphi}\,.$

\smallskip
 On the other hand, since for any pair of  $\N^{^\frak N}\!\-$morphisms $\varphi\,\colon R\to Q$ and $\psi\,\colon T\to R$
 we have chosen $x_\varphi$ in $\big(N_{\!\M^{^{\frak X}}\!} (U)\big)(Q,R)$ and $x_\psi$ in 
 $\big(N_{\!\M^{^{\frak X}}\!} (U)\big)(R,T)$, from equality~£6.22.6 we get
 $$\frak h_U (x_\varphi)\.\frak h_U (x_\psi) = \frak h_U (x_{\varphi\circ \psi})\. \frak h_U (k_{\varphi,\psi})
 \eqno £6.22.17\phantom{.}$$
 and therefore we still get
 $$\bar\frak h_U (\bar x_\varphi)\.\frak h_U (\bar x_\psi) = 
 \bar\frak h_U (\bar x_{\varphi\circ \psi})\. \bar\frak h_U (\bar k_{\varphi,\psi})
 \eqno £6.22.18,$$
 so that equalities~£6.22.5 force $\bar\frak h_U (\bar k_{\varphi,\psi}) = 1\,;$
 moreover,  it follows from~£6.11 and from isomorphism~£6.21.7 above that actually we have
 $$\bar\frak h_U (\bar k_{\varphi,\psi}) = (\overline{\frak z^{_{\frak X,\ad}}_U\! * \theta_U \!}\,)_T 
(\bar k_{\varphi,\psi}) 
\eqno £6.22.19.$$
That is to say,  the map (cf.~£6.11)
$$\Bbb C^2 (\tilde\frak i_U,\overline{\frak z^{_{\frak X,\ad}}_U\! * \theta_U \!}\,) : 
\Bbb C^2\big(\tilde\F^{^\frak X}\!,\Ker (\bar\rho^{_\frak X})\big)\too 
\Bbb C^2\big(\tilde\N^{^\frak N}\!,\Ker (\bar\rho^{_\frak N}_{_U})\big)
\eqno £6.22.20\phantom{.}$$
sends the {\it $\Ker (\bar\rho^{_\frak X})\-$valued $2\-$cocycle\/} $\bar k$ over $\tilde\F^{^\frak X}\!$ to the trivial
{\it $\Ker \big(\bar\rho^{_\frak N})\-$valued $2\-$cocycle\/} over $\tilde\N^{^\frak N}\!$; then, it follows from~£6.13 above that $\bar k$ is indeed a {\it $2\-$co-boundary\/}.
\eject

\smallskip
Thus, we have obtained a {\it functorial section\/} $\sigma^{_\frak X}\,\colon \F^{^\frak X}\too \bar\M^{^\frak X}$
of $\bar\rho^{^\frak X}\,;$ actually, $\sigma^{_\frak X}$ can be  modified  in order to get an {\it $\F^{^\frak X}\!\-$locality functorial section\/}. Indeed, for any 
$\F^{^\frak X}_{\! P}\-$morphism $\zeta\,\colon R\to Q\,,$ choosing $u_\zeta$ in $T_P (R,Q)$ lifting $\zeta\,,$
both $\bar\M^{^{\frak X}}\-$morphisms $\sigma^{_{\frak X}} (\zeta)$ and $\bar\upsilon^{^{\frak X}}_{_{Q,R}}(u_\zeta)$ 
 lift $\zeta\,;$ once again,  the {\it divisibility\/} of~$\bar\M^{^{\frak X}}$ guarantees the existence and the uniqueness 
 of~$m_\zeta\in {\rm Ker}(\bar\rho^{^{\frak X}}_{_R})$ fulfilling
$$\bar\upsilon^{^{\frak X}}_{_{Q,R}}(u_\zeta) = 
\sigma^{_{\frak X}} (\zeta)\.m_\zeta
\eqno £6.22.22;$$
with the notation above, note that for any $\F^{^{\frak X}}\!\-$morphism $\varphi\,\colon R\to Q$
and any $u\in Q$ we still have
$$\bar x_{\kappa^{_{\frak X}}_{_Q}(u)\circ\varphi}
\.(\bar\ell_{\tilde\varphi})^{-1} = \bar\upsilon^{_{\frak X}}_{_Q} (u)\.\bar x_\varphi\.(\bar\ell_{\tilde\varphi})^{-1}
\eqno £6.22.23,$$
so that we may assume that $m_\zeta$ only depends on 
$\tilde\zeta\in \tilde\F_{\!P}(Q,R)$ and, as above, we write $m_{\tilde\xi}$ instead of $m_\xi\,;$ moreover, for a second 
$\F^{^\frak X}_{\! P}\-$morphism $\xi\,\colon T\to R\,,$ we get
$$\eqalign{\sigma^{_{\frak X}} (\zeta\circ\xi)\.m_{\tilde\zeta\circ\tilde\xi} &= 
\bar\upsilon^{^{\frak X}}_{_{Q,T}}(u_{\zeta\circ\xi}) = \bar\upsilon^{^{\frak X}}_{_{Q,R}} (u_\zeta)\.\bar\upsilon^{^{\frak X}}_{_{R,T}}(u_\xi)\cr
&= \sigma^{_{\frak X}} (\zeta)\.m_{\tilde\zeta}\.\sigma^{_{\frak X}} (\xi)
\.m_{\tilde\xi}\cr
&= \sigma^{_{\frak X}} (\zeta\circ\xi)\.\big(\Ker (\bar\rho^{^{\frak X}})(\tilde\xi)\big)(m_{\tilde\zeta})\.m_{\tilde\xi}\cr}
\eqno £6.22.24.$$

\smallskip
The {\it divisibility\/} of $\bar\M^{^{\frak X}}$ always forces
$$m_{\tilde\zeta\circ\tilde\xi} = \big(\Ker (\bar\rho^{^{\frak X}})
(\tilde\xi)\big)(m_{\tilde\zeta})\.m_{\tilde\xi}
\eqno £6.22.25\phantom{.}$$
and, since ${\rm Ker}(\bar\rho^{^{\frak X}}_T)$ is Abelian, in the additive notation we obtain
$$0 = \big(\Ker (\bar\rho^{^{\frak X}})(\tilde\xi)\big)(m_{\tilde\zeta}) - 
m_{\tilde\zeta\circ\tilde\xi}+ m_{\tilde\xi}
\eqno £6.22.26;$$
that is to say,  denoting by $\frak i^{_\frak X}_P\,\colon \tilde\F^{^\frak X}_{\!P}\i \tilde\F^{^\frak X}$ the obvious inclusion functor, the correspondence $m$ sending any  $\tilde\F^{^\frak X}_{\! P}\-$morphism $\tilde\zeta\,\colon R\to Q$ 
to~$m_{\tilde\zeta}$ defines a {\it $1\-$cocycle\/} in $\Bbb C^1\big(\tilde\F^{^\frak X}_{\! P},
\Ker (\bar\rho^{^{\frak X}}) \circ\frak i^{_\frak X}_P\big)\,;$ but, since the category $\tilde\F^{^\frak X}_{\! P}$
has a final object, we actually have [11,~Corollary~A4.8]
$$\Bbb H^1\big(\tilde\F^{^\frak X}_{\! P},\Ker (\bar\rho^{^{\frak X}})
\circ\frak i^{_\frak X}_P\big) = \{0\}
\eqno £6.22.27;$$
consequently, we obtain $m = d^{^\frak X}_0 (w)$ for some element 
$w = (w_Q)_{Q\in \frak X}$ in 
$$\Bbb C^0\big(\tilde\F^{^\frak X}_{\! P},\Ker (\bar\rho^{^{\frak X}})
\circ\frak i^{_\frak X}_P\big) = \Bbb C^0\big(\tilde\F^{^\frak X},
\Ker (\bar\rho^{^{\frak X}})\big)
\eqno £6.22.28.$$
In conclusion, equality~£6.22.22 becomes
$$\bar\upsilon^{^{\frak X}}_{_{Q,R}}(u_\zeta) = 
\sigma^{_{\frak X}} (\zeta)\.\big(\Ker (\bar\rho^{^{\frak X}})
(\tilde\zeta)\big)(w_Q)\.w_R^{-1} = w_Q\.\sigma^{_{\frak X}} (\zeta)\.w_R^{-1}
\eqno £6.22.29\phantom{.}$$
and therefore  the new correspondence  sending 
$\varphi\in \F (Q,R)$ to  $w_Q\.\sigma^{_{\frak X}} (\varphi)\.w_R^{-1}$ 
defines  a {\it $\F^{^{\frak X}}\-$locality functorial section\/}   
of $\bar\rho^{^{\frak X}}\,.$
\eject

\smallskip
Finally, assume that $\sigma^{_{\frak X}}$ and $\sigma'^{_{\frak X}}$ are {\it $\F^{^{\frak X}}\!\-$locality functorial sections\/}  of~$\bar\rho^{^{\frak X}}$;  for any $\F^{^{\frak X}}\-$morphism $\varphi\,\colon R\to Q\,,$  since 
the $\bar\M^{^{\frak X}}\-$morphisms $\sigma^{_{\frak X}}(\varphi)$ and $\sigma'^{_{\frak X}}(\varphi)$ have the same image in $\F (Q,R)\,,$ the {\it divisibility\/} of $\bar\M^{^{\frak X}}$ forces  the existence of a unique  $\bar\ell_\varphi\in {\rm Ker}(\bar\rho^{_{\frak X}}_R)$ fulfilling 
$\sigma'^{_{\frak X}}(\varphi) = \sigma^{_{\frak X}}(\varphi)\.\bar\ell_\varphi\,;$ actually, since $\sigma^{_{\frak X}}$ 
and $\sigma'^{_{\frak X}}$ are {\it $\F^{^{\frak X}}\!\-$locality functors\/}, it is easily checked that  $\bar\ell_\varphi$ only depends on $\tilde\varphi\in \tilde\F(Q,R)$  and, as above, we write $\bar\ell_{\tilde\varphi}$ instead of $\bar\ell_\varphi\,.$ Denote by 
$\bar\ell$ the element of $\Bbb C^1\big(\tilde\F^{^{\frak X}}\!,\Ker (\bar\rho^{^{\frak X}})\big)$ defined by the correspondence sending $\tilde\varphi$ to~$\bar\ell_{\tilde\varphi}\,;$ as above, it is easily checked that this correspondence  defines a {\it $1\-$cocycle\/} in~$\Bbb C^1\big(\tilde\F^{^\frak X}\!,\Ker (\bar\rho^{^{\frak X}}) \big)\,.$

\smallskip
Once again, it is well-known that it suffices to prove that $\bar\ell$ is a {\it $1\-$coboun-dary\/}; moreover, it follows
again from~£6.13 and from~isomorphism~£6.21.7 that actually it suffices to prove that the  image 
$\Bbb C^1 (\tilde\frak i_U,\overline{\frak z^{_{\frak X,\ad}}_U\! * \theta_U \!}\,)(\bar\ell)$ of this 
{\it $1\-$cocycle\/} in~$\Bbb C^1 \big(\tilde\N^{^{\frak N}}\!\!,\Ker (\bar\rho^{^{\frak N}})\big)$
is a {\it $1\-$coboundary\/}. First of all, consider the group homomorphisms (cf.~£6.21)
$$\eqalign{\sigma^{_{\frak X}} (U) &: \F (U)\too \bar\M^{^{\frak X}}(U) = \big(N_{\bar\M^{^{\frak X}}} (U)\big)(U)\cr \sigma'^{_{\frak X}} (U) &: \F (U) \too \bar\M^{^{\frak X}}(U) = \big(N_{\bar\M^{^{\frak X}}} (U)\big)(U)\cr
\bar\frak h_U (U) &: \big(N_{\bar\M^{^{\frak X}}} (U)\big)(U)\too \bar\M^{^{\frak N}}(U)\cr}
\eqno £6.22.30;$$
in particular, we have the compositions (cf.~£6.18.1)
$$\bar\lambda_U = \bar\frak h_U (U)\circ \sigma^{_{\frak X}} (U) \circ \pi_U\qq 
\bar\lambda'_U = \bar\frak h_U (U)\circ \sigma'^{_{\frak X}} (U) \circ \pi_U
\eqno £6.22.31\phantom{.}$$
from $L$ to $\bar\M^{^{\frak N}}(U) = \M^{^{\frak N}}(U)\big/\upsilon^{_{\frak N}}_U \big(Z(U)\big)$ (cf.~£6.12).

\smallskip
As in Proposition~£6.19 above, we claim that we can apply twice Lemma 18.8 in~[11] to the groups $L$ and~$\M^{^{\frak N}}(U)\,,$ to the quotient $\bar\M^{^{\frak N}}(U)$ of 
$\M^{^{\frak N}}(U)$ endowed with the surjective group homomorphims $\bar\lambda_U$ or $\bar\lambda'_U\,,$ and to the group homomorphism 
$$\eta_U : \tau_{_U} (N)\too \M^{^{\frak N}}(U)
\eqno £6.22.32\phantom{.}$$
mapping $\tau_{_U}(u)$ on $\upsilon^{_{\frak N}}_{_U}( u)$ for any $u\in N\,,$ which makes sense since $\tau_{_U}$ is injective; indeed, we already know that $\tau_{_U}(N)$ is a Sylow $p\-$subgroup of $L$  and it is easily checked
that in both cases the condition in [11,~Lemma~18.8] holds. Consequently, according to [11,~Lemma~18.8], 
there are two group homomorphisms
$$\lambda_U : L\too \M^{^{\frak N}}(U)\qq \lambda'_U : L\too \M^{^{\frak N}}(U)
\eqno £6.22.33\phantom{.}$$
respectively lifting $\bar\lambda_U$ and $\bar\lambda'_U\,,$ and extending $\eta_U\,;$ actually, as in~£6.20 and~£6.21
above, these homomorphisms come from two {\it $\N^{^\frak N}\!\-$locality functors\/}
$$\frak l_U : \T_L^{^{\frak N}}\too \M^{^{\frak N}}\qq \frak l'_U : \T_L^{^{\frak N}}\too \M^{^{\frak N}}
\eqno £6.22.34\phantom{.}$$
and, since we clearly have
$$\rho^{_{\frak N}}_U\circ \lambda_U = \pi_U = \rho^{_{\frak N}}_U\circ \lambda'_U
\eqno £6.22.35,$$
they are {\it naturally $\N^{^\frak N}\!\-$isomorphic\/} (cf.~£2.9).
	
\smallskip
On the other hand, it is clear that any $\N^{^{\frak N}}\-$morphism $\varphi\,\colon R\to Q$ is already induced 
by some $y\in \T_L^{^{\frak N}}(Q,R)$ and, from a {\it natural $\N^{^\frak N}\!\-$isomorphism\/} 
$z\,\colon \frak l'_U\cong \frak l_U\,,$ we get
$$\frak l'_U (y) = (z_Q)^{-1}\.\frak l_U (y)\. z_R = \frak l_U (y)\.\big(\Ker (\rho^{_\frak N})
(\tilde\varphi)\big)(z_Q)^{-1}\.z_R
\eqno £6.22.36;$$
but, denoting by $\overline{\frak l'_U (y)\!}\,,$ $\overline{\frak l_U (y)\!}\,,$ $\bar z_Q$ and $\bar z_R$
the respective images in $\bar\M^{^{\frak N}} (Q,R)\,,$ $\bar\M^{^{\frak N}} (Q,R)\,,$ $\bar\M^{^{\frak N}}(Q)$ and $\bar\M^{^{\frak N}}(R)\,,$
it follows from equalities~£6.22.31 that
$$\overline{\frak l'_U (y)\!}\, = \bar\frak h_U\big(\sigma'^{_{\frak X}} (\varphi)\big)\qq
\overline{\frak l_U (y)\!}\, = \bar\frak h_U\big(\sigma^{_{\frak X}} (\varphi)\big)
\eqno £6.22.37;$$
hence, from equality~£6.22.36, we get
$$\bar\frak h_U\big(\sigma'^{_{\frak X}}(\varphi)\big) = \bar\frak h_U\big(\sigma^{_{\frak X}}(\varphi)\big)\.
\big(\Ker (\bar\rho^{_{\frak N}})(\tilde\varphi)\big)(\bar z_Q)^{-1}\.\bar z_R 
\eqno £6.22.38,$$
which forces $\bar\frak h_U(\bar\ell_{\tilde\varphi}) = \big(\Ker (\bar\rho^{_{\frak N}})(\tilde\varphi)\big)(\bar z_Q)^{-1}\.\bar z_R \,.$
We are done.

\bigskip
 \bigskip
\noindent
{\bf £7. The perfect $\F\-$locality extending the perfect $\F^{^{\rm sc}}\-$locality }

\bigskip
£7.1\phantom{.}  Let $P$ be a finite $p\-$group and $\F$ a Frobenius $P\-$category; from section~£6 
we already know the existence and the uniqueness of a {\it perfect $\F^{^{\rm sc}}\-$locality\/}  $\P^{^{\rm sc}}\!\,;$
as a matter of fact, in [11,~Chap.~20]{\footnote{\dag}{\cds The argument in
 [11,~20.16] has been scratched; below we develop a correct argument.} we already have proved that any
 {\it perfect $\F^{^{\rm sc}}\-$locality\/} $\P^{^{\rm sc}}$ can be extended to a unique  
 {\it perfect $\F\-$locality\/}~$\P$ (cf.~£2.8); in this section, we prove
the following more precise result which actually shows the existence and the uniqueness of a  
{\it perfect $\F\-$locality\/}.

\bigskip
\noindent
{\bf Theorem~£7.2.} {\it  Any perfect $\F^{^{\rm sc}}\-$locality  $\P^{^{\rm sc}}$ can be extended to a unique 
perfect $\F\-$locality~$\P\,.$ Moreover, for any {\it $p\-$coherent $\F\-$locality\/} $\L\,,$
 any $\F^{^{\rm sc}}\-$locality functor $\frak h^{^{\rm sc}}$ from  $\P^{^{\rm sc}}$ to $\L^{^{\rm sc}}$ can be extended to a unique $\F\-$locality functor
$$\frak h : \P\too \L
\eqno £7.2.1.$$\/}

\par
 £7.3.  Let us consider a set $\frak X$ of subgroups of $P$ containing
the set of $\F\-$selfcentralizing subgroups and any subgroup $Q$ of $P$ such that $\F(Q,R)\not= \emptyset$
for some $R\in \frak X\,.$ Arguing by induction on $\vert\frak X\vert\,,$
we will construct the {\it perfect $\F^{^{\frak X}}\!\-$locality\/}  $\P^{^{\frak X}}$ extending $\P^{^{\rm sc}}$ 
 and the $\F^{^{\frak X}}\!\-$locality functor $\frak h^{^{\frak X}}\!\,\colon \P^{^{\frak X}}\!\to \L^{^{\frak X}}$
 extending $\frak h^{^{\rm sc}}\,.$ We may
assume that $\frak X$ contains subgroups of $P$ which are not $\F\-$selfcentralizing,
choose a minimal one $U$ and set 
$$\frak Y = \frak X - \{\theta (U)\mid \theta\in \F(P,U)\}
\eqno £7.3.1;$$
\eject
\noindent
from now on, we assume that there exist a perfect 
$\F^{^{\frak Y}} \-$locality~$\P^{^{\frak Y}}$ and  a functor 
$\frak h^{^{\frak Y}}\!\,\colon \P^{^{\frak Y}}\!\to \L^{^{\frak Y}}\,,$
 and we denote by 
$$\eqalign{\tau^{_{\frak Y}} : \T_{\!P}^{^{\frak Y}} \too \P^{^{\frak Y}}&\qq \pi^{_{\frak Y}} : 
\P^{^{\frak Y}}\too \F^{^{\frak Y}}\cr 
\bar\tau^{_{\frak Y}} : \T_{\!P}^{^{\frak Y}} \too \L^{^{\frak Y}}&\qq \bar\pi^{_{\frak Y}} : 
\L^{^{\frak Y}}\too \F^{^{\frak Y}}\cr}
\eqno £7.3.2\phantom{.}$$
the corresponding structural functors (cf.~£2.7.1); for any pair of subgroups
$Q$ and $R$ in $\frak Y$ we set $\P^{^{\frak X}}\!(Q,R) = \P^{^{\frak Y}}(Q,R)$ and
$\frak h_{_{Q,R}}^{^{\frak X}}\! = \frak h_{_{Q,R}}^{^{\frak Y}}\! \,,$ and if $R\i Q$ then 
we set $i_{_R}^{^Q} = \tau^{_{\frak Y}}_{_{Q,R}}(1)$ and $\bar \imath_{_R}^{^Q} 
= \bar\tau^{_{\frak Y}}_{_{Q,R}}(1)$ for short.

\medskip
£7.4. If $V\in \frak X -\frak Y$ is {\it fully centralized\/} in $\F$ then we consider $\widehat V = V\.C_P (V)$ 
which is clearly $\F\-$selfcentralizing; in particular, $\P^{^{\frak X}}\! (\widehat  V) $ has been already defined above, and the structural functor $\tau^{_{\frak Y}}\,\colon \T_{\!P}^{^{\frak Y}}\to \P^{^{\frak Y}}$  determines a group
homomorphism $\tau^{_{\frak X}}_{_{\hat  V}} \,\colon N_P (\widehat  V)\to \P^{^{\frak X}}\! (\widehat  V)\,.$
Let $Q$ be a subgroup in $\frak Y$  which contains and normalizes~$V\,;$ thus, $Q$ normalizes 
$\widehat  V$ and we set $\widehat  Q =  Q\.\widehat  V$ which  coincides with the converse image of $\F_{\! Q}(V)$ 
in~$N_P (V)\,;$ since $V$ is also fully $\F_{\! Q}(V)\-$normalized in $\F$ [11,~2.10] and we have 
$N_P^{\F_{\! Q} (V)} (V) = \widehat  Q\,,$ we get  the  Frobenius $\widehat  Q\-$category $\F^{^{V,Q}}$ [11,~Proposition~2.16]  and the associated  perfect $\F^{^{V,Q}}\-$locality [11,~£17.4 and~£17.5]
$$\F^{^{V,Q}} = N_{\F}^{\F_{Q}(V)} (V)\qq \P^{^{V,Q}} = N_{\P^{^{\frak Y}}}^{\F_{\! Q}(V)} (V)
\eqno £7.4.1$$  
defined over the set $\frak Y^{^{V,Q}}$ of elements of $\frak Y$ contained in~$\widehat  Q\,;$ 
we identify $\F^{^{V,Q}}$ and  $\P^{^{V,Q}}$ with their canonical image in
$\F^{^{\frak Y}}$ and $\P^{^{\frak Y}}$ respectively.

\medskip
£7.5. Since the {\it hyperfocal subgroup\/} $H_{\!\F^{^{V,Q}}}$ (cf.~£2.4) is a 
$\F^{^{V,Q}}\!\-$stable subgroup  of~$\widehat  Q$ [11,~Lemma~13.3], it follows from [11,~Theorem~17.18]  that we have the {\it quotient\/}  perfect 
$\F^{^{V,Q}}\!/H_{\!\F^{ {V,Q}}}\-$locality $\P^{^{V,Q}}\!/H_{\!\F^{^{V,Q}}}\,;$ but, it is easily
checked from [11,~Lemma~13.3] that $\P^{^{V,Q}}\!/H_{\!\F^{^{V,Q}}}$ can be identified to the 
{\it full\/} subcategory of~$\T_{\hat Q/H_{\!\F^{^{V,Q}}} }$ [11,~17.2] over the set of images in 
$\widehat  Q/H_{\!\F^{^{V,Q}}}$ of the subgroups in~$\frak Y^{^{V,Q}};$ hence, we have a canonical functor
$$\frak t^{^{V,Q}} : \P^{^{V,Q}}\too \T_{\hat Q/H_{\!\F^{^{V,Q}}} }
\eqno £7.5.1\phantom{.}$$
compatible with the structural functors; in particular, we have a group homomorphism
$(\frak t^{^{V,Q}} )_{_{\hat  V}}\,\colon \P^{^{\frak X}} (\widehat  V)\to \widehat  Q/H_{\!\F^{^{V,Q}}}$
and we consider its kernel
$$O (V) = {\rm Ker}\big((\frak t^{^{V,Q}} )_{_{\hat V}}\big) = {\Bbb  O}^p\Big(C_{\P^{^{\frak X}} (\hat V)}\big(\tau^{_{\frak X}}_{_{\hat V}}(V)\big)
\Big)\. \tau^{_{\frak X}}_{_{\hat V}}(H_{C_{\F} (V)})
\eqno £7.5.2;$$
similarly, in $\L^{^{\frak X}}\! (\widehat  V)$ we set
$$\bar O (V) = {\Bbb  O}^p\Big(C_{\L^{^{\frak X}} (\hat V)}\big(\bar\tau^{_{\frak X}}_{_{\hat V}}(V)\big)
\Big)\. \bar\tau^{_{\frak X}}_{_{\hat V}}(H_{C_{\F} (V)})
\eqno £7.5.3;$$
note that $\frak h_{_{\hat V}}^{^{\frak Y}}\big( O(V)\big) \i \bar O (V)$ (cf.~£7.3).
\eject

\medskip
£7.6. If $V'\in \frak X -\frak Y$ is also fully centralized  in $\F\,,$ setting 
$\widehat  V' = V'\.C_P (V')$ it follows from [11,~statement~2.10.1] that any 
$\varphi\in \F (V',V)$ can be extended to a suitable $\hat\varphi\in 
\F  (\widehat  V',\widehat  V)\,;$ thus, the restriction map 
$$\frak f_{_{V',V}}^{^{\hat V',\hat V}} :\F (\widehat  V',\widehat  V)_{V',V}\too \F(V',V)
\eqno £7.6.1\phantom{.}$$
 is surjective. In this case, $\P^{^{\frak X}}\! (\widehat  V',\widehat  V)$ has been already defined above, and the groups 
 $\P^{^{\frak X}} (\widehat  V)$ and $\P^{^{\frak X}} (\widehat  V')$ act on this set by composition on the right-hand and on the left-hand respectively; moreover, it is quite clear that the respective subgroups $O (V)$ and $O (V')$ stabilize 
 $\P^{^{\frak X}}\! (\widehat  V',\widehat  V)_{V',V}\,,$
and that the corresponding quotient sets $\P^{^{\frak X}}\! (\widehat  V',\widehat  V)_{V',V}/O (V)$ and
$O (V')\backslash \P^{^{\frak X}}\! (\widehat  V',\widehat  V)_{V',V}$ coincide with each other. Thus, we can define
$$\P^{^{\frak X}}\! (V',V) = \P^{^{\frak X}}\! (\widehat  V',\widehat  V)_{V',V}/O (V) =
 O (V')\backslash \P^{^{\frak X}}\! (\widehat  V',\widehat  V)_{V',V}
\eqno £7.6.2\phantom{.}$$
and we denote by 
$$\frak g_{_{V',V}}^{^{\hat  V',\hat  V}} : \P^{^{\frak X}}\! (\widehat  V',\widehat  V)_{V',V} \too \P^{^{\frak X}}\! (V',V)
\eqno £7.6.3\phantom{.}$$
 the canonical map; moreover, since $\L$ is {\it $p\-$coherent\/}, the  image of 
 $\bar O( V)$ in $\L(V)$ is trivial, and therefore 
 $\frak h_{_{\hat V',\hat V}}^{^{\frak X}}\! $ induces a map
 $$\frak h_{_{V',V}}^{^{\frak X}}\! : \P (V',V)\too \L (V',V)
 \eqno £7.6.4\phantom{.}$$
 fulfilling $\,\bar\frak g_{_{V',V}}^{^{\hat V',\hat V}}\!\circ \frak h_{_{\hat V',\hat V}}^{^{\frak X}}\!
 = \frak h_{_{V',V}}^{^{\frak X}}\!\circ \frak g_{_{V',V}}^{^{\hat V',\hat V}}$ where
 $$\bar\frak g_{_{V',V}}^{^{\hat V',\hat V}} : \L^{^{\frak X}}\! (\widehat  V',\widehat  V)_{V',V} \too \L^{^{\frak X}}\! (V',V)
\eqno £7.6.5\phantom{.}$$
denotes the map determined by the {\it divisibility\/} in $\L\,.$

\medskip
£7.7.  It is clear that there is a unique map 
$$\pi_{_{V',V}} : \P^{^{\frak X}}\! (V',V)\too \F (V',V)
\eqno £7.7.1\phantom{.}$$ 
such that, for any $\hat x\in \P^{^{\frak X}}\! (\widehat  V',\widehat  V)_{V',V}\,,$ we have (cf.~£7.6.1)
$$\pi_{_{V',V}}^{_{\frak X}}\! \big(\frak g_{_{V',V}}^{^{\hat V',\hat V}}(\hat x)\big) 
= \frak f_{_{V',V}}^{^{\hat V',\hat V}}\big( \pi_{_{\hat V',\hat V}}(\hat x)\big)
\eqno £7.7.2.$$
Similarly, if $u\in P$ belongs to $\T_{\! P} (V',V)$ then it belongs to $\T_{\! P} (\widehat  V',\widehat  V)$ too, and we consider the map $\tau_{_{V',V}}^{_{\frak X}}\!\,\colon \T_{\! P} (V',V)\to\P^{^{\frak X}}\! (V',V)$ defined~by
$$\tau_{_{V',V}}^{_{\frak X}} (u) = \frak g_{_{V',V}}^{^{\hat V',\hat V}}
\big(\tau_{_{\hat V',\hat V}}^{_{\frak Y}}(u)\big)
\eqno £7.7.3.$$
Moreover, it is easily checked that
$$\frak h_{_{V',V}}^{^{\frak X}}\!\circ \tau_{_{V',V}}^{_{\frak X}} = \bar\tau_{_{V',V}}^{_{\frak X}}
\qq \bar \pi_{_{V',V}}\circ \frak h_{_{V',V}}^{^{\frak X}}\! = \pi_{_{V',V}}
\eqno £7.7.4.$$
\eject

\medskip
£7.8.  Further, it is clear that the composition in $\P^{^{\frak Y}}$ defines a
compatible {\it composition\/} between those morphisms; explicitly, if
$V''\in \frak X -\frak Y$ is fully centralized in~$\F\,,$ setting $\widehat  V'' = V''\.C_P (V'')$  we get a
{\it composition map\/} 
$$\P^{^{\frak X}}\! (V'',V')\times \P^{^{\frak X}}\! (V',V)\too \P^{^{\frak X}}\! (V'',V)
\eqno £7.8.1\phantom{.}$$
 such that, for any $\hat x\in \P^{^{\frak X}}\! (\widehat  V',\widehat  V)_{V',V}$ and 
 any~$\hat x'\in \P^{^{\frak X}}\! (\widehat  V'',\widehat  V')_{V'',V'}\,,$ we have
$$\frak g_{_{V'',V}}^{^{\hat V'',\hat V}}(\hat x'\. \hat x) =
\frak g_{_{V'',V'}}^{^{\hat V'',\hat V'}}(\hat x')\. \frak g_{_{V',V}} 
^{^{\hat V',\hat V}}(\hat x)
\eqno £7.8.2\phantom{.}$$
 and the associativity of the composition in $\P^{^{\frak Y}}\!$ forces the obvious {\it associa-tivity\/} here.
Once again, it is easily checked that for any $x\in \P^{^{\frak X}}\! (V',V)$ and any $x'\in \P^{^{\frak X}}\! (V'',V')$ we have
$$\frak h_{_{V'',V}}^{^{\frak X}}\!(x'\. x) = \frak h_{_{V'',V'}}^{^{\frak X}}\!(x')\.
\frak h_{_{V',V}}^{^{\frak X}}\!(x)
\eqno £7.8.3.$$

\bigskip
\noindent
{\bf Proposition £7.9.} {\it With the notation and the hypothesis above, let $Q'$ be a subgroup $\F\-$isomorphic to $Q$ which contains and normalizes $V'\,,$ and  set $\widehat  Q' = Q'\.\hat V'\,.$  For any element $x\in \P^{^{\frak X}}\! (Q',Q)_{V',V}$ there is $\hat a\in \P^{^{\frak X}}\! (\widehat  Q',\widehat  Q)_{V',V}$ such that 
$$\hat a^{-1}\. i_{_{Q'}}^{^{\hat Q'}}\. x\in \P^{^{V,Q}} (\widehat  Q, Q)\quad and
\quad \frak t^{^{V,Q}} (\hat a^{-1}\. i_{_{Q'}}^{^{\hat Q'}}\.x) = 1
\eqno £7.9.1$$
 Moreover, the correspondence sending $x$ to
$\frak g_{_{V',V}}^{^{\hat V',\hat V}} 
\big(\frak g_{_{\hat V',\hat V}}^{^{\hat Q',\hat Q}}(\hat a)\big)$ determines a 
map 
$$\frak g_{_{V',V}}^{^{Q',Q}} : \P^{^{\frak X}}\! (Q',Q)_{V',V}\too \P^{^{\frak X}}\!(V',V)
\eqno £7.9.2.$$\/}

\par
\noindent
{\bf Proof:} Denoting by $\chi\in\F (V',V)$ the element fulfilling  $\pi_{_{Q',Q}}^{_{\frak X}}(x)\circ \iota_{_V}^{_{Q}} 
= \iota_{_{V'}}^{_{Q'}}\circ \chi\,,$ we clearly have ${}^{\chi} \F_{\!Q} (V) = \F_{\!Q'} (V')\,;$ then, since $V$ is normal 
in $\widehat  Q$ and we have $\,\F_{\!\hat Q}(V) = \F_{\!Q} (V)\,,$ it follows from [11,~statement~2.10.1] that $\chi$ can be extended to an $\F\-$morphism $\alpha\,\colon\widehat  Q\to P$ and, since $\widehat  Q'$  coincides with the converse image of
$\F_{\! Q'}(V')$ in~$N_P(V')\,,$  we clearly have  $\alpha (\widehat   Q) =\widehat  Q'\,;$ in
particular, there is $\hat a\in \P (\widehat  Q',\widehat  Q)_{R',R}$ such that 
$\hat a^{-1}\. i_{_{Q'}}^{^{\hat Q'}}\. x$ belongs to $\P^{^{V,Q}}\! (\widehat  Q,Q)\,,$ 
and actually we can modify  our choice of $\hat a$ in such a way that 
$$\frak t^{^{V,Q}} (\hat a^{-1}\. i_{_{Q'}}^{^{\hat Q'}}\. x) = 1
\eqno £7.9.3.$$

\smallskip
 Moreover, if $\hat a'\in\P (\widehat  Q',\widehat  Q)_{V',V}$ is another
choice such that 
$$\hat a'^{-1}\. i_{_{Q'}}^{^{\hat Q'}}\.x\in \P^{^{V,Q}}(\widehat  Q,Q) 
\qq \frak t^{^{V,Q}} (\hat a'^{-1}\. i_{_{Q'}}^{^{\hat Q'}}\.x) = 1
\eqno £7.9.4,$$
 then the difference $\hat a^{-1}\.\hat a'$ belongs to~$\P^{^{V,Q}}(\widehat  Q)$ and we have
$\frak t^{^{V,Q}} (\hat a^{-1}\.\hat a') = 1\,;$ but, as in~£7.5.2 above, we have
$${\rm Ker}\big((\frak t^{^{V,Q}} )_{_{\hat Q}}\big) = 
{\Bbb O}^p\Big(C_{\P^{^{\frak X}}\! (\hat Q)} \big(\tau_{_{\hat Q}}^{_{\frak X}}\!(V)\big)\Big) \.
\tau_{_{\hat Q}}^{_{\frak X}}\!(H_{C_{\F} (V)})
\eqno £7.9.5;$$
\eject
\noindent
consequently, since we have
$$\eqalign{\frak g_{_{\hat V,\hat V}}^{^{\hat Q,\hat Q}}\bigg({\Bbb  O}^p\Big(C_{\P^{^{\frak X}}\! 
(\hat Q)}\big(\tau_{_{\hat Q}}^{_{\frak X}}\!(V) \big)\Big)\bigg)
&\i {\Bbb  O}^p\Big(C_{\P^{^{\frak X}}\! (\hat V)}\big(\tau_{_{\hat V}}^{_{\frak X}}\!(V) \big)\Big)\cr 
\frak g_{_{\hat V,\hat V}}^{^{\hat Q,\hat Q}}\big(\tau_{_{\hat Q}}^{_{\frak X}}\!(H_{C_{\F} (V)})\big)
&= \tau_{_{\hat V}}^{_{\frak X}}\!(H_{C_{\F} (V)})\cr}
\eqno £7.9.6,$$
 the element  $\frak g_{_{\hat V,\hat V}}^{^{\hat Q,\hat Q}}(\hat a^{-1}\.\hat a')$
belongs to $O (V)\,.$ We are done.

\medskip
£7.10. Now, we claim that the family of maps $\frak g_{_{V',V}}^{^{Q',Q}}$
obtained in Proposition~£7.9 is compatible with the {\it composition maps\/} defined
in~£7.8 and it fulfills the corresponding {\it transitivity condition\/}.

\bigskip
\noindent
{\bf Proposition £7.11.} {\it With the notation and hypothesis above, let $Q''$ be a
subgroup of~$P$ which is $\F\-$isomorphic to $Q$ and~$Q'\,,$ and $V''\in \frak X -\frak Y$ a 
normal subgroup of~$Q''$    fully centralized in~$\F\,.$ Then, for any 
$x'\in \P (Q'',Q')_{V'',V'}$ and any $x\in \P (Q',Q)_{V',V}\,,$ we have 
$$\frak g_{_{V'',V}}^{^{Q'',Q}}(x'\.x) = \frak g_{_{V'',V'}}^{^{Q'',Q'}}
(x')\.\frak g_{_{V',V}}^{^{Q',Q}}(x)
\eqno £7.11.1.$$\/}
\par
\noindent
{\bf Proof:}  With the notation above, set $\widehat  Q'' = Q''\.\widehat  V''$ and choose an
element $\hat a'$ in $\P  (\widehat  Q'',\widehat  Q')_{V'',V'}$ such that 
$$\hat a'^{-1}\.   i_{_{Q''}}^{^{\hat  Q''}} \. x'\in \P^{^{V',Q'}} (\widehat  Q',Q')
\qq \frak t^{^{V',Q'}}  (\hat a'^{-1}\.  i_{_{Q''}}^{^{\hat Q''}}\. x') = 1
\eqno £7.11.2;$$
then, setting $\hat a'' = \hat a'\.\hat a$ and $x'' = x'\. x\,,$ we claim that
$$\hat a''^{-1}\.i_{_{Q''}}^{^{\hat  Q''}}\. x''\in\P^{^{V,Q}} (\widehat  Q,Q)
\quad{\rm and}\quad \frak t^{^{V,Q}}  (\hat a''^{-1}\.  
i_{_{Q''}}^{^{\hat Q''}}\.x'') = 1
\eqno £7.11.3.$$

\smallskip
We argue by induction on the {\it length\/} $\ell$ of $\pi_{_{\hat Q',Q'}}^{_{\frak X}}\! 
(\hat a'^{-1}\. i_{_{Q''}}^{^{\hat Q''}}\. x')$ as an {\it $\F^{^{V',Q'}}\!\-$ morphism\/} 
[11,~5.15 and~20.8.2]; if $\ell = 0$ then we have
$$\hat a'^{-1}\. i_{_{ Q''}}^{^{\hat Q''}}\. x' = n'\.i_{_{\hat Q'}}^{^{\hat Q'}}
\eqno £7.11.4\phantom{.}$$
for a suitable $n' \in \P^{^{V',Q'}} (\widehat  Q')$ fulfilling 
$\,\frak t^{^{{V',Q'}}}\! (n') = 1$ since we already know that
$\frak t^{^{{V',Q'}}}\! (i_{_{\hat Q'}}^{^{\hat Q'}}) = 1\,,$ and therefore we easily get
$$\hat a''^{-1}\.i_{_{ Q''}}^{^{\hat Q''}}\. x'' = 
(\hat a^{-1}\. n'\.\hat a)\.(\hat a^{-1}\.i_{_{Q'}}^{^{\hat Q'}}\. x)
\eqno £7.11.5;$$
since $\hat a^{-1}\. n'\.\hat a$ belongs to $\P^{^{V,Q}} (\widehat  Q)$ and we have $\frak t^{^{V,Q}} (\hat a^{-1}\. n'\.\hat a) = 1\,,$ in this case we are done.
\eject

\smallskip
 If $\ell\ge 1$ then we have  
$$\hat a'^{-1}\. i_{_{ Q''}}^{^{\hat Q''}}\. x' = i_{_{E'}}^{^{\hat Q'}}\. y'\. v'
\eqno £7.11.6\phantom{.}$$
 for some {\it $\F^{^{V',Q'}}\!\!\-$essential\/} subgroup~$E'$ [11,~5.18] in $\hat Q'\,,$ some $p'\-$element $ y'$ 
 of the converse image of~$X_{\F^{^{V',Q'}}}(E')$ [11,~Corollary~5.13] in $\P^{^{V',Q'}}\! (E')$ and some element
$v'$ of $\P^{^{V',Q'}} \!(E',Q')$ in such a way that $\pi_{_{\hat Q', Q'}}^{_{\frak X}}\! 
(i_{_{E'}}^{^{\hat Q'}}\. v')$ has length $\,\ell -1$ [11,~5.15]. Note that we have 
$\pi_{v'}^{_{\frak X}}\!(V') = V'$ and that, setting 
$$Q''' = \pi_{v'}^{_{\frak X}}\!(Q')\i E'\i \widehat  Q'
\eqno £7.11.7$$
and denoting by $v'_*$ and  $y'_*$ the respective elements of $\P^{^{V',Q'}} \!(Q''',Q')$
and $\P^{^{V',Q'}} \!(Q',Q''')$ determined by $v'$ and $y'\,,$ we get
$$i_{_{E'}}^{^{\hat Q'}}\. v' = i_{_{Q'''}}^{^{\hat Q'}}\.v'_*\qq
\hat a'^{-1}\. i_{_{ Q''}}^{^{\hat Q''}}\. x' =  i_{_{Q'}}^{^{\hat Q'}}\.y'_*\.v'_*
\eqno £7.11.8;$$
moreover, since $y'$ is a $p'\-$element, we have $\frak t^{^{V',Q'}}(y') = 1$ which implies that
$\frak t^{^{V',Q'}}(y'_*) = 1$ and therefore, since $\frak t^{^{V',Q'}}  (\hat a'^{-1}\.  
i_{_{Q''}}^{^{\hat Q''}}\. x') = 1\,,$ we successively obtain $\frak t^{^{V',Q'}}(v'_*) = 1$
and $\frak t^{^{V',Q'}}(i_{_{E'}}^{^{\hat Q'}}\. v' ) = 1\,.$

\smallskip
Thus, by the induction hypothesis, we already know that 
$$\hat a^{-1}\.(i_{_{E'}}^{^{\hat Q'}}\. v')\. x\in \P^{^{V,Q}} (\widehat  Q, Q)
\qq \frak t^{^{V,Q}} \big(\hat a^{-1}\.(i_{_{E'}}^{^{\hat Q'}}\. v')\.  x\big) = 1
\eqno £7.11.9\phantom{.}$$
and therefore, setting $E = (\pi_{\hat a}^{_{\frak X}}\!)^{-1} (E')$ and denoting by 
$b\in \P^{^{\frak X}}\! (E',E)$ the element such that 
$\hat a^{-1}\.  i_{_{E'}}^{^{\hat Q'}} = i_{_{E}}^{^{\hat Q}}\.b^{-1}\,,$ the
divisibility in $\P^{^{V,Q}}$ implies that $ b^{-1}\.  v'\. x$ belongs to $\P^{^{V,Q}} (E,Q)$ and then
 we still have $\frak t^{^{V,Q}} \big( b^{-1}\. v'\.  x\big) = 1\,;$
consequently, we still get
$$\hat a''^{-1}\.i_{_{Q''}}^{^{\hat Q''}}\. x'' = \hat a^{-1}\.( i_{_{E'}}^{^{\hat Q'}}\. y'\. v')\. x = 
i_{_{E}}^{^{\hat Q}}\.( b^{-1} \. y'\. b)\.( b^{-1}\.  v' \. x)
\eqno £7.11.10\phantom{.}$$
and, since $ b^{-1}\. y\. b$ is a $p'\-$element of $\P^{^{V,Q}}\! (E)\,,$ we
have  $\frak t^{^{V,Q}}\! ( b^{-1}\. y'\. b) = 1\,,$ which proves our claim.

\smallskip
Now, according to~£7.8 and to~Proposition~£7.9, we have
$$\eqalign{\frak g_{_{V'',V}}^{^{Q'',Q}}(x'') & = \frak g_{_{V'',V}}^{^{\hat V'',\hat V}} 
\big(\frak g_{_{\hat V'',\hat V}}^{^{\hat Q'',\hat Q}}(\hat a'')\big)= \frak g_{_{V'',V}}^{^{\hat V'',\hat V}} 
\big(\frak g_{_{\hat V'',\hat V'}}^{^{\hat Q'',\hat Q'}}(\hat a')\.\frak g_{_{\hat
V',\hat V}}^{^{\hat Q',\hat Q}}(\hat a)\big)\cr
&= \frak g_{_{V'',V'}}^{^{\hat V'',\hat V'}} 
\big(\frak g_{_{\hat V'',\hat V'}}^{^{\hat Q'',\hat Q'}}(\hat a')\big)\.
\frak g_{_{V',V}}^{^{\hat V',\hat V}} \big(\frak
g_{_{\hat V',\hat V}}^{^{\hat Q',\hat Q}}(\hat a)\big)\cr
&= \frak g_{_{V'',V'}}^{^{Q'',Q'}}
(x')\.\frak g_{_{V',V}}^{^{Q',Q}}(x)\cr} 
\eqno £7.11.11.$$
We are done.
\eject

\bigskip
\noindent
{\bf Corollary~£7.12.} {\it With the notation and the hypothesis above, we have a group homomorphism
$$\frak g_{_{V}}^{^{Q}} : \P^{^{\frak X}}\! (Q)_{V}\too \P^{^{\frak X}}\!(V)
\eqno £7.12.1.$$
fulfilling $\,\Bbb O^p \big({\rm Ker}(\frak g_{_{V}}^{^{Q}})\big) = 
\Bbb O^p\Big(C_{\P^{^{\frak X}}\! (Q)}\big(\tau_{_Q}^{_{\frak X}}
\!(V)\big)\Big)\,,$ and moreover have the commutative diagram
$$\matrix{\P^{^{\frak X}}\! (Q',Q)_{V',V}&\buildrel \frak g_{_{V',V}}^{^{Q',Q}} \over\too 
&\P^{^{\frak X}}\!(V',V) \cr
\cr
\hskip-40pt{\scriptstyle \frak h_{_{Q',Q}}^{^{\frak X}}\!}\hskip5pt\big\downarrow&
&\big\downarrow\hskip5pt{\scriptstyle \frak h_{_{V',V}}^{^{\frak X}}\!}\hskip-30pt\cr
\L^{^{\frak X}}\! (Q',Q)_{V',V}&\buildrel \bar\frak g_{_{V',V}}^{^{Q',Q}} \over\too 
&\L^{^{\frak X}}\!(V',V) \cr}
\eqno £7.12.2.$$\/}

\par
\noindent
{\bf Proof:} The first statement is easily checked from proposition~£7.11 and from the following exact sequence (cf.~£2.12.2)
$$1\too H_{C_\F (V)}\too C_P(V)\too \P(V)\too \F (V)\too 1
\eqno £7.12.3.$$
Moreover, for any $x\in \P^{^{\frak X}}\! (Q',Q)_{V',V}$ and any $\hat a\in \P^{^{\frak X}}\! 
(\widehat  Q',\widehat  Q)_{V',V}$ fulfilling condition~£7.9.1 above, it follows from the very definition of the  {\it perfect $\F^{^{V,Q}}\!/H_{\!\F^{ {V,Q}}}\-$locality\/} 
$\P^{^{V,Q}}\!/H_{\!\F^{^{V,Q}}}$ that the equality 
$\frak t^{^{V,Q}} (\hat a^{-1}\. i_{_{Q'}}^{^{\hat Q'}}\.x) = 1$ implies that 
(cf.~£7.5.1)
 $$\hat a^{-1}\. i_{_{Q'}}^{^{\hat Q'}}\.x \in  i_{_{Q}}^{^{\hat Q}}\.
 {\rm Ker}\big((\frak t^{^{V,Q}} )_{_{ Q}}\big) 
 \eqno £7.12.4\phantom{.}$$
and, since $\L$ is {\it $p\-$coherent\/}, it is easily checked that 
$$\bar\frak g_{_V}^{^Q}\bigg(\frak h_{_Q}^{^{\frak X}}\!
\Big( {\rm Ker}\big((\frak t^{^{V,Q}} )_{_{ Q}}\big) \Big)\bigg) = \{\bar\imath_{_V}^{^V}\}
\eqno £7.12.5;$$
hence, we get
$$\eqalign{\bar\imath_{_{Q'}}^{^{\hat Q'}}\.\frak h_{_{Q',Q}}^{^{\frak X}}\!(x)
&= \frak h_{_{\hat Q',Q}}^{^{\frak X}}\!(i_{_{Q'}}^{^{\hat Q'}}\.x) 
= \frak h_{_{\hat Q',Q}}^{^{\frak X}}\!\big(\hat a\.(\hat a^{-1}\. i_{_{Q'}}^{^{\hat Q'}}\.x)\big) \cr
&= \frak h_{_{\hat Q',\hat Q}}^{^{\frak X}}\!(\hat a)\. \frak h_{_{\hat Q, Q}}^{^{\frak X}}\!(\hat a^{-1}\. i_{_{Q'}}^{^{\hat Q'}}\.x) = \frak h_{_{\hat Q',\hat Q}}^{^{\frak X}}\!(\hat a)\.\bar\imath_{_Q}^{^{\hat Q}}\cr
&=\bar\imath_{_{Q'}}^{^{\hat Q'}}\.\bar\frak g_{_{Q',Q}}^{^{\hat Q',\hat Q}}\big(\frak h_{_{\hat Q',\hat Q}}^{^{\frak X}}\!(\hat a)\big) \cr}
\eqno £7.12.6\phantom{.}$$
and therefore  from Proposition~£7.9 we still get
$$\eqalign{\bar\frak g_{_{V',V}}^{^{Q',Q}}\big(\frak h_{_{Q',Q}}^{^{\frak X}}\!(x)\big) &= 
\bar\frak g_{_{V',V}}^{^{\hat Q',\hat Q}}\big(\frak h_{_{\hat Q',\hat Q}}^{^{\frak X}}\!(\hat a)\big)
= \bar\frak g_{_{V',V}}^{^{\hat V',\hat V}}\Big(\bar\frak g_{_{\hat V',\hat V}}^{^{\hat Q',\hat Q}}\big(\frak h_{_{\hat Q',\hat Q}}^{^{\frak X}}\!(\hat a)\big)\Big)\cr
&= \bar\frak g_{_{V',V}}^{^{\hat V',\hat V}}
\Big(\frak h_{_{\hat V',\hat V}}^{^{\frak X}}\!
\big(\frak g_{_{\hat V',\hat V}}^{^{\hat Q',\hat Q}}(\hat a)\big)\Big)\cr
&=\frak h_{_{ V', V}}^{^{\frak X}}\!\Big(\frak g_{_{V',V}}^{^{\hat V',\hat V}}
\big(\frak g_{_{\hat V',\hat V}}^{^{\hat Q',\hat Q}}(\hat a)\big)\Big) =  
\frak h_{_{ V', V}}^{^{\frak X}}\!\big(\frak g_{_{V',V}}^{^{Q',Q}}(x)\big)\cr}
\eqno £7.12.7.$$
We are done.
\eject

\bigskip
\noindent
{\bf Proposition £7.13.} {\it With the notation and the hypothesis above, let  $\,R$
be a  subgroup of~$Q$~which contains~$V\,.$ Then,
for any  $x\in \P (Q',Q)_{V',V}\,,$ setting $R' = \pi_x^{_{\frak X}}\! (R)$ we have
$$\frak g_{_{V',V}}^{^{R',R}} \big(\frak g_{_{R',R}}^{^{Q',Q}}(x)\big)
= \frak g_{_{V',V}}^{^{Q',Q}}(x)
\eqno £7.13.1.$$ \/} 
\par
\noindent
{\bf Proof:} As above we set
$$\widehat  R  = R\.\widehat  V\i \widehat  Q \qq \widehat  R'  = R'\.\widehat  V' \i \widehat  Q'
\eqno £7.13.2\phantom{.}$$
and,  denoting by $\chi\in\F (V',V)$ the element fulfilling $\pi_{_{Q',Q}}^{^{\frak X}}(x)\circ \iota_{_V}^{^{Q}} 
= \iota_{_{V'}}^{^{Q'}}\circ \chi\,,$ we clearly have ${}^{\chi} \F_{\!R} (V) = \F_{\!R'} (V')\,;$ as above, considering an 
$\F\-$morphism $\alpha\,\colon\widehat  Q\to P$ extending $\chi\,,$  since $\widehat  R'$  coincides with the converse image of $\F_{\! R'}(V')$ in~$N_P(V')$  we clearly have  $\alpha (\widehat   R) =\widehat  R'\,.$ Moreover, up to suitable identifications, it is quite clear that (cf.~£7.4.1)
$$\F^{^{V,R}}\i \F^{^{V,Q}}\qq \P^{^{V,R}}\i \P^{^{V,Q}}
\eqno £7.13.3.$$
Consequently, for a choice of $\hat a\in \P^{^{\frak X}}\! (\widehat  Q',\widehat  Q)_{V',V}$ such that 
$$\hat a^{-1}\. i_{_{Q'}}^{^{\hat Q'}}\. x\in \P^{^{V,Q}} (\widehat  Q, Q)\qq 
\frak t^{^{V,Q}} (\hat a^{-1}\. i_{_{Q'}}^{^{\hat Q'}}\.x) = 1
\eqno £7.13.4,$$
the element $\frak g_{_{\hat R',\hat R}}^{^{\hat Q',\hat Q}}(\hat a)$ in $\P^{^{\frak X}}\!  (\widehat  R',\widehat  R)_{V',V}$ fulfills the analogous conditions with respect to $\frak g_{_{R',R}}^{^{Q',Q}}(x)\,,$
which proves the proposition.

\medskip 
£7.14. We are ready to define the set $\P^{^{\frak X}}\! (V',V)$ for any pair of
 subgroups $V$ and $V'$ in $\frak X -\frak Y\,;$ we clearly have $N = N_P (V)\not= V$ and it
follows from [11,~Proposition~2.7] that there is an $\F\-$morphism $\nu\, \,\colon N\to P$ such that  
$\nu (V)$ is fully centralized in~$\F\,;$ moreover, we choose 
$n\in \P^{^{\frak X}}\! \big(\nu (N),N\big)$
lifting~$\nu\,.$ That is to say, we may assume that
\smallskip
\noindent
£7.14.1\quad {\it There is a pair $(N,n)$ formed by a subgroup $N$ of $P$ which
strictly contains and normalizes $V\,,$ and by an element $n$ in 
$\P^{^{\frak X}}\!\big(\nu (N),N\big)$ lifting~$\nu$ for a suitable $\F\-$morphism $\nu\,\colon N\to P$ such that $\nu (V)$ is fully centralized in~$\F\,.$\/}
\smallskip
\noindent
 We denote by $\frak N (V)$ the set of such pairs and often we write $n$
instead of~$(N,n)\,,$ setting ${}^n \!N = \nu (N)$ and $\pi_n^{_{\frak X}}\! = 
\pi_{_{\nu (N),N}}^{_{\frak X}}\!(n)\,.$

\medskip
£7.15. For another pair $(\widehat  N,\hat n)$ in $\frak N (V)\,,$  denoting by $\hat\nu\,\colon \widehat  N\to P$  
the  $\F\-$morphism determined by $\hat n\,,$ setting $M = \langle N,\widehat  N\rangle$ and considering a new $\F\-$morphism 
$\mu\, \colon M\to P$  such that $\mu (V)$ is fully centralized in~$\F\,,$  we can obtain a third pair 
$(M,m)$ in $\frak N (V)\,;$ then,  $\frak g_{_{{}^m\! N, N}}^{^{{}^m M,M}}(m)\.n^{-1}$ 
and $\frak g_{_{{}^m\!\hat N,\hat N}}^{^{{}^m M,M}}(m)\.\hat n^{-1}$
respectively belong to $\P^{^{\frak X}}\! ({}^m\! N,{}^n\! N)$
and to $\P^{^{\frak X}}\!({}^m\!\widehat  N,{}^{\hat n}\! \widehat  N\big)\,;$
in particular, since ${}^n  V\,,$ ${}^{\hat n}  V$ and ${}^mV$ are fully
centralized in~$\F\,,$ the sets $\P^{^{\frak X}}\!({}^mV,{}^n  V) \,,$ 
$\P^{^{\frak X}}\! ({}^mV,{}^{\hat n}  V) $ and $\P^{^{\frak X}}\!  ({}^{\hat n}  V,{}^n  V)$ 
have been already defined above, and we consider the element
$$g_{\hat n,n} = \frak g_{_{{}^m V,{}^{\hat n}  V}}^{^{{}^m \! \hat N,{}^{\hat n}\! \hat N}}
\big(\frak g_{_{{}^m\!\hat N,\hat N}}^{^{{}^m M,M}}(m)\.\hat n^{-1}\big)^{-1}
\.\frak g_{_{{}^m V,{}^{n} V}}^{^{{}^m\! N,{}^{n}\! N}}
\big(\frak g_{_{{}^m\! N, N}}^{^{{}^m M,M}}(m)\.n^{-1}\big) 
\eqno £7.15.1\phantom{.}$$
in $\P^{^{\frak X}}\! ({}^{\hat n}  V,{}^n  V)\,,$ which actually does not depend on the choice of $m\,.$

\medskip
£7.16. Indeed, for another pair $(M,m')$ in $\frak N(V)$ we have
$$\eqalign{\frak g_{_{{}^{m'}\! N, N}}^{^{{}^{m'}\! M,M}}(m') &= 
\frak g_{_{{}^{m'}\! N, {}^m\! N}}^{^{{}^{m'}\! M,{}^m\!M}}(m'\.m^{-1})
\.\frak g_{_{{}^{m}\! N, N}}^{^{{}^{m}\! M,M}}(m)\cr
\frak g_{_{{}^{m'}\!\hat  N,\bar  N}}^{^{{}^{m'}\! M,M}}(m') &= 
\frak g_{_{{}^{m'}\!\hat  N, {}^m\!\hat  N}}^{^{{}^{m'}\! M,{}^m\!M}}(m'\.m^{-1})
\.\frak g_{_{{}^{m}\!\hat  N,\bar  N}}^{^{{}^{m}\! M,M}}(m)\cr}
\eqno £7.16.1\phantom{.}$$
and therefore it follows from Proposition~£7.13 that we get
$$\eqalign{\frak g_{_{{}^{m'} V,{}^{ n}  V}}^{^{{}^{m'} \!  N,{}^{ n}\!  N}}
&\big(\frak g_{_{{}^{m'}\! N, N}}^{^{{}^{m'} M,M}}(m')\. n^{-1}\big) \cr
&= \frak g_{_{{}^{m'} V,{}^{ n}  V}}^{^{{}^{m'} \!  N,{}^{ n}\!  N}}
\big(\frak g_{_{{}^{m'}\!  N, {}^m\!  N}}^{^{{}^{m'}\! M,{}^m\!M}}(m'\.m^{-1})
\.\frak g_{_{{}^{m}\!  N,  N}}^{^{{}^{m}\! M,M}}(m)\. n^{-1}\big) \cr
&= \frak g_{_{{}^{m'} V,{}^m  V}}^{^{{}^{m'}\! M,{}^m\!M}}(m'\.m^{-1})
\. \frak g_{_{{}^{m'} V,{}^{ n}  V}}^{^{{}^{m'} \!  N,{}^{ n}\!  N}}
\big(\frak g_{_{{}^{m}\!  N,  N}}^{^{{}^{m}\! M,M}}(m)\. n^{-1}\big)\cr
\frak g_{_{{}^{m'} V,{}^{\hat n}  V}}^{^{{}^{m'} \! \hat N,{}^{\hat n}\! \hat N}}
&\big(\frak g_{_{{}^{m'}\!\hat N,\bar N}}^{^{{}^{m'} M,M}}(m')\.\hat n^{-1}\big) \cr
&= \frak g_{_{{}^{m'} V,{}^{\hat n}  V}}^{^{{}^{m'} \! \hat N,{}^{\hat n}\! \hat N}}
\big(\frak g_{_{{}^{m'}\!\bar  N, {}^m\!\hat  N}}^{^{{}^{m'}\! M,{}^m\!M}}(m'\.m^{-1})
\.\frak g_{_{{}^{m}\!\hat  N,\hat  N}}^{^{{}^{m}\! M,M}}(m)\.\hat n^{-1}\big) \cr
&= \frak g_{_{{}^{m'} V,{}^m  V}}^{^{{}^{m'}\! M,{}^m\!M}}(m'\.m^{-1})
\. \frak g_{_{{}^{m'} V,{}^{\hat n}  V}}^{^{{}^{m'} \! \hat N,{}^{\hat n}\! \bar N}}
\big(\frak g_{_{{}^{m}\!\hat  N,\hat  N}}^{^{{}^{m}\! M,M}}(m)\.\hat n^{-1}\big)\cr}
\eqno £7.16.2,$$
which proves our claim. Similarly, for any triple of pairs $(N,n)\,,$ $(\widehat  N,\hat n)$ and 
$(\skew3\widehat {\widehat  N},\skew1\hat{\hat n})$  in $\frak N (V)\,,$ considering a pair 
$(\langle N, \widehat  N, \skew3\widehat {\widehat  N}\rangle,m)$ in $\frak N (V)\,,$ it follows from Propositions~£7.11 and~£7.13 that
$$g_{\skew1\hat{\hat n},\hat n}\. g_{\hat n,n}  = g_{\skew1\hat{\hat n},n}
\eqno £7.16.3.$$ 
Note that if $V$ is fully centralized in $\F$ then $N = N_P (V)$ is 
$\F\-$selfcentralizing, so that it is fully centralized too, and therefore the pair
 $(N,i_{_N}^{^N})$ belongs to~$\frak N (V)\,.$

\medskip
£7.17.  Then, for any pair of subgroups $V$ and $V'$ in $\frak X -\frak Y\,,$  since for any $(N,n)\in \frak N (V)$ and 
any $(N',n')\in \frak N (V')$ the set $\P^{^{\frak X}}\!  ({}^{n'}\! V',{}^{n} V)$ is already defined,
we denote by~$\P^{^{\frak X}}\!  (V',V)$ the subset of the product
$$ \prod_{n\in\frak N (V)}\,\prod_{n'\in\frak N (V')} \P^{^{\frak X}}\!  ({}^{n'}\! V',{}^{n} V)
\eqno £7.17.1\phantom{.}$$
formed by the families $\{x_{n',n}\}_{n\in\frak N (V),n'\in \frak N (V')}$ fulfilling
$$g_{\hat n',n'}\.x_{n',n} = x_{\hat n',\hat n}\. g_{\hat n,n} 
\eqno £7.17.2.$$ 
In other words, the set $\P^{^{\frak X}}\!  (V',V)$ is the {\it inverse limit\/} of the family formed by the sets 
$\P^{^{\frak X}}\! \big({}^{n'}\! V', {}^{n} V\big)$ and by the bijections between them induced by the 
$\P^{^{\frak X}}\!\-$morphisms $g_{\hat n,n}$ and~$g_{\hat n',n'}\,.$

\medskip
£7.18. Note that, according to equalities~£7.16.3, the {\it projection map\/} onto the factor labeled by the pair 
$\big((N,n),(N',n')\big)$ induces a bijection 
$$\frak n_{n',n} : \P^{^{\frak X}}\! (V',V)\cong \P^{^{\frak X}}\!\big({}^{n'}\! V',{}^{n} V\big)
\eqno £7.18.1;$$ 
in particular, if $V$ and $V'$ are fully centralized in $\F\,,$ setting $N = N_P (V)$ and 
$N' = N_P (V')\,,$  the pairs $(N,i_{_N}^{^N})$ and $(N',i_{_{N'}}^{^{N'}})$ respectively belong
 to~$\frak N (V)$ and to $\frak N (V')\,,$ and therefore we have a {\it canonical\/} bijection
 $$\frak n_{i_{_{N'}}^{^{N'}},i_{_N}^{^N}} : \P^{^{\frak X}}\! (V',V)\cong \P^{^{\frak X}}\!\big(\,{}^{i_{_{N'}}^{^{N'}}}\! V',{}^{i_{_N}^{^N}} V\big)
\eqno £7.18.2,$$ 
so that our notation is coherent. At this point, since the map $\frak h_{{}^{n'}\! V',{}^{n} V}$
is already defined, we can define a map
$$\frak h_{_{V',V}}^{^{\frak X}}\! : \P^{^{\frak X}}\! (V',V)\too \L^{^{\frak X}}\! (V',V)
\eqno £7.18.3\phantom{.}$$
sending $x\in \P^{^{\frak X}}\! (V',V)$ to
$$\frak h_{_{V',V}}^{^{\frak X}}\! (x)= \bar\frak g_{_{{}^{n'}\!V',V'}}^{^{{}^{n'}\!N',N'}}(\bar n')^{-1}\. 
\frak h_{{}^{n'}\! V',{}^{n} V}^{^{\frak X}}\!\big(\frak n_{n',n}(x)\big)\.
\bar\frak g_{_{{}^{n}\!V,V}}^{^{{}^{n}\!N,N}}(\bar n)
\eqno £7.18.4\phantom{.}$$
where we are setting $\bar n = \frak h_{_{{}^{n}\!N,N}}^{^{\frak X}}\!(n)$ and 
$\bar n' = \frak h_{_{{}^{n'}\!N',N'}}^{^{\frak X}}\!(n')\,.$ From Corollary~£7.12 and Proposition~£7.13
it is not difficult to check that this map does not depend on the choice of the pairs $(N,n)$ and $(N',n')\,.$

\medskip
£7.19. Moreover, we have a map
$$\pi_{_{V',V}}^{_{\frak X}}\! : \P^{^{\frak X}}\! (V',V)\too \F (V',V)
\eqno £7.19.1\phantom{.}$$
sending $x\in \P^{^{\frak X}}\! (V',V)$ to (cf.~£7.7) 
$$\pi_{_{V',V}}^{_{\frak X}}\! (x) = \frak f^{^{{}^{n'}\!\!N',N'}}_{_{{}^{n'}\! V',V'}}
(\pi_{n'}^{_{\frak X}}\!)^{-1} \circ \pi_{_{{}^{n'}\! V',{}^{n} V }}^{_{\frak X}}\! \big(\frak n_{n',n}(x) \big) \circ \frak f^{^{{}^{n}\! N,N}}_{_{{}^{n} V,V}}(\pi_n^{_{\frak X}}\!)
\eqno £7.19.2;$$
then, from condition~£7.17.2, it is not difficult to prove that this map does not depend
on the choice of the pairs $(N,n)$ and $(N',n')\,.$ Similarly, if $u$ belongs to $\T_{\! P}
(V',V)$ then we may assume that it belongs to $\T_{\! P} (N',N)$ too, and we\break
\eject
\noindent
 consider the map 
$\tau_{_{V',V}}^{_{\frak X}}\!\,\colon \T_{\! P} (V',V)\to\P^{^{\frak X}}\! (V',V)$ determined~by
$$\frak n_{n',n}\big(\tau_{_{V',V}}^{_{\frak X}} (u)\big) = 
\frak g_{_{{}^{n'}\!V',{}^n V}}^{^{{}^{n'}\!\!N',{}^n\! N}}
\big(n'\.\tau_{_{N',\hat N}}^{_{\frak Y}}(u)\.n^{-1}\big)
\eqno £7.19.3.$$
Now, it is easy to check that
$$\frak h_{_{V',V}}^{^{\frak X}}\!\circ \tau_{_{V',V}}^{_{\frak X}}\! = \bar\tau_{_{V',V}}^{_{\frak X}}\!
\qq \bar \pi_{_{V',V}}^{_{\frak X}}\!\circ \frak h_{_{V',V}}^{^{\frak X}}\! = \pi_{_{V',V}}^{_{\frak X}}\!
\eqno £7.19.4.$$

\medskip
£7.20. On the other hand, for any $V''\in \frak X -\frak Y\,,$ the {\it composition map\/} in~£7.8 can be extended to a new {\it composition map\/}
$$\P^{^{\frak X}}\! (V'',V')\times \P^{^{\frak X}}\! (V',V)\too \P^{^{\frak X}}\! (V'',V)
\eqno £7.20.1\phantom{.}$$
sending $(x',x)\in \P^{^{\frak X}}\! (V'',V')\times \P^{^{\frak X}}\! (V',V)$ to 
$$ x'\.x = (\frak n_{n'',n})^{-1}\big(\frak n_{n'',n'}(x')\.\frak n_{n',n}(x)\big)
\eqno £7.20.2\phantom{.}$$
for a choice of $(N,n)$ in $\frak N (V)\,,$ of $(N',n')$ in $\frak N (V')$ and of $(N'',n'')$ in
$\frak N (V'')\,.$ 
This  {\it composition map\/} does not depend on our choice; indeed, for another choice of pairs 
$(\widehat  N,\hat n)\in \frak N (V)\,,$ $(\widehat  N',\hat n')\in\frak N (V')$ and $(\widehat  N'',\hat n'')\in \frak N (V'')\,,$
we get (cf.~£7.17.2)
$$\eqalign{g_{n'',{\hat n}''} &\.\big(\frak n_{\hat n'',\hat n'}(x')\.
\frak n_{\hat n',\hat n} (x)\big) = \frak n_{n'',n'}(x')\.g_{n',\hat n}\.\frak n_{\hat n',\hat n}(x)\cr 
& = \frak n_{n'',n'}(x')\.\frak n_{n',n}(x)\.g_{n,\hat n} = \frak n_{n'',n}(x'\.x)\.g_{n,\hat n}\cr}
\eqno £7.20.3.$$
Moreover, it is compatible with the maps $\frak h_{_{V',V}}^{^{\frak X}}\!$ defined in equality~£7.18.4 above since, setting $\bar m = \bar\frak g_{_{{}^{n}\!V,V}}^{^{{}^{n}\!N,N}}
\big(\frak h_{_{{}^{n}\!N,N}}^{^{\frak X}}\!(n)\big)\,,$ 
$\bar m' =\bar\frak g_{_{{}^{n'}\!V',V'}}^{^{{}^{n'}\!N',N'}}
\big(\frak h_{_{{}^{n'}\!N',N'}}^{^{\frak X}}\!(n')\big) $ and 
$\bar m'' = \bar\frak g_{_{{}^{n''}\!V'',V''}}^{^{{}^{n''}\!N'',N''}}
\big(\frak h_{_{{}^{n''}\!N'',N''}}^{^{\frak X}}\!(n'')\big)\,,$ we get
$$\eqalign{&\frak h_{_{V'',V}}^{^{\frak X}}\! (x'\.x) = 
\bar m''^{-1} \. \frak h_{{}^{n''}\! V'',{}^{n} V}^{^{\frak X}}\!\big(\frak n_{n'',n}(x'\.x)\big)\.\bar m\cr
&=\bar m''^{-1} \. \frak h_{{}^{n''}\! V'',{}^{n} V}^{^{\frak X}}
\!\big(\frak n_{n'',n'}(x')\.\frak n_{n',n}(x)\big)\.\bar m\cr
&=\big(\bar m''^{-1} \. 
\frak h_{{}^{n''}\! V'',{}^{n'} V'}^{^{\frak X}}\!\big(\frak n_{n'',n'}(x')\big)\.\bar m'\big)\.\big(\bar m'^{-1}
\.\frak h_{{}^{n'}\! V',{}^{n} V}^{^{\frak X}}
\!\big(\frak n_{n',n}(x)\big)\.\bar m\big)\cr
&= \frak h_{_{V'',V'}}^{^{\frak X}}\! (x')\. \frak h_{_{V',V}}^{^{\frak X}}\! (x)\cr}
\eqno £7.20.4.$$
Finally, for any $V'''\in \frak X -\frak Y$ and any $x''\in \P^{^{\frak X}}\! (V''',V'')\,,$ it is quite clear that
$$(x''\.x')\. x = x''\.(x'\.x)
\eqno £7.20.5.$$

\medskip
£7.21. We are ready to complete our construction of the announced
{\it perfect $\F^{^{\frak X}}\!\-$locality\/} $\P^{^{\frak X}}\!$ and the {\it $\F^{^{\frak X}}\!\-$locality\/} functor 
$\frak h^{^{\frak X}}\!\,\colon \P^{^{\frak X}}\!\to \L^{^{\frak X}}\!\,;$
for any subgroups $V$ in $\frak X -\frak Y$ and $Q$ in $\frak Y$ we define
$$\P^{^{\frak X}}\! (V,Q) = \emptyset\qq 
\P^{^{\frak X}}\! (Q,V) = \bigsqcup_{V'} \P^{^{\frak X}}\!(V',V)
\eqno £7.21.1\phantom{.}$$
where $V'$ runs over the set of subgroups $V'\in \frak X -\frak Y$ contained in $Q\,,$
and the map 
$$\frak h_{_{Q,V}}^{^{\frak X}}\! : \P^{^{\frak X}}\! (Q,V)\too \L^{^{\frak X}}\! (Q,V)
\eqno £7.21.2$$
sends $x\in \P^{^{\frak X}}\!(V',V)\i \P^{^{\frak X}}\! (Q,V)$ to 
$\bar\imath_{_{V'}}^{^Q}\.\frak h_{_{V',V}}^{^{\frak X}}(x)\,.$
In order to define the composition of two  $\P^{^{\frak X}}\!\-$morphisms 
$x\,\colon R\to Q$ and $y\,\colon T\to R$ we already may assume that $T$ does not belong to $\frak Y\,;$ if $Q$ does not belong to~$\frak Y$ then the composition 
$x\.y$ is given by the map~£7.20.1
which is compatible with the maps $\frak h_{_{V',V}}^{^{\frak X}}\!$ defined above.
If $Q\in\frak Y$ but $R$ does not belong to~$\frak Y$ then, setting 
$R' = \big(\pi_{_{Q,R}}^{_{\frak X}}\!(x)\big)(R)\,,$ $x$ is actually an element 
of~$\P^{^{\frak X}}\!(R',R)$ and, according to 
definition~£7.21.1, the element $x\. y$ defined by the map~£7.20.1 belongs 
to~$\P^{^{\frak X}}\!(R',T)\i \P^{^{\frak X}}\!(Q,T)\,,$ so that we still have 
$$\frak h_{_{Q,T}}^{^{\frak X}}\!(x\.y) = \bar\imath_{_{R'}}^{^Q}\.\frak h_{_{R',R}}^{^{\frak X}}\!(x)\. 
\frak h_{_{R,T}}^{^{\frak X}}\!(y) = \frak h_{_{Q,R}}^{^{\frak X}}\!(x)\. 
\frak h_{_{R,T}}^{^{\frak X}}\!(y)
\eqno £7.21.3.$$

\medskip
£7.22. Finally, assume that $R$ belongs to $\frak Y$ and consider the respective subgroups of $R$ and $Q$
$$T' = \big(\pi_{_{R,T}}^{_{\frak X}}\!(y)\big)(T)\qq T'' = \big(\pi_{_{Q,R}}^{_{\frak X}}\!(x)\big)(T')
\eqno £7.22.1;$$
setting $N' = N_P (T')$ and $N'' = N_P (T'')\,,$ and  considering pairs $(N',n')$  
in  $\frak N (T')$ and $(N'',n'')$  in $\frak N (T'')\,,$ the {\it divisibility\/} of 
$ \P^{^{\frak Y}}\!$ forces the existence of a unique $ \P^{^{\frak Y}}\!\-$morphism $r\,\colon N'\to N''$ fulfilling $i^{^Q}_{_{N''}}\!\.r = x\. i^{^R}_{_{N'}}\,;$ then, we consider the element $s$ in $\P^{^{\frak X}}\! (T'', T') $ determined  by the equality (cf.~£7.18.1)
$$\frak n_{n'',n'}(s) = \frak g^{{}^{n''} \! N'', {}^{n'}\! N'}_{{}^{n''} 
\! T'', {}^{n''}\! T'}(n''\. r\.n'^{-1})
\eqno £7.22.2\phantom{.}$$
and, since $y\in \P^{^{\frak X}}\! (T', T)\i \P^{^{\frak X}}\! (R, T')\,,$ we can define $x\. y = s\. y\,.$

\medskip
£7.23. Once again, we have $\frak h_{_{Q,T}}^{^{\frak X}}\!(x\.y) 
= \bar\imath_{_{T''}}^{^Q}\.\frak h_{_{T'',T'}}^{^{\frak X}}\!(s)\. 
\frak h_{_{T',T}}^{^{\frak X}}\!(y)\,;$ but, according to our definition, we have
$\frak h_{_{R,T}}^{^{\frak X}}\!(y) 
= \bar\imath_{_{T'}}^{^R}\!\. \frak h_{_{T',T}}^{^{\frak X}}\!(y)$
and from  $i^{^Q}_{_{N''}}\!\.r = x\. i^{^R}_{_{N'}}$ we  get 
$$\frak h_{_{Q,R}}^{^{\frak X}}\!(x)\.\bar\imath^{^R}_{_{N'}} = 
\bar\imath^{^Q}_{_{N''}}\.\frak h_{_{N'',N'}}^{^{\frak X}}\! (r)
\eqno £7.23.1,$$
so that we obtain 
$$\frak h_{_{Q,R}}^{^{\frak X}}\!(x)\.\bar\imath^{^R}_{_{T'}} = 
\bar\imath^{^Q}_{_{N''}}\.\frak h_{_{N'',N'}}^{^{\frak X}}\! (r)
\.\bar\imath_{_{T'}}^{^{N'}} = 
\bar\imath^{^Q}_{_{T''}}\.\bar\frak g_{_{T'',T'}}^{^{N'',N'}}
\big(\frak h_{_{N'',N'}}^{^{\frak X}}\! (r)\big)
\eqno £7.23.2;$$
 thus, setting 
$\bar m' =\bar\frak g_{_{{}^{n'}\!V',V'}}^{^{{}^{n'}\!N',N'}}
\big(\frak h_{_{{}^{n'}\!N',N'}}^{^{\frak X}}\!(n')\big) $ and 
$\bar m'' = \bar\frak g_{_{{}^{n''}\!V'',V''}}^{^{{}^{n''}\!N'',N''}}
\big(\frak h_{_{{}^{n''}\!N'',N''}}^{^{\frak X}}\!(n'')\big)\,,$ we get
$$\eqalign{&\bar\frak g_{_{T'',T'}}^{^{N'',N'}}
\big(\frak h_{_{N'',N'}}^{^{\frak X}}\! (r)\big)\cr
& = \bar\frak g_{_{T'',T'}}^{^{N'','N'}}\big(\frak h_{_{{}^{n''}\!N'',N''}}^{^{\frak X}}\!(n''^{-1})\.\frak h_{_{{}^{n''}\!N'',{}^{n'}\!N'}}^{^{\frak X}}\!(n''\.r\.n'^{-1})\.
\frak h_{_{{}^{n'}\!N',N'}}^{^{\frak X}}\!(n')\big)\cr
&= \bar m''^{-1}\. \bar\frak g_{_{{}^{n''}\!T'',{}^{n'}\!T'}}^{^{{}^{n''}
\!N'',{}^{n'}\!N'}}\big(\frak h_{_{{}^{n''}\!N'',{}^{n'}\!N'}}^{^{\frak X}}
\!(n''\.r\.n'^{-1})\big)\.\bar m'\cr
&= \bar m''^{-1}\. \frak h_{_{{}^{n''}\!T'',{}^{n'}\!T'}}^{^{\frak X}}\!
\big(\frak g_{_{{}^{n''}\!T'',{}^{n'}\!T'}}^{^{{}^{n''}\!N'',{}^{n'}\!N'}}
(n''\.r\.n'^{-1})\big)\.\bar m'\cr
&= \bar m''^{-1}\. \frak h_{_{{}^{n''}\!T'',{}^{n'}\!T'}}^{^{\frak X}}\!
\big(\frak n_{n'',n'}(s)\big)\.\bar m' = \frak h_{_{T'',T'}}^{^{\frak X}}\! (s)\cr}
\eqno £7.23.3;$$
consequently, we have 
$$\eqalign{\frak h_{_{Q,R}}^{^{\frak X}}\!(x)\.\frak h_{_{R,T}}^{^{\frak X}}\!(y)
&= \frak h_{_{Q,R}}^{^{\frak X}}\!(x)\.\bar\imath_{_{T'}}^{^R}\!\. 
\frak h_{_{T',T}}^{^{\frak X}}\!(y)\cr
& = \bar\imath^{^Q}_{_{T''}}\.\bar\frak g_{_{T'',T'}}^{^{N'',N'}}
\big(\frak h_{_{N'',N'}}^{^{\frak X}}
\! (r)\big)\. \frak h_{_{T',T}}^{^{\frak X}}\!(y)\cr
&=\bar\imath^{^Q}_{_{T''}}\.\frak h_{_{T'',T'}}^{^{\frak X}}\! (s)\. \frak h_{_{T',T}}^{^{\frak X}}\!(y) = \frak h_{_{Q,T}}^{^{\frak X}}\!(x\.y)\cr}
\eqno £7.23.4.$$
This completes the definition of the composition in $\P^{^{\frak X}}\!$ and the compatibility of this composition with  
$\frak h^{^{\frak X}}\,.$

\medskip
£7.24. This composition $x\.y$ does not depend on our choice; indeed,
for another choice of pairs  $(N',\hat n')\in \frak N (T')$ and  $(N'',\hat n'')\in \frak N (T'')\,,$
it follows from Proposition~£7.11 and from equality~£7.17.2 that we have
$$\eqalign{\frak g^{{}^{\hat n''} \! N'', {}^{\hat n'}\! N'}_{{}^{\hat n''} \! T'', {}^{\hat n''}\! T'}
(\hat n''\.r\.\hat n'^{-1}) &= 
\frak g^{{}^{\hat n''} \! N'', {}^{\hat n'}\! N'}_{{}^{\hat n''} \! T'', {}^{\hat n''}\! T'}
\big((\hat n''\. n''^{-1})\.(n''\. r\.n')\.(n'\.\hat n'^{-1})\big)\cr
&= g_{\hat n'',  n''}\.\frak g^{{}^{ n''} \! N'', {}^{ n'}\! N'}_{{}^{ n''} \! T'', {}^{ n''}\! T'}
(n''\. r\.n')\.g_{n',\hat n'} \cr
&= \frak n_{\hat n'',\hat n'}(s) \cr}
\eqno £7.24.1\phantom{.}$$
Then, the {\it associativity\/} follows from equality~£7.20.5, the structural functors
are easily defined from~£7.19 and from the right-hand definition in~£7.21.1, and the
 equalities~£7.19.4 show that $\frak h^{^{\!\frak X}}\!$ is an
 {\it $\F^{^\frak X}\!\-$locality\/} functor. We are done.
\vfill
\eject

\bigskip
 \bigskip
\noindent
{\bf £8. A functor from the perfect $\F\-$locality to the basic $\F\-$locality 
$\L^{^{\rm b}}$ }

\bigskip
£8.1. Let $\F$ be a Frobenius $P\-$category. From section~£6 we already know the existence of a {\it perfect $\F^{^{\rm sc}}\-$locality\/} $\P^{^{\rm sc}}$ canonically contained in the {\it natural $\F^{^{\rm sc}}\-$locality\/} $\bar\L^{^{n,\rm sc}}$ and therefore $\P^{^{\rm sc}}$ is also contained in the corresponding quotient 
$\bar\L^{^{b,\rm sc}}$ (cf.~£4.13.3) of the {\it full\/} 
$\F^{^{\rm sc}}\!\-$sublocality $\L^{^{b,\rm sc}}$ of the 
{\it basic $\F\-$locality\/}~$\L^{^{\rm b}}$ (cf.~Corollary~£4.11); actually,
$\bar\L^{^{b,\rm sc}}$ is also the  {\it full\/} $\F^{^{\rm sc}}\!\-$sublocality of a quotient $\bar\L^{^{\rm b}}$ of the {\it basic\/} $\F\-$locality $\L^{^{\rm b}}$ (see £8.5 below) and therefore it follows from Theorem~£7.2
 that the inclusion $\P^{^{\rm sc}}\i \bar\L^{^{b,\rm sc}}$ can be extended to a unique {\it $\F\-$locality functor\/}  $\bar\frak h\,\colon \P\to \bar\L^{^{\rm b}}\,;$ the main purpose of this section is to prove that this functor can be lifted to an  {\it $\F\-$locality functor\/}   $\frak h\,\colon \P\to  \L^{^{b}}\,,$  in an essentially unique way.

\medskip
£8.2. From~£4.4 and £4.12 we have a {\it contravariant\/} functor 
$$\tilde\frak c^{^{\rm b}} : \tilde\F \too \Ab
\eqno £8.2.1\phantom{.}$$
mapping any subgroup $Q$ of $P$ on the Abelian group
$$\tilde\frak c^{^{\rm b}}(Q) =  \prod_{\tilde O\in \frak O_Q} \ab \big({\rm Aut}(O)\big)
\eqno £8.2.2,$$
where we denote by  $\frak O_Q$ the  set of isomorphism classes  of indecomposable $Q\times P\-$sets 
$(Q\times P)/\Delta_\theta (U)$ where $U$ is a subgroup of $P$ and $\theta\,\colon U\to Q$ an $\F\-$morphism,  
and mapping any $\tilde\F\-$morphism $\tilde\varphi\,\colon R\to Q$ on the group homomorphism
$$\tilde\frak c^{^{\rm b}}(\tilde\varphi) :   \prod_{\tilde O\in \frak O_Q} \ab \big({\rm Aut}(O)\big)\too  
 \prod_{\tilde O\in \frak O_R} \ab \big({\rm Aut}(O)\big)
\eqno £8.2.3\phantom{.}$$
 described in~Proposition~£4.6 above.

\medskip
£8.3.  For any $\tilde O\in \frak O_Q\,,$ note that the homomorphism $\tilde\frak c^{^{\rm b}} (\tilde\varphi)$ sends  an element of $\ab \big({\rm Aut}(O)\big)$ to a family of terms indexed by $R\times P\-$orbits with ``smaller'' stabilizers. More precisely, consider a set $\frak N$ of subgroups  of $P$ such that  any subgroup $U$ of $P$ fulfilling $\F (T,U)\not= \emptyset$ for  some $T\in \frak N$ belongs to~$\frak N\,,$ and for any subgroup $Q$ of $P$ denote by $\frak O^{^{\!\frak N}}_Q$ the subset of 
$\tilde O\in \frak O_Q$ such that  $(Q\times P)/\Delta_\eta (T)$ belongs to $\tilde O$ if and only if $T$ belongs to $\frak N\,.$

\bigskip
\noindent
{\bf Corollary~£8.4.} {\it With the notation above, the correspondence sending  any subgroup $Q$ of $P$ to 
$$\tilde\frak c^{^{\frak N}}\! (Q) = \prod_{\tilde O\in \frak O^{^{\!\frak N}}_Q} 
\ab \big({\rm Aut}(O)\big)
\eqno £8.4.1\phantom{.}$$
defines a contravariant subfunctor $\tilde\frak c^{^{\frak N}}\colon \tilde\F\to  \Ab$ of $\,\tilde\frak c^{^{\rm b}}\,.$\/}

\medskip
\noindent
{\bf Proof:} Straightforward.  
\eject

\medskip
£8.5. In particular, considering the set of subgroups  of $P$ which are {\it not\/}  $\F\-$selfcentralizing and denoting 
by $\tilde\frak c^{^{\rm nsc}}\colon \tilde\F\to  \Ab$ the corresponding subfunctor of~$\tilde\frak c^{^{\rm b}}\,,$ 
we get the quotient $\bar\L^{^{\rm b}} = \L^{^{\rm b}}/\tilde\frak c^{^{\rm nsc}}$ of the {\it basic\/} $\F\-$locality (cf.~£2.9); as mentioned above, it follows from~£4.13.3 and from  section 6 that the {\it perfect $\F^{^{\rm sc}}\!\-$locality\/} 
$\P^{^{\rm sc}}$ is contained in the {\it full\/} subcategory $\bar\L^{^{\rm b, sc}}$ of $\bar\L^{^{\rm b}}$ over the set of 
$\F\-$selfcentralizing subgroups of~$P\,,$ and then it follows from Theorem~£7.2 that this inclusion can be extended to a unique {\it $\F\-$locality functor\/}
$$\bar\frak h : \P\too \bar\L^{^{\rm b}}
\eqno £8.5.1\phantom{.}$$
where $\P$ denotes the  {\it perfect $\F\-$locality\/}. As a matter of fact, in this case it is not necessary to quote
Theorem~£7.2 since, for any pair of subgroups $Q$ and $R$ of $P\,,$ if $R$  is not  $\F\-$selfcentralizing then  we just have 
$\bar\L^{^{\rm b}}\!(Q,R) = \F (Q,R)\,.$

\medskip
£8.6. More generally, for any set $\frak N$ as in~£8.3 above, we consider the corres-ponding quotient --- denoted by 
$(\bar\tau^{_{\frak N,\rm b}},\bar\L^{^{\frak N,\rm b}},\bar\pi^{^{\frak N,\rm b}})$ --- of the {\it basic\/} $\F\-$locality; 
if all the subgroups in $\frak N$ are {\it not\/} $\F\-$selfcentralizing, then we claim that $\bar\frak h$ can be lifted to a  unique
{\it natural $\F\-$isomorphism\/} class of {\it $\F\-$locality functors\/} 
$$\bar\frak h^{^\frak N} : \P\too \bar\L^{^{\frak N,\rm b}}
\eqno £8.6.1.$$
The induction argument on the cardinal of the complement of $\frak N$ in  the set of all the subgroups  of $P$ justifies  
 the following general construction;  assume that this complement is not empty, choose on it a minimal element 
$U$  {\it fully normalized\/} in $\F$ and set
$$\frak M = \frak N \cup \{\theta (U)\mid \theta\in \F (P,U)\}
\eqno £8.6.2;$$
then,  we clearly have a canonical functor $\bar\frak l^{^{\frak M,\frak N}} \,\colon \bar\L^{^{\frak N,\rm b}}
\to  \bar\L^{^{\frak M,\rm b}}$ and, since we have $\tilde\P = \tilde\F\,,$
it makes sense to consider the quotient {\it contravariant\/}  functor  
$$\tilde\frak t^{^U }=\Ker (\bar\frak l^{^{\frak M,\frak N}}) =  \tilde\frak c^{^{\frak M}}/\tilde\frak c^{^{\frak N}} : 
\tilde\F \too \Ab
\eqno £8.6.3;$$
actually, we prove below that this  {\it contravariant\/} functor admits a {\it compatible complement\/}
in the sense of [12,~5.7] and therefore, denoting by $\Bbb H_*^n (-,-)$ the $n\-$th {\it stable cohomology groups\/}
[11,~A317], for any $n\ge 1$ we have [12,~Proposition~5.8]
$$\Bbb H^n_* (\tilde\F, \tilde\frak t^{^U }) = \{0\}
\eqno £8.6.4;$$
this fact will be also quoted in section 9.

\medskip
£8.7. Recall that, for any  subgroup $R$ of $P$ and any element $\eta\in \F (R,U)\,,$ we have (cf.~£3.6.1)
$${\rm Aut}\big((R\times P)/\Delta_\eta(U)\big) \cong \bar N_{R\times P}\big(\Delta_{\eta}(U)\big)
\eqno £8.7.1,$$
and, for any $\F\-$morphism $\varphi\,\colon R\to Q\,,$ denote by 
$$\bar\varphi_\eta : \bar N_{R\times P}\big(\Delta_{\eta}(U)\big)\too 
\bar N_{Q\times P}\big(\Delta_{\varphi\circ\eta}(U)\big)
\eqno £8.7.2\phantom{.}$$
the group homomorphism induced by $\varphi\times {\rm id}_P\,.$ Moreover, for any  $\tilde\F\-$mor-phisms
 $\tilde\varphi\,\colon R\to Q$ and  $\tilde\theta\,\colon U\to Q\,,$ the existence of an injective $R\times P\-$set homomorphism
 $$f : (R\times P)/\Delta_\eta (U)\too {\rm Res}_{\varphi\times {\rm id}_P}\big((Q\times P)/\Delta_\theta (U)\big)
 \eqno £8.7.3,$$
 for some $\varphi\in \tilde\varphi$ and some $\theta\in \tilde\theta\,,$ is equivalent to the existence of $(v,u)\in Q\times P$ such that
 $$(\varphi\times {\rm id}_P)\big(\Delta_\eta (U)\big) = \Delta_\theta (U)^{(v,u)}
 \eqno £8.7.4,$$
 which implies that $u$ belongs to $N_P (U)$ and that $\tilde\varphi\circ\tilde\eta  = \tilde\theta\circ\tilde\kappa_{_U}\!(u)$ (cf.~£2.7).

 \medskip
 £8.8. More precisely, as in~£4.5 above set
 $$ M = (R\times P)/\Delta_\eta (U)\qq O = (Q\times P)/\Delta_\theta (U)
 \eqno £8.8.1\phantom{.}$$
and denote by ${\rm Inj}_{R\times P}\big(M,{\rm Res}_{\varphi\times {\rm id}_P}(O)\big)$ the corresponding set of  injective 
$R\times P\-$set homomorphisms; assuming that this set is not empty, in~£8.7.3 let us denote by 
$(w_f,\hat w_f)\in Q\times P$ an element belonging to the class in $O$ which is the image by $f$ of the class of $(1,1)$ in $M\,;$  from [11,~Lemma~22.19]  it follows that
 \smallskip
 \noindent
£8.8.2\quad {\it The correspondence mapping $f$ on the double class  $\varphi (R)w_f\theta (U)$ determines a bijection from the set of  ${\rm Aut}(M)\-$orbits in  ${\rm Inj}_{R\times P}\big(M,{\rm Res}_{\varphi\times {\rm id}_P}(O)\big)$ 
onto the set of double classes $\varphi (R)v\theta (U)$ in  $\varphi (R)\backslash Q/\theta (U)$ admitting a re-presentative $v\in Q$ such that $(v,\hat v)$ normalizes  
$\Delta_\theta (U)$ for some $\hat v\in N_P (U)\,.$\/} 
\smallskip
 \noindent
In particular, since $(v,\hat v)$ belongs to $\bar N_{Q\times P}\big(\Delta_\theta (U)\big)\cong {\rm Aut}(O)\,,$ 
in the present situation with the notation in~£4.5 {\it we have $\vert \I^{\tilde O}_{\tilde M} \vert = 1\,,$ we may assume that 
$\tilde\varphi\circ\tilde\eta  = \tilde\theta\,,$ the homomorphism $\delta_f$ in~{\rm £4.5.2} is an isomorphism  and moreover we have\/}
 $$\ab(\delta_f)\circ \ab^\circ (\varepsilon_f) = \ab^\circ (\bar\varphi_\eta)
 \eqno £8.8.3.$$
 At this point, choosing   a set of representatives  $\frak O^{\tilde\varphi}_{\tilde\theta}$ for the set of classes in 
 $\tilde\F (R,U)/\tilde\F_P(U)$ having an element $\tilde\gamma$ such that  
 $\tilde\varphi\circ\tilde\gamma   = \tilde\theta\,,$ and setting 
 $M_{\tilde\gamma} = (R\times P)/\Delta_\gamma (U)$
 for some $\gamma\in \tilde\gamma\,,$ we actually have a bijection
  $$\bigsqcup_{\tilde\gamma\in \frak O^{\tilde\varphi}_{\tilde\theta}} M_{\tilde\gamma}/{\rm Aut}(M_{\tilde\gamma})
  \cong O/{\rm Aut}(O)
 \eqno £8.8.4.$$

\bigskip
\noindent
{\bf Proposition~£8.9.} {\it With the notation above, the contravariant functor $\tilde\frak t^{^U }$
maps any subgroup $Q$ of $P$ on 
$$\tilde\frak t^{^U } \!(Q) = \bigg(\prod_{\tilde\theta\in \tilde\F (Q,U)} 
\ab \Big(\bar N_{Q\times P}\big(\Delta_{\theta}(U)\big)\Big)\bigg)^{\tilde\F_P (U)}
\eqno £8.9.1$$
and any $\tilde\F\-$morphism $\tilde\varphi\,\colon R\to Q$ on the  homomorphism induced by the sum of
the group homomorphisms
$$\ab^\circ (\bar\varphi_\eta) : \ab\Big(\bar N_{Q\times P}\big(\Delta_{\varphi\circ\eta}(U)\big)\Big)
\too \ab\Big(\bar N_{R\times P}\big(\Delta_{\eta}(U)\big)\Big)
\eqno £8.9.2\phantom{.}$$
where $\tilde\eta$ runs over $\tilde\F (R,U)\,.$ In particular,   $\,\tilde\frak t^{^U }$
admits a compatible complement mapping $\tilde\varphi\,\colon R\to Q$ on 
the  homomorphism $\tilde\frak t^{^U }\! (\tilde\varphi)^\circ$ induced by the sum of the group homomorphisms
$$\ab (\bar\varphi_\eta) : \ab\Big(\bar N_{R\times P}\big(\Delta_{\eta}(U)\big)\Big)
\too \ab\Big(\bar N_{Q\times P}\big(\Delta_{\varphi\circ\eta}(U)\big)\Big)
\eqno £8.9.3 \phantom{.}$$
where $\tilde\eta$ runs over $\tilde\F (R,U)\,.$ \/}

\medskip
\noindent
{\bf Proof:} It is clear that 
$$\tilde\frak t^{^U }\! (Q) = \prod_{\tilde O\in \frak O^{^{\!\frak M}}_Q - \frak O^{^{\!\frak N}}_Q} 
\ab \big({\rm Aut}(O)\big)
\eqno £8.9.4;$$
but, for any $\tilde O\in \frak O^{^{\!\frak M}}_Q - \frak O^{^{\!\frak N}}_Q\,,$ we necessarily have 
 $O\cong (Q\times P)/\Delta_\theta (U)$ for some $\theta\in \F (Q,U)\,;$ moreover, it is clear that
 we have a $Q\times P\-$set isomorphism
 $$(Q\times P)/\Delta_\theta (U)\cong (Q\times P)/\Delta_{\theta'} (U)
 \eqno £8.9.5\phantom{.}$$
 if and only if $\Delta_\theta (U)$ and $\Delta_{\theta'} (U)$ are $Q\times P\-$conjugate to each other;
now, equality~£8.9.1 follows from~£3.6.1 up to  suitable identifications.

\smallskip
On the other hand, if $\tilde\varphi\,\colon R\to Q$ is an $\tilde\F\-$morphism then it follows from Proposition~£4.6
that, for any $\theta\in \F (Q,U)$ and any $a\in \ab\Big(\bar N_{Q\times P}\big(\Delta_{\theta}(U)\big)\Big)\,,$
setting $O = (Q\times P)/\Delta_\theta (U)$ we get
$$\big(\tilde\frak c^{^\frak M} (\tilde\varphi)\big)(a) = \sum_{\tilde M\in \frak O^{^\frak M}_R} 
\,\sum_{f\in \I_{\tilde M}^{\tilde O}(\tilde\varphi)}\big(\ab(\delta_f)\circ \ab^\circ (\varepsilon_f)\big)(a)
\eqno £8.9.6,$$
\eject
\noindent
so that, according to~equality~£8.9.4, we still get
$$\big(\tilde\frak t^{^U} \!(\tilde\varphi)\big)(a) = \sum_{\tilde M\in \frak O^{^\frak M}_R - \frak O^{^\frak N}_R} \,\sum_{f\in \I_{\tilde M}^{\tilde O}(\tilde\varphi)}\big(\ab(\delta_f)\circ \ab^\circ (\varepsilon_f)\big)(a)
\eqno £8.9.7;$$
once again,  for any $\tilde M\in \frak O^{^{\!\frak M}}_R - \frak O^{^{\!\frak N}}_R\,,$ we necessarily have 
 $M\cong (R\times P)/\Delta_\eta (U)$ for some $\eta\in \F (R,U)\,;$ consequently, from~£8.8 we obtain
 $$\tilde\varphi\circ\tilde\eta   = \tilde\theta\qq \big(\tilde\frak t^{^U} \!(\tilde\varphi)\big)(a) 
 = \sum_{\tilde \eta\in \frak O^{\tilde\varphi}_{\tilde\theta}} \big(\ab^\circ(\bar\varphi_\eta)\big)(a)
 \eqno £8.9.8.$$

\smallskip
In order to prove the last statement in the proposition, with the notation above we have to compute the following element of $\bar\frak t^{^U}\!(Q)$
$$\sum_{\tilde\eta\in \frak O^{\tilde\varphi}_{\tilde\theta}} 
\big(\ab (\bar\varphi_\eta)\big)\Big(\big(\ab^\circ (\bar\varphi_\eta)\big)(a)\Big)
= \sum_{\tilde\eta\in \frak O^{\tilde\varphi}_{\tilde\theta}}
{\big\vert \bar N_{Q\times P}\big(\Delta_\theta (U)\big) 
\big\vert\over \big\vert \bar N_{R\times P}\big(\Delta_\eta (U)\big)\big\vert}\. a
\eqno £8.9.9;$$
but, according to the bijection~£8.8.4, we get
$$\sum_{\tilde\eta\in \frak O^{\tilde\varphi}_{\tilde\theta}} {\vert R\times P\vert\over \big\vert  N_{R\times P}\big(\Delta_\eta (U)\big)\big\vert} = {\vert Q\times P\vert\over \big\vert  N_{Q\times P}\big(\Delta_\theta (U)\big)\big\vert}
\eqno £8.9.10;$$
consequently, we obtain $\tilde\frak t^{^U}\! (\tilde\varphi)^\circ\circ \tilde\frak t^{^U}\! (\tilde\varphi) = \vert Q\vert/
\vert R\vert\. {\rm id}_{\tilde\frak t^{^U}\!(Q)}\,.$

\smallskip
Finally, consider a  {\it special $\ad(\tilde\F)\-$square\/} [12,~5.1 and~5.2]
$$\matrix{&&Q\cr
&{\tilde\varphi\atop}\hskip-5pt\nearrow\hskip-5pt&&\hskip-5pt\nwarrow\hskip-5pt{\tilde\psi\atop}\cr
R\hskip-10pt&&&&\hskip-10pt T\cr
&{\atop\tilde\zeta}\hskip-5pt\nwarrow\hskip-5pt&&\hskip-5pt\nearrow\hskip-5pt{\atop\tilde\xi}\cr
&&S\cr}
\eqno £8.9.11;$$ 
in order to prove that $\tilde\frak t^{^U}$ admits a {\it compatible complement\/}, we may assume that
$Q\,,$ $R$ and $T$  are just subgroups of $P\,;$ more precisely, up to isomorphisms, we may assume
that $Q$ contains $R$ and $T\,,$ and that $\tilde\varphi = \tilde\iota_R^Q$ and $\tilde\psi = \tilde\iota_T^Q\,;$
in this case, by the very definition of the
{\it special $\ad(\tilde\F)\-$squares,\/} we have
 $$S = \bigoplus_{w\in W} S_w
\eqno £8.9.12\phantom{.}$$
where $W\i Q$ is  a set of representatives for the set of double classes $R\backslash Q/T$ and, 
for any $w\in W\,,$ we set $S_w = R^w\cap T$
and respectively denote by
$$\iota_w^R : S_w\too R\qq \iota_w^T : S_w\too T
\eqno £8.9.13\phantom{.}$$
the $\F\-$morphisms mapping $s\in S_w$ on $wsw^{-1}$ and on $s\,;$ moreover, the  $\ad(\tilde\F)\-$morphisms
$$\tilde\zeta :\bigoplus_{w\in W} S_w\too R\qq \tilde\xi :\bigoplus_{w\in W} S_w\too T
\eqno £8.9.14\phantom{.}$$
are respectively determined by the families $\{\tilde\iota_w^R\}_{w\in W}$ and $\{\tilde\iota_w^T\}_{w\in W}\,.$
\eject

\smallskip
Now, it suffices to prove the commutativity of the following diagram
$$\matrix{&&\tilde\frak t^{^U} \!(Q)\cr
&{\tilde\frak t^{^U} \!(\tilde\iota_R^Q)^\circ\atop}\hskip-5pt\nearrow\hskip-10pt&
&\hskip-10pt\searrow\hskip-5pt{\tilde\frak t^{^U}\!(\tilde\iota_T^Q)\atop}\cr
\tilde\frak t^{^U}\! (R)\hskip-40pt&&\phantom{\bigg\uparrow}&&\hskip-30pt \tilde\frak t^{^U}\! (T)\cr
&{\atop\tilde\frak t^{^U}\! (\tilde\zeta)}\hskip-5pt\searrow\hskip-15pt&
&\hskip-15pt\nearrow\hskip-5pt{\atop\tilde\frak t^{^U}\! (\tilde\xi)^\circ}\cr
&&{\displaystyle \prod_{w\in W}}\tilde\frak t^{^U}\!  (S_w)\cr}
\eqno £8.9.15;$$ 
that is to say, for any $\eta\in \F (R,U)$ and any $a\in \ab\Big(\bar N_{R\times P}\big(\Delta_{\eta}(U)\big)\Big)\,,$
it suffices to prove that
$$\tilde\frak t^{^U}\!(\tilde\iota_T^Q)\big(\tilde\frak t^{^U} \!(\tilde\iota_R^Q)^\circ (a)\big) = 
\tilde\frak t^{^U}\! (\tilde\xi)^\circ\big(\frak t^{^U}\! (\tilde\zeta) (a)\big)
\eqno £8.9.16\,.$$
According to our definition of $\tilde\frak t^{^U} \!(\tilde\iota_R^Q)^\circ$ (cf.~£8.9.3), we have
$$\tilde\frak t^{^U} \!(\tilde\iota_R^Q)^\circ (a) = \big(\ab ((\bar\iota_R^Q)_\eta)\big)(a)
\eqno £8.9.17\phantom{.}$$
and therefore, setting $\theta = \iota_R^Q\circ\eta\,,$ it follows from £8.9.8 that we get
$$\tilde\frak t^{^U}\!(\tilde\iota_T^Q)\big(\tilde\frak t^{^U} 
\!(\tilde\iota_R^Q)^\circ (a)\big) = 
\sum_{\tilde\gamma\in \frak O_{\!\tilde\iota_T^Q}^{\tilde\theta}} \big(\ab^\circ 
((\,\overline{\!\iota_T^Q\!}\,)_{\gamma})\circ \ab ((\,\overline{\!\iota_R^Q\!}\,)_\eta)\big)(a)
\eqno £8.9.18.$$
On the other hand, it follows from~£8.9.8 that
$$\frak t^{^U}\! (\tilde\zeta) (a) = \sum_{w\in W}\,\sum_{\tilde\upsilon\in \frak O^{\tilde\eta}_{\!\tilde\iota_w^R}}
\Big(\ab^\circ \big((\,\overline{\!\iota_w^R\!}\,)_\upsilon\big)\Big)(a)
\eqno £8.9.19\phantom{.}$$
and, according to our definition of $\,\tilde\frak t^{^U}\! (\tilde\xi)^\circ$ (cf.~£8.9.3), we obtain
$$\tilde\frak t^{^U}\! (\tilde\xi)^\circ\big(\frak t^{^U}\! (\tilde\zeta) (a)\big) = 
\sum_{w\in W}\,\sum_{\tilde\upsilon\in \frak O^{\tilde\eta}_{\tilde\iota_w^R}}
\Big(\ab \big((\,\overline{\!\iota_w^T\!}\,)_\upsilon\big)\circ \ab^\circ 
\big(((\,\overline{\!\iota_w^R\!}\,)_\upsilon\big)\Big)(a)
\eqno £8.9.20.$$

\smallskip
Thus, it suffices to prove the commutativity of the following diagram
$$\hskip-12pt\matrix{&&\ab\Big(\bar N_{Q\times P}\big(\Delta_\theta (U)\big)\Big)\cr
&{\ab ((\,\overline{\!\iota_R^Q\!}\,)_\eta)\atop}\hskip-5pt\nearrow\hskip-120pt&
&\hskip-100pt\searrow\hskip-5pt{\sum_{\tilde\gamma\in \frak O_{\!\tilde\iota_T^Q}^{\tilde\theta}}\ab^\circ 
((\,\overline{\!\iota_T^Q\!}\,)_{\gamma})\atop}\cr
\ab\Big(\bar N_{R\times P}\big(\Delta_\eta (U)\big)\Big)\hskip-180pt&&\phantom{\Big\uparrow}
&&\hskip-140pt \prod_{\tilde\gamma\in \frak O_{\!\tilde\iota_T^Q}^{\tilde\theta}} \ab\Big(\bar N_{T\times P} \big(\Delta_\gamma (U)\big)\Big)\cr
&{\atop \sum_{w\in W}\,\sum_{\tilde\upsilon\in \frak O^{\tilde\eta}_{\!\tilde\iota_w^R}} \ab^\circ 
(((\,\overline{\!\iota_w^R\!}\,)_\upsilon)}\searrow\hskip-50pt&
&\hskip-60pt\nearrow{\atop \sum_{w\in W}\,\sum_{\tilde\upsilon\in \frak O^{\tilde\eta}_{\!\tilde\iota_w^R}} \ab 
(((\,\overline{\!\iota_w^T\!}\,)_\upsilon)}\cr
&&{\displaystyle \prod_{w\in W}\,\prod_{\tilde\upsilon\in \frak O^{\tilde\eta}_{\!\tilde\iota_w^R}}}\ab\Big( \bar N_{S_w\times P}\big(\Delta_\upsilon (U)\big)\Big)\cr}
\eqno £8.9.21;$$ 
but, denoting by $Q'$ the converse image in $P$ of 
${}^{\theta^*}\! \F_P\big(\theta (U)\big)\cap \tilde\F_P (U)$ where $\theta^*\,\colon \theta (U)\cong U$
is the inverse of the isomorphism induced by $\theta\,,$ $N_{Q\times P}\big(\Delta_\theta (U)\big)$ is clearly contained in $Q'\times P$ and therefore we have
$$\ab\Big(\bar N_{Q\times P}\big(\Delta_\theta (U)\big)\Big) = 
\ab\Big(\bar N_{Q'\times P}\big(\Delta_\theta (U)\big)\Big)
\eqno £8.9.22;$$
 explicitly, it follows from [11,~statement~2.10.1] that the composition $\iota_U^P\circ \theta^*$ can be extended
to an $\F\-$morphism $\chi'\,\colon Q'\to P$ and then we actually have
$$N_{Q'\times P}\big(\Delta_\theta (U)\big) = \big(\{1\}\times C_P (U)\big)\rtimes 
N_{\Delta^{\!*}_{\chi'}( Q')}\big(\Delta_\theta (U)\big)
\eqno £8.9.23\phantom{.}$$
where we set $\Delta^{\!*}_{\chi'}( Q') = \Delta_{{\rm id}_{Q'},\chi'}(Q')\,.$

\smallskip
According to our choice of $\theta\,,$ we have $\Delta_\eta (U) = \Delta_\theta (U)$ and therefore, setting 
$R' = R\cap Q'\,,$ we similarly  have
$$\ab\Big(\bar N_{R\times P}\big(\Delta_\theta (U)\big)\Big) = 
\ab\Big(\bar N_{R'\times P}\big(\Delta_\theta (U)\big)\Big)
\eqno £8.9.24\phantom{.}$$
and, denoting by $\rho'\,\colon R'\to P$ the restriction of $\chi'\,,$
we still have
$$N_{R'\times P}\big(\Delta_\theta (U)\big) = \big(\{1\}\times C_P (U)\big)\rtimes 
N_{\Delta^*_{\rho'}( R')}\big(\Delta_\theta (U)\big)
\eqno £8.9.25.$$
On the other hand, for any $\tilde\gamma\in \frak O_{\!\tilde\iota_T^Q}^{\tilde\theta}\,,$ we have an injective
$T\times P\-$set homomorphism
$$(T\times P)/\Delta_\gamma (U)\too {\rm Res}_{\iota_T^Q\times {\rm id}_P}
\big((Q\times P)/\Delta_\theta (U)\big)
\eqno £8.9.26\phantom{.}$$
and therefore, setting $T' = T\cap Q'$ and  denoting by $\theta'\,\colon U\to Q'$ the restriction of $\theta\,,$ 
it follows from statement~£8.8.2 that, up to a modification of our choice of the set of representatives 
in $\tilde\F (T,U)/\tilde\F_P (U)$ and our choice of $\gamma$ in~$\tilde\gamma\,,$
we also may assume that $\Delta_\gamma (U)= \Delta_\theta (U)\,,$ and then  it is easily checked that
we get an $\F\-$morphism $\gamma'\,\colon U\to T'$ such that $\iota_{T'}^{T}\circ \gamma'$ is a representative 
of~$\tilde\gamma$ and that we have $\tilde\theta' = \tilde\iota_{T'}^{Q'}\circ \tilde\gamma'\,;$
that is to say, there is $v'\in Q'$ such that we have $\theta' (u) = {}^{v'}\! \gamma' (u)$ for any $u\in U$
and therefore we get (cf.~£2.7)
$$\theta' \circ \kappa_{_U}\!\big(\chi'(v')\big) = \iota_{T'}^{Q'}\circ \gamma'
\eqno £8.9.27;$$
since $\chi'(v')$ belongs to $N_P (U)\,,$ we actually obtain 
$\vert \frak O_{\!\tilde\iota_T^Q}^{\tilde\theta}\vert = 1$ (cf.~£8.8). Moreover,  denoting    by $\tau'\,\colon T'\to P$ 
the restriction of $\chi'\,,$ we still have
$$\eqalign{\ab\Big(\bar N_{T\times P}\big(\Delta_\theta (U)\big)\Big) &= 
\ab\Big(\bar N_{T'\times P}\big(\Delta_\theta (U)\big)\Big)\cr
N_{T'\times P}\big(\Delta_\theta (U)\big) &= \big(\{1\}\times C_P (U)\big)\rtimes 
N_{\Delta^*_{\tau'}( T')}\big(\Delta_\theta (U)\big)\cr}
\eqno £8.9.28.$$

\smallskip
Finally, for any $w\in W$ and any $\tilde\upsilon\in \frak O^{\tilde\eta}_{\!\tilde\iota_w^R}$  we have an injective
$S_w\times P\-$set homomorphism
$$(S_w\times P)/\Delta_\upsilon (U)\too {\rm Res}_{\iota_w^R\times {\rm id}_P}
\big((R\times P)/\Delta_\eta (U)\big)
\eqno £8.9.29\phantom{.}$$
and therefore, since  $\Delta_\eta (U) = \Delta_\theta (U)\,,$ setting $\iota_w^Q = \iota_R^Q\circ \iota_w^R$ 
we also have  an injective
$S_w\times P\-$set homomorphism
$$(S_w\times P)/\Delta_\upsilon (U)\too {\rm Res}_{\iota_w^Q\times {\rm id}_P}
\big((Q\times P)/\Delta_\theta (U)\big)
\eqno £8.9.30;$$
once again,  setting $S'_w = S_w\cap Q'\,,$ 
it follows from statement~£8.8.2 that, up to a modification of our choice of the set of representatives 
in $\tilde\F (S_w,U)/\tilde\F_P (U)$ and our choice of $\upsilon$ in~$\tilde\upsilon\,,$
we also may assume that ${}^{(w,1)}\Delta_\upsilon (U)= \Delta_\theta (U)\,;$ in this case,  up to a modification of our choice of $W\,,$ we may assume that $w$ belongs to $Q'\,.$

\smallskip
 Thus, denoting by $W'$ the set of $w\in W$ such that 
 $\frak O^{\tilde\eta}_{\!\tilde\iota_w^R}\not= \emptyset\,,$
 we may assume that $W'$ is contained in $Q'$ and then it is actually a set of representatives for the set of double classes $R'\backslash Q'/ T'\,.$ Moreover, denoting by $\eta'\,\colon U\to R'$ the restriction of $\eta$ and by 
 $\iota_w^{R'}\,\colon S'_w\to R'$ the restriction of $\iota_w^R$ for any $w\in W'\,,$
 the argument above also proves that,  we still have 
 $\vert\frak O^{\tilde\eta'}_{\!\tilde\iota_w^{R'}}\vert = 1$ (cf.~£8.8) and then we still obtain
 $$\eqalign{\ab\Big(\bar N_{S_w\times P}\big(\Delta_{\upsilon} (U)\big)\Big) &= 
\ab\Big(\bar N_{S'_w\times P}\big(\Delta_{\upsilon} (U)\big)\Big)\cr
N_{S'_w\times P}\big(\Delta_{\upsilon} (U)\big) &= \big(\{1\}\times C_P (U)\big)\rtimes 
N_{\Delta^*_{\sigma'}( S'_w)}\big(\Delta_{\upsilon} (U)\big)\cr}
\eqno £8.9.31\phantom{.}$$
where  $\sigma'\,\colon S'_w\to P$ is the respective restriction of $\tau'\,.$

 \smallskip
 In conclusion, respectively denoting by $\tilde\gamma$ and by $\tilde\upsilon_w$ the unique elements of  $\frak O_{\!\tilde\iota_T^Q}^{\tilde\theta}$ and of 
 $\frak O^{\tilde\eta}_{\!\tilde\iota_w^R}$ for any $w\in W'\,,$
 diagram~£8.9.21 becomes
 $$\matrix{&&\ab\Big(\bar N_{Q\times P}\big(\Delta_\theta (U)\big)\Big)\cr
&{\ab ((\,\overline{\!\iota_R^Q\!}\,)_\eta)\atop}\hskip-5pt\nearrow\hskip-80pt&
&\hskip-70pt\searrow\hskip-5pt{\ab^\circ 
((\,\overline{\!\iota_T^Q\!}\,)_{\gamma})\atop}\cr
\ab\Big(\bar N_{R\times P}\big(\Delta_\eta (U)\big)\Big)\hskip-100pt&&\phantom{\Big\uparrow}
&&\hskip-100pt  \ab\Big(\bar N_{T\times P} \big(\Delta_\gamma (U)\big)\Big)\cr
&{\atop \sum_{w\in W'} \ab^\circ 
((\,\overline{\!\iota_w^R\!}\,)_{\upsilon_w})}\searrow\hskip-30pt&
&\hskip-30pt\nearrow{\atop \sum_{w\in W'} \ab 
((\,\overline{\!\iota_w^T\!}\,)_{\upsilon_w})}\cr
&&{\displaystyle \prod_{w\in W'}}\ab\Big( \bar N_{S_w\times P}\big(\Delta_{\upsilon_w} (U)\big)\Big)\cr}
\eqno £8.9.32;$$
but, denoting by $\frak i\Gr$ the category of finite groups with injective group homomorphisms,  the {\it contravariant\/} functor determined by the {\it transfert\/}
$$\ab^\circ : \frak i\Gr\too \Ab
\eqno £8.9.33\phantom{.}$$
admits the functor $\ab\,\colon \frak i\Gr\to \Ab$ as a {\it Mackey complement\/}; moreover, the (bijective) image
of $W'$ in $\bar N_{Q\times P}\big(\Delta_\theta (U)\big)$ is a set of representatives for the set of double classes
$$\bar N_{R\times P}\big(\Delta_\eta (U)\big)\big\backslash\bar N_{Q\times P}\big(\Delta_\theta (U)\big)
\big/\bar N_{T\times P} \big(\Delta_\gamma (U)\big)
\eqno £8.9.34;$$
consequently, the commutativity of  diagram~£8.9.32 follows from the  {\it Mackey formula\/} applied to
the pair $(\ab^\circ,\ab)\,.$ We are done.

\bigskip
\noindent
{\bf Theorem~£8.10.} {\it The  $\F\-$locality functor $\bar\frak h\,\colon \P\to  \bar\L^{^{b}}$ can be lifted to a unique natural 
$\F\-$isomorphism class of  $\F\-$locality functors 
$\frak h\,\colon \P\to  \L^{^{b}}\,.$\/}

\medskip
\noindent
{\bf Proof:} As above,  consider a set $\frak N$ of subgroups  of $P$ in  such a way that  any subgroup $V$ of $P$ fulfilling $\F (T,V)\not= \emptyset$ for  some 
$T\in \frak N$  belongs to $\frak N\,;$ assume that  all the subgroups in $\frak N$ are not $\F\-$selfcentralizing; arguing by induction on the cardinal $\rm c$ of the complement of $\frak N$ in  this set, we will prove that   $\bar\frak h\,\colon \P\to  \bar\L^{^{b}}$ can be lifted to a unique  {\it natural 
$\F\-$isomorphism\/} class of {\it $\F\-$locality functors\/}
$$\bar\frak h^{^\frak N} : \P\too \bar\L^{^{\frak N,\rm b}}
\eqno £8.10.1.$$
We may assume that $\rm c \not= 0$ and then,  in the complement of $\frak N$ we choose a minimal element $U$  {\it fully normalized\/} in $\F\,,$ setting
$$\frak M = \frak N \cup \{\theta (U)\mid \theta\in \F (P,U)\}
\eqno £8.10.2;$$
according to the induction hypothesis, we have  a unique {\it  natural $\F\-$isomor-phism\/} class of {\it $\F\-$locality functors\/} 
$\bar\frak h^{^\frak M}\,\colon \P\to \bar\L^{^{\frak M,\rm b}}$ lifting~$\bar\frak h\,;$ then, it suffices to prove that
such a functor can be lifted to a unique  {\it natural  $\F\-$isomorphism\/} class of {\it $\F\-$locality functors\/} 
$\bar\frak h^{^\frak N}\!$ as in~£8.10.1.

\smallskip
As in~£2.16 above, let us consider  the functors [11,~18.20.3]
$$\loc_{\,\P} : \ch^*(\F)\too \widetilde{\Loc}\qq
\loc_{ \bar\L^{^{\frak M,\rm b}}} : \ch^*(\F )\too \widetilde{\Loc}
\eqno £8.10.3,$$
respectively mapping any {\it $\F\-$chain\/} $\frak q\,\colon \Delta_n\to \F$
on $\big(\P (\hat\frak q),{\rm Ker}(\pi_{\hat\frak q})\big)$ and on $\big(\bar\L^{^{\frak M,\rm b}} (\hat\frak q^{^\frak M}),
{\rm Ker}(\bar\pi^{^{\frak M,\rm b}}_{\hat\frak q^{^\frak M}})\big)$ where 
$$\hat\frak q : \Delta_n\too \P\qq \hat\frak q^{^\frak M} : \Delta_n\too \bar\L^{^{\frak M,\rm b}} 
 \eqno £8.10.4\phantom{.}$$
 are respective $\P\-$ and $\bar\L^{^{\frak M,\rm b}}\! \-$chains lifting $\frak q\,,$
  and consider the obvious {\it natural map\/}
$$\loc_{\,\bar\frak h^{^\frak M}} : \loc_{\,\P}\too 
\loc_{\bar\L^{^{\frak M,\rm b}} }
\eqno £8.10.5\phantom{.}$$
 determined by the {\it $\F\-$locality functor $\bar\frak h^{^\frak M}\,.$ \/}

 \smallskip
 Actually,  from the uniqueness  part of [11,~Proposition~18.19], it is clear that the functor $\loc_{\,\P}$ coincides with the {\it $\F\-$localizing functor\/} $\loc_\F$ [11,~18.12.1] and  the point is that, according to Proposition~£2.17 above, the {\it natural map\/}
$\loc_{\,\bar\frak h^{^\frak M}}$ already can be lifted to a  {\it unique natural map\/}
$$\lambda_{\bar\L^{^{\frak N,\rm b}}} : \loc_\F = \loc_{\,\P}\too \loc_{\bar\L^{^{\frak N,\rm b}} }
\eqno £8.10.6\phantom{.}$$
fulfilling the conditions there. That is to say, for any  {\it $\F\-$chain\/} $\frak q\,\colon\Delta_n\to\F\,,$ we have a group homomorphism
 $$\lambda^{^{\!\frak N}}_\frak q = (\lambda_{\bar\L^{^{\frak N,\rm b}}} )_\frak q : \P(\hat\frak q)\too 
\bar\L^{^{\frak N,\rm b}} (\hat\frak q^{^\frak N})
\eqno £8.10.7$$
lifting $(\loc_{\,\bar\frak h^{^\frak M}})_\frak q$
which is compatible with the corresponding  structural functors and, according to definition~£8.6.3, is {\it unique\/} up to 
$\tilde\frak t^{^U}\!\big(\frak q (n)\big)\-$conjugation; analogously, for a second
{\it $\F\-$chain\/} $\frak r\,\colon \Delta_m\to \F$ and any
$\ch^*(\F)\-$morphism $(\mu,\delta)\,\colon  (\frak r,\Delta_m)\to (\frak q,\Delta_n)$ [11,~A2.8],
the diagram
$$\matrix{\P(\hat\frak q)&\buildrel \lambda^{^{\!\frak N}}_\frak q\over\too
 &\bar\L^{^{\frak N,\rm b}} (\hat\frak q^{^\frak N})\cr
\big\uparrow&\phantom{\Big\uparrow}&\big\uparrow\cr
\P(\hat\frak r)&\buildrel \lambda^{^{\!\frak N}}_\frak r\over\too 
&\bar\L^{^{\frak N,\rm b}} (\hat\frak r^{^\frak N})\cr}
\eqno £8.10.8\phantom{.}$$
 is commutative up to $\tilde\frak t^{^U}\!\big(\frak q (n)\big)\-$conjugation (cf.~£2.15).

\smallskip
In particular, assume that $n= 0\,,$ $m=1$ and $\delta = \delta^0_1\,,$ and setting $Q = \frak r (1)\,,$ $R = \frak r (0)\,,$  
$x = \hat\frak r^{^\frak N} (0\bullet 1)$ and 
$\hat\varphi = \hat\frak r (0\bullet 1)\,,$ assume that $\frak q (0) = R$
and that (cf.~£2.15)
$$\loc_\P (\mu,\delta) = \big(\tau_{_R}(1),{\rm id}_{\Delta_0}\big) \qq
\loc_{\bar\L^{^{\frak N,\rm b}} } = 
\big(\bar\tau^{^{\frak N,\rm b}}_{_R} (1),{\rm id}_{\Delta_0}\big)
\eqno £8.10.9;$$ 
then,  $\P(\hat\frak r)$ coincides with the  
stabilizer~$\P(Q)_{\hat\varphi}$ of $\varphi (R)$ in 
 $\P(Q)\,,$ $\bar\L^{^{\frak N,\rm b}} (\hat\frak r^{^\frak N})$ 
 coincides with the stabilizer $\bar\L^{^{\frak N,\rm b}} \!(Q)_x$ of $\varphi (R)$ in $\bar\L^{^{\frak N,\rm b}} \!(Q)$
 and diagram~£8.10.8 becomes
$$\matrix{\P(R)&\buildrel \lambda^{^{\!\frak N}}_{_R}\over\too&\bar\L^{^{\frak N,\rm b}} \!(R)\cr
\big\uparrow&\phantom{\Big\uparrow}
&\hskip-10pt{\scriptstyle \mu_x}\big\uparrow\cr
\P(Q)_{\hat\varphi}&\buildrel \lambda^{^{\!\frak N}}_x\over\too & \bar\L^{^{\frak N,\rm b}} \!(Q)_x\cr}
\eqno £8.10.10\phantom{.}$$
where $\mu_x\,\colon \bar\L^{^{\frak N,\rm b}}\! (Q)_x\to \bar\L^{^{\frak N,\rm b}}\!(R)$ sends 
$a\in  \bar\L^{^{\frak N,\rm b}}\! (Q)_x$  to the unique $b\in \bar\L^{^{\frak N,\rm b}}\! (R)$  fulfilling $x\. b = a\.x\,;$ moreover, since $\P(Q)_{\varphi}$ and $\bar\L^{^{\frak N,\rm b}}\!(Q)_x$  are respectively contained in 
$\P(Q)$ and $\bar\L^{^{\frak N,\rm b}}\!(Q)\,,$ and since $\bar\L^{^{\frak N,\rm b}}\!(Q)_x$ contains ${\rm Ker}(\bar\pi^{^{\frak N,\rm b}}_{_Q})\,,$ we actually may assume that $\lambda^{^{\!\frak N}}_x$ is just the restriction of~$\lambda^{^{\!\frak N}}_Q\,;$ then note that, for some choice of $x$ lifting $\bar\frak h^{^\frak M}(\hat\varphi)\,,$ diagram~£8.10.10 becomes commutative.
\eject

\smallskip
Consider the actions of $\P(Q)\times\P (R)$ on $\bar\L^{^{\frak M,\rm b}}\! (Q,R)$ and on 
$\bar\L^{^{\frak N,\rm b}}\! (Q,R)$  defined  by the composition on the left- and  the right-hand respectively {\it via\/}
 the functor $\bar\frak h^{^\frak M}$ and  {\it via\/} the group homomorphisms
$$\lambda^{^{\!\frak N}}_Q : \P (Q)\too \bar\L^{^{\frak N,\rm b}} \!(Q)\qq 
\lambda^{^{\!\frak N}}_R : \P (R)\too \bar\L^{^{\frak N,\rm b}}\!(R)
\eqno £8.10.11;$$ 
 for any $\hat\varphi\in \P(Q,R)\,,$ choose a lifting  $x_{\hat\varphi}\in \bar\L^{^{\frak N,\rm b}}\! (Q,R)$
of $\bar\frak h^{^\frak M}(\hat\varphi)$
such that the corresponding diagram £8.10.10 is commutative; then, we claim that we have the equality of stabilizers
$$\big(\P(Q)\times\P (R)\big)_{x_{\hat\varphi}} =  \big(\P(Q)\times\P(R)\big)_{\bar\frak h^{^\frak M}(\hat\varphi)}
\eqno £8.10.12.$$
Indeed, since $x_{\hat\varphi}$ lifts $\bar\frak h^{^\frak M}(\hat\varphi)\,,$ the inclusion of the left-hand member
in the right-hand one is clear; conversely, for any pair $(\hat\alpha,\hat\beta)\in 
\big(\P(Q)\times\P (R)\big)_{\bar\frak h^{^\frak M}(\hat\varphi)}\,,$ we have 
$\hat\alpha\.\bar\frak h^{^\frak M}(\hat\varphi) = \bar\frak h^{^\frak M}(\hat\varphi) \.\hat\beta\,;$  in
particular, denoting by $\pi\,\colon \P\to \F$ the second structural functor, we get 
$$\pi_{_Q} (\hat\alpha)\circ\bar\pi^{^{\frak M,\rm b}}_{_{Q,R}}(\hat\varphi) = 
\bar\pi^{^{\frak M,\rm b}}_{_{Q,R}}(\hat\varphi)\circ \pi_{_R}(\hat\beta)
\eqno £8.10.13\phantom{.}$$
and therefore $\hat\alpha$ belongs to $\P (Q)_{\hat\varphi}\,;$ then, since we assume that  the  cor-responding diagram~£8.10.10 is commutative, we still get
$$\mu_{x_{\hat\varphi}}\big(\lambda^{^{\!\frak N}}_Q(\hat\alpha)\big) = \lambda^{^{\!\frak N}}_R(\hat\beta)
\eqno £8.10.14,$$
which amounts to saying that $x_{\hat\varphi}\.\lambda^{^{\!\frak N}}_R(\hat\beta) = 
\lambda^{^{\!\frak N}}_Q(\hat\alpha)\.x_{\hat\varphi}\,,$
so that $(\hat\alpha,\hat\beta)$ belongs to~$\big(\P(Q)\times\P (R)\big)_{x_{\hat\varphi}}\,.$

\smallskip
This allows us to choose a family of liftings $\{x_{\hat\varphi}\}_{\hat\varphi}\,,$ where $\hat\varphi$ runs 
over the set of $\P\-$morphisms, which is compatible with $\P\-$isomorphisms. Precisely,
 choose a set of representatives $\X$ for the set of $\P\-$isomorphism classes of subgroups of $P\,,$
 for any  pair of  subgroups $Q$ and $R$ in~$\X$ choose a set of  representatives $\P_{Q,R}$ in $\P (Q,R)$ 
 for the set of  $\P (Q)\times\P (R)\-$orbits, and for~any $\hat\varphi\in \P (Q,R)$ choose a lifting  
$x_{\hat\varphi}\in \bar\L^{^{\frak N,\rm b}}\! (Q,R)$ such that the corresponding equality~£8.10.12  holds;
thus, any subgroup $Q$ of $P$ determines a unique $\bar Q$ in $\X$ and, moreover, we choose a 
$\P\-$isomorphism $\hat\omega_Q\,\colon Q\cong \bar Q$ and a lifting 
$x_Q\in \bar\L^{^{\frak N,\rm b}}\! (\bar Q,Q)$ of $\hat\omega_Q$ such that the corresponding equality~£8.10.12  holds. Hence,  any  $\P\-$morphism $\hat\varphi\,\colon R\to Q$ determines $\bar Q,\bar R\in \X$ and 
$\skew3\hat{\bar\varphi}\in \P_{\bar Q,\bar R}$
fulfilling 
$$\hat\varphi = \hat\omega_Q^{-1}\. \hat\alpha\.\skew3\hat{\bar\varphi}\.\hat\beta\. \hat\omega_R
\eqno £8.10.15\phantom{.}$$
for suitable $\hat\alpha\in \P (\bar Q)$ and $\hat\beta\in \P (\bar R)\,,$ and then we define
$$x_{\hat\varphi }= x_Q^{-1}\. \lambda^{^{\!\frak N}}_{\bar Q}(\hat\alpha)\. x_{\hat\varphi}
\. \lambda^{^{\!\frak N}}_{\bar R}(\hat\beta)\. x_R
\eqno £8.10.16.$$
At this point, it is routine to check that 
\smallskip
\noindent
£8.10.17\quad {\it We have $x_{\hat\alpha\circ \hat\varphi\circ \hat\beta} =  
x_{\hat\alpha}\.x_{\hat\varphi}\. x_{\hat\beta}$ for any  $\P\-$isomorphisms $\hat\alpha\in \P (Q',Q)$ 
and  $\hat\beta\in \P (R,R')\,.$\/}\smallskip
\noindent
Note that, for any  subgroup $Q$  of $P$ and any $\hat\alpha\in \P (Q)\,,$ we may assume that 
$$x_{\hat\alpha} = \lambda^{^{\!\frak N}}_Q(\hat\alpha)
\eqno £8.10.18;$$
in particular, since the section $\lambda^{^{\!\frak N}}_Q$ is compatible with the first structural functor 
$\tau\,\colon \P\to \F\,,$ for any $u\in N_P (Q)$ we have $x_{\tau_{_Q}(u)} = \bar\tau^{^{\frak N,\rm b}}_{_Q} (u)\,.$

\smallskip
Then, for any triple of subgroups $Q\,,$ $R$ and $T$ of $P\,,$ and any pair of $\P\-$morphisms 
$\hat\psi\,\colon T\to R$ and $\hat\varphi\,\colon R\to Q\,,$ since $ x_{\hat\varphi}\.x_{\hat\psi}$ and 
$x_{\hat\varphi\circ\hat\psi}$ have the same image $\bar\frak h^{^\frak M} (\hat\varphi\.\hat\psi)$ in 
$\bar\L^{^{\frak M,\rm b}}\!(Q,T)\,,$ the {\it divisibility\/} of  $\bar\L^{^{\frak N,\rm b}}$ guarantees 
the existence and the uniqueness of $t_{\hat\varphi,\hat\psi}\in \tilde\frak t^{^U} (T)$ fulfilling
$$x_{\hat\varphi}\.x_{\hat\psi} = x_{\hat\varphi\.\hat\psi}\.t_{\varphi,\psi}
\eqno £8.10.19.$$
That is to say, we get a correspondence mapping any {\it $\P\-$chain\/} 
$\frak q\,\colon \Delta_2\to \P$  on~$ t_{\frak q (0\bullet 1),\frak q(1\bullet 2)}$ 
and,  considering the {\it contravariant\/} functor $\tilde\frak t^{^U}$ (cf.~£8.6.3)
 and setting
$$\Bbb C^n (\tilde\F,\tilde\frak t^{^U}) = 
\prod_{\tilde\frak q\in \Fct(\Delta_n,\tilde\F)}\tilde\frak t^{^U}(\tilde\frak q (0))
\eqno £8.10.20\phantom{.}$$
 for any $n\in \Bbb N\,,$  we claim that this correspondence determines a {\it stable\/} element $t$ of 
$ \Bbb C^2 (\tilde\F,\tilde\frak t^{^U}) $ [11,~A3.18].

\smallskip 
Indeed, for another isomorphic  {\it $\P\!\-$chain\/}  $\frak q'\,\colon \Delta_2\to \P$ and a {\it natural isomorphism\/} $\nu\,\colon \frak q\cong\frak q'\,,$
setting 
$$\eqalign{\hat\psi = \frak q (0\bullet 1)\!\!\quad,\quad \!\!\hat\varphi = \frak q(1\bullet 2)\!\!\quad,\quad  
\!\!&\hat\psi' = \frak q'(0\bullet 1)\!\!\quad,\quad \!\!\hat\varphi' = \frak q'(1\bullet 2)\cr
\nu_0 = \hat\gamma\quad,\quad  \nu_1 = \hat\beta&\qq \nu_2 = \hat\alpha\cr}
\eqno £8.10.21,$$
  from~statement~£8.10.17 we have
 $$x_{\hat\varphi'} = x_{\hat\alpha}\.x_{\hat\varphi}\.x_{\hat\beta}^{-1} \quad \!\! , \!\! \quad 
 x_{\hat\psi'} = x_{\hat\beta}\.x_{\hat\psi}\.x_{\hat\gamma}^{-1}\!\!\qq\!\! x_{\hat\varphi'\.\hat\psi'} = 
x_{\hat\alpha}\.x_{\hat\varphi\.\hat\psi}\.x_{\hat\gamma}^{-1}
 \eqno £8.10.22\phantom{.}$$
 and therefore we get (cf.~£8.6.3)
$$\eqalign{x_{\hat\varphi'\.\hat\psi'}\.  t_{\hat\varphi',\hat\psi'} 
&= x_{\hat\varphi'}\. x_{\hat\psi'} = (x_{\hat\alpha}\.x_{\hat\varphi}\.x_{\hat\beta}^{-1})\. 
(x_{\hat\beta}\.x_{\hat\psi}\.x_{\hat\gamma}^{-1})\cr
& = x_{\hat\alpha}\.(x_{\hat\varphi\.\hat\psi}\.  t_{\hat\varphi,\hat\psi})\.x_{\hat\gamma}^{-1} 
=  x_{\hat\varphi'\.\hat\psi'}\. \big(\tilde t^{^U}\!(\tilde{\hat\gamma}^{-1})\big)( t_{\hat\varphi,\hat\psi})\cr}
\eqno £8.10.23\phantom{.}$$ 
which proves  that the correspondence $t$ sending $(\skew3\tilde{\hat\varphi},\skew4\tilde{\hat\psi})$ to~$t_{\hat\varphi,\hat\psi}$ is {\it stable\/} and, in particular,  that $t_{\hat\varphi,\hat\psi}$ only depends on the corresponding 
$\tilde\F\-$morphisms $\skew3\tilde{\hat\varphi}$ and~$\skew4\tilde{\hat\psi}\,.$
\eject

\smallskip
Moreover, considering the usual differential map
$$d_2 : \Bbb C^2 (\tilde\F,\tilde\frak t^{^U})
\too \Bbb C^3 (\tilde\F,\tilde\frak t^{^U})
\eqno £8.10.24,$$
we claim that $d_2 (t) = 0\,;$ indeed, for a third $\P\-$morphism $\hat\eta\,\colon W\to T$
we get
$$\eqalign{(x_{\hat\varphi}\.x_{\hat\psi})\.x_{\hat\eta} 
&= (x_{\hat\varphi\.\hat\psi}\.t_{\hat\varphi,\hat\psi})\.x_{\hat\eta}
= (x_{\hat\varphi\.\hat\psi}\.x_{\hat\eta})\.\big(\tilde\frak t^{^U}(\skew3\tilde{\hat\eta})\big)
(t_{\hat\varphi,\hat\psi})\cr
&= x_{\hat\varphi\.\hat\psi\.\hat\eta}\.t_{\hat\varphi\.\hat\psi,\hat\eta}\.\big(\tilde\frak t^{^U}
(\skew3\tilde{\hat\eta})\big)(t_{\hat\varphi,\hat\psi})\cr
x_{\hat\varphi}\.(x_{\hat\psi}\.x_{\hat\eta}) &= x_{\hat\varphi}\.(x_{\hat\psi\.\hat\eta}\.t_{\hat\psi,\hat\eta}) = 
x_{\hat\varphi\.\hat\psi\.\hat\eta}\.t_{\hat\varphi,\hat\psi\.\hat\eta}\.t_{\hat\psi,\hat\eta}\cr}
\eqno £8.10.25\phantom{.}$$
and   the {\it divisibility\/} of $\bar\L^{^{\frak M,\rm b}}\!$ forces
$$t_{\hat\varphi\.\hat\psi,\hat\eta}\.\big(\tilde\frak t^{^U}
(\skew3\tilde{\hat\eta})\big)(t_{\hat\varphi,\hat\psi}) = t_{\hat\varphi,\hat\psi\.\hat\eta}\.t_{\hat\psi,\hat\eta}
\eqno £8.10.26;$$
since $\tilde\frak t^{^U}\!(W)$ is Abelian, in the additive notation we obtain
$$0 = \big(\tilde\frak t^{^U}
(\skew3\tilde{\hat\eta})\big)(t_{\hat\varphi,\hat\psi}) - t_{\hat\varphi,\hat\psi\.\hat\eta}
+ t_{\hat\varphi\.\hat\psi,\hat\eta} - t_{\hat\psi,\hat\eta}
\eqno £8.10.27,$$
proving our claim.

\smallskip
At this point, since $ \Bbb H_*^2 (\tilde\F,\tilde\frak t^{^U}) = 0$ (cf.~£8.6.4),   we have $t = d_1 (s)$ for some
element $s = (s_{\tilde\frak r})_{\tilde\frak r\in \Fct(\Delta_1,\tilde\F)}$ in 
$\Bbb C_*^1 (\tilde\F,\tilde\frak t^{^U})\,;$ that is to say, with the notation above we get
$$t_{\hat\varphi,\hat\psi} = \big(\tilde\frak t^{^U}(\skew3\tilde{\hat\psi})\big) (s_{\skew3\tilde{\hat\varphi}}) \.(s_{\skew3\tilde{\hat\varphi}\circ\skew3\tilde{\hat\psi}})^{-1}\.s_{\skew3\tilde{\hat\psi}}
\eqno £8.10.28\phantom{.}$$
where we  identify any $\tilde\F\-$morphism with the obvious {\it $\tilde\F\-$chain\/}
$\Delta_1\to \tilde\F\,;$ hence, from equality~£8.10.19 we obtain
$$\eqalign{\big(x_{\hat\varphi}\.(s_{\skew3\tilde{\hat\varphi}})^{-1}\big)\.\big(x_{\hat\psi} \.(s_{\skew3\tilde{\hat\psi}})^{-1}\big)
&=    (x_{\hat\varphi}\.x_{\hat\psi})\. \Big(\big(\tilde\frak t^{^U}(\skew3\tilde{\hat\psi})\big)
(s_{\skew3\tilde{\hat\varphi}})\.s_{\skew3\tilde{\hat\psi}}\Big)^{-1}\cr
&= x_{\hat\varphi\.\hat\psi}\.(s_{\skew3\tilde{\hat\varphi}\.\skew3\tilde{\hat\psi}})^{-1}\cr}
\eqno £8.10.29,$$
which amounts to saying that the correspondence sending $\hat\varphi\in \P (Q,R)$ to 
$x_{\hat\varphi}\.(s_{\skew3\tilde{\hat\varphi}})^{-1}\in \bar\L^{^{\frak N,\rm b}}(Q,R)$ defines  a 
functor $\bar\frak h^{^\frak N}\,\colon \P\to \bar\L^{^{\frak N,\rm b}}$ lifting $\bar\frak h^{^\frak M}\,.$  
Note that in the case that $Q = R = T$ and  that $\hat\varphi$ and $\hat\psi$ are both ``inner'' $\P\-$automorphisms then 
equality~£8.10.28 forces $s_{\skew3\tilde{\hat\varphi}} = 1$ and from~£8.10.18 we get 
$x_{\hat\varphi} = \bar\tau^{^{\frak N,\rm b}}_{_Q} (u)$ for some $u\in N_P (Q)\,.$

\smallskip
We can  modify this correspondence in order to get an {\it $\F\-$locality functor\/}; indeed, for any $u\in T_P (R,Q)\,,$ 
the $\bar\L^{^{\frak N,\rm b}}\!\-$morphisms $\bar\frak h^{^\frak N}\big(\tau_{_{Q,R}}(u)\big)$ and 
$\bar\tau^{^{\frak N,\rm b}}_{_{Q,R}}(u)$ both lift $\bar\tau^{^{\frak M,\rm b}}_{_{Q,R}}(u)\,;$ once again,  the 
{\it divisibility\/} of $\bar\L^{^{\frak N,\rm b}}$ guarantees the existence and the uniqueness of 
$\ell_u\in {\rm Ker}(\bar\frak l^{^{\frak N,\frak M}}_{_R})$ fulfilling
$$\bar\tau^{^{\frak N,\rm b}}_{_{Q,R}}(u) = \bar\frak h^{^\frak N}\big(\tau_{_{Q,R}}(u)\big)\.\ell_u
\eqno £8.10.30\phantom{.}$$
\eject
\noindent
and the remark above proves that $\ell_u$ only depends on  $\tilde u\in \tilde\F_{\!P}(Q,R)\,.$ Moreover, for a second 
$\T_{\! P}\-$morphism $v\,\colon T\to R\,,$ we get
$$\eqalign{\bar\frak h^{^\frak N}\big(\tau_{_{Q,T}}(u v)\big)\.\ell_{\tilde u \tilde v} 
&= \bar\tau^{^{\frak N,\rm b}}_{_{Q,T}}(uv) = 
\bar\tau^{^{\frak N,\rm b}}_{_{Q,R}}(u)\.\bar\tau^{^{\frak N,\rm b}}_{_{R,T}}(v)\cr
&= \bar\frak h^{^\frak N}\big(\tau_{_{Q,R}}(u)\big)\.\ell_{\tilde u}
\.\bar\frak h^{^\frak N}\big(\tau_{_{R,T}}(v)\big)\.\ell_{\tilde v}\cr
&= \bar\frak h^{^\frak N}\big(\tau_{_{Q,T}}(u v)\big)\.\big(\tilde\frak t^{^U}(\tilde v)\big)(\ell_{\tilde u})\.\ell_{\tilde v}\cr}
\eqno £8.10.31.$$

\smallskip
 The {\it divisibility\/} of $\bar\L^{^{\frak N,\rm b}}\!$ always  forces
$\ell_{\tilde u \tilde v}  = \big(\tilde\frak t^{^U}\!(\tilde v)\big)(\ell_{\tilde u})\.\ell_{\tilde v}$
and, since $\tilde\frak t^{^U}(T)$~is Abelian, in the additive notation we obtain
$$0 = \big(\tilde\frak t^{^U}\!(\tilde v)\big)(\ell_{\tilde u}) - \ell_{\tilde u\circ\tilde v}+ \ell_{\tilde v}
\eqno £8.10.32;$$
that is to say,  denoting by  $\frak i_P\,\colon \tilde\F_{\! P}\i \tilde\F$ the obvious inclusion functor, the correspondence $\ell$ sending any  $\tilde\F_{\! P}\-$morphism $\tilde u\colon R\to Q$ to~$\ell_{\tilde u}$ defines a {\it $1\-$cocycle\/} in 
$\Bbb C^1\big(\tilde\F_{\! P},\tilde\frak t^{^U}\!\circ\frak i\big)\,;$ but, since the category 
$\tilde\F_{\! P}$ has a final object, we actually have [11,~Corollary~A4.8]
$$\Bbb H^1\big(\tilde\F_{\! P},\tilde\frak t^{^U}\!\circ\frak i_P\big) = \{0\}
\eqno £8.10.33;$$
consequently, we obtain $\ell = d_0 (z)$ for some element $z = (z_Q)_{Q}$ in 
$$\Bbb C^0\big(\tilde\F_{\! P},\tilde\frak t^{^U}\!\circ\frak i_P\big) = \Bbb C^0\big(\tilde\F,\tilde\frak t^{^U}\big)
\eqno £8.10.34.$$
In conclusion, equality~£8.10.30 becomes
$$\bar\tau^{^{\frak N,\rm b}}_{_{Q,R}}(u)  = 
\bar\frak h^{^\frak N}\big(\tau_{_{Q,R}}(u)\big)\.\big(\tilde\frak t^{^U}\!(\tilde u)\big)(z_Q)\.z_R^{-1} 
= z_Q\.\bar\frak h^{^\frak N}\big(\tau_{_{Q,R}}(u)\big)\.z_R^{-1}
\eqno £8.10.35\phantom{.}$$
and therefore  the correspondence  sending $\varphi\in \F (Q,R)$ to  
$z_Q\. \bar\frak h^{^\frak N} (\varphi)\.z_R^{-1}$ defines  an {\it $\F\-$locality functor\/}   
 from $\P$ to $\bar\L^{^{\frak N,\rm b}}$  lifting~$\bar\frak h^{^\frak M}\,.$

\smallskip
Assume that we have two {\it $\F\-$locality functors\/} $\bar\frak h^{^\frak N}$ and
$\bar\frak h'^{^\frak N}$  from $\P$ to $\bar\L^{^{\frak N,\rm b}}$  lifting~$\bar\frak h^{^\frak M}\,;$ 
the uniqueness  of the {\it natural map\/} $\lambda_{\bar\L^{^{\frak N,\rm b}}}$ in~£8.10.6 already guarantees that, 
in order to prove that $\bar\frak h'^{^\frak N}$ is {\it naturally $\F\-$isomorphic\/} to $\bar\frak h^{^\frak N}\,,$ we may assume that 
$$\bar\frak h'^{^\frak N} (\hat\alpha) = \lambda^{^{\!\frak N}}_Q(\hat\alpha) = \bar\frak h^{^\frak N} (\hat\alpha)
\eqno £8.10.36\phantom{.}$$
for any subgroup $Q$ of $P$ and any $\hat\alpha\in \P (Q)\,;$ more precisely, we may assume that $\bar\frak h'^{^\frak N}\!(\alpha)  = \bar\frak h^{^\frak N}\! (\alpha)$ for any 
$\F^{^{\frak X}}\!\-$isomorphism $\alpha\in \F (Q',Q)\,.$ For any $\P\-$morphism $\hat\varphi\,\colon R\to Q\,,$
set $x_{\hat\varphi} =  \bar\frak h^{^\frak N}(\hat\varphi)$ and $x'_{\hat\varphi} =  
\bar\frak h'^{^\frak N}(\hat\varphi)$
for short; the {\it divisibility\/} of $\bar\L^{^{\frak N,\rm b}}$
forces again the existence of a unique $s_{\hat\varphi}\in \tilde\frak t^{^U}\! (R))$ fulfilling 
$x'_{\hat\varphi} = x_{\hat\varphi}\.s_{\hat\varphi}\,;$ note that equality~£8.10.36 forces $s_{\hat\alpha} = 1\,.$
That is to say, we get a correspondence mapping any 
{\it $\P\-$chain\/}  $\frak r\,\colon \Delta_1\to \P$   on~$s_{\frak r (0\bullet 1)}$ 
and  we claim that this correspondence determines a {\it stable\/} element $s$ of 
$ \Bbb C^1 (\tilde\F,\tilde\frak t^{^U}) $ [11,~A3.18].

\smallskip 
Indeed, for another isomorphic {\it $\P\-$chain\/} $\frak q'\,\colon \Delta_1\to \P$ 
and a {\it natural isomorphism\/} $\nu\,\colon \frak q\cong\frak q'\,,$
setting 
$$\eqalign{\hat\varphi = \frak q (0\bullet 1)\!\!\quad,\quad  \!\!\hat\varphi' = \frak q'(0\bullet 1)
\!\!\quad,\quad  \!\!\nu_0 = \hat\beta\qq \nu_1 = \hat\alpha\cr}
\eqno £8.10.37,$$
  from~our choice we have $s_{\hat\alpha} = 1$ and $s_{\hat\beta} = 1$ and therefore we get
$$\eqalign{x'_{\hat\varphi'} &= x_{\hat\varphi'}\.s_{\hat\varphi'} = 
(x_{\hat\alpha}\.x_{\hat\varphi}\.x_{\hat\beta}^{-1})\.s_{\hat\varphi'}
=  (x_{\hat\alpha}\.x'_{\hat\varphi}\.s_{\hat\varphi}^{-1}\.x_{\hat\beta}^{-1})\.s_{\hat\varphi'}\cr
& =  (x'_{\hat\alpha}\.x'_{\hat\varphi}\.x'^{-1}_{\hat\beta})\. 
\big(\tilde\frak t^{^U}\!(\skew3\tilde{\hat\beta}^{-1})\big)(s_{\hat\varphi}^{-1})\.s_{\hat\varphi'}\cr
& = x'_{\hat\varphi'}\. \big(\tilde\frak t^{^U}\!(\skew3\tilde{\hat\beta}^{-1})\big)
( s_{\hat\varphi}^{-1})\.s_{\hat\varphi'}\cr}
\eqno £8.10.38\phantom{.}$$ 
which proves  that the correspondence $s$ sending $\skew3\tilde{\hat\varphi}$ to~$s_{\hat\varphi}$ is {\it stable\/}
and, in particular,  that $s_{\hat\varphi}$ only depends on the corresponding 
$\tilde\F\-$morphism $\skew3\tilde{\hat\varphi}\,.$

\smallskip
Moreover, we also claim that $d_1 (s) = 0\,;$ indeed, for a second $\P\-$morphism 
$\hat\psi\,\colon T\to R$
we get
$$\eqalign{x'_{\hat\varphi\.\hat\psi}& = x'_{\hat\varphi}\.x'_{\hat\psi} = 
(x_{\hat\varphi}\.s_{\hat\varphi})\.(x_{\hat\psi}\.s_{\hat\psi})
= x_{\hat\varphi\.\hat\psi}\.\big(\tilde\frak t^{^U}\!(\skew4\tilde{\hat\psi})\big)(s_{\hat\varphi})\.s_{\hat\psi}\cr}
\eqno £8.10.39\phantom{.}$$
and   the {\it divisibility\/} of $\bar\L^{^{\frak N,\rm b}}$ forces $s_{\hat\varphi\.\hat\psi} 
=\big(\tilde\frak t^{^U}\!(\skew4\tilde{\hat\psi})\big)(s_{\hat\varphi})\.s_{\hat\psi}\,;$
since $\tilde\frak t^{^U}\!(T)$ is Abelian, in the additive notation we obtain
$$0 = \big(\tilde\frak t^{^U}\!(\skew4\tilde{\hat\psi})\big)(s_{\hat\varphi}) - s_{\hat\varphi\.\hat\psi}  + s_{\hat\psi}
\eqno £8.10.40,$$
proving our claim.

\smallskip
Finally, since $ \Bbb H_*^1 (\tilde\F,\tilde\frak t^{^U}) = 0$ (cf.~£8.6.4),   we have $t = d_0 (n)$ for some
element $n = (n_Q)_Q$ in  $\Bbb C^0 (\tilde\F,\tilde\frak t^{^U})$ where we  identify any 
subgroup $Q$ of $P$  with the obvious {\it $\tilde\F\-$chain\/}
$\Delta_0\to \tilde\F\,;$ that is to say, with the 
notation above we get
$$s_{\hat\varphi} = \big(\tilde\frak t^{^U}\!(\skew3\tilde{\hat\varphi})\big) (n_Q)\.n_R^{-1}
\eqno £8.10.41\phantom{.}$$
 hence,  we obtain
$$\bar\frak h'^{^\frak N}(\hat\varphi) = x'_{\hat\varphi} =   x_{\hat\varphi}\.\big(\tilde\frak t^{^U}\!(\skew3\tilde{\hat\varphi})\big)  (n_Q)\.n_R^{-1} = n_Q\.\bar\frak h^{^\frak N}(\hat\varphi)\.n_R^{-1}
\eqno £8.10.42,$$
which amounts to saying that the correspondence sending $Q$ to 
$n_Q$ defines  a {\it natural $\F\-$isomorphism\/} between $\bar\frak h^{^\frak N}$ and $\bar\frak h'^{^\frak N}\,.$ 
We are done.

\bigskip
 \bigskip
\noindent
{\bf £9. Functoriality of the perfect $\F\-$locality }

\medskip 
£9.1. It remains to dicuss the {\it functoriality\/} of the {\it perfect $\F\-$locality\/} $\P\,;$  let $P'$ be a second finite $p\-$group, $\F'$ a Frobenius $P'\-$category and $\P'$ the corresponding {\it perfect $\F'\-$locality\/}, and denote by
$$\tau' : \T_{P'}\too \P'\qq \pi' : \P'\too \F'
\eqno £9.1.1\phantom{.} $$
the structural functors; let $\alpha\,\colon P\to P'$ be an $(\F,\F')\-$functorial
group homomorphism [11,~12.1]; recall that we have a so-called {\it Frobenius functor $\frak f_\alpha\,\colon \F\to \F'$\/} [11,~12.1], and let us  denote by  $\frak t_\alpha\,\colon \T_{P}\to \T_{P'}$ the functor induced by~$\alpha\,.$ In this section, replacing $\P$ and $\P'$ by the quotients
$$\bar\P = \P/[\frak c_\F^\frak h,\frak c_\F^\frak h]\qq 
\bar\P' = \P'/[\frak c_{\F'}^\frak h,\frak c_{\F'}^\frak h]
\eqno £9.1.2,$$
we prove that  there is a unique isomorphism class of functors 
$\bar\frak g_\alpha\,\colon \bar\P\to \bar\P'$  fulfilling 
$$\bar\tau'\circ \frak t_\alpha = \bar\frak g_\alpha\circ \tau \qq 
 \pi'\circ \bar\frak g_\alpha = \frak f_\alpha\circ \pi
\eqno £9.1.3;$$
as a consequence, if $P''$ is a third finite $p\-$group, $\F''$ a Frobenius $P''\-$category, $\P''$ the {\it perfect 
$\F''\-$locality\/}  and $\alpha'\,\colon P'\to P''$ an $(\F',\F'')\-$functorial group homomorphism, then  the functors $\bar\frak g_{\alpha'}\circ \bar\frak g_{\alpha}$ and $\bar\frak g_{\alpha'\circ\alpha}$ from $\bar\P$ to $\bar\P''$ are 
{\it naturally isomorphic\/}.

\medskip
£9.2. As a matter  of fact, assuming the existence of $\P\,,$ we already have 
proved in [11,~Theorem~17.18] the existence of all the possible 
{\it perfect quotients $\skew4\bar{\bar\P}$\/} of $\P\,;$  presently, this simplifies our work. Let us recall the construction  of~$\skew4\bar{\bar \P}\,;$ let $U$ be an {\it $\F\-$stable\/} subgroup of $P$ (cf.~£2.5), set $\skew4\bar{\bar P} = P/U$ and denote by $\skew4\bar{\bar\F}$ the {\it quotient Frobenius $\skew4\bar{\bar P}\-$category\/} $\F/U$ [11,~Proposition~12.3]; for any subgroup $Q$ of $P\,,$ denote by $\skew4\bar{\bar Q}$ the image of $Q$  in $\skew4\bar{\bar P}$ and by $U_{\F} (Q)$ the kernel of the canonical group homomorphism 
$\F (Q)\to \skew4\bar{\bar\F} (\skew4\bar{\bar Q})\,;$ moreover, if $Q$ is fully normalized in~$\F\,,$ for short we set
$$P^{^Q} = N_P^{U_\F (Q)} (Q)\qq   \F^{^Q} = N_\F^{U_\F (Q)} (Q)
\eqno £9.2.1,$$
so that $\F^{^Q}$ is a {\it Frobenius $P^{^Q}\-$category\/}; then,  in the group $\P (Q)\,,$ we set (cf.~£2.4)
$$U_\P (Q) = \Bbb O^p\Big(\pi_{_Q}^{-1}\big(U_\F (Q)\big)\Big)\.\tau_{_Q}
\big(N_U (Q)\. H_{\F^{^Q}}\big)
\eqno £9.2.2;$$
actually, {\it via\/} $\P\-$isomorphisms we can extend the definition of $U_\P (Q)$ to any subgroup $Q$ of $P\,.$ 
Then, $\skew4\bar{\bar\P} = \P/U$ is the {\it perfect  
$\skew4\bar{\bar\F}\-$locality\/} fulfilling [11,~17.15-17.17]
$$\skew4\bar{\bar\P} (\skew4\bar{\bar Q},\skew4\bar{\bar R}) = \P (Q,R)\big/U_\P (R)
\eqno £9.2.3\phantom{.}$$
for any pair of subgroups $Q$ and $R$ of $P\,.$
\eject

\medskip
£9.3. In the general case, setting $U = {\rm Ker}(\alpha)$ and $\skew4\bar{\bar P} = P/U\,,$ we have an injective group homomorphism  $\skew2\bar{\bar\alpha}\,\colon \skew4\bar{\bar P}\to P'$ and from [11,~Proposition~12.3] it is easily checked that we have a   {\it faithful Frobenius functor  $\frak f_{\skew2\bar{\bar\alpha}}\,\colon 
\skew4\bar{\bar\F}\to \F'$\/}; moreover, from the description above, it is easily checked that any functor
$\bar\frak g_\alpha\,\colon \bar\P\to \bar\P'$  fulfilling condition~£9.1.3 factorizes
throughout the quotient~$\skew4\bar{\bar\P}/[\frak c_{\skew4\bar{\bar\F}}^\frak h,
\frak c_{\skew4\bar{\bar\F}}^\frak h]\,;$ consequently, in order to prove the  existence and the uniqueness of 
$\bar\frak g_\alpha\,,$ it suffices to prove the existence and the uniqueness of a suitable  functor  
$\bar\frak g_{\skew2\bar{\bar\alpha}}\,\colon \skew4\bar{\bar\P}/[\frak c_{\skew4\bar{\bar\F}}^\frak h,
\frak c_{\skew4\bar{\bar\F}}^\frak h]\to \bar\P'\,.$ 
 Thus, we may
 assume that $\alpha$ is injective and, in this case, let us identify $P$ and $\F$ with their respective images in~$P'$ and $\F'\,.$ In order to relate $\P$ and $\P'\,,$
 we start by getting a relationship between the 
 {\it natural $\F^{^{\rm sc}}\-$locality\/} $\bar\L^{^{\rm n,sc}}$ and the 
  {\it basic  $\F'\-$locality\/}~$\L'^{^{\rm b}}\,;$  more explicitly, the converse image ${\rm Res}_{\F^{^{\rm sc}}}(\L'^{^{\rm b}})$ of  $\F^{^{\rm sc}}$ in~$\L'^{^{\rm b}}$  is clearly a {\it $p\-$coherent $\F^{^{\rm sc}}\-$locality\/} and we will exhibit  a {\it canonical  $\F^{^{\rm sc}}\-$locality functor\/}  from 
  $\bar\L^{^{\rm n,sc}}$ to a suitable {\it quotient $\F^{^{\rm sc}}\-$locality\/} of 
  ${\rm Res}_{\F^{^{\rm sc}}}(\L'^{^{\rm b}})\,.$

 \medskip
£9.4. Choose a {\it  natural $\F\-$basic $P\times P\-$set\/} $\Omega$ (cf.~£3.3 and~£3.5), denote by~$G$ the group 
of $\{1\}\times P\-$set automorphisms of~${\rm Res}_{\{1\}\times P}(\Omega)$  and identify the $p\-$group $P$ with the image of $P\times \{1\}$ in $G\,,$ so that for any pair of $\F\-$selfcentralizing subgroups $Q$ and $R$ of $P$ we have 
(cf.~£4.12 and~£4.14)
$$\L^{^{\Omega,\rm sc}}\! (Q,R) = T_{G}(R,Q)/\frak S_{\Omega}^1 (R)\qq
\bar\L^{^{\rm n,sc}}\! = \L^{^{\Omega,\rm sc}}\! /
\tilde\frak c^{^{\rm nsc}}
\eqno £9.4.1.$$
It is clear that $G$ acts faithfully on the $P\times P'\-$set $\Omega\times_{P} P'$  centralizing the action 
of~$\{1\}\times P'\,;$ we claim that 
\smallskip
\noindent
£9.4.2.\quad {\it There is a  {\it thick $\F'\-$basic  $P'\times P'\-$set\/}~$\Omega'$ such  that ${\rm Res}_{P\times P'}(\Omega')$ contains~$\Omega\times_{P} P'\,.$\/}
\smallskip
\noindent 
Indeed, any $P\times P'\-$orbit~$O''$ of $\Omega\times_{P} P'$ is isomorphic to  $(P\times P')/\Delta_{\varphi}(Q)$ for some subgroup $Q$ 
of~$P$ and some $\varphi\in \F (P,Q)\,,$  and therefore it is isomorphic to a $P\times P'\-$orbit of  {\it thick $\F'\-$basic  $P'\times P'\-$set\/}~$\Omega'\,;$ 
hence, up to replacing $\Omega'$ by the disjoint union of $k$ copies of $\Omega'$ for a suitable $k$ prime to~$p\,,$ we may assume that ${\rm Res}_{P\times P'}(\Omega')$ contains~$\Omega\times_{P} P'\,.$  Similarly, denote  by~$G'$ the group of $\{1\}\times P'\-$set automorphisms of ${\rm Res}_{\{1\}\!\times\!P'}(\Omega')$  and identify the $p\-$group~$P'$ with the image of $P'\times \{1\}$ in $G'\,;$ once again, for any pair of~subgroups~$Q'$ and $R'$ of $P'$ we have (cf.~£4.14)
$$\L'^{^{\rm b}} (Q',R') = T_{G'}(R',Q')/\frak S_{\Omega'} (R')
\eqno £9.4.3.$$
In particular,  for any subgroup $Q$ of $P\,,$  
$\big({\rm Res}_{\F}(\L'^{^{\rm b}})\big)(P,Q)$ is the set 
of classes $f'\frak S_{\Omega'} (Q)$ in $\L'^{^{\rm b}} (P,Q) $ where 
$$f' : {\rm Res}_{Q\times P'}(\Omega') \cong {\rm Res}_{\varphi\times {\rm id}_{P'}}
\big({\rm Res}_{P\times P'}(\Omega')\big)
\eqno £9.4.4\phantom{.}$$
is a $Q\times P'\-$set isomorphism for some $\varphi\in \F (P,Q)\,.$

\medskip
£9.5.  As in section~£6 above, for induction purposes we have to consider  a nonempty set~$\frak X$  of $\F\-$selfcentralizing subgroups of $P$ which contains any subgroup of $P$ admitting an $\F\-$morphism from some subgroup 
in~$\frak X\,;$  recall that $\Omega$ contains  the {\it natural $\F^{^\frak X}\!\!\-$basic $P\times P\-$set\/}~$\Omega^{^\frak X}$ (cf.~£3.5) 
and, for any $Q\in \frak X\,,$ denote by  $\,\Omega^{^\frak X}_{Q}$ the union of all the  $Q\times P\-$orbits of 
$\Omega$ isomorphic to $(Q\times P)/\Delta_\eta (T)$ for some $T\in \frak X$ and some $\eta\in \F (Q,T)$ (cf.~£5.18); then, for any 
$\varphi\in \F(P,Q)$ note that any   $Q\times P\-$set isomorphism
$$f : {\rm Res}_{Q\times P}(\Omega) \cong {\rm Res}_{\varphi\times {\rm id}_{P}}
\big({\rm Res}_{P\times P}(\Omega)\big)
\eqno £9.5.1\phantom{.}$$
maps $\,\Omega^{^\frak X}_{Q}$ onto $\,\Omega^{^\frak X}_{\varphi (Q)}\,.$
 Moreover,  let us denote by  $\hat\Omega^{^\frak X}_Q\i {\rm Res}_{Q\times P'}(\Omega')$ the union of all the  
 $Q\times P'\-$orbits of   $\Omega'$ which are isomorphic to some  $Q\times P'\-$orbit 
 of~$\Omega^{^\frak X}_Q\!\times_P P'\,.$
  
\bigskip
\noindent
{\bf Proposition~£9.6.} {\it For any  $Q,R\in \frak X\,,$ any  $\varphi\in \F(Q,R)$
and any $R\times P'\-$set isomorphism
$$f' : {\rm Res}_{R\times P'}(\Omega') \cong 
{\rm Res}_{\varphi\times {\rm id}_{P'}}\big({\rm Res}_{Q\times P'}(\Omega')\big)
\eqno £9.6.1\phantom{.}$$
we have $f'(\hat\Omega^{^\frak X}_{R}) \i \hat\Omega^{^\frak X}_{Q}\,.$\/}

\medskip
\noindent
{\bf Proof:} For any $T,U\in \frak X\,,$ any $\theta\in \F (R,U)$ and any 
$\eta'\in \F' (Q,T)\,,$ it suffices to prove that,
if there is an injective $R\times P'\-$set homomorphism
$$(R\times P')/\Delta_{\theta}(U) \too {\rm Res}_{\varphi\times 
{\rm id}_{P'}}\big( (Q\times P')/\Delta_{\eta'}(T)\big)
\eqno £9.6.2,$$
then there is $u'\in P'$ such that $T^{u'}$ is contained in $P$ and that 
$\eta'\circ \kappa_{_{T,T^{u'}}}(u')$\break
 belongs to $\F (Q,T^{u'})$ (cf.~£2.6.4); we argue by induction on 
$\vert P\colon R\vert$ and may assume that~$R\not= P\,.$ 
Since we have an injective $Q\times P'\-$set homomorphism
$$(Q\times P')/\Delta_{\eta'}(T)\too {\rm Res}_{\iota_Q^P\times 
{\rm id}_{P'}}\big( (P\times P')/\Delta_{\iota_Q^P\circ \eta'}(T)\big)
\eqno £9.6.3,$$
we may assume that $Q = P\,;$ since we still have  an injective 
$U\times P'\-$set homomorphism
$$(U\times P')/\Delta_{{\rm id}_U}(U) \too {\rm Res}_{\theta\times 
{\rm id}_{P'}}\big( (R\times P')/\Delta_{\theta}(U)\big)
\eqno £9.6.4,$$
we actually may assume that $U = R$ and $\theta = {\rm id}_R\,.$

\smallskip
In this case, we also may assume that 
$$(\varphi \times {\rm id}_{P'})\big(\Delta (R)\big)\i \Delta_{\eta'}(T)
\eqno £9.6.5,$$
so that $R$ is contained in $T$ and that $\eta'$ extends $\varphi\,;$ in particular, we still may assume that $R\not= T\,,$ 
setting $\widehat R = N_T (R)$ and  $\widehat  P = N_P (R)\,;$ then, it follows from [11,~2.10] that $\varphi$ can be extended to 
an $\F\-$morphism $\hat\varphi\,\colon \widehat  R\to P$ and that $R$ is actually {\it fully $\F_{\!P} (R)\-$normalized\/} in $\F\,.$
\eject

\smallskip
Since $\eta'$ and $\hat\varphi$ coincide over $R$ and we have 
$C_P (R) = Z(R)\,,$  we still have $\eta' (\widehat  R) = \hat\varphi (\widehat  R)\,;$ 
that is to say, there is $\sigma'\in \F' (\widehat  R)$ such that 
the restriction of $\eta'$ to $\widehat  R$ coincides with $\hat\varphi\circ \sigma'\,;$ thus, $\sigma'$ acts trivially on $R$ and therefore, denoting by $\F' (\widehat  R)_R$ the stabilizer of $R$ in $\F' (\widehat  R)\,,$  $\sigma'$ belongs to 
$\Bbb O_p\big(\F' (\widehat  R)_R\big)$ [4,~Ch.~5, Theorem~3.4].

\smallskip
On the other hand, it follows from [11,~Proposition~2.7] that there is 
$\zeta'\in \F' (P',\widehat  P)$ such that $\overline{\!R} = \zeta' (R)$ is {\it  fully normalized\/} in $\F'$ and, moreover,
 that~$\skew2\widehat {\overline{\!R}} = \zeta' (\widehat  R)$  is {\it  fully normalized\/} in  $N_{\F'}(\,\overline{\!R}\,)\,;$ then 
$\F_{\!P'} (\,\skew2\widehat {\overline{\!R}}\,)_{\bar R}$ is a Sylow $p\-$subgroup of $\F' (\,\skew2\widehat {\overline{\!R}}\,)_{\bar R}$ [11,~Proposition~2.11] and therefore it contains 
$\Bbb O_p\big(\F' (\,\skew2\widehat {\overline{\!R}}\,)_{\bar R}\big)\,;$ set~$\,\skew2\widehat {\overline{\!P}}\, 
= \zeta' (\widehat P)$  and respectively denote by  $\skew3\hat{\bar\varphi}\,\colon 
\,\skew2\widehat {\overline{\!R}}\,\to \,\skew2\widehat {\overline{\!P}}\,$ 
  and $\bar\sigma'\in \F' (\,\skew2\widehat {\overline{\!R}}\,)$  the items determined by $\hat\varphi$ and
 $\sigma'$ {\it via\/}  the group isomorphism $\zeta'_*\,\colon \widehat P\cong \skew2\widehat {\overline{\!P}}$ 
 induced  by~$\zeta'\,.$

 \smallskip
 In particular, $\bar\sigma'$ belongs to $\F_{\!P'}(\,\skew2\widehat {\overline{\!R}}\,)$ and  we have an injective 
 $\,\skew2\widehat {\overline{\!R}}\,\times P'\-$set homomorphism
$$(\,\skew2\widehat {\overline{\!R}}\,\times P')/\Delta (\,\skew2\widehat {\overline{\!R}}\,) \too 
{\rm Res}_{\skew3\hat{\bar\varphi}\times {\rm id}_{P'}}\big( (\,\skew2\widehat {\overline{\!P}}\,\times P')/
\Delta_{\skew3\hat{\bar\varphi}\circ \bar\sigma'}(\,\skew2\widehat {\overline{\!R}}\,)\big)
\eqno £9.6.6;$$
then, still denoting by $\hat\varphi\,\colon \widehat R\to \widehat P$ the $\F\-$morphism determined by $\hat\varphi$
above, the inverse of the group  isomorphism $\zeta'_*\times {\rm id}_{P'}$ determines an  injective 
$\widehat R\times P'\-$set homomorphism
$$(\widehat R\times P')/\Delta (\widehat R) \too  {\rm Res}_{\hat\varphi\times {\rm id}_{P'}}\big( (\widehat P\times P')/
\Delta_{\hat\varphi\circ \sigma'}(\widehat R)\big)
\eqno £9.6.7;$$
moreover, since the restriction of $\eta'$ to $\widehat R$ coincides with  $\hat\varphi\circ \sigma'\,,$ we have 
$$\Delta_{\hat\varphi\circ \sigma'}(\widehat R) =  (\widehat P\times P')\cap \Delta_{\eta'}(R)
\eqno £9.6.8\phantom{.}$$
and therefore we still have  an  injective $\widehat R\times P'\-$set homomorphism
$$(\widehat R\times P')/\Delta (\widehat R) \too {\rm Res}_{\hat\varphi\times 
{\rm id}_{P'}}\big( ( P\times P')/\Delta_{\eta'}(T)\big)
\eqno £9.6.9;$$
finally, it suffices to apply the induction hypothesis.

\bigskip
\noindent
{\bf Proposition~£9.7.} {\it For any $Q\in \frak X$ the inclusion of the set of
$Q\times P'\-$orbits in $\Omega^{^\frak X}_Q\times_P P'$ in the set of
$Q\times P'\-$orbits in $\hat\Omega^{^\frak X}_Q$ admits a section 
$s^{_\frak X}_Q$ such that,  for any $Q\times P'\-$orbit $\hat O$ in 
$\hat\Omega^{^\frak X}_Q\,,$ we have a $Q\times P'\-$set isomorphism
$s^{_\frak X}_Q (\hat O)\cong \hat O$ and that,  for any $Q\times P\-$orbit $O$ in 
$\Omega^{^\frak X}_Q\,,$ any $\varphi\in \F (P,Q)$ and any $Q\times P\-$set isomorphism
$$f : {\rm Res}_{Q\times P}(\Omega) \cong {\rm Res}_{\varphi\times {\rm id}_{P}}
\big({\rm Res}_{P\times P}(\Omega)\big)
\eqno £9.7.1,$$
setting $\bar Q = \varphi (Q)$ and $\bar O = f(O)\,,$ we have 
$$\big\vert (s^{_\frak X}_{\bar Q})^{-1}(\bar O\times_P P')\big\vert = 
\big\vert (s^{_\frak X}_{ Q})^{-1}( O\times_P P')\big\vert
\eqno £9.7.2.$$\/}
\eject

\par
\noindent
{\bf Proof:} Let $\hat O$ be a $Q\times P'\-$orbit in $\hat\Omega^{^\frak X}_Q$
and $\hat\omega$ an element of $\hat O\,;$ since $\hat O$ is isomorphic to
$O\times_P P'$ for some $Q\times P\-$orbit $O$ in $\Omega^{^\frak X}_Q\,,$
there is $u'\in P'$ such that the stabilizer of $\hat\omega\.u'^{-1}$ in $Q\times P'$
coincides with the stabilizer $(Q\times P)_\omega$ in $Q\times P$ of some 
$\omega\in O\,;$ but, we have $(Q\times P)_\omega = \Delta_\eta (T)$ for
some $T\in \frak X$ and some $\eta\in \F(Q,T)$ and therefore, according to 
Proposition~£5.17, this stabilizer determines the set $Z(T)\.\omega\,.$
Moreover, if $v'$ is another element of~$P'$ such that the stabilizer 
of~$\hat\omega\.v'^{-1}$ coincides with the stabilizer of some element of $O$
then, since this stabilizer coincides with the stabilizer of $\omega\.(u'v'^{-1})\,,$
the element $u'v'^{-1}$ belongs to $P$ and thus, in $\Omega^{^\frak X}_Q\times_P P'\,,$
we have $(\omega,u') = \big(\omega\.(u'v'^{-1}),v'\big)\,.$

\smallskip 
In conclusion, the correspondence mapping $\hat O$ on $O\times_P P'$ defines a map $s^{_\frak X}_Q$ from  the set 
of $Q\times P'\-$orbits in $\hat\Omega^{^\frak X}_Q$ to  the set of $Q\times P'\-$orbits in 
$\Omega^{^\frak X}_Q\times_P P'\,,$ and it is clear that 
$s^{_\frak X}_Q (O\times_P P') = O\times_P P'\,.$ Finally, for  any $\varphi\in \F (P,Q)$ and any $Q\times P\-$set isomorphism
$$f : {\rm Res}_{Q\times P}(\Omega) \cong {\rm Res}_{\varphi\times {\rm id}_{P}}
\big({\rm Res}_{P\times P}(\Omega)\big)
\eqno £9.7.3,$$
setting $\bar Q = \varphi (Q)$ and denoting by $\varphi_*\,\colon Q\cong \bar Q$ the isomorphism induced by~$\varphi\,,$ we have an obvious  $Q\times P'\-$set isomorphism
$$f\times_P {\rm id_{P'}} : {\rm Res}_{Q\times P'}(\Omega\times_P P') \cong 
{\rm Res}_{\varphi_*\times {\rm id}_{P'}}\big({\rm Res}_{\bar Q\times P'}(\Omega\times_P P')\big)
\eqno £9.7.4;$$
but, since $\F(P,Q)\i \F'(P,Q)\,,$ we also have a $Q\times P'\-$set isomorphism
$${\rm Res}_{Q\times P'}(\Omega') \cong {\rm Res}_{\varphi_*\times {\rm id}_{P}}
\big({\rm Res}_{\bar Q\times P'}(\Omega')\big)
\eqno £9.7.5;$$
hence, we get a $Q\times P'\-$set isomorphism
$${\rm Res}_{Q\times P'}(\Omega' - \Omega\times_P P') \cong {\rm Res}_{\varphi_*\times {\rm id}_{P}}
\big({\rm Res}_{\bar Q\times P'}(\Omega' - \Omega\times_P P')\big)
\eqno £9.7.6\phantom{.}$$
and therefore $f\times_P {\rm id_{P'}}$ can be extended to a  $Q\times P'\-$set isomorphism
$$f' : {\rm Res}_{Q\times P'}(\Omega') \cong {\rm Res}_{\varphi_*\times {\rm id}_{P}}
\big({\rm Res}_{\bar Q\times P'}(\Omega')\big)
\eqno £9.7.7.$$
At this point, it follows from Proposition~£9.6 that $f'(\hat\Omega^{^\frak X}_Q)
= \hat\Omega^{^\frak X}_{\bar Q}$ and then it is quite clear that
$$s^{_\frak X}_{\bar Q}\big(f'(\hat O)\big) = f'\big(s^{_\frak X}_Q (\hat O)\big)
\eqno £9.7.8,$$
which proves equality~£9.7.2.

\medskip

£9.8. Now, for any $Q\in \frak X$ and any $Q\times P'\-$orbit $O'$ in 
$\Omega^{^\frak X}_Q\times_P P'$ we clearly can choose a subset 
$\s^{^\frak X}_{Q,O'}$ of $\frak S_{\Omega'}(Q)$ (cf.~£4.3) containing the trivial element and fulfilling
$$(s^{_\frak X}_Q)^{-1} (O') = \{s'(O')\}_{s'\in \s^{^\frak X}_{Q,O'}}\qq
\vert (s^{_\frak X}_Q)^{-1} (O') \vert = \vert \s^{^\frak X}_{Q,O'}\vert
\eqno £9.8.1;$$
thus, according to Proposition~£9.7, for  any $\varphi\in \F (P,Q)$ and any $Q\times P\-$set isomorphism
$$f : {\rm Res}_{Q\times P}(\Omega) \cong {\rm Res}_{\varphi\times {\rm id}_{P}}
\big({\rm Res}_{P\times P}(\Omega)\big)
\eqno £9.8.2,$$
setting $\bar Q = \varphi (Q)$ and $\bar O' = (f\times_P {\rm id}_{P'})(O')\,,$ 
we can choose a bijection 
$$\sigma^{_\frak X}_{Q,O',\varphi,f}   : \s^{^\frak X}_{Q,O'}\cong 
\s^{^\frak X}_{\bar Q,\bar O'}
\eqno £9.8.3\phantom{.}$$
preserving the trivial element.

\medskip
£9.9. As we mention in~£9.3 above, we have to consider a suitable quotient of the $\F^{^\frak X}\-$locality
${\rm Res}_{\F^{^\frak X}}(\L'^{^{\rm b}})\,;$ let us denote by 
$\tilde\frak k^{^{\rm b,\frak X}}\colon \tilde\F^{^\frak X}\to \Ab$ the  {\it contravariant\/} functor mapping 
$Q\in \frak X$ on (cf.~Corollary~£8.4)
$$\tilde\frak k^{^{\rm b,\frak X}} (Q) = \prod_{O} \ab \big({\rm Aut}(O)\big)
\eqno £9.9.1,$$
where $O$ runs over a set of representatives for the set of isomorphism classes of 
$Q\times P'\-$sets $(Q\times P')/\Delta_{\eta'} (T')$ where $T'$ is a subgroup of $P'$ such that any subgroup $U'$ of $P'$ fulfilling $\F' (T',U')\not= \emptyset$ does not belong to $\frak X\,,$ and $\eta'$ belongs to~$\F' (Q,T')\,;$
actually, this definition still makes sense for any subgroup $Q$ of $P\,,$ thus defining a {\it contravariant\/} functor
$\tilde\F\to \Ab$ that we still denote by $\tilde\frak k^{^{\rm b,\frak X}}\,;$  then, the quotient 
{\it $\F\-$locality\/} $ {\rm Res}_{\F}(\L'^{^{\rm b}})/\tilde\frak k^{^{\rm b,\frak X}}$ (cf.~£2.10) just maps any pair of subgroups $Q$ and $R$ of $P$ such that $R\not\in \frak X$ on $\F (Q,R)\,;$ this remark will be useful in~£9.11 below.

\bigskip
\noindent
{\bf Theorem~£9.10.} {\it With the notation and the choice above, there is a unique 
$\F^{^\frak X}\!\-$locality functor
$$\tilde\frak l^{^\frak X} : \bar\L^{^{\rm n,\frak X}}\too 
{\rm Res}_{\F^{^\frak X}}(\L'^{^{\rm b}})/\tilde\frak k^{^{\rm b,\frak X}}
\eqno £9.10.1\phantom{.}$$
in such a way that, for any $Q,R\in \frak X\,,$ any $\varphi\in \F (Q,R)$ and
any $R\times P\-$set isomorphism
$$f : {\rm Res}_{R\times P}(\Omega) \cong {\rm Res}_{\varphi\times {\rm id}_{P}}
\big({\rm Res}_{Q\times P}(\Omega)\big)
\eqno £9.10.2,$$
$\tilde\frak l^{^\frak X}$ maps the image $\bar f$ of $f$ in  $\bar\L^{^{\rm n,\frak X}}\!(Q,R)$ on the class $\tilde f'$ in 
$\L'^{^{\rm b,\frak X}}\!(Q,R)\big/\tilde\frak k^{^{\rm b,\frak X}}(R)$
of a $R\times P'\-$set isomorphism
$$f' : {\rm Res}_{R\times P'}(\Omega') \cong 
{\rm Res}_{\varphi\times {\rm id}_{P'}}\big({\rm Res}_{Q\times P'}(\Omega')\big)
\eqno £9.10.3\phantom{.}$$
which, for any  $R\times P'\-$orbit $O'$ in $\Omega^{^\frak X}_R\times_P P'\,,$ any $\omega'\in O'$ and any $s'\in \s^{^\frak X}_{R,O'}\,,$ fulfills  
$$f'\big(s'(\omega')\big) = \bar s' \big((f\times_P {\rm id}_{P'})(\omega')\big)
\eqno £9.10.4\phantom{.}$$
where we set $\bar s' = \sigma^{_\frak X}_{R,O',\varphi^{_P},f}(s')$ for $\varphi^{_P} = \iota_{_Q}^{_P}\circ \varphi\,.$\/}
\eject

\medskip
\noindent
{\bf Proof:} Setting $\bar R = \varphi (R)$ and denoting by $\varphi_*\,\colon R\cong \bar R$ the isomorphism determined by $\varphi,$ we already know that 
$f(\Omega^{^\frak X}_R) = \Omega^{^\frak X}_{\bar R}$ (cf.~£9.5) and therefore $f$ induces an $R\times P'\-$set isomorphism
$$f^{^\frak X}_R = f\times_P {\rm id}_{P'} : \Omega^{^\frak X}_R\times_P P' \cong {\rm Res}_{\varphi_*\times {\rm id}_{P'}} (\Omega^{^\frak X}_{\bar R}\times_P P')
\eqno £9.10.5;$$
then, it is quite clear that this $R\times P'\-$set isomorphism can be uniquely extended to a  $R\times P'\-$set isomorphism
$$\hat f^{^\frak X}_R : \hat\Omega^{^\frak X}_R  \cong {\rm Res}_{\varphi_*\times {\rm id}_{P'}} (\hat\Omega^{^\frak X}_{\bar R})
\eqno £9.10.6\phantom{.}$$
fulfilling condition~£9.10.4; once again, since  we also have an $R\times P'\-$set isomorphism
$${\rm Res}_{R\times P'}(\Omega') \cong {\rm Res}_{\varphi_*\times {\rm id}_{P}}
\big({\rm Res}_{\bar R\times P'}(\Omega')\big)
\eqno £9.10.7,$$
$\hat f^{^\frak X}_R$ can be extended to a $R\times P'\-$set isomorphism
$$\eqalign{f' : {\rm Res}_{R\times P'}(\Omega') \cong 
&\,{\rm Res}_{\varphi_*\times {\rm id}_{P'}}\big({\rm Res}_{\bar R\times P'}(\Omega')\big)\cr
&\hskip50pt\Vert\cr
&\, {\rm Res}_{\varphi\times {\rm id}_{P'}}\big({\rm Res}_{Q\times P'}(\Omega')\big)\cr}
\eqno £9.10.8.$$

\smallskip
Moreover, if
$$g' : {\rm Res}_{R\times P'}(\Omega') \cong 
{\rm Res}_{\varphi_*\times {\rm id}_{P'}}\big({\rm Res}_{\bar R\times P'}(\Omega')\big)
\eqno £9.10.9\phantom{.}$$
is another $R\times P'\-$set isomorphism extending $\hat f^{^\frak X}_R$
then the composition
$$f'^{-1}\circ g' : {\rm Res}_{R\times P'}(\Omega')\cong 
{\rm Res}_{R\times P'}(\Omega')
\eqno £9.10.10\phantom{.}$$
acts trivially on $\hat\Omega^{^\frak X}_R$ and in particular  the image of 
$f'^{-1}\circ g'$ in $C_{G'}(R)/\frak S_{\Omega'} (R)$ belongs to 
$\tilde\frak k^{^{\rm b,\frak X}}(R)\,;$
that is to say, $f$ determines a unique class $\tilde f'$ in the quotient set
$\L'^{^{\rm b,\frak X}}(Q,R)\big/\tilde\frak k^{^{\rm b,\frak X}}(R)\,.$

\smallskip
On the other hand,  the class $\tilde f'$ does not depend on our choice of the sets $\s^{^\frak X}_{R,O'}$ and the bijections 
$\sigma^{_\frak X}_{R,O',\varphi^{_P},f}$ in~£9.8; indeed, another choice of them determines another  $R\times P'\-$set isomorphism
$$\hat g^{^\frak X}_R : \hat\Omega^{^\frak X}_R  \cong {\rm Res}_{\varphi_*\times {\rm id}_{P'}} (\hat\Omega^{^\frak X}_{\bar R})
\eqno £9.10.11\phantom{.}$$
extending $f^{^\frak X}_R$ and it is easily checked that the extension  to $\Omega'$
by the identity on $\Omega' -\hat\Omega^{^\frak X}_R$ of the difference
$(\hat f^{^\frak X}_R)^{-1}\circ \hat g^{^\frak X}_R$ 
 belongs to 
$\frak S_{\Omega'}(R)\,.$

\smallskip
Now, we claim that the correspondence mapping $f\in T_G (R,Q)$ on the class 
$\tilde f'\in \L'^{^{\rm b,\frak X}}(Q,R)\big/\tilde\frak k^{^{\rm b,\frak X}}(R)$ is functorial;
namely we claim that, for any $T\in \frak X\,,$ any $\psi\in \F (R,T)\,,$  any $T\times P\-$set isomorphism
$$g : {\rm Res}_{T\times P}(\Omega) \cong {\rm Res}_{\psi\times {\rm id}_{P}}
\big({\rm Res}_{R\times P}(\Omega)\big)
\eqno £9.10.12\phantom{.}$$
and any $T\times P'\-$set isomorphism
$$g' : {\rm Res}_{T\times P'}(\Omega') \cong 
{\rm Res}_{\psi\times {\rm id}_{P'}}\big({\rm Res}_{R\times P'}(\Omega')\big)
\eqno £9.10.13\phantom{.}$$
\eject
\noindent
which for any  $T\times P'\-$orbit $M'$ in $\Omega^{^\frak X}_T\times_P P'\,,$ 
any $\omega'\in M'$ and any $s'\in \s^{^\frak X}_{T,M'}\,,$ fulfills  
$$g'\big(s'(\omega')\big) = \bar s' \big((g\times_P {\rm id}_{P'})(\omega')\big)
\eqno £9.10.14\phantom{.}$$
where  $\bar s' = \sigma^{_\frak X}_{T,M',\psi^{_P},g}(s')\,,$ 
the above correspondence maps the composition $f\circ g$ on the composition
of the classes 
$$\tilde g'\in \L'^{^{\rm b,\frak X}}(R,T)\big/
\tilde\frak k^{^{\rm b,\frak X}}\!(T)\qq \tilde f'\in \L'^{^{\rm b,\frak X}}(Q,R)\big/
\tilde\frak k^{^{\rm b,\frak X}}\!(R)
\eqno £9.10.15.$$

\smallskip
Indeed, setting $\bar T = \psi (T)$ and $\skew4\bar{\bar T}= \varphi (\bar T)\,,$
and denoting by $\psi_*\,\colon T\cong \bar T$ and  $\bar\varphi_*\,\colon \bar T\cong \skew4\bar{\bar T}$ the respective isomorphisms induced by $\psi$  and  $\varphi_*\,,$ it is clear that $g'$ and $f'$ induce a $T\times P'\-$ and a  $\bar T\times P'\-$set isomorphisms
$$\eqalign{g' : {\rm Res}_{T\times P'}(\Omega') &\cong {\rm Res}_{\psi_*\times {\rm id}_{P'}}\big({\rm Res}_{\bar T\times P'}(\Omega')\big)\cr
\bar f' : {\rm Res}_{\bar T\times P'}(\Omega') &\cong \,{\rm Res}_{\bar\varphi_*\times {\rm id}_{P'}}\big({\rm Res}_{\skew4\bar{\bar T}\times P'}(\Omega')\big)\cr}
\eqno £9.10.16$$
and, according to Proposition~£9.6, we have 
$\bar f' (\hat\Omega^{^\frak X}_{\bar T}) = 
\hat\Omega^{^\frak X}_{\skew4\bar{\bar T}}\,;$ moreover, since $f'$ extends 
$f^{^\frak X}_R\,,$ it is clear that $\bar f'$ extends the corresponding 
$f^{^\frak X}_{\bar T}\,;$ in particular, denoting by $O'$ the $R\times P'\-$orbit
in $\Omega^{^\frak X}_R\times_P P'$ containing $\bar M' = g'(M')$ and by 
$\bar O'$ the $\bar R\times P'\-$orbit in $\Omega^{^\frak X}_{\bar R}\times_P P'$ containing $\skew3\bar{\bar M}' = f'(\bar M')\,,$ it follows again from 
Proposition~£5.17 that we can do our choice in~£9.8 above in such a way that we have
$\s^{^\frak X}_{\bar T,\bar M'}\i \s^{^\frak X}_{R,O'}$ and
 $\s^{^\frak X}_{\skew4\bar{\bar T},\skew3\bar{\bar M}'}\i 
 \s^{^\frak X}_{\bar R,\bar O'}\,,$  that $\sigma^{_\frak X}_{R,O',\varphi^{_P},f}$
 maps $\s^{^\frak X}_{\bar T,\bar M'}$ onto $\s^{^\frak X}_{\skew4\bar{\bar T},\skew3\bar{\bar M}'}\,,$ and that 
 $\sigma^{_\frak X}_{\bar T,\bar M',\bar\varphi_*^{_P},\bar f'}$ coincides with the restriction
 of $\sigma^{_\frak X}_{R,O',\varphi^{_P},f}\,.$

 \smallskip
 In this situation, considering the $ T\times P'\-$set isomorphism
$$\bar f'\circ g' : {\rm Res}_{ T\times P'}(\Omega') \cong 
\,{\rm Res}_{\bar\varphi_*\circ \psi_*\times {\rm id}_{P'}}\big({\rm Res}_{\skew4\bar{\bar T}\times P'}(\Omega')\big)
\eqno £9.10.17,$$
for any  $T\times P'\-$orbit $M'$ in $\Omega^{^\frak X}_T\times_P P'\,,$ 
any $\omega'\in M'$ and any $s'\in \s^{^\frak X}_{T,M'}\,,$ we get  
$$\eqalign{(\bar f'\circ g')\big(s'(\omega')\big) &= \bar f'\Big(\bar s' \big((g\times_P {\rm id}_{P'})(\omega')\big)\big)\cr
&= \skew2\bar{\bar s}'\Big(\big((f\circ g)\times _P {\rm id}_{P"}\big)(\omega')\Big)\cr} 
\eqno £9.10.18\phantom{.}$$
where $\bar s' = \sigma^{_\frak X}_{T,M',\psi^{_P},g}(s')$ and 
$\skew2\bar{\bar s}' = \sigma^{_\frak X}_{\bar T,\bar M',\varphi_*^{_P},\bar f'}(\bar s')\,;$
thus, the composition $f'\circ g'$ also fulfills the corresponding condition~£9.10.4
and therefore our correspondence above maps $f\circ g$ on the class of $f'\circ g'$
in $\L'^{^{\rm b,\frak X}}(Q,T)\big/\tilde\frak k^{^{\rm b,\frak X}}\!(T)\,,$ proving our claim.

\smallskip
Finally, recall that (cf.~£4.14)
$$\bar\L^{^{\rm n,\frak X}}(Q,R) = \L^{^{\Omega,\frak X}}(Q,R)\big/\big(\tilde\frak k^{^\frak X}(R)\times \tilde\frak c^{^{\rm nsc}}(R)\big)
\eqno £9.10.19$$
and that $\L^{^{\Omega,\frak X}}(Q,R) = T_G (R,Q)/\frak S^1_{\Omega}(R)$ (cf.~£4.7);
thus, if $Q = R\,,$  $\varphi = {\rm id}_R$ and $f$ belongs to the converse image
in $T_G (R,Q)$ of $\tilde\frak k^{^\frak X}(R)\times \tilde\frak c^{^{\rm nsc}}(R)\,,$ the action\break
\eject
\noindent
 of $f$ on $\Omega^{^\frak X}_R$ coincides with the action of some element in $\frak S^1_{\Omega}(R)\,;$
then, it is easily checked that  the action of the  uniquely extended  
 $R\times P'\-$set isomorphism
$$\hat f^{^\frak X}_R : \hat\Omega^{^\frak X}_R  \cong {\rm Res}_{{\rm id}_R\times {\rm id}_{P'}} (\hat\Omega^{^\frak X}_{ R})
\eqno £9.10.20\phantom{.}$$
fulfilling condition~£9.10.4 also coincides with the action of some element 
belonging to  $\frak S_{\Omega'}(R)\,;$ in this case, it is not difficult to check that
the class in~$\L'^{^{\rm b,\frak X}}(R)\big/\tilde\frak k^{^{\rm b,\frak X}}\!(R)$ of a 
$R\times P'\-$set isomorphism $f'$ extending $\hat f^{^\frak X}_R$ is trivial.
This proves the existence and the uniqueness of the functor  $\tilde\frak l^{^\frak X}$
in~£9.10.1; the compatibility with the corresponding structural functors is easily checked, proving that it is actually 
an {\it $\F^{^\frak X}\!\-$locality functor\/}. We are done.

\medskip
£9.11. It follows from section 6 that the {\it $\F^{^\frak X}\!\-$locality\/} 
$\bar\L^{^{\rm n,\frak X}}$ contains a {\it perfect $\F^{^\frak X}\!\-$locality\/}
$\P^{^\frak X} =\bar\P^{^\frak X}$ (cf.~£9.1.2) and therefore from Theorem~£9.10 we get a 
{\it $\F^{^\frak X}\!\-$locality functor\/}  
$$\frak h^{^{\!\frak X}} : \P^{^\frak X} \too {\rm Res}_{\F^{^\frak X}}
(\L'^{^{\rm b}})/\tilde\frak k^{^{\rm b,\frak X}}
\eqno £9.11.1;$$
as a matter of fact, $\frak h^{^{\!\frak X}}$ can be extended to a {\it $\F^{^\frak X}\!\-$locality\/} functor (cf.~£9.9)
$$\P\too {\rm Res}_\F (\L'^{^{\rm b}})/\tilde\frak k^{^{\rm b,\frak X}}
\eqno £9.11.2,$$
that we still denote by $\frak h^{^{\!\frak X}}\,,$
mapping any pair of subgroups $Q$ and $R$ of $P$ such that $R\not\in \frak X$ on  the structural map
$\P(Q,R)\to \F(Q,R)\,;$ this remark will be useful in the proof of the theorem below.

\bigskip
\noindent
{\bf Theorem~£9.12.} {\it Any $\F^{^\frak X}\!\-$locality functor from 
$\P^{^\frak X}$ to ${\rm Res}_{\F^{^\frak X}}(\L'^{^{\rm b}})/
\tilde\frak k^{^{\rm b,\frak X}}$ is naturally $\F^{^\frak X}\!\-$isomorphic to 
$\frak h^{^\frak X}\,.$\/}

\medskip
\noindent
{\bf Proof:} We argue by induction on $\vert \frak X\vert\,;$ if $\frak X = \{P\}$
then the statement follows from Proposition~£2.17; thus,  assume that 
$\frak X\not= \{P\}\,,$  choose a minimal element $U$ in $\frak X$ {\it fully normalized\/} in $\F$ and set 
$$\frak Y = \frak X - \{\theta(U)\mid \theta\in \F(P,U)\}
\eqno £9.12.1;$$
we may assume that the {\it full\/} subcategory of $\P^{^\frak X}$ over $\frak Y$ coincides
with $\P^{^\frak Y}\,;$ on the other hand, we have a canonical functor
$$\big({\rm Res}_{\F^{^\frak X}}(\L'^{^{\rm b}})/
\tilde\frak k^{^{\rm b,\frak X}})^\frak Y\too {\rm Res}_{\F^{^\frak Y}}
(\L'^{^{\rm b}})/\tilde\frak k^{^{\rm b,\frak Y}}
\eqno £9.12.2\phantom{.}$$
which composed with the restriction of $\frak h^{^{\!\frak X}}$ to $\P^{^\frak Y}$
clearly coincides with $\frak h^{^{\!\frak Y}}\,.$

\smallskip
In particular, by the induction hypothesis the restriction to $\P^{^\frak Y}$ of any {\it $\F^{^\frak X}\!\-$locality functor\/}
$$\frak f^{^{\frak X}} : \P^{^\frak X} \too {\rm Res}_{\F^{^\frak X}}
(\L'^{^{\rm b}})/\tilde\frak k^{^{\rm b,\frak X}}
\eqno £9.12.3\phantom{.}$$
 composed with the canonical functor above is {\it naturally $\F^{^\frak Y}\-$isomorphic\/} to~$\frak h^{^{\!\frak Y}}\,;$ 
 thus, up to a modification of $\frak f^{^{\frak X}}$ by conjugation with a suitable 
element\break
\eject
\noindent
 in~$\big({\rm Res}_{\F^{^{\frak X}}}(\L'^{^{\rm b}})\big)\!(P)\big/\tilde\frak k^{^{\rm b,\frak X}}\!(P)$ 
(cf.~£2.9.3), we may assume that the restriction of~$\frak f^{^{\frak X}}$ to $\P^{^\frak Y}$ composed with the canonical functor above coincides with~$\frak h^{^{\!\frak Y}}\,.$
In this situation, the converse image $\widehat{\frak h^{^\frak Y}(\P^{^\frak Y})}$
in $\big({\rm Res}_{\F^{^\frak X}}(\L'^{^{\rm b}})/\tilde\frak k^{^{\rm b,\frak X}})^\frak Y$ {\it via\/} 
homomorphism~£9.12.2 of the image $\frak h^{^\frak Y}\!(\P^{^\frak Y})$ of $\P^{^\frak Y}$ by
$\frak h^{^{\frak Y}}$ contains the image $\frak f^{^\frak X} \!(\P^{^\frak Y})$ of~$\P^{^\frak Y}\,;$ 
more explicitly, denoting by $\widehat{\frak h^{^\frak X}(\P^{^\frak X})}$ the subcategory of 
${\rm Res}_{\F}(\L'^{^{\rm b}})/\tilde\frak k^{^{\rm b,\frak X}}$ (cf.~£9.9)
which coincides with $\widehat{\frak h^{^\frak Y}(\P^{^\frak Y})}$ over $\frak Y$
and maps on
$$\widehat{\frak h^{^\frak X}(\P^{^\frak X})}(Q,R) = 
\big({\rm Res}_{\F}(\L'^{^{\rm b}})\big)(Q,R)\big/
\tilde\frak k^{^{\rm b,\frak X}}(R)
\eqno £9.12.4\phantom{.}$$
 any pair of subgroups $Q$ and $R$ of $P$ such that $R\not\in \frak Y\,,$
we get two {\it $\F\-$locality functors\/}  (cf.~£9.11.2)
$$\frak f^{^{\frak X}} : \P \too \widehat{\frak h^{^\frak X}(\P^{^\frak X})}\qq 
\frak h^{^{\frak X}} : \P \too \widehat{\frak h^{^\frak X}(\P^{^\frak X})}
\eqno £9.12.5.$$
where we still denote by $\frak f^{^\frak X}$ the obvious extension of the functor~£9.12.3
mapping any pair of subgroups $Q$ and $R$ of $P$ such that $R\not\in \frak X$ on  the structural map
$\P(Q,R)\to \F(Q,R)\,.$

\smallskip
Set $\M^{^\frak X} = \widehat{\frak h^{^\frak X}(\P^{^\frak X})}$ for short, and 
denote by $\rho^{_\frak X}\,\colon \M^{^\frak X}\to \F$ the second structural 
functor; first of all note that, according to~£2.16, we have  functors
$$\loc_{\P} : \ch^*(\F)\too \widetilde{\Loc}\qq
\loc_{\M^{^\frak X}} : \ch^*(\F)\too \widetilde{\Loc}
\eqno £9.12.6\phantom{.}$$
and it is clear that the functors $\frak f^{^{\frak X}}$ and $\frak h^{^{\frak X}}$ determine {\it natural\/} maps
$\loc_{\frak f}$ and $\loc_{\frak h}$ from $\loc_{\P}$ to 
$\loc_{\M^{^\frak X}}\,;$ moreover, we know that 
$\loc_{\P} = \loc_\F$ (cf.~£2.15) and,
since $\frak f^{^{\frak X}}$ and $\frak h^{^{\frak X}}$ are {\it $\F\-$locality functors\/}, it is easily checked that 
$\loc_{\frak f^{^{\frak X}}}$ and $\loc_{\frak h^{^{\frak X}}}$ fulfill the conditions in~Proposition~£2.17; consequently, it follows from the uniqueness part of this proposition that we actually have the equality $\loc_{\frak f^{^{\frak X}}} = \loc_{\frak h^{^{\frak X}}}\,.$ In particular, for any  {\it $\P\-$chain\/} $\frak q\,\colon \Delta_n \to \P\,,$ denoting by 
$\bar\frak q\,\colon \Delta_n \to \F$ the corresponding {\it $\F\-$chain\/}, we have (cf.~£2.16)
$$\eqalign{\loc_{\P} (\bar\frak q) &= \big(\P (\frak q),{\rm Ker}(\pi_\frak q)\big)\cr
\loc_{\M^{^\frak X}} (\bar\frak q)  &= \big(\M^{^\frak X}\! (\frak f^{^{\frak X}}\!\!\circ\frak q),
{\rm Ker}(\rho_{\frak f^{^{\frak X}}\!\circ\frak q})\big) = \big(\M^{^\frak X} \!(\frak h^{^{\frak X}}\!\!\circ\frak q),{\rm Ker}(\rho_{\frak h^{^{\frak X}}\!\circ\frak q})\big)\cr}
\eqno £9.12.7\phantom{.}$$
and it is not difficult to check that, up to replacing $\frak f^{^{\frak X}}$ by a {\it naturally
$\F\-$iso-morphic\/} functor, we may assume that we have
$$\M^{^\frak X} \!(\frak f^{^{\frak X}}\!\circ\frak q) = \M^{^\frak X} \!(\frak h^{^{\frak X}}\!\circ\frak q)
\eqno £9.12.8\phantom{.}$$
and that the group homomorphisms mapping any $\P\-$automorphism 
$\alpha\,\colon\frak q\cong \frak q$ on 
$$\frak f^{^{\frak X}}\! \!*\alpha : \frak f^{^{\frak X}}\!\circ\frak q\cong \frak f^{^{\frak X}}\!\circ\frak q\qq
\frak h^{^{\frak X}}\! \!*\alpha : \frak h^{^{\frak X}}\!\circ\frak q\cong \frak h^{^{\frak X}}\!\circ\frak q
\eqno £9.12.9\phantom{.}$$
coincide with each other.

\smallskip
On the other hand, since $\frak f^{^{\frak X}}$ and $\frak h^{^{\frak X}}$ coincide over $\P^{^\frak Y}$ and they are 
{\it $\F\-$locality functors\/}, for any $\P\-$morphism $\varphi\,\colon R\to Q$ we have
$$\frak f^{^{\frak X}} (\varphi) = \frak h^{^{\frak X}} (\varphi)\.\ell_\varphi
\eqno £9.12.10\phantom{.}$$
\eject
\noindent
for some element $\ell_\varphi\in \M^{^\frak X}(R)$
belonging either to the kernel of the canonical homomorphism (cf.~£9.12.2)
$$\L'^{^{\rm b}}(R)/\tilde\frak k^{^{\rm b,\frak X}}(R)\too \L'^{^{\rm b}}(R)/\tilde\frak k^{^{\rm b,\frak Y}}(R)
\eqno £9.12.11\phantom{.}$$
whenever $R$ belongs to $\frak Y\,,$ or to ${\rm Ker}(\rho^{_\frak X}_{_R})$ otherwise. Moreover, in the situation above, for any $\alpha\in \P(R)$ we get $\frak f^{^{\frak X}}\! (\alpha) = \frak h^{^{\frak X}} \!(\alpha)\,,$ so that 
$\ell_\alpha = 1\,;$ actually, up to replacing $\frak f^{^{\frak X}}$ by a {\it naturally $\F\-$isomorphic\/} functor, we may assume that we have $\frak f^{^{\frak X}} \!(\alpha) = \frak h^{^{\frak X}}\! (\alpha)$ for any $\P\-$isomorphism 
$\alpha\in \P (Q',Q)\,.$

\smallskip
 But, we already know that the kernel of the structural homomorphism
$\pi'^{^{\rm b}}_{_R}\,\colon \L'^{^{\rm b}}(R)\to \F'(R)$ is given by (cf.~£4.3.2)
$${\rm Ker}(\pi'^{^{\rm b}}_{_R}) = \prod_{\tilde O'\in \frak O'_R} 
\ab \big({\rm Aut} (O')\big)
\eqno £9.12.12\phantom{.}$$
where we denote by  $\frak O'_R$ the  set of isomorphism classes  of  the $R\times P'\-$sets $(R\times P')/\Delta_{\theta'} (T')$ for any subgroup $T'$ of $P'$ and any element $\theta'\in \F' (R,T')\,;$ consequently, if $R$ belongs to $\frak X\,,$ the kernel above involves the isomorphism classes of the $R\times P'\-$orbits in 
$\hat\Omega^{^\frak X}_R$ and therefore we can choose a set
of representatives in~$\Omega^{^\frak X}_R\times_P P'\,;$ explicitly, considering
the {\it contravariant\/} functor $\tilde\frak t^{^U}\colon \tilde\F\to \Ab$
introduced in Proposition~£8.9 above, since $\tilde\frak  t^{^U}\! (R) = \{0\}$ whenever
$R$ does not belong to $\frak X\,,$ in all the cases it is easily checked that $\ell_\varphi$ above belongs to (cf.~£3.6.1)
$$\eqalign{\tilde\frak t^{^{U}}\! (R) \cong \bigg(\prod_{\tilde\theta\in \tilde\F (R,U)} \ab\Big( {\rm Aut}\big((R\times P')/\Delta_{\theta} (U)\big)\Big)\bigg)^{\F_P(U)}\cr
\cong \bigg(\prod_{\tilde\theta\in \tilde\F (R,U)} \ab\Big(\bar N_{R\times P'}\big(\Delta_{\theta} (U)\big)\Big)\bigg)^{\F_P(U)}\cr}
\eqno £9.12.13.$$

\smallskip
That is to say, we have obtained a correspondence mapping any {\it $\P\-$chain\/}  $\frak q\,\colon \Delta_1\to \P$  
on~$\ell_{\frak q (0\bullet 1)}\in \tilde\frak t^{^{U}}\! \big(\frak q (0)\big)$ and we may assume that it vanish over the $\P\-$isomorphisms; in this case,  we claim that it determines a {\it stable\/} element $\ell$ 
of~$ \Bbb C^1 (\tilde\F,\tilde\frak  t^{^U)}$ [11,~A3.17]. Indeed, for another isomorphic {\it $\P\-$chain\/} $\frak q'\,\colon \Delta_1\to \P$ 
and a {\it natural isomorphism\/} $\nu\,\colon \frak q\cong\frak q'\,,$
setting 
$$\eqalign{\varphi = \frak q (0\bullet 1)\!\!\quad,\quad  \!\!\varphi' = \frak q'(0\bullet 1)
\!\!\quad,\quad  \!\!\nu_0 = \beta\qq \nu_1 = \alpha\cr}
\eqno £9.12.14,$$
  from~our choice we have $\ell_\alpha = 1$ and $\ell_\beta = 1$ and therefore we get
$$\eqalign{\frak f^{^{\frak X}}(\varphi') &= \frak h^{^{\frak X}} (\varphi')\.\ell_{\varphi'} = \big(\frak h^{^{\frak X}} (\alpha)\.\frak h^{^{\frak X}}(\varphi)\.\frak h^{^{\frak X}} (\beta)^{-1}\big)\.\ell_{\varphi'}\cr
&=  \big(\frak f^{^{\frak X}} (\alpha)\.\frak f^{^{\frak X}} (\varphi)\.\ell_\varphi^{-1}\.\frak h^{^{\frak X}}(\beta)^{-1}\big)\.\ell_{\varphi'}\cr
& = \frak f^{^{\frak X}} (\varphi')\. \big(\tilde\frak t^{^U}\!(\tilde\beta^{-1})\big)( \ell_{\varphi}^{-1})\.\ell_{\varphi'}\cr}
\eqno £9.12.15\phantom{.}$$ 
which proves that the correspondence $\ell$
sending $\tilde\varphi$ to~$\ell_{\varphi}$ is {\it stable\/} and, in particular, that $\ell_{\varphi}$ only depends on the corresponding $\tilde\P\-$morphism $\tilde\varphi\,.$
\eject

\smallskip
Moreover, we also claim that $d_1 (\ell) = 0\,;$ indeed, for a second $\P\-$mor-phism 
$\psi\,\colon T\to R$
we get
$$\eqalign{\frak f^{^{\frak X}} (\varphi\circ\psi)& = \frak f^{^{\frak X}} (\varphi)\.\frak f^{^{\frak X}} (\psi) = 
\big(\frak h^{^{\frak X}} (\varphi)\.\ell_\varphi\big)\.\big(\frak h^{^{\frak X}} (\psi)\.\ell_\psi\big)\cr
&= \frak h^{^{\frak X}}(\varphi\circ\psi)\.\big(\tilde\frak t^{^U}\!(\tilde\psi)\big)(\ell_{\varphi})\.\ell_\psi\cr}
\eqno £9.12.16\phantom{.}$$
and   the {\it divisibility\/} of $\M^{^{\frak X}}$ forces
$$\ell_{\varphi\circ\psi} =\big(\tilde\frak t^{^U}\!(\tilde\psi)\big)(\ell_{\varphi})\.\ell_\psi
\eqno £9.12.17;$$
since $\tilde\frak t^{^U}\! (T)$ is Abelian, in  the additive notation we obtain
$$0 = \big(\tilde\frak t^{^U}\!(\tilde\psi)\big)(\ell_{\varphi}) - \ell_{\varphi\circ\psi} + \ell_\psi
\eqno £9.12.18,$$
proving our claim.

\smallskip
Finally, since $\Bbb H_*^1(\tilde\F,\tilde\frak t^{^U}) = \{0\}$ (cf.~£8.6.4), we have $\ell = d_0 (n)$ for some
element $n = (n_Q)_Q$ of 
$\Bbb C^0 \big(\tilde\F,\tilde\frak t^{^U}\big)\,;$ that is to say, with the 
notation above we get
$$\ell_{\varphi} = \big(\tilde\frak t^{^U}\!(\tilde\varphi)\big) (n_Q)\.n_R^{-1}
\eqno £9.12.19\phantom{.}$$
where we  identify any $\tilde\F\-$object with the obvious {\it $\tilde\F\-$chain\/}
$\Delta_0\to \tilde\F\,;$ hence,  we obtain
$$\frak f^{^{\frak X}}(\varphi) =    \frak h^{^{\frak X}} (\varphi)\.\big(\tilde\frak t^{^U}\!
(\tilde\varphi)\big) (n_Q)\.n_R^{-1} = n_Q\.\frak h^{^{\frak X}}(\varphi)\.n_R^{-1}
\eqno £9.12.20,$$
which amounts to saying that the correspondence  sending $Q$ to 
$n_Q$ defines  a {\it natural isomorphism\/} between $\frak h$ and $\frak f\,;$ 
the compatibility with the cor-responding structural functors is easily checked,
proving that $n$ defines a {\it natural $\F\-$isomorphism\/}.
We are done.

\medskip
£9.13. Recall that, according to Theorem~£8.10 above, we have a canonical
{\it $\F'\-$locality\/} functor
$$\frak h' : \P' \too \L'^{^{\rm b}}
\eqno £9.13.1;$$
thus, we also have an {\it $\F\-$locality\/} functor
$${\rm Res}_\F(\frak h') : {\rm Res}_\F(\P') \too {\rm Res}_\F(\L'^{^{\rm b}})
\eqno £9.13.2\phantom{.}$$
and we consider the induced {\it $\F^{^{\rm sc}}\-$locality functor\/}  (cf.~£9.1.2)
$${\rm Res}_{\F^{^{\rm sc}}}(\bar\frak h') : {\rm Res}_{\F^{^{\rm sc}}}(\bar\P') \too {\rm Res}_{\F^{^{\rm sc}}}(\L'^{^{\rm b}})/\frak k^{^{\rm b,sc}}
\eqno £9.13.3\phantom{.}$$
which is actually {\it faithful\/} as it proves  the following lemma; recall that the {\it kernel\/} of the structural functor 
$\bar\P'\to \F'$ is given by the {\it contravariant\/} functor [11,~Proposition~13.14]
$$\frak c^\frak f_{\F'} = \frak c^\frak h_{\F'}/[\frak c^\frak h_{\F'},\frak c^\frak h_{\F'}]
\eqno £9.13.4.$$
\eject

\bigskip
\noindent
{\bf Lemma~£9.14.} {\it For any  subgroup $Q$ in~$\frak X$ 
the group homomorphism 
$$\frak c^\frak f_{\F'} (Q)\too {\rm Ker}(\pi'^{^{\rm b}}_{_Q})/
\tilde\frak k^{^{\rm b,\frak X}} (Q)
\eqno £9.14.1\phantom{.}$$
determined by the  $\F^{^{\frak X}}\!\-$locality functor ${\rm Res}_{\F^{^{\frak X}}}(\bar\frak h')$ admits 
a section $\sigma'^{_{\frak X}}_Q$ which is stable  through
$\F\-$isomorphisms.\/}

\medskip
\noindent
{\bf Proof:} Choose an $\F'\-$morphism $\varphi'\,\colon Q\to P'$ such that 
$Q' = \varphi'(Q)$ is {\it fully centralized\/} in $\F'\,;$ then, we know that
$\frak c^\frak f_{\F'} (Q')$ is the {\it direct limit\/}  of the cano-nical functor from 
the {\it Frobenius $C_{P'}(Q')\-$category\/} $C_{\F'}(Q')$ to $\Gr$ [11,~13.1
and~Proposition~£13.14]
and therefore we have a canonical group homomorphism
$$\rho_{Q'} : \ab\big(C_{P'}(Q')\big)\too \frak c^\frak f_{\F'} (Q')
\eqno £9.14.2.$$

\smallskip
On the other hand, we know that (cf.~£9.9.1)
$${\rm Ker}(\pi'^{^{\rm b}}_{_Q})/\tilde\frak k^{^{\rm b,\frak X}} (Q)\cong \prod_{O'} \ab \big({\rm Aut}(O)\big)
\eqno £9.14.3,$$
where $O'$ runs over a set of representatives for the isomorphism classes of $Q\times P'\-$sets $(Q\times P')/\Delta_{\eta'} (T')$ where $T'$ is a subgroup of $P'$ such that for some  subgroup $U$ in $\frak X$ of $P$ we have 
$\F' (T',U)\not= \emptyset\,,$ and $\eta'$ belongs to~$\F' (Q,T')\,;$ in particular, for $T' = Q'$ and $\eta' $ equal to the inverse $\varphi'^*$ of the isomorphism induced by $\varphi'\,,$  in the 
right-hand member of equality~£9.14.3 we have the factor (cf.~£3.6.1)
$$\ab\Big(\bar N_{Q\times P'}\big(\Delta_{\varphi'^*} (Q')\big)\Big) \cong 
\ab\big(C_{P'}(Q')\big) 
\eqno £9.14.4.$$

\smallskip
At this point, we denote by $\sigma'^{_{\frak X}}_Q$ the composition
$${\rm Ker}(\pi'^{^{\rm b}}_{_Q})/\tilde\frak k^{^{\rm b,\frak X}} (Q)\too 
\ab\big(C_{P'}(Q')\big)\buildrel \rho_{Q'}\over\too \frak c^\frak f_{\F'} (Q')
\buildrel \frak c^\frak f_{\F'}(\varphi'_*)\over{\hbox to 30pt {\rightarrowfill}} \frak c^\frak f_{\F'} (Q)
\eqno £9.14.5\phantom{.}$$
which is clearly stable  through $\F\-$isomorphisms and, since the 
{\it $\F'\-$locality\/} functor $\bar\frak h' \,\colon \bar\P' \to \L'^{^{\rm b}}$
maps $\bar\tau'_{_{Q'}}(u')$ on $\tau'^{_{\rm b}}_{_{Q'}}(u')$ for any 
$u'\in C_{P'}(Q')\,,$ it is easily checked that $\sigma'^{_{\frak X}}_Q$ is a section of the homomorphism~£9.14.1 above.

\bigskip
\noindent
{\bf Theorem~£9.15.} {\it There is a unique natural $\F\-$isomorphism class
of  $\F\-$locality functors $\,\bar\frak g\,\colon \bar\P\too {\rm Res}_{\F}(\bar\P')\,.$\/}

\medskip
\noindent
{\bf Proof:} This statement follows from Theorem~£7.2 applied to the {\it $p\-$coherent $\F\-$locality\/} 
${\rm Res}_{\F}(\bar\P')$ provided we prove that  there is a unique natural $\F^{^{\rm sc}}\!\-$iso-morphism class of  
$\F^{^{\rm sc}}\!\-$locality functors $\frak g^{_{\rm sc}}$ from 
$\P^{^{\rm sc}}\!$ to ${\rm Res}_{\F^{^{\rm sc}}\!}(\bar\P')\,;$   we actually
prove that, for any set $\frak X$ as in~£9.5 above, there is a unique natural 
$\F^{^{\frak X}}\!\-$isomorphism class of  $\F^{^{\frak X}}\!\-$locality functors 
$\frak g^{_{\frak X}}$ from $\P^{^{\frak X}}\!$ to ${\rm Res}_{\F^{^{\frak X}}\!}(\bar\P')\,.$

\smallskip
Arguing by induction on $\vert\frak X\vert\,,$ if $\frak X = \{P\}$ then the statement follows from Proposition~£2.17 above; 
thus,  assume that $\frak X\not= \{P\}\,,$  choose a minimal element $U$ in $\frak X$ {\it fully normalized\/} in $\F$ and set 
$$\frak Y = \frak X - \{\theta(U)\mid \theta\in \F(P,U)\}
\eqno £9.15.1;$$
\eject
\noindent
it follows from the induction hypothesis that, with the same notation as above,  there is a unique natural 
$\F^{^{\frak Y}}\!\-$isomorphism class of  {\it $\F^{^{\frak Y}}\!\-$locality functors\/} 
$\frak g^{_{\frak Y}}\colon\P^{^{\frak Y}}\!\to {\rm Res}_{\F^{^{\frak Y}}\!}(\bar\P')\,;$ then, considering the composition ${\rm Res}_{\F^{^{\frak Y}}}(\bar\frak h')\circ \frak g^{_{\frak Y}}\,,$ it follows from Theorem~£9.12 that,  for a suitable choice of 
a {\it $\F^{^{\frak Y}}\-$locality functor\/}~$\frak h^{^{\frak Y}}$, the image of this functor is contained in the image
of ${\rm Res}_{\F^{^{\frak Y}}}(\bar\frak h')$ (cf.~£2.9.3).

\smallskip
But, always according to this theorem,  $\frak h^{^{\frak Y}}$ can be extended to
an {\it $\F^{^{\frak X}}\-$lo-cality\/} functor $\frak h^{^{\frak X}}$ over 
$\P^{^{\frak X}}\,;$ thus, the image of this functor is contained in the
{\it $\F^{^{\frak X}}\-$sublocality\/} $\M^{^{\frak X}}$ of 
${\rm Res}_{\F^{^{\frak X}}}(\L'^{^{\rm b}})/\frak k^{^{\rm b,\frak X}}$ which coincides over $\frak Y$ 
with the image  of 
${\rm Res}_{\F^{^{\frak Y}}}(\bar\frak h')$ and maps any  $Q\in \frak X$ and
any $V\in \frak X -\frak Y$ on
$$\M^{^{\frak X}}(Q,V) = \big({\rm Res}_{\F}(\L'^{^{\rm b}})\big)(Q,V)\big/
\tilde\frak k^{^{\rm b,\frak X}}(V)
\eqno £9.15.2;$$
in particular, $\frak h^{^\frak X}$ and $\bar\frak h'$ induce two functors
$$\P^{^\frak X}\too \M^{^\frak X}\longleftarrow 
{\rm Res}_{\F^{^{\frak X}}}(\bar\P')
\eqno £9.15.3.$$

\smallskip
Then, denoting by $\rho^{_{\frak X}}\colon \M^{^\frak X}\to \F^{^\frak X}$
the second structural functor of $\M^{^\frak X}\,,$ it follows from Lemma~£9.14 that the correspondence
$\frak d^{^{\!\frak X}}\colon \tilde\F^{^\frak X}\to \Ab$ mapping any $Q\in \frak Y$ on $\{0\}$ and 
any $V\in \frak X -\frak Y$ on ${\rm Ker}(\sigma'^{_\frak X}_Q)$ actually defines a subfunctor of
$$\Ker (\rho^{_{\frak X}}) : \tilde\F^{^\frak X}\too \Ab
\eqno £9.15.4$$
and therefore we get a quotient {\it $\F^{^\frak X}\!\-$locality\/} $\M^{^\frak X}\!/\frak d^{^\frak X}$ (cf.~£2.9). 
At this point, it is easily checked that the composition of the right-hand functor in~£9.15.3 with the canonical functor
$\M^{^\frak X}\!\to\M^{^\frak X}\!/\frak d^{^\frak X}$  induces an {\it $\F^{^\frak X}\!\-$locality  isomorphism\/}
$${\rm Res}_{\F^{^{\frak X}}}(\bar\P')\cong \M^{^\frak X}\!/\frak d^{^\frak X}
\eqno £9.15.5;$$
thus,  the composition of the left-hand functor in~£9.15.3 with the  functor
$\M^{^\frak X}\!\to  {\rm Res}_{\F^{^{\frak X}}}(\bar\P')$ obtained from isomorphism~£9.15.5 supplies the announced {\it $\F^{^\frak X}\!\-$locality\/} functor
$$\frak g^{_{\frak X}} : \P^{^{\frak X}}\too {\rm Res}_{\F^{^{\frak X}}}(\bar\P')
\eqno £9.15.6.$$

\smallskip
Moreover, if we have another  {\it $\F^{^\frak X}\!\-$locality\/} functor 
$\hat\frak g^{_{\frak X}} \colon \P^{^{\frak X}}\!\to {\rm Res}_{\F^{^{\frak X}}}(\bar\P')$ then, from the induction 
hypothesis, we may assume that the restriction of this functor to $\P^{^\frak Y}$ coincides with $\frak g^{_\frak Y}$ and
may choose the same  {\it $\F^{^{\frak X}}\-$locality functor\/} $\frak h^{^{\frak X}}$ over $\P^{^{\frak X}}$; now, the images 
of~$\frak h^{^{\frak X}}$ and of the composition of~$\hat\frak g^{_{\frak X}}$ with the right-hand functor in~£9.15.3 are contained in $\M^{^\frak X}$ and it is not difficult to check from Theorem~£9.12 that the corresponding functors from  
$\P^{^\frak X}$ to~$\M^{^\frak X}$ still are {\it naturally $\F\-$isomorphic\/}. Finally, the compositions of these functors  with  the  functor $\M^{^\frak X}\!\to  {\rm Res}_{\F^{^{\frak X}}}(\bar\P')$ obtained from isomorphism~£9.15.5 
remain {\it naturally $\F\-$isomorphic\/} to each other and, on the other hand, respectively coincide with $\frak g^{^\frak X}$ and  with $\hat\frak g^{^\frak X}\,.$ We are done.
\eject

\vfill
\eject

\bigskip
\bigskip

\bigskip
\centerline{\large References}
 
\bigskip

\smallskip\noindent
[1]\phantom{.} Dave Benson, personal letter 1994
\smallskip\noindent
[2]\phantom{.} Carles Broto, Ran Levi and Bob Oliver,  {\it The homotopy theory
of fusion systems\/}, Journal of Amer. Math. Soc. 16(2003), 779-856.
\smallskip\noindent
[3]\phantom{.} Andrew Chermak. {\it Fusion systems and localities\/}, 
Acta Mathematica, 211(2013), 47-139.
\smallskip\noindent
[4]\phantom{.} Daniel Gorenstein, {\it ``Finite groups''\/} Harper's Series,
1968, Harper and Row.
\smallskip\noindent
[5]\phantom{.} Burkhard K\"ulshammer and Llu\'\i s Puig, {\it Extensions of
nilpotent blocks}, Inventiones math., 102(1990), 17-71.
\smallskip\noindent
[6] Bob Oliver. {\it Existence and Uniqueness of Linking Systems: Chermak's proof via obstruction theory\/}, 
Acta Mathematica, 211(2013), 141-175.
\smallskip\noindent
[7]\phantom{.} Llu\'\i s Puig, {\it Structure locale dans les groupes
finis\/}, Bulletin Soc. Math. France, M\'emoire 47(1976). 
\smallskip\noindent
[8]\phantom{.}  Llu\'\i s Puig, {\it Brauer-Frobenius categories\/}, Manuscript notes 1993
\smallskip\noindent
[9]\phantom{.} Llu\'\i s Puig, {\it Frobenius categories and localities: history and survey},
Workshop ``Topology, Representation theory and Cohomology'' 2005, EPFL 
\smallskip\noindent
[10]\phantom{.} Llu\'\i s Puig, {\it Frobenius categories},
Journal of Algebra, 303(2006), 309-357.
\smallskip\noindent
[11]\phantom{.} Llu\'\i s Puig, {\it ``Frobenius categories versus Brauer blocks''\/}, Progress in Math. 
274(2009), Birkh\"auser, Basel.
\smallskip\noindent
[12]\phantom{.} Llu\'\i s Puig, {\it A criterion on vanishing cohomology\/}, submitted to Journal of Algebra, arxiv.org/abs/1308.3765.
\smallskip\noindent
[13]\phantom{.} Llu\'\i s Puig, {\it The Hecke algebra of a Frobenius $P\-$category\/}, Algebra Colloquium 21(2014), 1-52.

\end